\documentclass[11pt,reqno]{amsart}
\usepackage[a4paper,margin=27mm]{geometry}
\usepackage{amsmath,amssymb,mathtools,mathrsfs}
\usepackage{enumitem,booktabs,array}
\usepackage{needspace}
\IfFileExists{lmodern.sty}{\usepackage{lmodern}\usepackage{microtype}}{}
\usepackage[english]{babel}
\usepackage{graphicx,tikz}
\usepackage{capt-of,placeins}
\usetikzlibrary{arrows.meta}
\usepackage[colorlinks=true,linkcolor=blue,citecolor=blue,urlcolor=blue]{hyperref}

\numberwithin{equation}{section}

\theoremstyle{plain}
\newtheorem{theorem}{Theorem}[section]
\newtheorem{lemma}[theorem]{Lemma}
\newtheorem{proposition}[theorem]{Proposition}
\newtheorem{corollary}[theorem]{Corollary}

\newtheorem{fact}[theorem]{Fact}

\theoremstyle{definition}
\newtheorem{definition}[theorem]{Definition}
\newtheorem{example}[theorem]{Example}
\theoremstyle{remark}
\newtheorem{remark}[theorem]{Remark}
\newtheorem*{idea}{Idea of the proof}

\newcommand{\C}{\mathbb C}\newcommand{\Z}{\mathbb Z}
\newcommand{\N}{\mathbb N}\newcommand{\R}{\mathbb R}
\newcommand{\Q}{\mathbb Q}\newcommand{\A}{\mathbb A}
\newcommand{\PP}{\mathbb P}
\newcommand{\GL}{\mathrm{GL}}\newcommand{\SL}{\mathrm{SL}}
\newcommand{\PGL}{\mathrm{PGL}}\newcommand{\SO}{\mathrm{SO}}
\newcommand{\cD}{\mathcal D}\newcommand{\cC}{\mathcal C}\newcommand{\cR}{\mathcal R}
\newcommand{\cM}{\mathcal M}
\newcommand{\der}{\mathrm{der}}\newcommand{\ad}{\mathrm{ad}}
\newcommand{\ab}{\mathrm{ab}}
\newcommand{\bs}{{\boldsymbol s}}\newcommand{\bcC}{{\boldsymbol{\cC}}}

\newcommand{\red}{\mathrm{red}}

\DeclareMathOperator{\Spec}{Spec}\DeclareMathOperator{\Proj}{Proj}
\DeclareMathOperator{\Hom}{Hom}\DeclareMathOperator{\End}{End}
\DeclareMathOperator{\tr}{tr}\DeclareMathOperator{\diag}{diag}
\DeclareMathOperator{\relint}{relint}\DeclareMathOperator{\cone}{cone}
\DeclareMathOperator{\gr}{gr}\DeclareMathOperator{\Frac}{Frac}
\DeclareMathOperator{\id}{id}\DeclareMathOperator{\Ad}{Ad}
\DeclareMathOperator{\gp}{gp}\DeclareMathOperator{\ord}{ord}
\DeclareMathOperator{\inn}{in}\DeclareMathOperator{\Lie}{Lie}
\DeclareMathOperator{\divi}{div}\DeclareMathOperator{\Span}{span}
\DeclareMathOperator{\Ext}{Ext}\DeclareMathOperator{\rank}{rank}
\DeclareMathOperator{\Pic}{Pic}\DeclareMathOperator{\Cl}{Cl}

\hypersetup{pdftitle={Log Calabi--Yau compactifications of relative character varieties},pdfauthor={Tao Su}}
\hypersetup{bookmarksdepth=2}

\makeatletter
\def\l@part{\@tocline{-1}{4pt plus1pt}{0pt}{}{\bfseries}}
\makeatother

\theoremstyle{plain}
\newtheorem{mainresult}{Theorem}

\begin{document}
\title[Log Calabi--Yau relative character varieties]{Log Calabi--Yau compactifications of relative character varieties}
\author{Tao Su}
\email{sutao@bimsa.cn}
\address{Beijing Institute of Mathematical Sciences and Applications, Beijing, China}

\date{}

\subjclass[2020]{Primary 14M35, 14M27;
  Secondary 14E30, 14B05, 14M25, 53D17}
\keywords{Relative character varieties for reductive groups,
  log Calabi--Yau compactifications, log canonical pairs,
  conjugacy-class closures, ideal triangulations,
  toric degenerations, Poisson geometry, normality}

\begin{abstract}
A folklore conjecture predicts that smooth relative character varieties admit
log Calabi--Yau compactifications. Let $G$ be a connected complex reductive
group and $\Sigma$ a punctured compact Riemann surface. We study relative
$G$-character varieties of $\Sigma$ with prescribed conjugacy-class closures
at the punctures, including the spectral case of fixed conjugation invariants.
Assuming that the prescribed semisimple parts satisfy an intrinsic genericity
condition, we prove that every irreducible component $U$ of the normalized
relative character variety has a normal projective compactification $X$ with
reduced boundary $D$, such that $(X,D)$ is log canonical and $K_X+D\sim0$
integrally. Moreover, $X$ can be chosen of Fano type.

Normalization is unnecessary when all prescribed class closures are normal;
for nonempty varieties in positive genus, the converse holds. In genus zero,
we give an exact local normality criterion. In the spectral case, $U$ is
$\mathbb{Q}$-factorial with canonical singularities, and $(X,D)$ admits a
projective $\mathbb{Q}$-factorial dlt modification that is an isomorphism over
$U$.

For simply connected semisimple groups, ideal triangulations and Lusztig's
canonical bases degenerate spectral fibers to reduced Gorenstein unions of
toric varieties through a canonical monomial identity; central descent yields
the general spectral compactifications. For arbitrary class closures, a
Poisson compactification reduces log canonicity to a nonzero Goldman volume
period, obtained from a Lagrangian Hitchin fiber via nonabelian Hodge theory.
Examples show that a normal relative character variety need not admit a
$\mathbb{Q}$-factorial dlt compactification.
\end{abstract}

\maketitle
{\small\tableofcontents}

\section{Introduction}\label{sec:intro}

\subsection{Relative character varieties}
Let $\overline\Sigma$ be a compact Riemann surface of genus $g$, let
$k\ge1$, and put $\Sigma^\circ=\overline\Sigma\setminus\{p_1,\dots,p_k\}$.
Let $G$ be a connected complex reductive group, e.g.\ $\GL_n$, $\SL_n$ or a
torus. A $G$-local system on $\Sigma^\circ$ is given by handle matrices
$A_j,B_j\in G$ and puncture monodromies $M_i\in G$ satisfying
\begin{equation}\label{eq:surface-word}
 \mu(A_j,B_j,M_i):=\prod_{j=1}^g[A_j,B_j]\prod_{i=1}^kM_i=1,\qquad
 [A,B]=ABA^{-1}B^{-1},
\end{equation}
where the loop around $p_i$ keeps the surface on its left. Changing the
trivialization conjugates all matrices simultaneously. For an affine variety
$Z$ with a $G$-action, the affine \emph{geometric invariant theory} (GIT)
quotient is $Z/\!\!/G=\Spec\C[Z]^G$; the points of
$Z/\!\!/G$ are the closed $G$-orbits in $Z$ (Fact~\ref{fact:git}).

\smallskip\noindent\emph{Prescribed class closures.} For conjugacy classes
$\mathbf C=(C_1,\dots,C_k)$ of $G$ let $Y_i=\overline C_i$, with the reduced
structure, and put
\begin{equation}\label{paper:eq-relative}
 R_{\mathbf Y}=\Bigl\{(A_j,B_j,M_i)\in G^{2g}\times\textstyle\prod_iY_i:
 \mu=1\Bigr\},\qquad
 U_{\mathbf Y}=(R_{\mathbf Y})_{\mathrm{red}}/\!\!/G .
\end{equation}
On the generic locus considered below, the representations with
$M_i\in C_i$ for all $i$ form the \emph{exact-class locus}
\begin{equation}\label{eq:intro-MB}
 \cM_B(G,\mathbf C)=\Bigl\{(A_j,B_j,M_i)\in G^{2g}\times\textstyle\prod_iC_i:
 \mu=1\Bigr\}/G .
\end{equation}
We write $U^\nu_{\mathbf Y}$ for the \emph{normalization} of $U_{\mathbf Y}$:
$U^\nu_{\mathbf Y}=\coprod_U\Spec\widetilde{\C[U]}$, the union over the
irreducible components $U$ of $U_{\mathbf Y}$, where $\widetilde{\C[U]}$ is the
integral closure of $\C[U]$ in its fraction field. For $\GL_n$, such varieties
were introduced by Letellier, who studied their irreducibility, resolutions and
intersection cohomology \cite[\S3.1, Thms.~3.5, 3.11]{Letellier2015}.

\smallskip\noindent\emph{Prescribed conjugation invariants.} Let
$\chi_G:G\to G/\!\!/G$ be the adjoint quotient, i.e.\ the morphism recording all
class functions; for $\GL_n$ it records the characteristic polynomial. For
$\mathbf a=(a_1,\dots,a_k)\in(G/\!\!/G)^k$ put
$S_i=\chi_G^{-1}(a_i)_{\mathrm{red}}$. Replacing $Y_i$ by $S_i$ in
\eqref{paper:eq-relative} defines the \emph{spectral character variety}
\begin{equation}\label{eq:spectral-definition}
 \cM_B^{\mathrm{sp}}(G,\mathbf a)
 =\Bigl\{(A_j,B_j,M_i)\in G^{2g}\times\textstyle\prod_iS_i:\mu=1\Bigr\}_{\mathrm{red}}
 /\!\!/G .
\end{equation}
A point $a\in G/\!\!/G$ fixes the conjugacy class of the semisimple part $x_s$
of the Jordan decomposition $x=x_sx_u$, and allows every compatible unipotent
part $x_u$ (Lemma~\ref{lem:spectral-fiber-local}). For example, the $\GL_2$-fiber
of the characteristic polynomial $(t-\lambda)^2$, $\lambda\ne0$, is
\begin{equation}\label{eq:intro-gl2-fiber}
 \Bigl\{\lambda I+N:\ N=\begin{pmatrix}a&b\\c&-a\end{pmatrix},\ a^2+bc=0\Bigr\},
\end{equation}
a quadric cone with vertex $\lambda I$; the other points form the class of the
Jordan block. If every $a_i$ is the class of a regular semisimple element
(e.g.\ a diagonalizable matrix with distinct eigenvalues), then $S_i$ is a single
conjugacy class $C_i$ and $\cM_B^{\mathrm{sp}}(G,\mathbf a)=\cM_B(G,\mathbf C)$.

Prescribing conjugation invariants is therefore a special case of
prescribing conjugacy-class closures: each $S_i$ is the reduced closure
of its unique regular conjugacy class
(Corollary~\ref{cor:adjoint-fiber-closure}).
Here regular means that the centralizer has dimension $\rank G$;
a regular element need not be semisimple.
In \eqref{eq:intro-gl2-fiber}, the Jordan-block class is dense in the
fiber, while the scalar class $\{\lambda I\}$ is a proper class closure
in the same fiber.
This identification concerns the reduced fibers in
\eqref{eq:spectral-definition}; scheme-theoretic fibers can differ
(Example~\ref{ex:pgl2-spectral-thickening}).

\subsection{History and motivation}
Relative character varieties are the Betti spaces of nonabelian Hodge theory on
punctured curves. For $\GL_n$, Simpson established a correspondence between
stable filtered local systems of degree zero, stable filtered regular singular
connections of degree zero, and stable parabolic Higgs bundles of degree zero
\cite[Thm., p.~718, and table, p.~720]{Simpson1990}; for reductive groups with fixed
monodromy classes see \cite[Thm.~7.10]{BiquardGarciaPradaMundet2020} and the
parahoric formulation \cite{HKSZ2025}. Goldman's construction gives the
symplectic form on character varieties of closed surfaces
\cite[\S\S1.7--1.8]{Goldman1984}; quasi-Hamiltonian reduction gives the
relative form \cite{AMM1998,Boalch2007}. The cohomology of character varieties
is studied by point counts \cite[Thm.~1.2.3]{HauselLetellierRodriguezVillegas2011}
and is related to the Hitchin fibration by the $P=W$ conjecture
\cite[\S1.1]{deCataldoHauselMigliorini2012}.

A \emph{log canonical compactification} of a normal variety $U$ is a normal
projective variety $X\supset U$ such that
\[
 D=(X\setminus U)_{\mathrm{red}}\ \text{is a divisor},\qquad
 (X,D)\ \text{is log canonical},\qquad K_X+D\sim0 ,
\]
where $\sim$ denotes integral linear equivalence. Log canonicity is a condition
on discrepancies (Section~\ref{subsec:bg-pairs}). When $K_X+D\sim0$, the condition
says that a rational volume form $\Omega$ on $X$ with $\divi(\Omega)=-D$ acquires
at most simple poles along every exceptional divisor of every resolution of $X$.
Thus log canonical pairs with $K_X+D\sim0$ are the singular analog of smooth
compactifications whose boundary has simple normal crossings and is
anticanonical. A folklore conjecture predicts that smooth relative character
varieties admit such compactifications; see \cite[Conj.~1.1]{ArguzBousseau2026}
for a precise formulation and further references.

For $G=\SL_2$, several compactification constructions are known. Whang proved
that every relative $\SL_2$-character variety with prescribed boundary traces,
on a genus-$g$ surface with $k\ge1$ boundary components and $3g+k-3>0$, has a normal
projective compactification with $K_X+D\sim0$ \cite[Thm.~1.1]{Whang2020}. Komyo
compactified $\SL_2$-character varieties of punctured spheres
\cite{Komyo2015}. Farajzadeh-Tehrani and Frohman built compactifications of
punctured surfaces from ideal triangulations, with toric boundary divisors
\cite{FarajzadehTehraniFrohman2026}, and Ayilliath-Kutteri,
Farajzadeh-Tehrani and Frohman treated closed surfaces of genus $g>1$
\cite{AyilliathKutteriEtAl2025}. Arg\"uz and Bousseau proved that these
filtration compactifications are log canonical with $K_X+D\sim0$, using a
general criterion for compactifications defined by filtrations
\cite[Thms.~3.7--3.8]{ArguzBousseau2026}. For punctured surfaces with
$2g-2+k>0$, this gives log canonical log Calabi--Yau compactifications for all
boundary traces; generic regular semisimple boundary data additionally give
dlt compactifications \cite[Thm.~A]{ArguzBousseau2026}.
For the $\GL_n$- and $\SL_n$-character varieties of the closed torus,
Mauri, Mazzon and Stevenson proved the existence of dlt compactifications with
$K_X+D\sim_{\Q}0$ \cite[Thms.~D--E]{MauriMazzonStevenson2022}.
For $\SL_3$ on punctured surfaces with $2g-2+k>0$, Kim constructed compactifications by skein
algebras, for arbitrary prescribed boundary characteristic polynomials,
with $K_X+D\sim_{\Q}0$ \cite[Thm.~A]{Kim2026}.

The boundary also carries topological information. For a smooth variety $U$
with a smooth compactification whose boundary has simple normal crossings, the
\emph{dual boundary complex} has a vertex for each boundary component and a
$j$-cell for each component of an intersection of $j+1$ boundary components;
its homotopy type depends only on $U$ \cite[Thm.~1.1]{Payne2013}. The geometric
$P=W$ conjecture predicts that the dual boundary complex of a smooth character
variety of dimension $d$ is homotopy equivalent to $S^{d-1}$
\cite{KNPS2015}, \cite[Conj.~11.1]{Simpson2016}; see
\cite[Thm.~1.1]{Simpson2016} and \cite[Thm.~0.4]{Su2023} for known cases
of this sphere prediction. Kim proved an analogous sphere statement for the
boundary incidence complex of his $\SL_3$ compactifications
\cite[Thm.~B]{Kim2026}. In the closed-surface setting, Mauri, Mazzon and
Stevenson conjectured that the dual complex of a dlt compactification with
$K_X+D\sim_{\Q}0$ is homeomorphic to a sphere
\cite[Conj.~1.3.1]{MauriMazzonStevenson2022}. For
a proper log canonical pair with $K_X+D\sim_{\Q}0$, Kontsevich,
Koll\'ar and Xu asked whether the dual complex is a finite orthogonal
quotient of a sphere \cite[Question~4]{KollarXu2016}; see
\cite[p.~321]{Simpson2016} for the attribution. A compactification alone does not
identify these complexes, and we make no such identification here.

\subsection{Main results}
Write $G_\der=[G,G]$ and $G_\ab=G/G_\der$, and let $s_i$ be the semisimple part
of an element of $C_i$ (respectively, a semisimple element with
$\chi_G(s_i)=a_i$). The data are \emph{compatible} if
$\prod_i\bar s_i=1$ in $G_\ab$, where $\bar s_i$ is the image of $s_i$. The data are
\emph{intrinsically generic} if, for every proper parabolic subgroup
$Q\subsetneq G$ (Section~\ref{subsec:bg-groups}) and all
$x_i\in(G\cdot s_i)\cap Q$,
\begin{equation}\label{eq:intro-generic}
 x_1\cdots x_k\notin[Q,Q].
\end{equation}
This condition excludes reducible representations for all Jordan types
(Lemma~\ref{lem:spectral-genericity}). For compatible $\GL_n$ data, this says that no product
obtained by choosing $d$ eigenvalues of each $M_i$, counted with multiplicity,
equals $1$, for $0<d<n$
\cite[Def.~2.1.1]{HauselLetellierRodriguezVillegas2011},
\cite[Def.~3.1, Lem.~3.3]{Letellier2015}; for reductive groups see
\cite{GianniniEtAl2026}. A projective variety $X$ is of \emph{Fano type} if
$(X,\Delta)$ is klt and $-(K_X+\Delta)$ is ample for some effective $\Q$-divisor
$\Delta$ (Section~\ref{subsec:bg-pairs}).

\begin{mainresult}[Section~\ref{r2:sec-full-goal}]\label{paper:thm-normalized}
Assume that $\mathbf C$ is compatible and intrinsically generic. Every nonempty
irreducible component $U$ of $U_{\mathbf Y}^{\nu}$ has a normal projective
log canonical compactification $(X,D)$ with $K_X+D\sim0$. The variety $X$ may be
chosen of Fano type; in particular $X$ has rational singularities and is
Cohen--Macaulay.
\end{mainresult}

Theorem~\ref{r2:thm-full-normalized-logcy} proves a slightly stronger statement:
$X$ may be taken to be a Poisson compactification whose boundary supports an
effective ample Cartier divisor. No condition on the Jordan types of the
prescribed classes, such as regular semisimplicity, is imposed. The normalization cannot be omitted in general: an
open subset of a normal variety is normal. The next result says exactly when
normalization is needed.

\begin{mainresult}[Section~\ref{r2:norm-section}]\label{paper:thm-normality}
Assume that $\mathbf C$ is compatible and intrinsically generic. If $g\ge1$, then
for every $r=(A_j,B_j,M_i)\in R_{\mathbf Y}$ with image $p\in U_{\mathbf Y}$,
\[
 U_{\mathbf Y}\ \text{is normal at}\ p
 \iff Y_i\ \text{is normal at}\ M_i\ \text{for every}\ i .
\]
Every tuple of puncture points occurs, so $U_{\mathbf Y}\ne\varnothing$, and
$U_{\mathbf Y}$ is normal if and only if every $Y_i$ is normal. In every genus,
the normal locus and the Cohen--Macaulay locus of $U_{\mathbf Y}$ coincide.
\end{mainresult}

In genus zero a finite stabilizer can hide a non-normal puncture: at a point
$p$, the completed local ring of $U_{\mathbf Y}$ is the ring of invariants of a
finite group $F$ acting on a power series ring over the completed
puncture-transversal rings, and normality at $p$ is decided by the subgroup of $F$ acting trivially on
the smooth deformation directions (Theorem~\ref{r2:norm-kernel-criterion}).

For the spectral subfamily we obtain more precise information.

\begin{mainresult}[Section~\ref{subsec:proofA}]
\label{thm:main}\label{thm:spectral-main}
Let $\mathbf a\in(G/\!\!/G)^k$ be compatible and intrinsically generic, with
$\cM_B^{\mathrm{sp}}(G,\mathbf a)\ne\varnothing$. Then every connected component
$U$ of $\cM_B^{\mathrm{sp}}(G,\mathbf a)$ has a normal Cohen--Macaulay projective
log canonical compactification $(X,D)$ with $K_X+D\sim0$. The variety $U$ is
$\Q$-factorial and Gorenstein with canonical singularities, and $(X,D)$ has a
projective $\Q$-factorial dlt modification which is an isomorphism over $U$. If
$G_\der$ is simply connected, then $\cM_B^{\mathrm{sp}}(G,\mathbf a)$ is connected
and equals the scheme-theoretic fiber of the peripheral morphism
$G^{2g+k-1}/\!\!/G\to(G/\!\!/G)^k$, $\varphi\mapsto(\chi_G(M_i))_i$.
\end{mainresult}

Here a \emph{dlt modification} is a projective birational morphism $X'\to X$ such
that the pair consisting of $X'$ and the reduced preimage of $D$ together with the
exceptional divisors is dlt with the same log canonical class
(Fact~\ref{fact:pairs}(i)), and $\Q$-factorial means that every Weil divisor has
a Cartier multiple.

\begin{corollary}[Regular semisimple monodromy;
Section~\ref{subsec:proof-corollaries}]\label{cor:regular-semisimple}
Let $\mathbf C$ be a compatible and intrinsically generic tuple of regular
semisimple classes with $\cM_B(G,\mathbf C)\ne\varnothing$. Then all conclusions
of Theorem~\ref{thm:spectral-main} hold for $\cM_B(G,\mathbf C)$. Moreover,
locally in the \'etale topology, $\cM_B(G,\mathbf C)$ is the quotient of a smooth
symplectic variety by a finite group preserving the symplectic form.
\end{corollary}

If the simply connected cover of $G_\der$ is a product of special linear groups,
the compactification conclusion of Theorem~\ref{thm:spectral-main} holds for all
compatible $\mathbf a$, without genericity (Corollary~\ref{cor:spectral-type-a}),
by the normality theorem for multiplicative quiver varieties
\cite[Cor.~5.5]{KaplanSchedler2023} and the identification
\cite[Thm.~3.6, Rem.~3.7]{SchedlerTirelli2022}.

The spectral compactifications come from an explicit degeneration. Let $H$ be
semisimple and simply connected of rank $r$, put $m=2g+k-1$ and
$A_0=\C[H^m]^H$, and for semisimple $s_i\in H$ let $A_\bs$ be the quotient of
$A_0$ prescribing $\chi_H(M_i)=\chi_H(s_i)$. Filtrations, Rees algebras and
toric-face rings are recalled in Sections~\ref{subsec:bg-filtrations}
and~\ref{subsec:bg-toric}, and ideal triangulations in
Section~\ref{subsec:ribbon}.

\begin{mainresult}[Theorems~\ref{thm:special-fiber} and~\ref{thm:AB},
Proposition~\ref{prop:spectral-flatness}]\label{paper:thm-degeneration}
Assume $2g+k\ge3$. Every ideal triangulation of $\Sigma^\circ$ determines an
increasing filtration of $A_0$ by finite-dimensional subspaces, with finitely
generated Rees algebra, such that for all $s_1,\dots,s_k$ the induced filtration
$F$ of $A_\bs$ satisfies:
\begin{enumerate}[label=\textup{(\roman*)}]
\item $\gr_FA_\bs$ is the quotient of a normal semigroup ring by a regular
sequence of $kr$ monomials independent of the $s_i$; it is a reduced Gorenstein
toric-face ring whose canonical module is generated in degree $0$;
\item the fan of $\gr_FA_\bs$ is homeomorphic, cone by cone, to a complete fan in
$\R^d$, $d=\dim\Spec A_\bs$; hence the positive slice of the fan is a sphere
$S^{d-1}$ (empty if $d=0$);
\item the peripheral morphism $\Spec A_0\to(H/\!\!/H)^k$ is faithfully flat, and
every fiber $\Spec A_\bs$ is reduced, connected and Gorenstein of dimension
$d=(2g+k-2)\dim H-kr$;
\item if $A_\bs$ is normal, then $X=\Proj\bigoplus_nF_nA_\bs$ with
$D=\Proj\gr_FA_\bs$ is a log canonical compactification of $\Spec A_\bs$ with
$K_X+D\sim0$.
\end{enumerate}
\end{mainresult}

In (iii), \emph{faithfully flat} means that $A_0$ is a flat module over the
coordinate ring of $(H/\!\!/H)^k$ and every fiber is nonempty; in fact $A_0$ is a
free module. The sphere in (ii) has the dimension predicted by the geometric
$P=W$ conjecture, but we do not identify it with any dual complex.

\subsection{The model calculation}\label{subsec:intro-example}
For $G=\SL_2$ and $(g,k)=(1,1)$, put $x=\tr A$, $y=\tr B$, $z=\tr AB$ and
$\kappa=\tr[A,B]$. Fricke's trace identities give
\cite[\S2, pp.~S.7--S.9]{Vogt1889}, \cite[\S2.1, (2.1.4)]{Goldman2003}
\[
 A_0=\C[x,y,z],\qquad \kappa=x^2+y^2+z^2-xyz-2,\qquad
 A_\tau=A_0/(\kappa-\tau),
\]
where prescribing the conjugation invariant of the puncture monodromy means
prescribing $\tau=\tr M$. Let $F$ be the total-degree filtration and $\bar F$ the quotient
filtration. Then
\[
 \gr_FA_0=\C[\N^3],\qquad \inn_F(\kappa-\tau)=-xyz,\qquad
 \gr_{\bar F}A_\tau=\C[x,y,z]/(xyz).
\]
The exponent $(1,1,1)$ of $xyz$ translates $\N^3$ onto the interior lattice points
of the octant, i.e.\ $xyz$ is the canonical monomial (Fact~\ref{fact:semigroup}):
\[
 xyz\,\frac{dx}{x}\wedge\frac{dy}{y}\wedge\frac{dz}{z}=dx\wedge dy\wedge dz .
\]
By adjunction, the canonical module of $\C[x,y,z]/(xyz)$ is generated in degree
$3-3=0$. The Rees construction gives the cubic surface and its boundary,
\[
 X_\tau=\{t(x^2+y^2+z^2)-xyz=(\tau+2)t^3\}\subset\PP^3_{[x:y:z:t]},\qquad
 D_\tau=X_\tau\cap\{t=0\}=\{t=xyz=0\},
\]
a triangle of lines with $K_{X_\tau}+D_\tau\sim0$. For $\tau=\pm2$ the
surface acquires canonical nodes away from $D_\tau$: one node for $\tau=-2$,
where the puncture fiber $-I+\mathcal N(\mathfrak{sl}_2)$ is a quadric cone and
the data are generic (Theorems~\ref{paper:thm-normality}
and~\ref{thm:spectral-main}); four nodes for $\tau=2$, where the data are not
generic (Corollary~\ref{cor:spectral-type-a}). Examples~\ref{ex:rees-cubic},
\ref{ex:cubic-special-fiber}, \ref{ex:cubic-flatness}
and~\ref{ex:exceptional-cubics} verify these statements directly.

\subsection{Ideas of the proofs}\label{subsec:strategy}
\emph{The spectral route.} An ideal triangulation replaces $\C[x,y,z]$ by a
normal semigroup ring $\C[Q_\Gamma]$, and each puncture equation $f_{i,a}$ has a
monomial initial form. The local corner identity (Theorem~\ref{thm:corner}) and
the fact that the boundary walks of the triangulation visit each corner once
give
\begin{equation}\label{eq:intro-sum}
 \sum_{i,a}v_\Gamma(f_{i,a})=\Omega_\Gamma,
\end{equation}
where $x^{\Omega_\Gamma}$ generates the canonical module of $\C[Q_\Gamma]$. This
identity forces the initial equations to form a radical regular sequence and
produces a toric-face special fiber with canonical module generated in degree
$0$. A canonical partition of the facets gives a basis over the peripheral
monomials, which lifts to a free basis over the peripheral characters. The Rees
algebra then gives $K_X+D\sim0$, and inversion of adjunction gives log
canonicity. The cancellation of central characters in \eqref{eq:intro-sum}
preserves the integral trivialization under the finite central cover of $G$.

\emph{The Poisson route.} For a smaller class closure, restricting a spectral
compactification gives neither adjunction nor log canonicity. Instead we use a
filtration of degree $0$ for the Poisson bracket. The resulting projective model
has Poisson boundary, and a logarithmic volume form which does not vanish forces
every boundary pole to be simple and proves log canonicity
(Theorem~\ref{thm:poisson-nonvanishing}). The Goldman class lies in $F^2$ of
Deligne's Hodge filtration (Proposition~\ref{prop:goldman-f2-all-types}), so a
compact cycle with nonzero Goldman-volume period supplies this form. The cycle is
a Lagrangian fiber of a projective Hitchin map on a Dolbeault moduli space with
fixed exact residues. On this stratum the symplectic class is a linear
combination of Chern classes; Hodge theory makes projective fibers isotropic,
Poisson commutativity gives the opposite dimension bound, and hyperk\"ahler
rotation turns the Lagrangian fiber into the required period. Only the
correspondence for the exact stratum is used.

\emph{How the two routes compare.} Spectral varieties are class-closure
varieties (Corollary~\ref{cor:adjoint-fiber-closure}), so the Poisson route also
compactifies them. The spectral route gives more: an explicit toric-face
degeneration with a free peripheral map (Theorem~\ref{paper:thm-degeneration}),
an interior-preserving $\Q$-factorial dlt modification for generic spectra
(Theorem~\ref{thm:spectral-main}), and compactifications for all compatible
spectra in type~$A$. The Poisson route handles arbitrary generic class closures
after normalization and gives a model of Fano type. The two routes are
independent: the Poisson route uses a matrix-coefficient filtration
(Proposition~\ref{prop:matrix-poisson-compactification}) and never the string
filtration of Sections~\ref{sec:graph}--\ref{sec:main-proofs}. Both use a finite central cover with simply connected derived group.
The spectral route descends a compactified pair and its invariant logarithmic
volume form; the Poisson route transfers the Goldman volume through the finite
quotient and constructs a Poisson compactification downstairs.

\emph{Normality.} A computation in the group directions gives the formal model
\[
 \widehat{\mathcal O}_{U_{\mathbf Y},p}\simeq(A[[W_r^*]])^F,
\]
where $A$ is the completed product of the puncture-transversal rings, $W_r$ is the smooth deformation
representation and $F$ is the finite stabilizer. Polynomial functions on a
faithful representation of a finite group detect every nonzero normalization
defect. In genus at least two the representation $W_r$ is faithful; genus one
needs a separate central-isogeny argument. This gives
Theorem~\ref{paper:thm-normality} and explains the genus-zero criterion.

\subsection{Organization and conventions}\label{subsec:guide}
After this introduction, the main text has four parts, with continuous
section numbering.
Part~\ref{part:background} (Sections~\ref{sec:background}--\ref{new:sec-normalization})
collects background and develops the interior geometry.
Section~\ref{sec:background} can be consulted as needed.
Section~\ref{sec:interior} proves the common genericity and surface-relation
arguments and the stronger spectral interior statements;
Sections~\ref{new:sec-base-change}--\ref{new:sec-normalization} treat arbitrary
class closures and their normalizations.

Part~\ref{part:spectral} (Sections~\ref{sec:graph}--\ref{sec:main-proofs})
develops the spectral construction and proves
Theorems~\ref{thm:spectral-main} and~\ref{paper:thm-degeneration}.
Readers following this route need only Section~\ref{sec:interior} from
the interior analysis.
Part~\ref{part:poisson} (Sections~\ref{sec:poisson-bridge}--\ref{r2:sec-full-goal})
uses the normalized class-closure interiors to develop the Poisson
construction and prove Theorem~\ref{paper:thm-normalized}.

\begin{figure}[htbp]
\centering
\begin{tikzpicture}[>=stealth,font=\small,
 box/.style={draw,rounded corners,align=center,inner sep=4pt,minimum width=3.3cm}]
 \node[box] (bg) at (4.2,2.7) {Part~\ref{part:background}: \S\ref{sec:background} Background};
 \node[box] (int) at (4.2,1.35) {Part~\ref{part:background}: \S\ref{sec:interior}\\Genericity, local geometry};
 \node[box] (sp) at (0,0) {Part~\ref{part:spectral}: \S\S\ref{sec:graph}--\ref{sec:main-proofs}\\Spectral construction};
 \node[box] (cl) at (4.2,0) {Part~\ref{part:background}: \S\S\ref{new:sec-base-change}--\ref{new:sec-normalization}\\Class closures, normalization};
 \node[box] (po) at (8.4,0) {Part~\ref{part:poisson}: \S\S\ref{sec:poisson-bridge}--\ref{r2:sec-full-goal}\\Poisson construction};
 \node[box] (norm) at (4.2,-1.4) {Part~\ref{part:normality}: \S\ref{r2:norm-section}\\Normality criterion};
 \node[box] (thCD) at (0,-2.8) {Theorems~\ref{thm:spectral-main}, \ref{paper:thm-degeneration}};
 \node[box] (thB) at (4.2,-2.8) {Theorem~\ref{paper:thm-normality}};
 \node[box] (thA) at (8.4,-2.8) {Theorem~\ref{paper:thm-normalized}};
 \draw[->] (bg)--(int);
 \draw[->] (int)--(sp); \draw[->] (int)--(cl);
 \draw[->] (cl)--(po); \draw[->] (cl)--(norm);
 \draw[->] (sp)--(thCD); \draw[->] (po)--(thA);
 \draw[->] (norm)--(thB);
\end{tikzpicture}
\caption{The four parts and their proof dependencies. The branching inside
Part~\ref{part:background} distinguishes the spectral interior analysis
from the geometry of arbitrary class closures. Theorems~\ref{paper:thm-normalized}
and~\ref{paper:thm-normality} together determine when the original reduced
interior can be retained.}\label{fig:roadmap}
\end{figure}
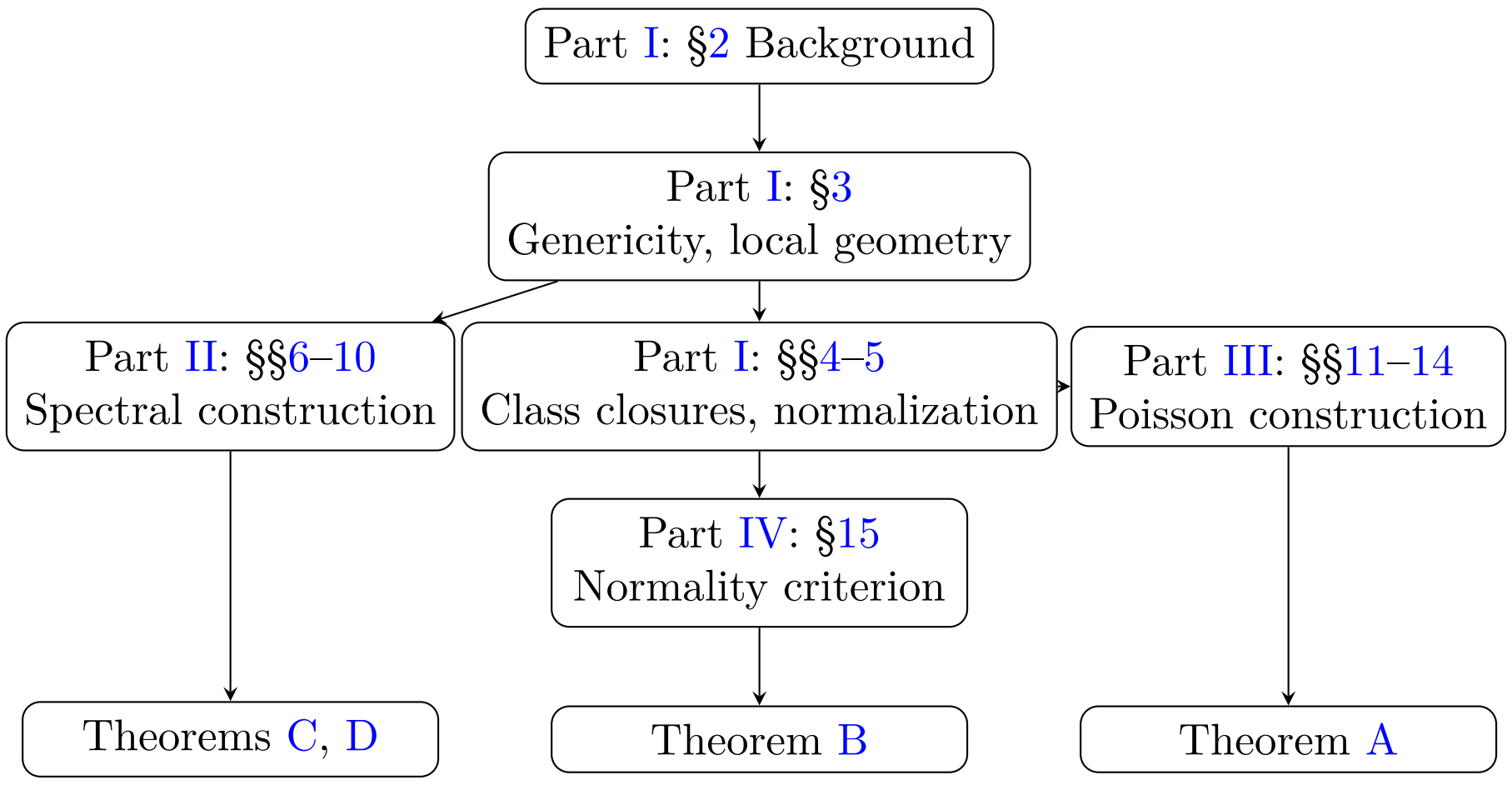

Part~\ref{part:normality} (Sections~\ref{r2:norm-section}--\ref{sec:examples})
proves the normality classification, Theorem~\ref{paper:thm-normality},
and gives examples. Neither compactification construction uses the
classification in Section~\ref{r2:norm-section}; that classification and
Theorem~\ref{paper:thm-normalized} together decide when the given reduced
interior can be retained. Figure~\ref{fig:roadmap} records these dependencies.

The three appendices lie outside the four parts and separate the
supporting arguments by their role.
Appendix~\ref{app:spectral-inputs} contains the algebraic inputs to the
spectral construction: the local string valuation, scalarization and
auxiliary geometric facts.
Appendix~\ref{app:poisson-inputs} contains the inputs to the Poisson
construction: extension of forms, Poisson normalization and Hitchin
projectivity. The final subsection of this appendix records an optional cohomological
nonvanishing criterion, which the main period proof does not use.
Appendix~\ref{new:sec-counterexample} gives the realization and obstruction
proofs summarized in Section~\ref{subsec:closure-obstructions}, explaining
the limits of the compactification results.

All schemes are of finite type over $\C$ unless completed local rings are
displayed, and varieties are reduced. We write $\N=\{0,1,2,\dots\}$,
$W^\vee=\Hom(W,\C)$ and $Z_{\mathrm{red}}$ for the reduction of a scheme $Z$. A
dual complex belongs to a pair $(X,D)$; a dual boundary complex belongs to a
smooth open variety. The letters $F$ and $\bar F$ denote a filtration and the
induced quotient filtration.

\part{Background and interior geometry}\label{part:background}

\section{Background}\label{sec:background}

We fix the conventions and record the background used in the interior
analysis and the two compactification constructions. The supporting local proofs that supplement the cited
theorems are collected in Appendices~\ref{app:spectral-inputs}
and~\ref{app:poisson-inputs}.

\subsection{Reductive groups and their representations}\label{subsec:bg-groups}
Let $G$ be a connected complex reductive group and fix $T\subset B=TU$.
Here $B$ is a Borel subgroup, a maximal connected solvable subgroup;
$T$ is a maximal torus, and $U$ is the unipotent radical of $B$.
Over $\C$, a torus is a group isomorphic to $(\C^*)^\ell$.
Characters and
cocharacters pair by $\lambda(\alpha^\vee(z))=z^{\langle\lambda,\alpha^\vee\rangle}$.
Write $\Phi^+$ for the positive roots, $\alpha_i,\alpha_i^\vee$ for the simple
roots and coroots, and $W=N_G(T)/T$, with simple reflections $s_i$ and longest
element $w_0$ \cite[\S7.4, Cor.~8.1.2]{SpringerLAG1998}. A reduced word
$w_0=s_{i_1}\cdots s_{i_N}$ has $N=|\Phi^+|$ letters. Any two reduced words
of an element of $W$ are connected by braid moves \cite{Matsumoto1964}; see
\cite[Prop.~8.3.3]{SpringerLAG1998}. We use the Bruhat decomposition
$G=\bigsqcup_{w\in W}B\dot wB$. Parabolic subgroups contain a conjugate of $B$;
in particular
$P(\lambda)=\{g:\lim_{t\to0}\lambda(t)g\lambda(t)^{-1}\text{ exists}\}$
is parabolic for a cocharacter $\lambda$, and equals $G$ precisely when
$\lambda$ is central \cite[Prop.~8.4.5]{SpringerLAG1998}.
A semisimple $s\in T$ is regular when $\alpha(s)\ne1$ for every root,
equivalently $Z_G(s)^\circ=T$.

\begin{example}[$G=\SL_n$]\label{ex:sln}
Take $T$ diagonal, $B$ upper triangular and $U$ upper unitriangular. Writing
$\epsilon_i(\diag(t_1,\dots,t_n))=t_i$ (additively, $\epsilon_i$), the roots
are $\epsilon_i-\epsilon_j$ ($i\ne j$), the simple roots are
$\alpha_i=\epsilon_i-\epsilon_{i+1}$, $W=S_n$ permutes the diagonal entries, and
$w_0$ is the permutation $i\mapsto n+1-i$. The parabolic subgroups are the
stabilizers of partial flags $0\subsetneq V_1\subsetneq\dots\subsetneq\C^n$.
A semisimple element is regular if and only if its eigenvalues are distinct.
\end{example}

\begin{fact}[Central isogeny {\cite[Cor.~8.1.6, \S10.1, Ex.~10.1.4(1)]{SpringerLAG1998}}]
\label{fact:isogeny}
Let $G_\der^{\mathrm{sc}}$ be the simply connected cover of the derived group
$G_\der=[G,G]$. Multiplication defines a surjective homomorphism
\[
 p:\widetilde G=Z(G)^\circ\times G_\der^{\mathrm{sc}}\longrightarrow G
\]
whose kernel $F_p$ is finite and central.
\end{fact}

For example, for $G=\GL_n$ the map is $\C^*\times\SL_n\to\GL_n$,
$(z,h)\mapsto zh$, with kernel $\{(\zeta^{-1},\zeta I):\zeta^n=1\}$.

For the rest of this subsection let $H$ be simply connected semisimple of
rank $r$, with $N=|\Phi^+|$. The fundamental weights of $H$ satisfy
$\langle\varpi_a,\alpha_b^\vee\rangle=\delta_{ab}$; put
$\rho=\sum_a\varpi_a=\frac12\sum_{\alpha>0}\alpha$ and
$X^+(T)=\bigoplus_a\N\varpi_a$. The dominance order is
$\nu\le\lambda$ when $\lambda-\nu$ is a nonnegative integral sum of simple
roots. The function $\ell(\lambda)=\langle\lambda,2\rho^\vee\rangle$, with
$2\rho^\vee$ the sum of the positive coroots, satisfies $\ell(\alpha_i)=2$;
thus $\nu<\lambda$ implies $\ell(\nu)<\ell(\lambda)$.

\begin{fact}[Highest weights {\cite[\S14.3, \S\S31.3--31.6]{HumphreysLAG1975}}]\label{fact:hw}
Every finite-dimensional rational representation of $H$ is completely
reducible. The irreducible representations are the modules $V(\lambda)$,
$\lambda\in X^+(T)$, characterized as follows: the $U$-invariant vectors of
$V(\lambda)$ form a line $\C v_\lambda$ (the \emph{highest weight line}), $T$
acts on $v_\lambda$ by $\lambda$, and every weight $\nu$ of $V(\lambda)$
satisfies $\nu\le\lambda$. Dually, $V(\lambda)$ has a \emph{lowest weight
line}, of weight $w_0\lambda$. Moreover $V(\lambda)^\vee\simeq V(\lambda^*)$
with $\lambda^*=-w_0\lambda$. We write $\chi_\lambda(g)=\tr(g\mid V(\lambda))$
for the character of $V(\lambda)$.
\end{fact}

\begin{example}\label{ex:sln-reps}
For $H=\SL_2$ there is one fundamental weight $\varpi_1=\rho$, and
$V(n\varpi_1)=\operatorname{Sym}^n\C^2$; every representation is self-dual. For
$H=\SL_n$, $\varpi_a=\epsilon_1+\dots+\epsilon_a$ and
$V(\varpi_a)=\bigwedge^a\C^n$, with highest weight vector
$e_1\wedge\dots\wedge e_a$; duality is $\varpi_a^*=\varpi_{n-a}$, reflecting
$(\bigwedge^a\C^n)^\vee\simeq\bigwedge^{n-a}\C^n$. For $\SL_3$, $V(\varpi_1)=\C^3$
and $V(\varpi_2)=(\C^3)^\vee$.
\end{example}

\begin{fact}[Cartan component]\label{fact:cartan}
For $\lambda,\mu\in X^+(T)$, the module $V(\lambda)\otimes V(\mu)$ contains
$V(\lambda+\mu)$ with multiplicity one, generated by $v_\lambda\otimes v_\mu$;
every other irreducible constituent has highest weight $<\lambda+\mu$. We call
$V(\lambda+\mu)$ the \emph{Cartan component}, and the $H$-equivariant projection
$V(\lambda)\otimes V(\mu)\to V(\lambda+\mu)$ with
$v_\lambda\otimes v_\mu\mapsto v_{\lambda+\mu}$ the \emph{Cartan product}.
\end{fact}

Fact~\ref{fact:cartan} follows from Fact~\ref{fact:hw}, because the weight $\lambda+\mu$ occurs
in $V(\lambda)\otimes V(\mu)$ only on the line $\C v_\lambda\otimes v_\mu$. For
$\SL_2$, the Clebsch--Gordan rule
$\operatorname{Sym}^a\otimes\operatorname{Sym}^b=\operatorname{Sym}^{a+b}
\oplus\operatorname{Sym}^{a+b-2}\oplus\cdots$ exhibits the Cartan component
$\operatorname{Sym}^{a+b}$, and the Cartan product is multiplication of
polynomials in two variables.

\begin{fact}[Peter--Weyl]\label{fact:PW}
As a representation of $H\times H$ acting by $((h_1,h_2)f)(g)=f(h_1^{-1}gh_2)$,
\[
 \C[H]=\bigoplus_{\lambda\in X^+(T)}V(\lambda)^\vee\otimes V(\lambda),
\]
where $\ell\otimes w$ corresponds to the \emph{matrix coefficient}
$c_{\ell,w}(g)=\ell(gw)$. Products of matrix coefficients are matrix
coefficients of tensor products: $c_{\ell,w}\,c_{\ell',w'}=c_{\ell\otimes\ell',
w\otimes w'}$ \textup{(}see \cite[Thm.~4.2.7]{GoodmanWallach2009}\textup{)}.
\end{fact}

For instance, the matrix entries $g_{ij}$ of $g\in\SL_n$ are the matrix
coefficients of $V(\varpi_1)=\C^n$.

\begin{fact}[The basic affine space {\cite[\S\S2--3]{VinbergPopov1972}}]\label{fact:basic-affine}
Let $\C[H/U]=\C[H]^U$ be the ring of functions on $H$ invariant under right
translation by $U$, and $\mathcal Z_H=\Spec\C[H/U]$.
\begin{enumerate}[label=\textup{(\roman*)}]
\item By Fact~\ref{fact:PW},
$\C[H/U]=\bigoplus_{\lambda}V(\lambda)^\vee\otimes\C v_\lambda\simeq
\bigoplus_\lambda V(\lambda^*)$, graded by the action of $T$ by right
translation. The product of the summands indexed by $\lambda$ and $\mu$ lies in
the summand indexed by $\lambda+\mu$, and is the Cartan product; in particular
the product of two nonzero homogeneous elements is nonzero.
\item The map $H/U\to\mathcal Z_H$ is an open embedding whose complement has
codimension at least $2$.
\item The ring $\C[H/U]$ is a factorial normal domain whose units are the
constants.
\end{enumerate}
\end{fact}

\begin{example}\label{ex:basic-affine}
For $\SL_2$, the first column of $g$ identifies $\SL_2/U$ with $\C^2\setminus0$,
because $\SL_2$ acts transitively on nonzero vectors and the stabilizer of
the first standard basis vector is $U$. Regular functions extend across
the omitted origin, so $\mathcal Z_{\SL_2}=\C^2$.

For $\SL_3$, let $v$ be the first column of $g$ and $\phi$ the last row of
$g^{-1}$. Both are unchanged by right multiplication of $g$ by $U$, and
$\phi(v)=0$ follows from $g^{-1}g=I$. Conversely, given nonzero $v,\phi$
with $\phi(v)=0$, choose $w$ so that $(v,w)$ is a basis of $\ker\phi$, and
choose $z$ with $\phi(z)=1$. Rescale $w$ so that $\det(v,w,z)=1$.
Then $g=(v,w,z)\in\SL_3$ has the prescribed first column and last inverse
row. Two such matrices differ by right multiplication by $U$.
Thus $\SL_3/U$ is the locus $v\ne0$, $\phi\ne0$ in
\[
 Z=\{(v,\phi)\in\C^3\times(\C^3)^\vee:\phi(v)=0\}.
\]
The hypersurface $Z$ is normal: its only singular point is the origin,
so Serre's criterion applies (Section~\ref{subsec:bg-algebra}).
The omitted loci $v=0$ and $\phi=0$ each have dimension $3$, hence
codimension $2$ in $Z$. Regular functions therefore extend across these
loci, giving $\C[\SL_3/U]=\C[Z]$ and $\mathcal Z_{\SL_3}=Z$.
\end{example}

\begin{fact}[Steinberg {\cite[Thms.~6.1, 6.11, 8.1, Rem.~6.15(c)]{Steinberg1965}}]\label{fact:steinberg}\mbox{}
\begin{enumerate}[label=\textup{(\roman*)}]
\item The ring $\C[H]^H$ of class functions is the polynomial ring
$\C[\chi_{\varpi_1},\dots,\chi_{\varpi_r}]$.
\item If $s\in H$ is regular semisimple, the subscheme of $H$ defined by the
equations $\chi_{\varpi_a}(g)=\chi_{\varpi_a}(s)$, $1\le a\le r$, is the
conjugacy class $H\cdot s$ with its reduced structure; in particular this
subscheme is smooth of codimension $r$.
\item For every $a\in H/\!/H$, the scheme fiber $\chi_H^{-1}(a)$ is
an irreducible reduced normal complete intersection of codimension $r$.
The fiber has a unique open regular conjugacy class, whose complement has
codimension at least two, and is a finite union of conjugacy classes.
Here \emph{regular} means that the centralizer has dimension $r$;
the element need not be semisimple.
\end{enumerate}
\end{fact}

For $H=\SL_n$, the character $\chi_{\varpi_a}=\tr\bigwedge^a$ is, up to sign,
the coefficient of $t^{n-a}$ in the characteristic polynomial $\det(t-g)$. So
(i) says that class functions are polynomials in the coefficients of the
characteristic polynomial, and (ii) says that a matrix with distinct eigenvalues
is determined up to conjugacy by its characteristic polynomial, scheme
theoretically. 

\begin{fact}[Reductive GIT {\cite[Ch.~I, \S\S1--2, pp.~26--29; Ch.~II, \S1, pp.~50--55]{MumfordFogarty1982}}]
\label{fact:git}
Let a reductive group $L$ act on an affine variety $Y$.
\begin{enumerate}[label=\textup{(\roman*)}]
\item Taking $L$-invariants is exact on rational $L$-modules, and there is an
$\C[Y]^L$-linear projection $\C[Y]\to\C[Y]^L$ (the \emph{Reynolds operator}).
Consequently $\C[Y]^L$ is finitely generated, normal if $\C[Y]$ is normal, and,
if $f_1,\dots,f_n\in\C[Y]^L$, then
$(\C[Y]/(f_1,\dots,f_n))^L=\C[Y]^L/(f_1,\dots,f_n)$.
\item (Hilbert--Mumford.) A point $y\in Y$ has closed orbit and finite
stabilizer if and only if $\lim_{t\to0}\lambda(t)y$ does not exist for every
nontrivial cocharacter $\lambda$ of $L$.
\item (Luna slice theorem \cite[III, \S1, pp.~97--98]{Luna1973}.)
If $L\cdot y$ is closed and $F=L_y$ is finite, there is an affine
$F$-stable locally closed slice $S\ni y$ such that
$L\times^F S\to Y$ and $S/F\to Y/\!/L$ are \'etale. The first map has
a saturated affine open image $Y'$: $Y'$ is the preimage of an open set
in $Y/\!/L$. Over this image the quotient diagram is Cartesian,
$L\times^F S\simeq Y'\times_{Y'/\!/L}(S/F)$.
Here $L\times^F S$ is the quotient of $L\times S$ by
$f\cdot(l,s)=(lf^{-1},fs)$.
Consequently the analytic germ of $Y$ at $y$ is the product of the germ
of $L$ at $1$ and the germ of $S$ at $y$, and the quotient germ is the
germ of $S/F$. If $Y$ is smooth, $S$ is smooth and transverse to the orbit.
\end{enumerate}
\end{fact}

An \'etale morphism is locally an isomorphism in the complex analytic
topology. Thus the slice statement allows us to study singularities of the
quotient through $S/F$, including when $Y$ and $S$ are singular.

\subsection{Canonical modules}\label{subsec:bg-algebra}
All rings are Noetherian. A regular sequence consists of successive
nonzerodivisors and generates a proper ideal. We use the Cohen--Macaulay
condition and Serre's criterion: a reduced ring is normal if and only if it
is $R_1$ and $S_2$. On a normal variety, reflexive sheaves and rational
sections of them are determined in codimension one. In particular, a
rational top form $\beta$ determines a canonical divisor
$K_X=\divi(\beta)$ and $\omega_X\simeq\mathcal O_X(K_X)$.
For a pure $d$-dimensional Cohen--Macaulay algebra
$R=\C[x_1,\dots,x_n]/I$, our canonical module is the dualizing module
relative to $\C$,
\[
 \omega_R=\Ext^{n-d}_{\C[x_1,\dots,x_n]}
 (R,\C[x_1,\dots,x_n]\,dx_1\wedge\cdots\wedge dx_n),
\]
independent of the presentation \cite[\S3.3]{BrunsHerzog1998}.
For normal $R$ it agrees with the sections of $\omega_{\Spec R}$.
The ring is Gorenstein when $\omega_R$ is locally free of rank one.
A module $M$ is maximal Cohen--Macaulay when
$\operatorname{depth}_{R_{\mathfrak p}}M_{\mathfrak p}=\dim R_{\mathfrak p}$
at every prime in its support.

\begin{fact}[{\cite[Cor.~3.6.14, Thm.~3.3.10, Thm.~4.4.5]{BrunsHerzog1998}}]
\label{fact:canonical}
Let $R$ be a Cohen--Macaulay ring graded by a lattice $L$, with $R_0=\C$,
finite-dimensional graded pieces, and a linear form on $L$ which is positive on
the degrees of all nonconstant homogeneous elements. Then $R$ has a graded
canonical module $\omega_R$, and:
\begin{enumerate}[label=\textup{(\roman*)}]
\item (Adjunction.) If $f\in R$ is a homogeneous nonzerodivisor of degree $p$,
then $R/(f)$ is Cohen--Macaulay and
$\omega_{R/(f)}\simeq(\omega_R/f\omega_R)(p)$, where $M(p)_q=M_{p+q}$.
\item (Duality.) If $M$ is a maximal Cohen--Macaulay $R$-module, then
$\Hom_R(M,\omega_R)$ is maximal Cohen--Macaulay and $\Ext^i_R(M,\omega_R)=0$ for
$i>0$.
\item (Reciprocity \cite[Thm.~6.1(ii)]{Stanley1978}.) The Hilbert series
$H_M(\boldsymbol z)=\sum_q\dim M_q\,\boldsymbol z^q$ satisfy the identity of
rational functions $H_{\omega_R}(\boldsymbol z)=(-1)^{\dim R}H_R(\boldsymbol
z^{-1})$.
\end{enumerate}
\end{fact}

Bruns--Herzog treat $\Z$-gradings; the lattice-graded versions are in
\cite[Props.~4.2, 4.4, Thm.~4.5]{BarileBruns2025}. Thus quotienting a
Gorenstein ring with canonical generator of degree $\vartheta$ by a
homogeneous regular sequence of degrees $p_i$ gives canonical degree
$\vartheta-\sum_i p_i$. For a degree-$d$ hypersurface in $\PP^n$, this is
the familiar $K_Z=\mathcal O_Z(d-n-1)$. The shift convention can also be
checked on $\C[x]$: $H_{\omega_{\C[x]}}=z/(1-z)$, with $dx$ in degree one.

\begin{fact}[Veronese subrings {\cite[\S3.6]{BrunsHerzog1998}, \cite[Thm.~3.1.1, Cor.~3.1.3, \S5.1(5.1.8)]{GotoWatanabe1978}}]
\label{fact:veronese}
Let $R=\bigoplus_{d\ge0}R_d$ be a finitely generated graded domain with $R_0=\C$
and $X=\Proj R$. There is $e\ge1$ such that $R^{(e)}=\bigoplus_jR_{ej}$ is
generated by $R_e$; then $X=\Proj R^{(e)}$, and a finitely generated graded $R$-module $M$ defines
the coherent sheaf $\widetilde M=(M^{(e)})^\sim$ on $X$. If $R$ is
Cohen--Macaulay and normal, then $R^{(e)}$ is Cohen--Macaulay with canonical
module $(\omega_R)^{(e)}$, and $\omega_X\simeq\widetilde{\omega_R}$.
\end{fact}

\begin{fact}[Finite quotients preserve Cohen--Macaulayness]
\label{fact:finite-quotient-cm}
Let $R$ be a finitely generated Cohen--Macaulay $\C$-algebra, and let a
finite group $\Gamma$ act on $R$ by $\C$-algebra automorphisms. Then
$R^\Gamma$ is Cohen--Macaulay. Consequently a finite-group quotient of
a Cohen--Macaulay complex variety is Cohen--Macaulay.
\end{fact}

\begin{proof}
Put $S=R^\Gamma$. The extension $S\subset R$ is finite. Fix a maximal
ideal $\mathfrak m$ of $S$, put $d=\dim S_{\mathfrak m}$, and set
$B=R\otimes_S S_{\mathfrak m}$. The maximal ideals of $B$ form one
$\Gamma$-orbit, so their local rings have the same dimension. Since
$S_{\mathfrak m}\subset B$ is integral, $\dim B=d$; hence every such
local ring has dimension $d$ and is Cohen--Macaulay.

Choose a system of parameters $x_1,\ldots,x_d$ of $S_{\mathfrak m}$.
The ring $B/(x_1,\ldots,x_d)B$ is finite over the zero-dimensional
ring $S_{\mathfrak m}/(x_1,\ldots,x_d)$, so the $x_i$ form a system
of parameters, and therefore a regular sequence, in every maximal
localization of $B$. They are consequently a regular sequence on $B$.
The Reynolds projection
\[
 B\longrightarrow S_{\mathfrak m},\qquad
 b\longmapsto\frac1{|\Gamma|}\sum_{\gamma\in\Gamma}\gamma(b),
\]
is $S_{\mathfrak m}$-linear and descends after quotienting by any
initial segment of this sequence. Thus $S_{\mathfrak m}$ and each
of these successive quotients are direct summands of the corresponding
quotients of $B$. The sequence is therefore regular on
$S_{\mathfrak m}$, proving that this local ring is Cohen--Macaulay.
Checking maximal localizations suffices, since Cohen--Macaulayness
localizes. The variety assertion follows on affine quotient charts.
\end{proof}

\subsection{Affine semigroups and toric-face rings}\label{subsec:bg-toric}
Let $L$ be a lattice and $C\subset L_\R$ a full-dimensional pointed rational
polyhedral cone. For each facet $F$ of $C$, let $u_F\in L^\vee$ be the primitive
inward normal of $F$, so that $C=\{x:u_F(x)\ge0\ \forall F\}$ and $u_F(L)=\Z$; we
call $u_F(x)$ the \emph{height} of $x$ over $F$. Put $Q=C\cap L$ and let
$\relint Q$ be the set of points of $Q$ in the interior of $C$.

\begin{fact}[{\cite[Prop.~6.1.2, Thm.~6.1.4, Thm.~6.3.5]{BrunsHerzog1998}}]
\label{fact:semigroup}
The semigroup $Q$ is finitely generated (Gordan's lemma), and the semigroup ring
$\C[Q]=\bigoplus_{q\in Q}\C x^q$ is a normal Cohen--Macaulay domain (Hochster \cite[Thm.~1]{Hochster1972}).
The canonical module of $\C[Q]$ is the ideal spanned by the interior monomials
(Danilov--Stanley \cite[\S4.6]{Danilov1978}, \cite[Thm.~6.7]{Stanley1978}):
\[
 \omega_{\C[Q]}=\textstyle\bigoplus_{q\in\relint Q}\C x^q ,
\]
the monomial $x^q$ corresponding to the form $x^q\cdot\eta$, where
$\eta=\frac{dx^{m_1}}{x^{m_1}}\wedge\dots\wedge\frac{dx^{m_d}}{x^{m_d}}$ is the
logarithmic volume form of the torus $\Spec\C[L]$, for a basis
$m_1,\dots,m_d$ of $L$.
\end{fact}

For instance, for $Q=\N^n$ the interior monomials are the multiples of
$x_1\cdots x_n$, and $x_1\cdots x_n\cdot\frac{dx_1}{x_1}\wedge\dots\wedge
\frac{dx_n}{x_n}=dx_1\wedge\dots\wedge dx_n$. The logarithmic form $\eta$ has a
simple pole along every toric divisor, and the monomial $x^q$ vanishes along the
toric divisor of a facet $F$ to order $u_F(q)$; so $x^q\eta$ is regular exactly
when $u_F(q)\ge1$ for all $F$.

\begin{definition}\label{def:translator}
A \emph{canonical translator} of $Q$ is a point $\vartheta\in L$ such that
$\relint Q=\vartheta+Q$.
\end{definition}

By Fact~\ref{fact:semigroup}, $Q$ has a canonical translator $\vartheta$ if and
only if $\omega_{\C[Q]}=x^\vartheta\C[Q]$ is free, i.e.\ if and only if $\C[Q]$ is
Gorenstein; the translator is then unique.

\begin{lemma}\label{lem:height-one}
A point $\vartheta\in L$ is the canonical translator of $Q$ if and only if
$u_F(\vartheta)=1$ for every facet $F$ of $C$.
\end{lemma}

\begin{proof}
If $u_F(\vartheta)=1$ for every $F$, then for $q\in L$ the strict inequalities
$u_F(q)>0$, which for integers mean $u_F(q)\ge1$, are equivalent to
$u_F(q-\vartheta)\ge0$. Conversely, let $\vartheta$ be the translator and fix a
facet $F$. Choose $z\in L$ with $u_F(z)=1$ and $y\in L$ in the relative interior
of $F$. For $n\gg0$ the point $z+ny$ lies in the interior of $C$ and has height
$u_F(z+ny)=1$. Writing $z+ny=\vartheta+q$ with $q\in Q$ gives
$1=u_F(\vartheta)+u_F(q)$, with $u_F(\vartheta)\ge1$ and $u_F(q)\ge0$.
\end{proof}

\begin{figure}[htbp]
\centering
\begin{tikzpicture}[scale=0.62,font=\small]
 \fill[blue!7] (0,0)--(3.6,0)--(3.6,3.6)--(1.8,3.6)--cycle;
 \draw[thick] (0,0)--(3.6,0);\draw[thick] (0,0)--(1.8,3.6);
 \foreach \x in {0,...,3}\foreach \y in {0,...,3}{
   \pgfmathparse{(\y<=2*\x)?1:0}\ifnum\pgfmathresult=1
     \fill (\x,\y) circle (2pt);\fi}
 \foreach \x/\y in {1/1,2/1,2/2,2/3,3/1,3/2,3/3}\draw[red,thick] (\x,\y) circle (4pt);
 \node[below] at (1,1) {$\vartheta$};
 \node at (1.8,-0.8) {$\mathrm{cone}((1,0),(1,2))$};
 \begin{scope}[xshift=6cm]
 \fill[blue!7] (0,0)--(3.6,0)--(3.6,3.6)--(1.2,3.6)--cycle;
 \draw[thick] (0,0)--(3.6,0);\draw[thick] (0,0)--(1.2,3.6);
 \foreach \x in {0,...,3}\foreach \y in {0,...,3}{
   \pgfmathparse{(\y<=3*\x)?1:0}\ifnum\pgfmathresult=1
     \fill (\x,\y) circle (2pt);\fi}
 \foreach \x/\y in {1/1,1/2,2/1,2/2,2/3,3/1,3/2,3/3}\draw[red,thick] (\x,\y) circle (4pt);
 \node at (1.8,-0.8) {$\mathrm{cone}((1,0),(1,3))$};
 \end{scope}
\end{tikzpicture}
\caption{Interior lattice points (circled). On the left they form the translate
$(1,1)+Q$, so $\C[Q]=\C[u,v,w]/(uw-v^2)$ is Gorenstein with canonical
translator $\vartheta=(1,1)$. On the right, the interior points $(1,1)$ and
$(1,2)$ are both minimal, so there is no translator, and the cone over the
twisted cubic $\C[Q]$ is not Gorenstein.}\label{fig:translator}
\end{figure}
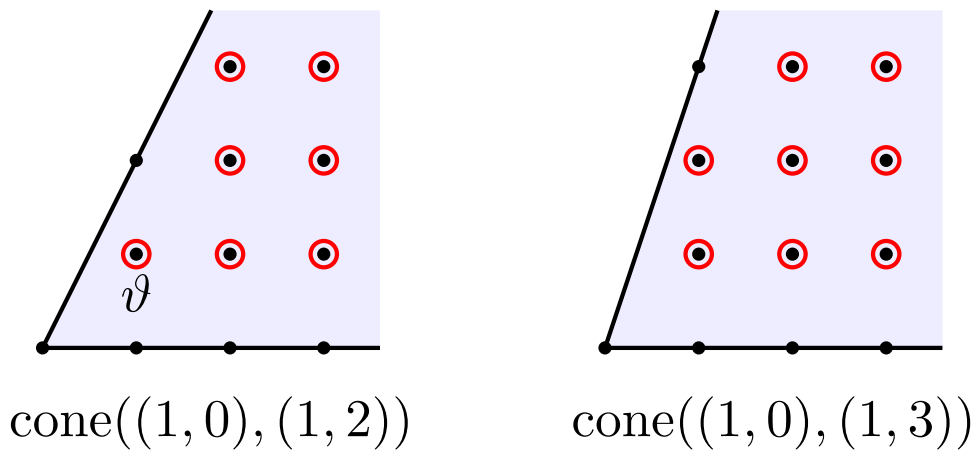

In Figure~\ref{fig:translator}, the left cone has facet normals $u_1(x,y)=y$
and $u_2(x,y)=2x-y$, and $(1,1)$ has height one over both, as
Lemma~\ref{lem:height-one} predicts. On the right, the facet normals are $y$ and
$3x-y$, and no lattice point has height one over both.

For a finite rational fan $\mathscr G$, the toric-face ring
$\C[\mathscr G]$ has basis $x^q$ for the lattice points in
$|\mathscr G|$. The product in $\C[\mathscr G]$ is $x^qx^{q'}=x^{q+q'}$ when $q,q'$ lie in a
common cone, and zero otherwise. The affine spectrum of $\C[\mathscr G]$ glues the toric
varieties of maximal cones along their common toric faces. For the fan
of proper faces of $\R^3_{\ge0}$ this gives $\C[x,y,z]/(xyz)$.

\begin{fact}[Toric varieties {\cite[Thms.~8.2.3, 9.2.3, 11.1.9]{CoxLittleSchenck2011}}]
\label{fact:toric}
Let $Y$ be a normal projective toric variety with toric boundary $\partial Y$
(the complement of the dense torus). Then $H^i(Y,\mathcal O_Y)=0$ for $i>0$,
the pair $(Y,\partial Y)$ is log canonical (see \S\ref{subsec:bg-pairs}), and the
logarithmic volume form of the torus has divisor $-\partial Y$.
\end{fact}
For the log canonicity assertion, take a smooth toric resolution
$f:Y'\to Y$ by subdividing the fan. The same logarithmic form gives
$K_{Y'}+\partial Y'=f^*(K_Y+\partial Y)=0$; therefore every exceptional
divisor has discrepancy $-1$, as required.

\subsection{Filtrations, valuations and Rees algebras}\label{subsec:bg-filtrations}
Let $A$ be a finitely generated $\C$-domain. A positive filtration is
exhaustive, increasing, multiplicative and finite dimensional in each degree,
with $F_0A=\C$. Set $F_{-1}A=0$ and
\[
 \gr_FA=\bigoplus_{d\ge0}F_dA/F_{d-1}A,\qquad
 \cR_F(A)=\bigoplus_{d\ge0}F_dA\,t^d\subset A[t].
\]
For $f\ne0$, write $\delta_F(f)=\min\{d:f\in F_dA\}$ and
$\inn_F(f)\in(\gr_FA)_{\delta_F(f)}$ for its initial form. Since
$\cR_F(A)/(t)=\gr_FA$ and $\cR_F(A)[t^{-1}]=A[t^{\pm1}]$,
$\Spec\cR_F(A)\to\A^1$ is a flat degeneration from $A$ to $\gr_FA$.
We call $F$ projective when $\cR_F(A)$ is finitely generated. Then
$\Proj\cR_F(A)$ compactifies $\Spec A=D_+(t)$ with boundary
$\Proj\gr_FA$.

\begin{example}\label{ex:rees-cubic}
Let $A_0=\C[x,y,z]$ have the total-degree filtration $F$, and give
$A=A_0/(f)$, where $f=x^2+y^2+z^2-xyz-c$, the quotient filtration $\bar F$.
The leading form $\inn_Ff=-xyz$ is a nonzerodivisor in $\gr_FA_0$.
Lemma~\ref{lem:filtered-regular} therefore gives
$\gr_{\bar F}A=\C[x,y,z]/(xyz)$ and, with $X=xt$, $Y=yt$, $Z=zt$,
\[
 \cR_{\bar F}(A)=\C[t,X,Y,Z]/\bigl(t(X^2+Y^2+Z^2)-XYZ-ct^3\bigr),
\]
with all four generators of degree one. For $c=\tau+2$, the Proj is
$X_\tau$ and the Proj of the quotient by $t$ gives the boundary $D_\tau$ of
\hyperref[subsec:intro-example]{the introductory model calculation}.

The compatibility between quotienting and taking associated graded rings
is expressed by the commuting square of ring homomorphisms
\[
 \begin{array}{ccc}
 \cR_F(A_0)&\longrightarrow&\cR_{\bar F}(A)\\
 \big\downarrow&&\big\downarrow\\
 \gr_FA_0&\longrightarrow&\gr_{\bar F}A .
 \end{array}
\]
The vertical maps quotient by $t$. The top map quotients by
$ft^3=t(X^2+Y^2+Z^2)-XYZ-ct^3$, and the bottom map quotients by
$\inn_Ff=-xyz$. Taking the fiber at $t=1$ instead recovers
$A_0\to A_0/(f)$. Thus the Rees algebra relates the original equation
to its leading form through actual ring maps; taking initial forms
alone is not an additive map $A_0\to\gr_FA_0$.
\end{example}

\begin{lemma}\label{lem:filtered-regular}
Let $F$ be a positive filtration of $A$ and let $f_1,\dots,f_s\in A$ be
elements whose initial forms form a regular sequence in $\gr_FA$. Give
$\bar A=A/(f_1,\dots,f_s)$ the quotient filtration $\bar F$. Then
\[
 \gr_{\bar F}\bar A=\gr_FA/(\inn_Ff_1,\dots,\inn_Ff_s),\qquad
 \cR_{\bar F}(\bar A)=\cR_F(A)/(f_1t^{d_1},\dots,f_st^{d_s}),\quad
 d_i=\delta_F(f_i).
\]
\end{lemma}

\begin{proof}
By induction on $s$ it suffices to treat $s=1$; write $f=f_1$ and $d=d_1$. Since
$\inn_Ff$ is a nonzerodivisor, $\inn_F(fg)=\inn_F(f)\inn_F(g)\ne0$, hence
$\delta_F(fg)=d+\delta_F(g)$ for every $g\ne0$. Therefore
$(f)\cap F_nA=f\cdot F_{n-d}A$ for every $n$. In each degree $n$, this equality
identifies $\bar F_n\bar A=F_nA/((f)\cap F_nA)$ with $F_nA/fF_{n-d}A$, which
gives both formulas.
\end{proof}

Without the regularity assumption the first formula can fail: the initial forms
of the elements of an ideal need not be generated by the initial forms of the
generators.

We use largest-term valuations
$v:A\setminus\{0\}\to\Z^M$ for the lexicographic order:
\[
 v(fg)=v(f)+v(g),\qquad v(f+g)\le\max(v(f),v(g)),\qquad v(\C^*)=0.
\]
The leaf at $q$ is $\{f:v(f)\le q\}/\{f:v(f)<q\}$ (with zero included).
Write $Q=v(A\setminus\{0\})$. One-dimensional leaves give a graded domain
$\gr_vA$ with one-dimensional nonzero pieces and hence a graded isomorphism
$\gr_vA\simeq\C[Q]$ \cite[\S2]{KavehManon2019}.
Here $\inn_v(f)$ denotes the class of $f$ in its leaf. A finite set whose
initial forms generate $\gr_vA$ is a Khovanskii basis.

\begin{proof}[The graded-algebra identification]
Let $B=\bigoplus_{q\in Q}B_q$ be a commutative domain graded by a
subsemigroup $Q$ of a lattice, with $B_0=\C$ and $\dim_\C B_q=1$ for every $q\in Q$.
Then $B\simeq\C[Q]$ as a graded algebra. The domain hypothesis is essential:
localizing at the nonzero homogeneous elements gives a graded ring with
one-dimensional pieces indexed by $\gp(Q)$. Choose a homogeneous unit for
each element of a basis of this free abelian group. Products of these units
identify the localization with $\C[\gp(Q)]$, and its subring $B$ with
$\C[Q]$.
\end{proof}

The following lemma, proved in Section~\ref{app:scalar} following
\cite[\S2.2, \S8.2]{KavehManon2019}, replaces a valuation of higher rank by a
filtration indexed by integers while preserving the associated graded ring
and the initial forms of a prescribed finite collection of elements.
The idea is to replace the lexicographic comparisons by a linear weight:
large positional weights reproduce any finite set of comparisons, and a
Gr\"obner basis reduces the problem to finitely many comparisons.

\begin{lemma}[Scalarization]\label{lem:scalarization}
Let $v:A\setminus0\to\N^M$ be a valuation with one-dimensional leaves whose value
semigroup $Q$ is finitely generated, and let $g_1,\dots,g_n\in A$ with
$\inn_v(g_i)=c_ix^{q_i}$, $c_i\in\C^*$, $q_i\ne0$. Then there are a positive
integral linear form $h$ on $\Z^M$ and a projective positive filtration $F$ of
$A$ such that $\gr_FA\simeq\C[Q]$, the monomial $x^q$ having degree $h(q)$, and
$\inn_F(g_i)=c_ix^{q_i}$ for all $i$. If a finite group $\Delta$ acts on $A$
preserving $v$ and acting on each leaf through a character, and if every $g_i$
is a $\Delta$-eigenvector, then $F$ may be chosen $\Delta$-stable.
\end{lemma}

We apply scalarization to the peripheral traces
$\chi_{\varpi_a}(\gamma_i)$, preserving the monomial initial equations needed
for Theorem~\ref{thm:special-fiber}.

\subsection{Singularities of pairs}\label{subsec:bg-pairs}
A pair $(X,D)$ consists of a normal variety and a reduced Weil divisor
with $K_X+D$ $\Q$-Cartier. On a log resolution $f:Y\to X$, the union of
the strict transform and exceptional divisors has simple normal crossings,
and we write
\[
 K_Y+f_*^{-1}D=f^*(K_X+D)+\sum_{E\ {\rm exceptional}}a(E)E.
\]
The pair is log canonical (lc) if all $a(E)\ge-1$, and dlt if some log
resolution has all exceptional $a(E)>-1$
\cite[Defs.~2.34, 2.37, Prop.~2.40, Thm.~2.44]{KollarMori1998}.
For an effective rational boundary $B=\sum b_iB_i$, lc also requires
$b_i\le1$; klt requires $b_i<1$ and all exceptional discrepancies $>-1$
\cite[Def.~2.34]{KollarMori1998}. Canonical singularities mean that $K_X$
is $\Q$-Cartier and the discrepancies of $(X,0)$ are nonnegative.
We call $X$ $\Q$-Gorenstein if $K_X$ is $\Q$-Cartier, and $\Q$-factorial
if every Weil divisor has a positive Cartier multiple.

When $\divi(\beta)=-D$, the pair is lc precisely when the pullback of
the logarithmic volume form $\beta$ to every resolution has at most
simple poles. For $D=\{xy=0\}\subset\A^2$, blowing up with $y=xs$ gives
$f^*(dx\wedge dy/(xy))=dx/x\wedge ds/s$, so the exceptional pole is simple.
For three concurrent lines the corresponding pullback is
$dx\wedge ds/(x^2s(1-s))$, with a double pole; that pair is not lc.

\begin{fact}[{\cite[Cor.~1.36]{KollarSingMMP2013}, \cite[Prop.~5.20(4)]{KollarMori1998}}]
\label{fact:pairs}
\begin{enumerate}[label=\textup{(\roman*)}]
\item (Dlt modification.) A projective lc pair $(X,D)$ admits a projective
birational morphism $f:X'\to X$ with $X'$ $\Q$-factorial, such that
$D'=f_*^{-1}D+\sum_{E\ {\rm exceptional}}E$ is reduced, $(X',D')$ is dlt and
$K_{X'}+D'=f^*(K_X+D)$.
\item (Finite morphisms.) Let $q:Y\to X$ be a finite surjective morphism of
normal varieties, let $D$, $E$ be effective rational divisors, and assume
$K_X+D$ is $\Q$-Cartier and $K_Y+E=q^*(K_X+D)$. Then $(Y,E)$ is lc
(respectively klt) if and only if $(X,D)$ is lc (respectively klt).
\item (Products.) The product of two lc pairs, with the sum of the pulled back
boundaries, is lc.
\end{enumerate}
\end{fact}
For (iii), take the product of log resolutions of the two pairs. The
boundary in the crepant pullback is the sum of the pulled-back boundaries
of those resolutions, so every coefficient is at most $1$ and the support
has simple normal crossings. This proves log canonicity of the product.

Log canonicity is local, is preserved by automorphisms, and improves when
boundary components are removed. For a smooth curve $C$, the pair
$(X\times C,D\times C)$ is lc if and only if $(X,D)$ is lc: multiply a
log resolution by $C$ to obtain the same discrepancies.

For a reduced divisor $S$ on a normal variety $\mathcal Y$, a rational
top form with at most a simple pole can be written generically as
$\beta=dg/g\wedge\gamma+\beta'$, with $g=0$ defining $S$ and $\beta'$
regular along $S$. The residue of $\beta$ is $\gamma|_S$, regarded as a rational
top form on $S^\nu$. For example,
$\operatorname{Res}_{y=0}(dx\wedge dy/y)=-dx$.

\begin{lemma}[Adjunction]\label{lem:adjunction}
Let $\mathcal Y$ be a Gorenstein variety and $S\subset\mathcal Y$ a reduced
effective Cartier divisor. Then $\omega_S\simeq\omega_{\mathcal Y}(S)|_S$, the
isomorphism being given by the residue. In particular, if $S$ is locally
isomorphic, at a point of codimension one in $S$, to the node $\{ab=0\}$ (two
smooth branches meeting transversally), then a local generator of $\omega_S$
restricts on each branch to a form with a simple pole along the intersection
of the branches.
\end{lemma}

\begin{proof}
The dualizing sheaf of the Cartier divisor $S$ is
$\mathcal Ext^1(\mathcal O_S,\omega_{\mathcal Y})$, and the exact sequence
$0\to\mathcal O(-S)\to\mathcal O\to\mathcal O_S\to0$ identifies this sheaf
with $\omega_{\mathcal Y}(S)/\omega_{\mathcal Y}=\omega_{\mathcal Y}(S)|_S$; on
the smooth locus the identification is the residue. For the node, compute in
the ambient plane: $\omega_{\{ab=0\}}$ is generated by the residue of
$\frac{da\wedge db}{ab}$, which is $\frac{db}{b}$ on the branch $\{a=0\}$ and
$-\frac{da}{a}$ on the branch $\{b=0\}$.
\end{proof}

Suppose $K_{\mathcal Y}+S+\Delta$ is Cartier, with $\Delta$ reduced and
having no component in common with $S$. For a local rational top form
$\beta$ with $\divi(\beta)=-S-\Delta$, define the different on $S^\nu$ by
\[
 \divi(\operatorname{Res}_S\beta)=-\operatorname{Diff},\qquad
 K_{S^\nu}+\operatorname{Diff}=\nu^*((K_{\mathcal Y}+S+\Delta)|_S).
\]
The different is independent of the choice of $\beta$; the right side denotes the
pullback of the restricted Cartier class. If $\mathcal Y$ and $S$ are smooth
and $S+\Delta$ has simple normal crossings, then
$\operatorname{Diff}=\Delta|_S$.

\begin{fact}[Inversion of adjunction {\cite[Theorem]{Kawakita2007}, \cite[Thm.~1.1]{Fujino2023}}]
\label{fact:inversion}
Let $\mathcal Y$ be a normal variety, $S$ a reduced divisor, and $\Delta$ a
reduced divisor without common components with $S$, such that
$K_{\mathcal Y}+S+\Delta$ is Cartier. Let $\nu:S^\nu\to S$ be the normalization.
Then $(\mathcal Y,S+\Delta)$ is log canonical in a neighborhood of $S$ if and
only if $(S^\nu,\operatorname{Diff})$ is log canonical.
\end{fact}

For $S=\{y=0\}$ and $\Delta=\{y=x^2\}$ in $\A^2$, the form
$dx\wedge dy/(y(y-x^2))$ has residue $dx/x^2$, giving
$\operatorname{Diff}=2[0]$; both pairs fail to be lc. Replacing $\Delta$
by the transverse line $x=0$ gives residue $-dx/x$ and an lc pair.

\subsection{Quasi-Hamiltonian spaces}\label{subsec:bg-qh}
Quasi-Hamiltonian reduction supplies symplectic local models for individual
conjugacy classes in Corollary~\ref{cor:symplectic-coarse}. No such structure
is assumed on a possibly singular full reduced adjoint fiber; that case uses
Section~\ref{subsec:spectral-interior}.
Let $L$ be connected reductive with a nondegenerate invariant symmetric form
$\kappa$ on $\mathfrak l$. The left and right Maurer--Cartan forms are
$\theta=g^{-1}dg$ and $\bar\theta=dg\,g^{-1}$, respectively; these are
$\mathfrak l$-valued one-forms on $L$. Pairing their coefficients with
$\xi\in\mathfrak l$ by $\kappa$ gives the ordinary one-form
$\kappa(\theta+\bar\theta,\xi)$ on $L$.
Put $\eta=\frac1{12}\kappa(\theta,[\theta,\theta])$. Our fundamental vector
field $v_\xi$ generates $\exp(-t\xi)$. A quasi-Hamiltonian $L$-space is a
smooth $L$-variety $M$ with an invariant two-form $\omega$ and an equivariant
map $\mu:M\to L$ (for conjugation) satisfying
\[
 d\omega=\mu^*\eta,\qquad
 \iota_{v_\xi}\omega=\tfrac12\mu^*\kappa(\theta+\bar\theta,\xi),\qquad
 \ker\omega_x=\{v_\xi(x):\Ad_{\mu(x)}\xi=-\xi\}.
\]
For example, $L=\C^*$ acts trivially on $M=(\C^*)^2$:
$\lambda\cdot(a,b)=(a,b)$. With $\mu(a,b)=1$ and
$\omega=da/a\wedge db/b$, the form is closed and nondegenerate; the
contraction identity is zero on both sides and $\Ad_1\xi=-\xi$ forces
$\xi=0$.

\begin{fact}[{\cite[Props.~3.1--3.2, Thm.~6.1]{AMM1998}, \cite[\S2, Examples~3--4, Thm.~4]{Boalch2007}}]\label{fact:qh}\mbox{}
\begin{enumerate}[label=\textup{(\roman*)}]
\item A conjugacy class $\cC\subset L$, with $\mu$ the inclusion, is a
quasi-Hamiltonian $L$-space.
\item The double $L\times L$, with action
$(g_1,g_2)\cdot(a,b)=(g_1ag_2^{-1},g_2bg_1^{-1})$ and moment map
$(ab,a^{-1}b^{-1})$, is a quasi-Hamiltonian $L\times L$-space.
\item (Fusion.) If $M$ is a quasi-Hamiltonian $L\times L\times K$-space with
moment map $(\mu_1,\mu_2,\mu_3)$, then $M$ with the diagonal $L\times K$-action
and moment map $(\mu_1\mu_2,\mu_3)$ is a quasi-Hamiltonian $L\times K$-space
with the modified two-form
\[
 \omega_{\mathrm{fus}}
 =\omega-\tfrac12\kappa(\mu_1^*\theta,\mu_2^*\bar\theta),
\]
where $\kappa$ pairs the Lie-algebra coefficients and wedges the one-forms
\cite[\S2, Thm.~4]{Boalch2007}. The sign agrees with the convention
$v_\xi(x)=\left.\frac{d}{dt}\right|_{t=0}(\exp(-t\xi)\cdot x)$ above.
\end{enumerate}
Consequently, first fusing the two $L$-actions on each double and then fusing
the $g$ doubles with the classes $\cC_1,\dots,\cC_k$ makes
$L^{2g}\times\prod_i\cC_i$, with the diagonal conjugation action, a
quasi-Hamiltonian $L$-space with moment map $\mu=\prod_j[A_j,B_j]\prod_iM_i$.
\end{fact}

\begin{lemma}[Reduction]\label{lem:qh-reduction}
Let $(M,\omega,\mu)$ be a quasi-Hamiltonian $L$-space on which $Z(L)$ acts
trivially and such that $\mu(M)\subset L_\der$, and suppose $1$ is a regular
value of $\mu:M\to L_\der$. Let $R=\mu^{-1}(1)$. Then the restriction of $\omega$
to $R$ is closed, $L$-invariant, and has kernel exactly the tangent spaces to the
$L$-orbits.
\end{lemma}

\begin{proof}
On $R$ we have $d\omega=\mu^*\eta=0$. At $x\in R$ we have $\mu(x)=1$, so the
kernel axiom says that $\omega_x$ is nondegenerate on $T_xM$, and the contraction
axiom says that $\iota_{v_\xi}\omega_x=\kappa(d\mu_x(\cdot),\xi)$. Since $Z(L)$ acts
trivially, $v_\xi=0$ for $\xi\in\mathfrak z(\mathfrak l)$, and
$\mathfrak l=\mathfrak z(\mathfrak l)\oplus\mathfrak l_\der$ is a
$\kappa$-orthogonal decomposition. Hence the $\omega_x$-orthogonal of
$T_xR=\ker d\mu_x$ is $\{v_\xi(x):\xi\in\mathfrak l_\der\}=T_x(L\cdot x)$, and the
kernel of $\omega_x|_{T_xR}$ is $T_xR\cap T_x(L\cdot x)=T_x(L\cdot x)$.
\end{proof}

\begin{fact}[Finite symplectic quotients {\cite[Props.~1.3, 2.4]{Beauville2000}, \cite[Lem.~5.16]{KollarMori1998}}]
\label{fact:symplectic-quotient}
Let a finite group $F$ act on a smooth variety $S$ preserving an algebraic
symplectic form $\sigma$. Then $S/F$ is normal, $\Q$-factorial,
Cohen--Macaulay and Gorenstein with canonical singularities, and the invariant
volume form $\sigma^n/n!$ ($2n=\dim S$) descends to a generator of
$\omega_{S/F}$.
\end{fact}

For the final assertion, a nontrivial finite-order symplectic automorphism
has fixed loci of even codimension. The effective quotient is therefore
unramified in codimension one, and the invariant volume descends without
zeros or poles there. For the action $\varepsilon(x,y)=(\varepsilon x,
\varepsilon y)$ of $F=\{\pm1\}$ on $\C^2$,
\[
 \C[x,y]^F=\C[u,v,w]/(uw-v^2),\qquad u=x^2,\ v=xy,\ w=y^2.
\]
The action of $F$ preserves $dx\wedge dy$, which induces a canonical generator on the
quadric cone, also the toric example on the left of
Figure~\ref{fig:translator}.

\subsection{Nilpotent cones}\label{subsec:bg-nilpotent}
Reduced puncture fibers are locally products of smooth varieties and
nilpotent cones (Lemma~\ref{lem:spectral-fiber-local}). For a semisimple
Lie algebra $\mathfrak l$, $\mathcal N(\mathfrak l)$ consists of elements
acting nilpotently in a faithful representation. For example,
$\mathcal N(\mathfrak{sl}_2)=\{(\begin{smallmatrix}a&b\\c&-a\end{smallmatrix}):
a^2+bc=0\}$ is the quadric cone of \eqref{eq:intro-gl2-fiber}, isomorphic
to $\C^2/\{\pm1\}$. Analytic local $\Q$-factoriality means that every
Weil divisor of each analytic germ has a positive Cartier multiple.

\begin{fact}[Nilpotent cones {\cite[\S3, Prop.~4.1, \S5]{LehnNamikawaSorger2012},
\cite[\S\S3.2--3.3, Thms.~6.7.2--6.7.4]{ChrissGinzburg1997}}]
\label{fact:nilpotent-cones}
The reduced, normal complete-intersection structure in (i) goes back to Kostant
\cite{Kostant1963}. The cited account of Chriss--Ginzburg also describes
the Springer resolution; analytic local $\Q$-factoriality is the assertion of
Lehn--Namikawa--Sorger. Let $L$ be a connected semisimple group with Lie algebra $\mathfrak l$ and Borel
subgroup $B_L$.
\begin{enumerate}[label=\textup{(\roman*)}]
\item The nilpotent cone $\mathcal N(\mathfrak l)$ is an irreducible normal
complete intersection of dimension $\dim\mathfrak l-\rank\mathfrak l$, and
a finite union of $L$-orbits, each of even dimension.
\item The \emph{Springer resolution}, a morphism
$\pi:T^*(L/B_L)\to\mathcal N(\mathfrak l)$, is proper and birational, and $K_{T^*(L/B_L)}=\pi^*K_{\mathcal N(\mathfrak l)}$.
Consequently $\mathcal N(\mathfrak l)$ is Gorenstein with canonical
singularities.
\item $\mathcal N(\mathfrak l)$ is analytically locally $\Q$-factorial.
\item Finite products of nilpotent cones, and products of such with smooth
varieties, are again irreducible normal local complete intersections with
canonical singularities which are analytically locally $\Q$-factorial.
\end{enumerate}
\end{fact}

For $\mathfrak l=\mathfrak{sl}_2$ the Springer resolution is the blow-up of the
quadric cone at the vertex, i.e.\ $T^*\PP^1$, whose exceptional curve has
self-intersection $-2$. The product assertion, especially analytic local $\Q$-factoriality, is
justified in Section~\ref{app:geometric-facts}.

\subsection{Poisson brackets and differential forms}\label{paper:bg-poisson}
A \emph{Poisson bracket} on an algebra $A$ is a skew bilinear operation
$\{-,-\}:A\times A\to A$ satisfying the Jacobi identity and the Leibniz
rule $\{a,bc\}=\{a,b\}c+b\{a,c\}$. The same definition applies to a
structure sheaf. A closed nondegenerate two-form on a smooth variety
defines a Poisson bracket by inversion. An ideal $I$ is Poisson when
$\{A,I\}\subset I$.

A normal variety $Z$ has \emph{symplectic singularities} if $Z_{\rm reg}$
has a closed nondegenerate algebraic two-form whose pullback extends
regularly over a resolution. Extension over one resolution implies
extension over every resolution by a common resolution
\cite[Def.~1.1]{Beauville2000}.

\begin{fact}[Normalization and form extension]\label{paper:fact-poisson-forms}
\begin{enumerate}[label=\textup{(\roman*)}]
\item A finitely generated integral Poisson $\C$-algebra has a unique
Poisson bracket on its normalization extending the given bracket.
Localization, reduction and passage to an irreducible component preserve
Poisson structures in characteristic zero
(Lemma~\ref{paper:lem-poisson-normalization}; see \cite[Thm.~1.5]{KaledinNormalization}).
\item If $Z$ is a complex quasiprojective variety with klt singularities
and $\pi:\widetilde Z\to Z$ is a resolution, every regular differential
form on $Z_{\rm reg}$ extends regularly to $\widetilde Z$
\cite[Thm.~1.4]{GrebKebekusKovacsPeternell2011}.
In particular this applies to finite quotient singularities.
Lemma~\ref{paper:lem-finite-quotient-forms} gives a direct proof of the
finite quotient case used below.
\end{enumerate}
\end{fact}
The finite-type algebras in (i) are excellent, as required in the
normalization theorem. Section~\ref{app:poisson-algebra} gives a
short proof of the normalization assertion.

\subsection{The Hodge-theoretic inputs}\label{paper:bg-hodge}
For a smooth complex algebraic variety $V$, write $W_\bullet$ and
$F^\bullet$ for the weight and Hodge filtrations on cohomology.
For $V=\overline V\setminus B$ with $\overline V$ smooth projective
and $B$ a simple-normal-crossings divisor, $F^\bullet$ comes from the
logarithmic de Rham complex.

\begin{fact}[Hodge theory {\cite[\S2.2, Thms.~2.3.5, 3.2.5,
Cor.~3.2.13]{Deligne1971}}]\label{paper:fact-hodge}
Algebraic pullbacks preserve both filtrations and are strict:
$f(F^pH)=\operatorname{im}f\cap F^pH'$.
Cup products preserve $F^\bullet$. For $d=\dim V$, logarithmic
Hodge--de Rham degeneration gives
\[
 F^dH^d(V,\C)=H^0(\overline V,\Omega_{\overline V}^d(\log B)).
\]
On a smooth projective variety, cohomology has the pure Hodge
decomposition; holomorphic two-forms inject into $H^2$ with type
$(2,0)$, whereas first Chern classes of algebraic line bundles have
type $(1,1)$.
\end{fact}
These statements are used on resolutions and smooth exact strata.
The comparison between Betti and Dolbeault spaces will only transport
differential forms and compact cycles.

\subsection{Parahoric Higgs bundles and the Hitchin map}
\label{subsec:bg-higgs}

Let $C$ be a smooth projective complex curve, let
$D=p_1+\cdots+p_k$ have distinct support, and let $G$ be connected
semisimple. Write $\mathfrak g=\operatorname{Lie}G$, choose a maximal
torus $T$ with Lie algebra $\mathfrak t$, and choose real weights
$\alpha_i\in X_*(T)\otimes\mathbb R$ in a closed affine Weyl alcove.
Here an alcove is a chamber cut out by the affine root hyperplanes
$r(\alpha)\in\mathbb Z$; its faces are called facets.

\paragraph{The local lattices.}
With a coordinate $z$ at $p_i$, put $\mathcal O=\mathbb C[[z]]$ and
$\mathcal K=\mathbb C((z))$. The parahoric subgroup of $G(\mathcal K)$ is generated by
$T(\mathcal O)$ and the root subgroups
$U_r(z^{\lceil-r(\alpha_i)\rceil}\mathcal O)$. The Lie lattice of this parahoric group is
\[
 \mathfrak p_{\alpha_i}
 =\mathfrak t\otimes\mathcal O\ \oplus\!
   \bigoplus_r z^{\lceil-r(\alpha_i)\rceil}
                    \mathfrak g_r\otimes\mathcal O.
\]
The pro-unipotent radical has Lie lattice
\[
 \mathfrak p_{\alpha_i}^{+}
 =z\mathfrak t\otimes\mathcal O\ \oplus\!
   \bigoplus_r z^{\lfloor-r(\alpha_i)\rfloor+1}
                    \mathfrak g_r\otimes\mathcal O.
\]
``Pro-unipotent'' means an inverse limit of unipotent groups.
The finite-dimensional quotient
$\mathfrak l_i=\mathfrak p_{\alpha_i}/\mathfrak p_{\alpha_i}^{+}$
is the Lie algebra of the reductive quotient $L_i$ of the parahoric
group. For $\SL_2$ and $0<r(\alpha_i)<1$, the positive and negative
root exponents are respectively $0$ and $1$; reduction records the
diagonal part, so $L_i\simeq\mathbb G_m$.
These groups glue with the constant group $G$ away from $D$ to a
smooth affine group scheme $\mathcal G_\alpha$ on $C$
\cite[arXiv v3, Defs.~2.1, 2.4]{KSZ2024}.

A \emph{parahoric torsor} $E$ is a principal $\mathcal G_\alpha$-bundle,
locally trivial in the \'etale topology. The adjoint bundle of $E$
$\operatorname{ad}(E)$ has these Lie lattices locally. A
\emph{weak logarithmic Higgs field} is a section, with $K_C$ the
cotangent line bundle,
\[
 \varphi\in H^0\bigl(C,\operatorname{ad}(E)\otimes K_C(D)\bigr);
 \qquad \varphi=A_i(z)\,dz/z,\quad A_i(z)\in\mathfrak p_{\alpha_i}.
\]
On a curve the Higgs integrability condition $[\varphi,\varphi]=0$
is automatic. The \emph{graded residue} of $\varphi$ at $p_i$ is the image of $A_i(z)$ in
$\mathfrak l_i$, defined up to $L_i$-conjugacy. Weak fields allow an
arbitrary graded residue; strongly parahoric fields require it to vanish.
Prescribing an exact orbit $\mathcal O_i\subset\mathfrak l_i$ is
different from prescribing its closure $Y_i=\overline{\mathcal O_i}$:
the latter permits smaller residue orbits. This is the logarithmic
Higgs notion of \cite[arXiv v3, Def.~3.1]{KSZ2024}.

\paragraph{The stability convention.}
A reduction $\sigma$ to a parabolic subgroup is \emph{Higgs-compatible}
if $\varphi$ takes values in the corresponding parabolic Lie subbundle.
For its parahoric extension $P_\alpha$ and a character
$\chi:P_\alpha\to\mathbb G_m$, let
$L(\sigma,\chi)=E_\sigma\times^{P_\alpha}\mathbb C_\chi$.
The parabolic degree used throughout is
\[
 \operatorname{pardeg}_\alpha(E;\sigma,\chi)
 =\deg L(\sigma,\chi)+\sum_i\langle\alpha_i,\chi\rangle.
\]
At a marked point the character is transported according to the
relative position of the reduction and the prescribed parahoric data.
Ramanathan stability, respectively semistability, requires this number
to be positive, respectively nonnegative, for every Higgs-compatible
proper parabolic reduction and every nontrivial anti-dominant
character trivial on $Z(G)$. Anti-dominance means
$\langle\chi,\beta^\vee\rangle\leq0$ for the positive simple
coroots of a Borel contained in that parabolic. Taking positive powers
removes the restriction on the finite center without changing any
inequality. These are exactly the signs of
\cite[arXiv v3, Defs.~4.2, 4.5 and Rem.~4.4]{KSZ2024}.
For an ordinary $\SL_2$-bundle and an invariant line subbundle $F$,
the anti-dominant character gives $F^{-1}$, so the test is
$-\deg F>0$, the usual slope-stability inequality.

Fix also the discrete topological type. The degree-zero condition of
tame nonabelian Hodge theory imposes zero parabolic degree for global
characters of $G$; for semisimple $G$ these characters are trivial.
A polystable object is induced from a stable Higgs torsor for a Levi
subgroup, with zero degree for each character of the corresponding
parabolic trivial on $Z(G)$. Every semistable object has a polystable
graded representative. Two objects are \emph{S-equivalent} when these
representatives are isomorphic. Their moduli space is quasiprojective, with a stable open
parametrizing isomorphism classes
\cite[arXiv v3, Thm.~6.1]{KSZ2024}; the cover reduction in
Lemma~\ref{r2:lem-auxiliary-parahoric-cover} supplies the same
construction in every genus. For real weights the algebraic
construction uses rational weights in the same facet and stability
cell, preserving all inequalities and equalities; see
Lemma~\ref{r2:lem-proper-parahoric} and
\cite[arXiv v3, Rem.~6.11]{KSZ2024}.

\paragraph{The Hitchin morphism and exact local data.}
The invariant-polynomial construction originates with Hitchin for compact
curves \cite[\S1, p.~91]{Hitchin1987Integrable}.
Choose homogeneous generators $p_a\in\mathbb C[\mathfrak g]^G$,
of degrees $d_a$. The Hitchin morphism records their values:
\[
 h(E,\varphi)=(p_a(\varphi))_a
 \in\mathcal A_G(C,D):=
       \bigoplus_a H^0\bigl(C,K_C(D)^{\otimes d_a}\bigr).
\]
This is a finite-dimensional affine space. Invariant polynomials are
regular on the parahoric lattices, so the displayed pole bounds hold.
Closed invariant residue conditions $Y_i$ define a reduced closed
moduli subspace $M_{\alpha,\mathbf Y}$; its normalization is the
Dolbeault closure used below. The Hitchin map can have image smaller
than this ambient base. Projectivity of the Hitchin map for the actual stability
condition is proved in Lemma~\ref{r2:lem-proper-parahoric}.

On each exact residue stratum, the tame correspondence matches the
weights and residues with filtered local-system data. At zero Betti
weight these are ordinary representations with the corresponding
exact conjugacy classes
\cite[Table~1, Prop.~7.9 and Thm.~7.10]{BiquardGarciaPradaMundet2020}.
The local dictionary is spelled out in
Lemma~\ref{r2:lem-local-dol-dictionary}; no correspondence between
entire orbit closures is assumed.

\paragraph{Stacks, filtrations, and Hamiltonian reduction.}
A stack retains families and their automorphisms. For a group $A$
acting on $Q$, an object of $[Q/A]$ over a scheme $S$ is an
$A$-torsor on $S$ with an equivariant map to $Q$; $BA=[\mathrm{pt}/A]$
classifies $A$-torsors. On the stable locus, a coarse moduli space
records isomorphism classes of geometric points. \emph{Finite inertia} means that the
family of automorphism groups is finite over the stack. The separated
stable loci here have Deligne--Mumford local charts $[S/F]$ with
$F$ finite; their coarse charts are $S/F$.

The stack $\Theta=[\mathbb A^1/\mathbb G_m]$, for the scaling action,
encodes a filtered object: evaluation at $1$ gives the object, and
evaluation at $0$ its graded limit. A Harder--Narasimhan reduction is
the canonical destabilizing parabolic reduction of an unstable object.
The resulting $\Theta$-stratification organizes these filtrations into
locally closed strata, formally closed in successive open substacks;
the definitions and well-ordering condition are
\cite[\S2.3, Def.~6.1]{AHLH2023}. For principal Higgs bundles in
characteristic zero, the required well-ordered stratification is
\cite[Lemma~6.19]{AHLH2023}.

Finally, a Hamiltonian action of $K$ on a smooth symplectic variety
$(P,\omega)$ has an equivariant moment map
$\mu:P\to\mathfrak k^*$ satisfying
$\iota_{\xi^\#}\omega=-d\langle\mu,\xi\rangle$.
At a regular zero level with free action, the reduced form on
$\mu^{-1}(0)/K$ is determined by pullback to the restricted form
on $\mu^{-1}(0)$. Finite stabilizers give the corresponding quotient
charts. For smooth $Q$ and $P=T^*Q$, take $\omega=d\lambda$ and
$\langle\mu,\xi\rangle=\lambda(\xi^\#)$, where $\lambda$ is the
tautological one-form. The classical cotangent presentation of
$[Q/A]$ is $[\mu_A^{-1}(0)/A]$. Additional puncture framings have
graded residue as moment map, as proved in
Proposition~\ref{r2:prop-dol-closure}.

\section{The interior}\label{sec:interior}\label{sec:charvar}\label{sec:spectral}

This section studies the spectral character variety itself, independently of
any compactification. The aim is Proposition~\ref{prop:spectral-interior}:
under intrinsic genericity, the spectral representation variety is reduced,
normal, canonical and $\Q$-factorial, and the spectral character variety is
normal, $\Q$-factorial and klt. These properties supply the normality needed
by the construction at infinity, and the interior conclusions of Theorem~\ref{thm:spectral-main}.
The argument has three steps. Genericity forces every representation to be
irreducible (Section~\ref{subsec:irreducibility}). A reduced adjoint fiber is
locally a smooth factor times a nilpotent cone
(Section~\ref{subsec:spectral-local}). At an irreducible representation, the
surface relation is transverse to the group directions, so the representation
variety inherits the local structure of the puncture fibers
(Section~\ref{subsec:spectral-interior}). For regular semisimple classes the
puncture fibers are smooth, and the classical quasi-Hamiltonian argument gives
a finer statement (Section~\ref{subsec:regular-interior}).

Throughout, $L$ denotes a connected reductive group and
\[
 \pi_1(\Sigma^\circ)=\bigl\langle a_1,b_1,\dots,a_g,b_g,c_1,\dots,c_k\ \big|\
 \textstyle\prod_j[a_j,b_j]\prod_ic_i=1\bigr\rangle .
\]
A representation $\varphi$ of a group in $L$ is \emph{$L$-irreducible} if the
image of $\varphi$ is contained in no proper parabolic subgroup of $L$. For
$L=\GL_n$ this means that the matrices have no common invariant subspace other
than $0$ and $\C^n$.

\subsection{Irreducibility}\label{subsec:irreducibility}

\begin{lemma}\label{lem:schur}
Let $\varphi:\Pi\to L$ be a representation of a finitely generated group. If
$\varphi$ is $L$-irreducible, then the $L$-orbit of $\varphi$ in
$\Hom(\Pi,L)$ is closed, and the stabilizer of $\varphi$ in $L_\ad=L/Z(L)$ is
finite.
\end{lemma}

\begin{proof}
Apply the Hilbert--Mumford criterion (Fact~\ref{fact:git}(ii)) to the action of
$L_\ad$. A nontrivial cocharacter of $L_\ad$ lifts, after replacing the
cocharacter by a positive multiple, to a noncentral cocharacter $\lambda$ of
$L$. The limit $\lim_{t\to0}\lambda(t)\varphi\lambda(t)^{-1}$ exists if and only
if every matrix of the tuple $\varphi$ lies in $P(\lambda)$, a proper parabolic
subgroup. Irreducibility excludes this.
\end{proof}

Intrinsic genericity is a condition on the semisimple classes
$\cC^{\mathrm{ss}}(a_i)$ only. The next lemma shows that the condition
nevertheless controls all elements of the reduced fibers, whatever their
unipotent parts.

\begin{lemma}[Genericity depends only on semisimple parts]
\label{lem:spectral-genericity}
Let $\mathbf a\in(L/\!\!/L)^k$ and $S_i=(\chi_L^{-1}(a_i))_{\red}$. The
data $\mathbf a$ are intrinsically $L$-generic if and only if, for every proper
parabolic subgroup $Q$ and all $x_i\in S_i\cap Q$, one has
$x_1\cdots x_k\notin[Q,Q]$. Under this condition, every point of
$\cR_B^{\mathrm{sp}}(L,\mathbf a)$ is $L$-irreducible.
\end{lemma}

\begin{proof}
If $x\in Q$, then the Jordan decomposition of $x$ takes place inside $Q$, and
the elements $x$ and $x_s\in\cC^{\mathrm{ss}}(a_i)\cap Q$ have the same image in
the torus $Q/[Q,Q]$, because this torus kills unipotent elements. Hence the
products $x_1\cdots x_k$ and $(x_1)_s\cdots(x_k)_s$ have the same image in
$Q/[Q,Q]$, which proves one implication; the converse follows by taking the
$x_i$ semisimple. If the image of a representation lay in a proper parabolic
subgroup $Q$, the surface relation would give
$M_1\cdots M_k=(\prod_j[A_j,B_j])^{-1}\in[Q,Q]$ with $M_i\in S_i\cap Q$, a
contradiction.
\end{proof}

\begin{example}\label{ex:generic-sl2}
For $L=\SL_2$ the proper parabolic subgroups are the Borel subgroups, and an
element of a Borel subgroup lies in the commutator subgroup $U$ exactly when
its diagonal entries are $1$. If $\cC^{\mathrm{ss}}(a_i)$ has eigenvalues
$\mu_i^{\pm1}$, the data are therefore intrinsically generic if and only if
$\mu_1^{\pm1}\cdots\mu_k^{\pm1}\ne1$ for all choices of signs. For the
once-punctured torus ($k=1$), genericity says $\mu_1\ne1$, i.e.\ $\tr M\ne2$,
whereas regular semisimplicity says $\tr M\ne\pm2$. Thus for $k=1$ every
regular semisimple class is generic, and the value $\tr M=-2$, whose fiber
$-\mathcal U(\SL_2)$ is singular, is generic as well.
\end{example}

\subsection{Reduced adjoint fibers}\label{subsec:spectral-local}
In \eqref{eq:intro-gl2-fiber}, the fiber of $(t-\lambda)^2$ is
$\lambda\cdot\mathcal U(\GL_2)$, where $\mathcal U$ denotes the variety of
unipotent elements, and the logarithm identifies $\mathcal U(\GL_2)$ with the
nilpotent cone $\mathcal N(\mathfrak{sl}_2)$. In general the semisimple part of
an element of a fiber runs over a single conjugacy class, and the unipotent part
runs over the unipotent variety of the centralizer of the semisimple part.

\begin{lemma}\label{lem:spectral-fiber-local}
Let $a\in L/\!\!/L$ and $S=(\chi_L^{-1}(a))_{\red}$. Choose a semisimple
element $s$ with $\chi_L(s)=a$, and put $C=Z_L(s)$ and $\mathfrak c=\Lie C$.
Then $S$ consists of the elements whose semisimple part is conjugate to $s$,
and there is an isomorphism
\begin{equation}\label{eq:spectral-associated-bundle}
 L\times^C\mathcal U(C^\circ)\xrightarrow{\ \sim\ }S,
 \qquad [h,u]\longmapsto hsuh^{-1},
\end{equation}
where $\mathcal U(C^\circ)$ is the variety of unipotent elements of the
identity component $C^\circ$, and $L\times^CY$ denotes the quotient of
$L\times Y$ by $c\cdot(h,y)=(hc^{-1},cyc^{-1})$. Consequently the analytic germ
of $S$ at a point $hsuh^{-1}$ is the product of a smooth germ and the germ of
$\mathcal N([\mathfrak c,\mathfrak c])$ at $\log u$. In particular $S$ is an
irreducible normal local complete intersection with canonical singularities,
analytically locally $\Q$-factorial, of dimension $\dim L-\rank L$.
\end{lemma}

\begin{idea}
The map $x\mapsto x_s$ to the semisimple part is a polynomial in $x$ on the
fiber, so the fiber maps equivariantly onto the class $L\cdot s=L/C$; the fiber
of this map over $s$ is $s$ times the unipotent variety of $C^\circ$. The only
point requiring care is that these identifications hold for reduced schemes,
not only for sets of points.
\end{idea}

\begin{proof}
The adjoint quotient records the conjugacy class of the semisimple part: for
semisimple groups this is \cite[Thm.~6.11(d)]{Steinberg1965}, and the reductive
case follows through the central isogeny of Fact~\ref{fact:isogeny}.

\emph{The semisimple part is a regular function on $S$.} Embed $L$ faithfully in
$\GL(V)$, put $n=\dim V$, and let $\lambda_1,\dots,\lambda_b$ be the distinct
eigenvalues of $s$. By the Chinese remainder theorem there is a polynomial
$P(z)$ with $P(z)\equiv\lambda_j\bmod(z-\lambda_j)^n$ for $1\le j\le b$, and the
Jordan normal form shows $P(x)=x_s$ for every $x\in S$. Hence
$x\mapsto x_s$ is an $L$-equivariant morphism $S\to L\cdot s=L/C$.

\emph{The fiber over $s$.} Let $F$ be the scheme-theoretic fiber of $S\to L/C$
over $s$. The closed points of $F$ are the elements $su$ with $u$ unipotent and
commuting with $s$, i.e.\ $u\in\mathcal U(C^\circ)$; here we use that a unipotent
element lies in the identity component of every algebraic group containing
it. Equivariance gives an isomorphism $L\times F\simeq L\times_{L/C}S$. The
right-hand side is reduced, because its projection to the reduced scheme $S$ is
smooth; hence $F$ is reduced, and $F=s\,\mathcal U(C^\circ)$ as schemes. Descent
along $L\to L/C$ gives \eqref{eq:spectral-associated-bundle}; no connectedness
of $C$ is needed.

\emph{Local structure.} The logarithm, which is a polynomial on unipotent
elements, identifies $\mathcal U(C^\circ)$ with the nilpotent cone of the
semisimple Lie algebra $[\mathfrak c,\mathfrak c]$. The bundle $L\to L/C$ has
local analytic sections, so the germ of $S$ is a smooth germ times a germ of
this nilpotent cone, and the singularity properties follow from
Fact~\ref{fact:nilpotent-cones}. Finally $\rank C=\rank L$ and
$\dim\mathcal U(C^\circ)=\dim C-\rank C$, so
$\dim S=\dim L-\dim C+\dim C-\rank C=\dim L-\rank L$.
\end{proof}

For semisimple groups, the following description is
\cite[Thm.~6.11]{Steinberg1965}. The preceding local model gives a proof
for every connected reductive group.

\begin{corollary}[Reduced fibers as class closures]
\label{cor:adjoint-fiber-closure}
With the notation of Lemma~\ref{lem:spectral-fiber-local}, let
$u_{\mathrm{reg}}$ be a regular unipotent element of $C^\circ$. Then
\[
 S=\overline{L\cdot(su_{\mathrm{reg}})}
\]
as reduced closed subvarieties of $L$. The displayed class is the unique
regular class in $S$, where regular means
$\dim Z_L(x)=\rank L$. Thus spectral character varieties are the special
case of \eqref{paper:eq-relative} in which each prescribed class is regular.
\end{corollary}

\begin{proof}
Under logarithm, the unipotent variety of $C^\circ$ is the nilpotent
cone of $[\mathfrak c,\mathfrak c]$. By
Fact~\ref{fact:nilpotent-cones}(i), this cone is irreducible with
finitely many orbits of dimension at most
$\dim[\mathfrak c,\mathfrak c]-\rank[\mathfrak c,\mathfrak c]$.
The nilpotent cone therefore has a unique dense orbit. The inverse image of
this orbit under logarithm is the dense regular-unipotent class
$\mathcal O_{\mathrm{reg}}$ of $C^\circ$; the centralizer dimension is
$\rank C^\circ$.
Conjugation by $C$ preserves this class, so
$L\times^C\mathcal O_{\mathrm{reg}}$ maps to a single dense $L$-class
in the associated bundle \eqref{eq:spectral-associated-bundle}.
Since $S$ is reduced, $S$ is its reduced closure.

Jordan decomposition gives $Z_L(su)=Z_C(u)$, and
$\dim Z_C(u)=\dim Z_{C^\circ}(u)$.
A maximal torus of $L$ containing $s$ lies in $C^\circ$, so
$\rank C^\circ=\rank L$.
Consequently $su$ is regular in $L$ exactly when $u$ is regular in
$C^\circ$, proving uniqueness even when $C$ is disconnected.
\end{proof}

For $\GL_n$, the dense regular class can be written explicitly in Jordan
form. This identifies spectral character varieties with Letellier's
character varieties \cite[\S3.1]{Letellier2015}.

\begin{lemma}[Full spectral fibers in $\GL_n$]\label{lem:gl-spectral-closure}
Let $P(t)=\prod_\lambda(t-\lambda)^{m_\lambda}$ be monic of degree $n$ with
nonzero constant term, and let $C_P\subset\GL_n$ be the conjugacy class with one
Jordan block of size $m_\lambda$ for each eigenvalue $\lambda$. Then
\[
 \{M\in\GL_n:\det(tI-M)=P(t)\}_{\red}=\overline{C_P}.
\]
Thus the spectral character variety of $\GL_n$ is Letellier's character
variety with the closures $\overline{C_{P_i}}$ at the punctures.
\end{lemma}

\begin{proof}
The characteristic polynomial is constant on $\overline{C_P}$, which gives one
inclusion. Conversely, fix any Jordan matrix with characteristic polynomial
$P$. For each eigenvalue, put $1$ on the superdiagonal inside the given Jordan
blocks and a parameter $z$ at the breaks between consecutive blocks. For
$z\ne0$, every superdiagonal entry in the block of the eigenvalue is nonzero, so
the matrix lies in $C_P$; for $z=0$ it is the given Jordan matrix. Hence every
point of the fiber lies in $\overline{C_P}$. Both sides are reduced and have
the same points. Imposing the surface relation, reducing, and taking
invariants gives the statement about character varieties.
\end{proof}

\begin{remark}[Characteristic polynomials versus the adjoint quotient]
\label{rem:charpoly-versus-adjoint}
For $\GL_n$ and $\SL_n$, the characteristic polynomial in the standard
representation determines the point of $G/\!\!/G$. For a general group a
single faithful representation need not have this property. In the standard
representation of $G=\SO(q)\subset\GL_4$, with $q(v)=v_1v_3+v_2v_4$, the elements
$s=\diag(2,3,1/2,1/3)$ and $s'=\diag(1/2,3,2,1/3)$ have the same characteristic
polynomial. The orthogonal transformation interchanging $e_1$ and $e_3$
conjugates $s$ to $s'$ and has determinant $-1$; since the eigenvalues are
distinct, every orthogonal element commuting with $s$ is of the form
$\diag(a,b,a^{-1},b^{-1})$, of determinant $1$. Hence every orthogonal
conjugator from $s$ to $s'$ has determinant $-1$, the elements $s,s'$ are not
conjugate in $\SO(q)$, and they define distinct points of $G/\!\!/G$.
\end{remark}

\subsection{The surface relation}\label{subsec:spectral-interior}
For each $i$ let $S_i=(\chi_L^{-1}(a_i))_{\red}$, assume the data
$\mathbf a$ are compatible, and put
\[
 M=L^{2g}\times\textstyle\prod_iS_i,\qquad
 \mu=\prod_j[A_j,B_j]\prod_iM_i:M\to L_\der,\qquad R=\mu^{-1}(1),
\]
so that $R=\cR_B^{\mathrm{sp}}(L,\mathbf a)$ (compatibility guarantees that
$\mu$ takes values in $L_\der$). The puncture fibers may be singular, so $M$
need not be smooth. The key point is that, at an irreducible representation,
the group directions alone already make $\mu$ a submersion. We first compute
the differential of $\mu$ in these directions.

\begin{lemma}\label{lem:differential}
Let $x=(A_j,B_j,M_i)\in R$. Identify the tangent spaces of $L$ with
$\mathfrak l$ by right translation, and $\mathfrak l$ with $\mathfrak l^\vee$ by
a nondegenerate invariant symmetric bilinear form. Consider the variations of $x$ obtained by multiplying each
$A_j$ and each $B_j$ on the left by an element of $L$ near $1$, and conjugating
each $M_i$ by an element of $L$ near $1$. Then the annihilator of the set of
derivatives of $\mu$ along these variations is the centralizer
$\mathfrak z_{\mathfrak l}(A_j,B_j,M_i)$ of the tuple $x$.
\end{lemma}

\begin{proof}
Write $d^Rg=(dg)g^{-1}$, so that $d^R(gh)=d^Rg+\Ad_g(d^Rh)$. For one handle,
with $C=[A,B]$ and $K=ABA^{-1}$, a direct computation gives
\[
 d^RC=(1-\Ad_K)d^RA+(\Ad_A-\Ad_C)d^RB,
\]
with $d^RA,d^RB\in\mathfrak l$ arbitrary; conjugating a factor $M$ by
$\exp(tX)$ gives $d^RM=(1-\Ad_M)X$, with $X\in\mathfrak l$ arbitrary. Suppose
that $\xi$ annihilates all these derivatives of $\mu$, and let
$\zeta=\Ad_{P^{-1}}\xi$, where $P$ is the product of the factors of the surface
word preceding a given handle. Annihilating the contribution of that handle
means $\Ad_K\zeta=\zeta$ and $\Ad_{A^{-1}}\zeta=\Ad_B\zeta$, which together
imply that $AB$, then $A$, then $B$ fix $\zeta$; in particular $C$ fixes
$\zeta$, so $\zeta$ is also the transported vector for the next factor. For a
puncture factor the condition is that $M$ fixes $\zeta$. By induction along
the word, every holonomy fixes $\xi$. The converse is immediate.
\end{proof}

\begin{lemma}[A local product for the surface relation]
\label{lem:spectral-relation-local}
Let $x\in R$ be $L$-irreducible. Then there are a smooth variety $Z$ of
dimension $\dim L_\der$, a point $e\in Z$, and a morphism $Z\times R\to M$
which is \'etale at $(e,x)$ and maps $(e,x)$ to $x$. In particular there is an
isomorphism of analytic germs
\begin{equation}\label{eq:spectral-local-product}
 (M,x)\simeq(R,x)\times(\C^{\dim L_\der},0).
\end{equation}
\end{lemma}

\begin{idea}
If $M$ were smooth, this would be the implicit function theorem for the
submersion $\mu$. Since $M$ may be singular, we use the group action instead:
moving $x$ by group elements realizes all directions transverse to $R$, and
these moves exist on $M$ whether or not $M$ is smooth.
\end{idea}

\begin{proof}
The group $K=L^{2g}\times L^k$ acts on $M$ by left translation of the handle
entries and conjugation of the puncture entries. By Lemma~\ref{lem:differential},
the differential at the identity of $h\mapsto\mu(h\cdot x)$ has image with
annihilator $\mathfrak z_{\mathfrak l}(x)$, which is $\mathfrak z(\mathfrak l)$
because $x$ is irreducible (Lemma~\ref{lem:schur}). Hence this differential
maps onto $\mathfrak l_\der$. Choose a smooth locally closed subvariety
$Z\subset K$ through the identity $e$, of dimension $\dim L_\der$, on which the
differential is an isomorphism, and consider
\[
 \Phi:Z\times M\longrightarrow L_\der\times M,\qquad
 (q,z)\longmapsto(\mu(q\cdot z),z).
\]
Both sides are smooth over $M$ of the same relative dimension, and the relative
differential is invertible at $(e,x)$; hence $\Phi$ is \'etale at $(e,x)$, a
criterion which does not require $M$ to be smooth. Base change along
$\{1\}\times M$ shows that $E=\{(q,z):\mu(q\cdot z)=1\}\to M$ is \'etale near
$(e,x)$. The map $(q,z)\mapsto(q,q\cdot z)$ identifies $E$ with $Z\times R$, with
inverse $(q,r)\mapsto(q,q^{-1}\cdot r)$, so the morphism
$Z\times R\to M$, $(q,r)\mapsto q^{-1}\cdot r$, is \'etale at $(e,x)$. An \'etale
morphism over $\C$ induces isomorphisms of analytic germs, which proves
\eqref{eq:spectral-local-product}.
\end{proof}

\begin{proposition}[The generic spectral interior]
\label{prop:spectral-interior}
Let $\mathbf a\in(L/\!\!/L)^k$ be compatible and intrinsically $L$-generic, and
suppose $R=\cR_B^{\mathrm{sp}}(L,\mathbf a)$ is nonempty. Then:
\begin{enumerate}[label=\textup{(\roman*)}]
\item $R$ is reduced, normal, analytically locally $\Q$-factorial, with
canonical singularities;
\item every $L$-orbit in $R$ is closed, with finite stabilizer in $L_\ad$;
\item the spectral character variety $\cM_B^{\mathrm{sp}}(L,\mathbf a)=R/\!\!/L$
is normal, $\Q$-factorial and klt, of dimension
\[
 (2g+k-2)\dim L+2\dim Z(L)-k\rank L.
\]
\end{enumerate}
\end{proposition}

Here \emph{klt} (Kawamata log terminal) means that all discrepancies are
greater than $-1$ (Section~\ref{subsec:bg-pairs}). We do not claim canonical
singularities at this stage: a finite quotient of a canonical singularity need
not be canonical. The integral relation $K_X+D\sim0$ will upgrade klt to
canonical in the proof of Theorem~\ref{thm:spectral-main}.

\begin{proof}
(ii) Every point of $R$ is irreducible (Lemma~\ref{lem:spectral-genericity}),
so Lemma~\ref{lem:schur} applies.

(i) By Lemma~\ref{lem:spectral-fiber-local}, the germs of $M$ are products of
nilpotent-cone germs and smooth germs, so they are reduced, normal, canonical
and analytically $\Q$-factorial (Fact~\ref{fact:nilpotent-cones}). By
\eqref{eq:spectral-local-product}, the germ of $R$ times a smooth germ has the
same properties, hence so does the germ of $R$; for $\Q$-factoriality, a Weil
divisor on the germ of $R$ pulls back to the product, a multiple of the
pullback is Cartier, and the restriction of this multiple along
$\{0\}\times(R,x)$ is a Cartier multiple of the original divisor.

(iii) By (ii) and the Luna slice theorem (Fact~\ref{fact:git}(iii)), every point
of $R/\!\!/L$ has an \'etale neighborhood of the form $S/F$, where $F\subset
L_\ad$ is finite and $S$ is a slice, possibly singular, such that the germ of $R$
is a smooth germ times the germ of $S$. Thus $S$ is normal, canonical and
analytically $\Q$-factorial. The quotient $S/F$ is normal. A finite quotient of
a $\Q$-factorial germ is $\Q$-factorial: pull back a Weil divisor, take a
Cartier multiple, and take its norm; in particular $S/F$ is $\Q$-Gorenstein.
For the quotient map $q:S\to S/F$, the ramification formula reads
\[
 K_S=q^*(K_{S/F}+B),\qquad B=\textstyle\sum_P(1-1/e_P)P\ge0,
\]
the sum running over the prime divisors $P$ of $S/F$ with ramification index
$e_P$. Since $S$ is canonical, the pair $(S,0)$ is klt, so $(S/F,B)$ is klt by
Fact~\ref{fact:pairs}(ii), and $S/F$ is klt because $B\ge0$. Analytic local
$\Q$-factoriality implies $\Q$-factoriality, by faithfully flat descent of the
reflexive powers of rank-one reflexive sheaves. Finally, $\dim M=2g\dim
L+k(\dim L-\rank L)$; forming $R$ subtracts $\dim L_\der$, and forming the
quotient subtracts $\dim L_\ad=\dim L_\der$, which gives the formula.
\end{proof}

\subsection{Regular semisimple classes}\label{subsec:regular-interior}
When all prescribed classes are regular semisimple, the puncture fibers are the
classes themselves, which are smooth, and the quasi-Hamiltonian formalism
yields symplectic local models.

\begin{proposition}\label{prop:generic-geometry}
Let $\bcC$ be a compatible and intrinsically $L$-generic tuple of semisimple
classes with $\cR_B(L,\bcC)\ne\varnothing$. Then every point of $\cR_B(L,\bcC)$ is
$L$-irreducible, every orbit is closed with finite stabilizer in $L_\ad$, and
$\cR_B(L,\bcC)$ is smooth of dimension $2g\dim L+\sum_i\dim\cC_i-\dim L_\der$.
Consequently
\[
 \dim\cM_B(L,\bcC)=(2g-2)\dim L+2\dim Z(L)+\textstyle\sum_i\dim\cC_i .
\]
\end{proposition}

\begin{proof}
Irreducibility follows from \eqref{eq:intro-generic} exactly as in the proof
of Lemma~\ref{lem:spectral-genericity}, and Lemma~\ref{lem:schur} gives the
orbit assertions. Since the classes are smooth and, by
Lemma~\ref{lem:differential} and irreducibility, the differential of
$\mu:L^{2g}\times\prod_i\cC_i\to L_\der$ is surjective at every point of
$\mu^{-1}(1)$, the scheme $\cR_B(L,\bcC)$ is smooth. The dimension formulas follow
by subtracting $\dim L_\der$ and then $\dim L_\ad$.
\end{proof}

\begin{corollary}\label{cor:symplectic-coarse}
Let $\bcC$ be semisimple classes of $L$ such that the surface word maps
$L^{2g}\times\prod_i\cC_i$ into $L_\der$ and every point of the nonempty scheme
$\cR_B(L,\bcC)$ is $L$-irreducible. Then $\cM_B(L,\bcC)$ is \'etale locally the
quotient of a smooth symplectic variety by a finite group preserving the
symplectic form. Every connected component of $\cM_B(L,\bcC)$ is a normal,
$\Q$-factorial, Gorenstein variety with canonical singularities and trivial
canonical sheaf.
\end{corollary}

\begin{proof}
As in Proposition~\ref{prop:generic-geometry}, the scheme $R=\cR_B(L,\bcC)$ is
smooth, and every orbit is closed with finite stabilizer in $L_\ad$. By
Fact~\ref{fact:qh}, $L^{2g}\times\prod_i\cC_i$ is quasi-Hamiltonian with moment
map the surface word, so by Lemma~\ref{lem:qh-reduction} the restriction
$\omega_R$ of the two-form to $R$ is closed and has kernel equal to the orbit
directions. By the Luna slice theorem (Fact~\ref{fact:git}(iii)), every point
of $\cM_B(L,\bcC)$ has an \'etale neighborhood of the form $S/F$, where
$S\subset R$ is a smooth slice transverse to an orbit and $F$ is a finite
stabilizer. The restriction of $\omega_R$ to $S$ is an $F$-invariant
symplectic form, and Fact~\ref{fact:symplectic-quotient} gives all the local
assertions. The local generators of the canonical sheaf obtained in this way
are restrictions of the top power of the global form $\omega_R$, so they glue
to a global trivialization.
\end{proof}

\begin{remark}
For closed surfaces of genus at least two, Bellamy--Schedler prove rational
Gorenstein singularities and $\Q$-factoriality for the identity component when
all simple quotients of the connected reductive group have type~$A$
\cite[Cor.~1.2, Thm.~1.3]{BellamySchedler2023}.
Herbig--Schwarz--Seaton establish symplectic singularities for connected
reductive groups \cite[Thms.~7.1, 7.4]{HerbigSchwarzSeaton2024};
Shu additionally proves $\Q$-factoriality for the full character variety
\cite[Thm.~A]{ShuSingularities2025}.
\end{remark}

\subsection{From log canonical to dlt}
Canonical singularities of the interior allow us to upgrade a log canonical
compactification to a dlt one without changing the interior.

\begin{corollary}\label{cor:dlt-open}
Let $U$ be a normal $\Q$-factorial variety with canonical singularities, and let
$(X,D)$ be a projective lc log Calabi--Yau compactification of $U$. Then $U$ has
a projective $\Q$-factorial dlt log Calabi--Yau compactification $(X',D')$ with
a birational morphism $X'\to X$ which is an isomorphism over $U$.
\end{corollary}

\begin{proof}
Let $f:X'\to X$ be a dlt modification (Fact~\ref{fact:pairs}(i)). Every
$f$-exceptional divisor $E$ has discrepancy $a(E)=-1$, while every divisor over
$U$ has discrepancy $\ge0$ because $U$ is canonical. Hence no exceptional
divisor maps into $U$, and $f^{-1}(U)\to U$ is a small projective birational
morphism. On a $\Q$-factorial variety such a morphism is an isomorphism: if
$A'$ is $f$-ample on $f^{-1}(U)$, a multiple of $f_*A'$ is Cartier, and
$f^*f_*A'-A'$ is supported on the (empty) exceptional divisor, so $A'$ is
$f$-trivial on every fiber, and the fibers are finite. Finally
$K_{X'}+D'=f^*(K_X+D)\sim0$.
\end{proof}

\section{Prescribed conjugacy closures and smooth boundary evaluation}
\label{new:sec-base-change}

This section and Section~\ref{new:sec-normalization} study the interior
$U_{\mathbf Y}$ for prescribed class closures $Y_i=\overline{\cC_i}$. The new
feature, compared with the spectral case, is that a class closure can be
non-normal, whereas a full reduced adjoint fiber is always normal
(Lemma~\ref{lem:spectral-fiber-local}). For example, in $\SL_2$ the closure of
the class of a nontrivial unipotent element is the whole unipotent variety, a
normal quadric cone; but in type $G_2$ the closure of the eight-dimensional
nilpotent orbit is not normal (Fact~\ref{new:fact-g2-orbit}). This section
describes how the singularities of the $Y_i$ enter the relation variety
$R_{\mathbf Y}$ and the finite stabilizer quotients.
Section~\ref{new:sec-normalization} constructs the normalization, and
Section~\ref{r2:norm-section} decides when normalization is needed. A
compactification of the given reduced variety means an open immersion of that
variety; compactifying its normalization is a separate assertion.

\subsection{Conjugacy-class closures}
\label{new:subsec-question}
Let $G$ be connected complex reductive, let $g\geq0$, $k\geq1$, and let
$\cC_1,\ldots,\cC_k$ be conjugacy classes.  Put
\[
 Y_i=\overline{\cC_i}_{\red},\qquad U_{\mathbf Y}=(R_{\mathbf Y})_{\red}/\!/G,
\]
\[
 R_{\mathbf Y}=\left\{(A_j,B_j,M_i)\in
 G^{2g}\times\prod_iY_i:\prod_j[A_j,B_j]\prod_iM_i=1\right\}.
\]
Write $s_i$ for the semisimple part of an element of $\cC_i$ and
$S_i=(\chi_G^{-1}(\chi_G(s_i)))_{\red}$.  Compatibility and intrinsic
genericity are the conditions of \hyperref[sec:intro]{the Introduction} on the
semisimple classes $G\cdot s_i$.

\subsection{The representation square and its quotient-stack form}
\label{new:subsec-square}
There is a scheme-theoretic Cartesian square
\begin{equation}\label{new:eq-representation-square}
\begin{array}{ccc}
 R_{\mathbf Y}&\longrightarrow&R_{\mathbf S}\\
 \downarrow&&\downarrow\\
 \prod_iY_i&\longrightarrow&\prod_iS_i,
\end{array}
\end{equation}
where the vertical maps record the puncture matrices.  On quotient
stacks this becomes evaluation in the conjugacy stacks.

Set
\[
 G_{\ad}=G/Z(G),\qquad d=\dim G,\qquad z=\dim Z(G),\qquad h=d-z.
\]
The actions by conjugation factor through $G_{\ad}$.  For a $G_{\ad}$-variety $V$,
$[V/G_{\ad}]$ denotes the quotient stack, which retains stabilizers as well as
orbits.  Let $R_{\mathbf S}^{\mathrm{irr}}$ denote the open locus of
$G$-irreducible tuples.  The boundary evaluation is the morphism
\begin{equation}\label{new:eq-boundary-evaluation}
 f:[R_{\mathbf S}^{\mathrm{irr}}/G_{\ad}]
       \longrightarrow\prod_i[S_i/G_{\ad}].
\end{equation}
The source uses simultaneous conjugation, and the target allows an
independent conjugation at each puncture.

\begin{theorem}[Smooth boundary evaluation]\label{new:thm-smooth-evaluation}
More generally, let the $S_i\subset G$ be any reduced invariant affine
subvarieties with constant compatible images in $G_{\ab}$.  The morphism
\eqref{new:eq-boundary-evaluation} is representable and smooth, of relative
dimension
\begin{equation}\label{new:eq-evaluation-dimension}
 q=2gd+(k-2)h.
\end{equation}
No normality or smoothness of the puncture varieties is required.
If $Y_i\subset S_i$ are invariant closed reduced subvarieties, then
\[
\begin{array}{ccc}
 [R_{\mathbf Y}^{\mathrm{irr}}/G_{\ad}]&\longrightarrow&
 [R_{\mathbf S}^{\mathrm{irr}}/G_{\ad}]\\
 \downarrow&&\downarrow f\\
 \prod_i[Y_i/G_{\ad}]&\longrightarrow&\prod_i[S_i/G_{\ad}]
\end{array}
\]
is Cartesian, and both vertical arrows are smooth.
\end{theorem}

Theorem~\ref{new:thm-smooth-evaluation} says that, at irreducible
representations, the puncture matrices can be deformed independently and
without obstruction. In the running example ($G=\SL_2$, $(g,k)=(1,1)$) one has
$d=h=3$ and $q=6-3=3$: the source of $f$ is the two-dimensional cubic surface
regarded as a stack, and the target $[S_1/\PGL_2]$ has dimension $2-3=-1$.

\begin{proof}
Pull back the atlas $B=\prod_iS_i\to\prod_i[S_i/G_{\ad}]$.  The fiber product of this atlas
with the source of $f$ is $Q/G_{\ad}$, where $Q$ is the open irreducible-tuple
locus of
\[
 \left\{(A_j,B_j,a_i,b_i)\in G^{2g}\times G_{\ad}^k\times B:
 \prod_j[A_j,B_j]\prod_i a_i b_i a_i^{-1}=1\right\}.
\]
Conjugation by an element of $G_{\ad}$ on $G$ is well defined.  The diagonal
$G_{\ad}$-action conjugates the handles and multiplies each $a_i$ on the left;
it fixes the $b_i$.  It is free, and $k\geq1$ allows the unique gauge
$a_1=1$.  Thus $Q/G_{\ad}$ is a scheme and the asserted map is representable.

The ambient scheme $G^{2g}\times G_{\ad}^k\times B$ is smooth over $B$.
Compatibility makes the displayed word take values in the smooth group
$G_{\der}$ of dimension $h$.  At a tuple over $1$, vary the handles and
the $a_i$, holding the $b_i$ fixed.  The differential computation in
Lemma~\ref{lem:differential} says that the annihilator of these word
derivatives in $\mathfrak g$ is the common infinitesimal centralizer.
Irreducibility identifies this centralizer with $\mathfrak z(\mathfrak g)$,
so the derivatives surject onto $\mathfrak g_{\der}$.
The relative Jacobian criterion, applied to smooth schemes over the
possibly singular base $B$, therefore proves that $Q\to B$ is smooth
of relative dimension $2gd+(k-1)h$.  The free quotient removes $h$
dimensions, giving \eqref{new:eq-evaluation-dimension}.
Smoothness can be checked after the atlas $B\to\prod_i[S_i/G_{\ad}]$.
The Cartesian assertion follows either from the moduli description or
from the same presentation with $b_i\in Y_i$; smoothness is preserved
by base change.
\end{proof}

\subsection{Local singularities and orbit strata}
\label{new:subsec-local-transfer}
For the rest of this section, assume the data are compatible and
intrinsically generic.  Every tuple is irreducible by
Lemma~\ref{lem:spectral-genericity}.

\begin{proposition}[Stratum-preserving local comparison]
\label{new:prop-local-transfer}
Let $M=G^{2g}\times\prod_iY_i$ and $r\in R_{\mathbf Y}$.
There is a smooth germ $(Z,e)$ of dimension $h$ and an \'etale morphism
\[
 (Z\times R_{\mathbf Y},(e,r))\longrightarrow(M,r).
\]
This morphism preserves the condition that each puncture belongs to any
specified conjugacy orbit.  Consequently:
\begin{enumerate}[label=\textup{(\roman*)}]
\item $R_{\mathbf Y}$ is reduced, so the reduction in its definition is
unnecessary under these hypotheses;
\item $R_{\mathbf Y}$ is normal at $r$ if and only if each $Y_i$ is normal
at $M_i$;
\item $R_{\mathbf Y}$ is pure of dimension
$2gd+\sum_i\dim Y_i-h$, and $U_{\mathbf Y}$ is pure of dimension
\begin{equation}\label{new:eq-closure-dimension}
 \dim U_{\mathbf Y}=2gd+\sum_i\dim Y_i-2h;
\end{equation}
\item the exact-class locus, where $M_i\in\cC_i$ for every $i$, is dense
in every component.  The complement of the exact-class locus has codimension at least two in
both $R_{\mathbf Y}$ and $U_{\mathbf Y}$.
\end{enumerate}
At a point with trivial stabilizer in $G_{\ad}$, the normality equivalence
in \textup{(ii)} holds for $U_{\mathbf Y}$ as well.
\end{proposition}

In words: after taking smooth factors, the germ of $R_{\mathbf Y}$ at $r$
is \'etale-locally isomorphic to the product of the puncture germs $(Y_i,M_i)$.
Thus local properties preserved and detected by \'etale maps and smooth
factors pass between these germs. For
example, if $G=\SL_2$, $(g,k)=(1,1)$ and $\cC_1$ is the class of
$-\begin{psmallmatrix}1&1\\0&1\end{psmallmatrix}$, then
$Y_1=-I\cdot\mathcal U(\SL_2)$ is a quadric cone, and $U_{\mathbf Y}$ is the
cubic surface $U_{-2}$ of \hyperref[subsec:intro-example]{the introductory model calculation}, which has a node.

\begin{proof}
Use the action of $G^{2g}\times G^k$ by left translation on handles and
independent conjugation on punctures.  The proof of
Lemma~\ref{lem:spectral-relation-local} only uses invariance of the
puncture conditions and the word differential.  It therefore gives
the stated \'etale map without alteration.  Independent conjugation
preserves each orbit stratum.

Reducedness, normality and dimension are preserved and detected under
this local \'etale comparison and a smooth factor.  A product of reduced
complex varieties is reduced, and a product germ is normal exactly when
each factor germ is normal.  This proves the representation assertions.
All $G_{\ad}$-orbits are closed and have dimension $h$, so the good quotient
is geometric and subtracts $h$ from each invariant stratum's dimension.

Each $Y_i$ contains finitely many conjugacy orbits: fixing the
semisimple part reduces this to nilpotent orbits in its reductive
centralizer, with its finite component-group action.  Orbit dimensions
are even, so each orbit in
$Y_i\setminus\cC_i$ therefore has codimension at least two in $Y_i$.
The stratum-preserving comparison gives the same dimension drop after
the relation and quotient.  Since the full spaces are pure of the
dimensions just computed, these strata cannot contain any component.
This proves density of the exact-class locus.

Finally, at a trivial effective stabilizer Luna's slice gives an
\'etale-local principal $G_{\ad}$-bundle from the representation space to
its quotient, so normality is detected on the quotient germ as well.
For a nontrivial finite stabilizer the quotient of the slice must be
retained, as in Section~\ref{r2:norm-section}.
\end{proof}

\subsection{What coarse quotients do and do not preserve}
The invariant ideal of $R_{\mathbf Y}\subset R_{\mathbf S}$ induces a
closed immersion $U_{\mathbf Y}\hookrightarrow U_{\mathbf S}$ by
exactness of reductive invariants.  The stack square need not remain
Cartesian on coarse spaces.

For example, let $\mu_2$ act on $\A^1_x$ by $x\mapsto-x$ and restrict
to its invariant reduced point $\{0\}$.  The quotient coordinate is
$t=x^2$.  Pulling back the quotient point $t=0$ gives
$\Spec\C[x]/(x^2)$, rather than the reduced point.
Thus the scheme-theoretic preimage of a reduced subvariety under a quotient map
need not be reduced. This is why $U_{\mathbf Y}$ is defined from the reduced
scheme $(R_{\mathbf Y})_{\mathrm{red}}$, and why normality of $U_{\mathbf Y}$
cannot be read off from $U_{\mathbf S}$.

\section{Normalization and symplectic singularities}
\label{new:sec-normalization}

Theorem~\ref{paper:thm-normalized} concerns the normalization
$U_{\mathbf Y}^\nu$, which in positive genus differs from $U_{\mathbf Y}$ exactly
when some puncture closure is not normal (Theorem~\ref{paper:thm-normality}).
The idea is to normalize the punctures first. The normalization $Y_i^\nu$ of a
class closure has rational Gorenstein singularities
(Lemma~\ref{new:lem-normalized-puncture}); the relation variety built from the
$Y_i^\nu$ inherits these properties, and its quotient is the normalization of
$U_{\mathbf Y}$. The Goldman form then extends to a symplectic form on the
regular locus, and the normalization has symplectic singularities.
Throughout, the data are compatible and intrinsically generic as in
Section~\ref{new:sec-base-change}.

\subsection{The normalization of a conjugacy-class closure}

\begin{lemma}\label{new:lem-normalized-puncture}
Let $Y$ be a reduced conjugacy-class closure in a connected complex
reductive group $G$, and let $\nu:Y^\nu\to Y$ be its normalization.
The conjugation action lifts uniquely to $Y^\nu$.  The variety $Y^\nu$
has rational Gorenstein, hence canonical, singularities.  The complement in $Y^\nu$
of the preimage of the open conjugacy class has codimension at least two.
\end{lemma}

\begin{proof}
The action lifts by the universal property applied to the normal
variety $G\times Y^\nu$; the action laws hold by density.

Fix the semisimple part $s$ and set $C=Z_G(s)$, which is allowed to be
disconnected.  On the fixed adjoint fiber, taking the semisimple part is
a regular morphism to $G\cdot s\simeq G/C$.  One way to see regularity
is to choose a faithful representation: on a fiber with fixed
eigenvalues, Hermite interpolation gives one polynomial in the matrix
whose value is its semisimple part.  This construction agrees with the
intrinsic Jordan decomposition and factors through $G\cdot s$.

The fiber in $Y$ over $s$ is $sZ$, where $Z$ is a closed union of
unipotent orbits in $C^\circ$, stable under $C$; equivalently, polynomial
logarithm identifies it with a finite union of nilpotent orbit closures
in $\Lie C^\circ$.  The irreducible components of $Z$ are permuted by
$C/C^\circ$.  Associated-bundle descent gives
\[
 Y\simeq G\times^C sZ,\qquad
 Y^\nu\simeq G\times^C sZ^\nu.
\]
Here $Z^\nu$ normalizes and separates the components of the
\emph{whole reduced union}, retaining the $C/C^\circ$-action.
The formulas are product statements after the smooth cover $G\to G/C$.

The normalization of a nilpotent orbit closure is a symplectic
singularity: the Kostant--Kirillov form extends by
\cite[Lem.~2]{Panyushev1991}; see also \cite[\S2.6]{Beauville2000}.  In particular it has rational
Gorenstein singularities \cite[Prop.~1.3]{Beauville2000}.
Smooth ascent and descent give these properties on $Y^\nu$; rational
Gorenstein singularities are canonical.  The other conjugacy orbits
have even, strictly smaller dimension, and finite normalization
preserves the resulting codimension-two bound.
\end{proof}

For $G=\GL_n$ every class closure is normal, since the closure is locally a
product of a smooth factor and nilpotent orbit closures in $\mathfrak{gl}_m$,
which are normal by Kraft--Procesi \cite{KraftProcesi1979}; then $Y^\nu=Y$. For
other groups class closures can be non-normal: an example in type $G_2$ is
Fact~\ref{new:fact-g2-orbit}, and Proposition~\ref{new:gl-pgl4-nonnormal} gives one in
$\PGL_4$, where a centralizer is disconnected.

\subsection{A concrete normalization formula}

Let $\nu_i:Y_i^\nu\to Y_i$ be the
puncture normalizations and define
\begin{equation}\label{new:eq-normalized-relation}
 \widehat R=
 \left\{(A_j,B_j,\widetilde M_i)\in
 G^{2g}\times\prod_iY_i^\nu:
 \prod_j[A_j,B_j]\prod_i\nu_i(\widetilde M_i)=1\right\},
 \qquad \widehat U=\widehat R/\!/G.
\end{equation}
The relation in this definition is scheme-theoretic.

\begin{theorem}[Normalization commutes with the generic relation and quotient]
\label{new:thm-normalization-commutes}
The finite maps induced by the $\nu_i$ identify
\[
 \widehat R=R_{\mathbf Y}^{\nu},\qquad
 \widehat U=U_{\mathbf Y}^{\nu}
\]
component by component.  The scheme $\widehat R$ is reduced and normal
and has rational Gorenstein canonical singularities.  The exact-class
locus is dense in every component, with complement of codimension at
least two.
\end{theorem}

\begin{proof}
The lifted conjugation actions allow exactly the group-direction proof
of Proposition~\ref{new:prop-local-transfer}: the derivatives of the
surface word in these directions are the derivatives of the underlying
tuple in $G$.  Genericity makes them surject onto $\mathfrak g_{\der}$.
Thus a germ of $G^{2g}\times\prod_iY_i^\nu$ is \`etale-locally a germ
of $\widehat R$ times a smooth factor of dimension $h$.
Lemma~\ref{new:lem-normalized-puncture} implies all the stated
singularity properties for $\widehat R$.  The finite map from $\widehat R$ to
$R_{\mathbf Y}$ is an isomorphism over the exact-class locus.  The same
stratum-preserving argument proves that the complementary locus has
codimension at least two and that the exact-class locus is dense in
every component of $\widehat R$ and $R_{\mathbf Y}$.  The finite map is
therefore birational on components and is the normalization.

A finite equivariant morphism of affine varieties induces a finite
morphism of affine reductive quotients.  The invariant ring of a normal
ring is normal, so $\widehat U$ is normal.  Since every orbit is closed,
the image of the invariant exact-class open set is open and the induced
map is an isomorphism there.  The complement of this open set has codimension at least
two, as all effective orbits have the same dimension $h$.
Consequently $\widehat U\to U_{\mathbf Y}$ is finite and birational on
each component.  This is the defining characterization of normalization.
\end{proof}

\subsection{The symplectic form on the normalization}

\begin{theorem}[The generic normalized interior]
\label{new:thm-normalized-interior}
Every connected component $V$ of $U_{\mathbf Y}^{\nu}$ is a symplectic
singularity.  In particular $V$ is rational Gorenstein and canonical,
and its canonical sheaf has an integral trivialization
\[
 \omega_V\simeq\mathcal O_V,\qquad K_V\sim0.
\]
If every $Y_i$ is normal, these conclusions hold for the original
$U_{\mathbf Y}$.  No assertion of $\Q$-factoriality is made.
\end{theorem}

\begin{proof}
At a closed orbit of $\widehat R$, its stabilizer
$F\subset G_{\ad}$ is finite.  Luna's slice theorem gives an \`etale-local
description by a finite quotient $S/F$.  The slice $S$ is normal and
has rational Gorenstein singularities, since locally
$G_{\ad}\times S\to\widehat R$ is \`etale and
Theorem~\ref{new:thm-normalization-commutes} applies.  We may work in an
$F$-invariant affine neighborhood of the point on $S$.

On the exact-class open part $S^\circ$, the quasi-Hamiltonian reduction
form is a regular nondegenerate closed two-form $\omega$;
Lemma~\ref{lem:qh-reduction} applies to the product of smooth conjugacy
classes.  The construction makes this form $F$-invariant.  The slice there is
smooth, and $S\setminus S^\circ$ has codimension at least two by
Theorem~\ref{new:thm-normalization-commutes}.  If $\dim S=2n$, the form
\[
 \Omega=\frac{\omega^{\wedge n}}{n!}
\]
trivializes $\omega_S$ on $S^\circ$.  Normality and invertibility of
$\omega_S$ extend both $\Omega$ and its inverse, giving a
nowhere-vanishing generator on $S$.

On $S_{\mathrm{reg}}$, the coefficients of $\omega$ extend by Hartogs
across $S_{\mathrm{reg}}\setminus S^\circ$; differential forms are
locally free on that smooth locus.  The extension remains closed and
its top wedge equals $n!\Omega$, hence it is nondegenerate everywhere
on $S_{\mathrm{reg}}$.  Namikawa's criterion
\cite[Thm.~6]{Namikawa2001} now applies: a normal variety with rational
Gorenstein singularities and a symplectic form on its regular locus is
a symplectic singularity.  The form is still $F$-invariant, by density.
Beauville's finite symplectic quotient theorem
\cite[Prop.~2.4]{Beauville2000} therefore makes $S/F$ a symplectic
singularity.  This property is local in the \`etale topology, proving
the assertion for $V$.

The local extensions agree on overlaps by their equality with the
Goldman form on the dense exact-class locus.  Their global top wedge
trivializes the reflexive canonical sheaf, which is invertible by the
Gorenstein property.  When all $Y_i$ are normal,
Theorem~\ref{new:thm-normalization-commutes} identifies $V$ with a
component of the original variety.
\end{proof}

\begin{remark}
For $\GL_n$, Schedler and Tirelli realized character varieties with
prescribed class closures as open subsets of multiplicative quiver varieties
and proved, under hypotheses on the classes, that their normalizations are
symplectic singularities \cite{SchedlerTirelli2022}.
Theorem~\ref{new:thm-normalized-interior} gives the analogous statement for every
reductive group under intrinsic genericity.
\end{remark}

\begin{remark}
Theorem~\ref{new:sl3-main} shows that these normal symplectic interiors
need not be $\Q$-factorial.  Controlling the poles of the volume form
at infinity requires the Poisson compactification argument in Part~\ref{part:poisson}.
\end{remark}

\part{The spectral construction}\label{part:spectral}

\section{A toric degeneration from an ideal triangulation}\label{sec:graph}

Throughout Sections~\ref{sec:graph}--\ref{sec:special}, $H$ is a simply
connected semisimple group of rank $r$ with $N$ positive roots, and we assume
$2g+k\ge3$. Put $m=2g+k-1$ and
\[
 A_0=\C[\Hom(F_m,H)]^H=\C[H^m]^H,
\]
the coordinate ring of the free-group character variety. We construct its
semigroup valuation (Theorem~\ref{thm:graph-valuation}) and compute the
initial forms of the peripheral equations
(Proposition~\ref{prop:peripheral-symbols}).

\subsection{Ribbon graphs}\label{subsec:ribbon}
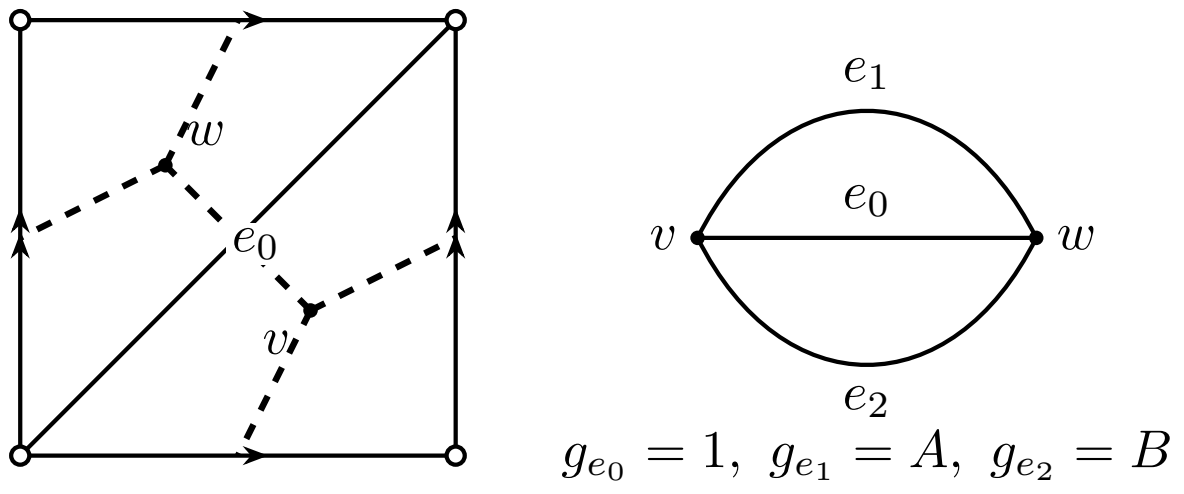
\begin{figure}[htbp]
\centering
\begin{tikzpicture}[scale=0.8,>=Stealth,font=\small]
 \draw[thick] (0,0) rectangle (3.6,3.6);
 \draw[thick] (0,0)--(3.6,3.6);
 \draw[->] (1.45,0)--(2.05,0);\draw[->] (1.45,3.6)--(2.05,3.6);
 \draw[->>] (0,1.45)--(0,2.05);\draw[->>] (3.6,1.45)--(3.6,2.05);
 \draw[dashed,very thick] (2.4,1.2)--(1.2,2.4);
 \draw[dashed,very thick] (2.4,1.2)--(1.8,0);
 \draw[dashed,very thick] (1.2,2.4)--(1.8,3.6);
 \draw[dashed,very thick] (2.4,1.2)--(3.6,1.8);
 \draw[dashed,very thick] (1.2,2.4)--(0,1.8);
 \foreach \x/\y in {0/0,3.6/0,0/3.6,3.6/3.6}
   \draw[fill=white,thick] (\x,\y) circle (0.075);
 \fill (2.4,1.2) circle (0.06) node[below left] {$v$};
 \fill (1.2,2.4) circle (0.06) node[above right] {$w$};
 \node[fill=white,inner sep=1pt] at (1.95,1.75) {$e_0$};
 \begin{scope}[xshift=5.6cm]
  \draw[thick] (0,1.8)..controls (0.7,3.2) and (2.1,3.2)..(2.8,1.8) node[pos=0.5,above] {$e_1$};
  \draw[thick] (0,1.8)--(2.8,1.8) node[midway,above] {$e_0$};
  \draw[thick] (0,1.8)..controls (0.7,0.4) and (2.1,0.4)..(2.8,1.8) node[pos=0.5,below] {$e_2$};
  \fill (0,1.8) circle (0.06) node[left] {$v$};
  \fill (2.8,1.8) circle (0.06) node[right] {$w$};
  \node at (1.4,0.0) {$g_{e_0}=1,\ g_{e_1}=A,\ g_{e_2}=B$};
 \end{scope}
\end{tikzpicture}
\caption{The once-punctured torus. The diagonal of the square gives an ideal
triangulation with two triangles; the dashed segments form the dual theta
graph $\Gamma$. Setting the transport along $e_0$ equal to $1$ leaves two
matrices $A,B$. The cyclic orders at $v$ and $w$ come from the square, not
from the planar drawing on the right; with these orders, $\Gamma$ has a
single boundary walk, which visits each vertex three times.}\label{fig:torus}
\end{figure}

\begin{figure}[htbp]
\centering
\begin{tikzpicture}[>=Stealth,font=\small,scale=0.9]
 \fill (0,0) circle (2pt) node[below=3pt] {$v$};
 \draw[thick] (0,0)--(90:2) node[above] {$h_1$};
 \draw[thick] (0,0)--(210:2) node[below left] {$h_2$};
 \draw[thick] (0,0)--(330:2) node[below right] {$h_3$};
 \node at (90:2.55) {};
 \draw[dashed,gray] (150:1.9)--(270:1.9)--(30:1.9)--cycle;
 \node[gray] at (270:2.25) {triangle};
 \draw[red,very thick,->] (105:1.5) to[out=-90,in=30] (195:1.5);
 \node[red] at (150:1.25) {$c$};
 \begin{scope}[xshift=6.3cm]
  \fill (0,0) circle (2pt);
  \draw[thick] (0,0)--(90:1.6) node[above] {$V(\varpi_a)^\vee$};
  \draw[thick] (0,0)--(210:1.6) node[below left] {$V(\varpi_a)$};
  \draw[thick,dotted] (0,0)--(330:1.6) node[below right] {$\C$};
  \node[align=center] at (0,-2.3) {corner pairing at $c$:\\ $\id_{V(\varpi_a)}\in V(\varpi_a)^\vee\otimes V(\varpi_a)\otimes\C$};
 \end{scope}
\end{tikzpicture}
\caption{A vertex $v$ of $\Gamma$ with half-edges $h_1,h_2,h_3$ in
counterclockwise order, inside the dual triangle (dashed). The red boundary
walk turns from $h_1$ to $h_2$, keeping the puncture corner $c$ on the right.
For $W=V(\varpi_a)$, the incoming half-edge $h_1$ carries $W^\vee$, the
outgoing half-edge $h_2$ carries $W$, and the third carries the trivial
representation; the turn contributes the identity tensor of $W$.}\label{fig:vertex}
\end{figure}
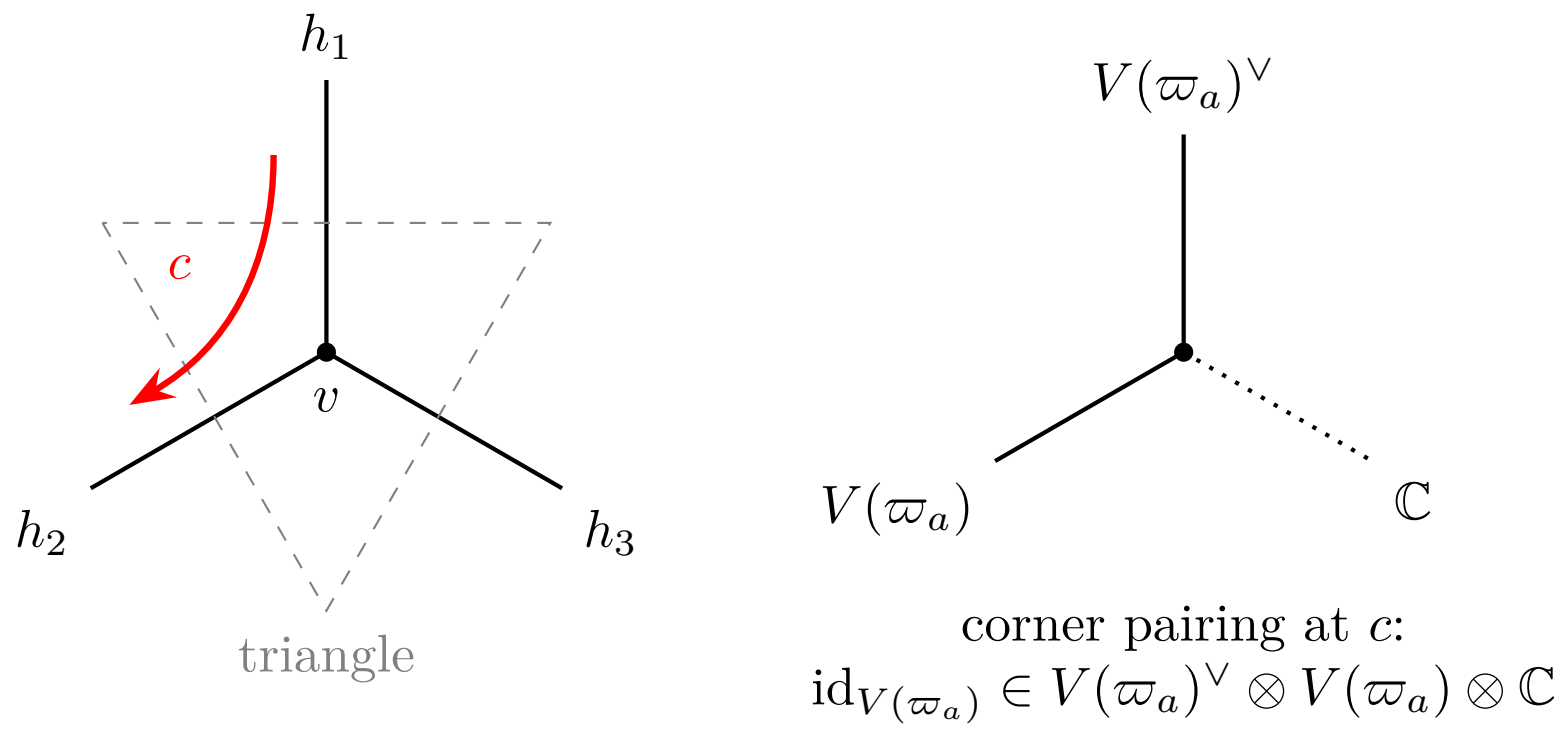

An \emph{ideal triangulation} of $\Sigma^\circ$ is a maximal collection of
pairwise disjoint and pairwise non-isotopic simple arcs with endpoints in $P$;
the complementary regions are ideal triangles, possibly self-folded, and ideal
triangulations exist because $2g+k\ge3$. The \emph{dual graph} $\Gamma$ has a
vertex in each triangle and an edge crossing each arc. The graph $\Gamma$ is
trivalent, embedded in $\Sigma^\circ$, and a deformation retract of
$\Sigma^\circ$; in particular the first Betti number of $\Gamma$ is $m$. The
orientation of $\Sigma^\circ$ gives a counterclockwise cyclic order on the
three outgoing half-edge rays at each vertex, in a positively oriented local
coordinate; a graph with such cyclic orders is called a \emph{ribbon graph}. We use the following dictionary:
\[
\begin{array}{@{}cccc@{}}
\toprule
\text{triangle}&\text{arc}&\text{puncture}&\text{corner of a triangle}\\
\midrule
\text{vertex of }\Gamma&\text{edge of }\Gamma&\text{boundary walk of }\Gamma
&\text{pair of half-edges at a vertex}\\
\bottomrule
\end{array}
\]
Here the \emph{boundary walk} of a puncture is the closed walk in $\Gamma$
obtained by retracting the oriented boundary circle around the puncture to
$\Gamma$. The orientation keeps the surface on the left and the puncture on
the right, so the circle is clockwise in a positively oriented coordinate
around the puncture. At each vertex the walk enters along one half-edge and
leaves along the next half-edge in the counterclockwise cyclic order
(Figure~\ref{fig:vertex}). A boundary walk may traverse an edge twice,
or visit a vertex several times; in the once-punctured torus
(Figure~\ref{fig:torus}) the unique boundary walk traverses each edge twice
and visits each vertex three times. We always count such passages separately,
as \emph{occurrences}. With this convention, every corner of every triangle
occurs exactly once in exactly one boundary walk.

Orient each edge $e$ from $s(e)$ to $t(e)$. A local system on $\Gamma$ with
fiber $H$ is a tuple of transports $(g_e)\in H^{E(\Gamma)}$, and changing the
trivializations at the vertices by $(h_v)\in H^{V(\Gamma)}$ replaces $g_e$ by
$h_{t(e)}g_eh_{s(e)}^{-1}$.

\begin{lemma}\label{lem:gauge}
Choose a spanning tree of $\Gamma$. Restriction to the transports along the
$m$ edges not in the tree identifies
$\C[H^{E(\Gamma)}]^{H^{V(\Gamma)}}$ with $A_0=\C[H^m]^H$. Under this
identification, the holonomy of the puncture loop $\gamma_i$ is the product of
the transports along the boundary walk of $p_i$.
\end{lemma}

\begin{proof}
Fix a root vertex. Moving away from the root along the tree, the gauge
transformations can be chosen uniquely (with $h_{\mathrm{root}}=1$) to make all
transports along tree edges equal to $1$. The remaining gauge freedom consists
of the elements of $H^{V(\Gamma)}$ preserving these conditions, i.e.\ of
constant tuples $(h,\dots,h)$, which act by simultaneous conjugation on the
$m$ remaining transports. An invariant function is thus determined by the
restriction of the function to the gauge-fixed locus, and conversely every conjugation
invariant function of the $m$ remaining transports extends uniquely to a
gauge-invariant function. The last statement holds because $\Gamma$ is a
spine of $\Sigma^\circ$.
\end{proof}

\subsection{The highest-weight filtration}\label{subsec:cartan}
By the Peter--Weyl theorem (Fact~\ref{fact:PW}), applied to each edge, and by
regrouping the tensor factors according to vertices, we obtain
\begin{equation}\label{eq:graph-tensors}
 A_0\simeq\bigoplus_{(\lambda_e)\in X^+(T)^{E(\Gamma)}}\ \bigotimes_{v\in V(\Gamma)}
 \Bigl(\bigotimes_{h\text{ at }v}V(\lambda_h)\Bigr)^H .
\end{equation}
Here the half-edges $h$ of an edge $e$ are labeled $\lambda_h=\lambda_e$ at
the source of $e$ and $\lambda_h=\lambda_e^*$ at the target of $e$, because the
coefficient of the matrix coefficient $c_{\ell,w}(g_e)=\ell(g_ew)$ in
$V(\lambda_e)$ is acted on by the gauge group at $s(e)$, and the coefficient
in $V(\lambda_e)^\vee=V(\lambda_e^*)$ by the gauge group at $t(e)$. Thus an
element of $A_0$ is a linear combination of \emph{tensor networks}: a choice of
a representation on each edge and an invariant tensor at each vertex.

Decomposition \eqref{eq:graph-tensors} is a decomposition of vector spaces,
not of rings: by Fact~\ref{fact:PW}, the product of two matrix coefficients of
$V(\lambda)$ and $V(\mu)$ is a matrix coefficient of $V(\lambda)\otimes
V(\mu)$, which decomposes into the Cartan component $V(\lambda+\mu)$ and
components $V(\nu)$ with $\nu<\lambda+\mu$ (Fact~\ref{fact:cartan}). The
latter have smaller \emph{degree} $\ell(\nu)<\ell(\lambda+\mu)$, where
$\ell(\lambda)=\langle\lambda,2\rho^\vee\rangle$. Filtering by degree therefore
replaces multiplication by the Cartan product.

\begin{example}\label{ex:cartan-sl3}
Let $H=\SL_3$, let $v_1,v_2,v_3$ be the standard basis of $V=\C^3$, and let
$v_1^*,v_2^*,v_3^*$ be the dual basis. The dual action is
$(g\phi)(v)=\phi(g^{-1}v)$. Thus
\[
 f(g)=g_{11}=v_1^*(gv_1),\qquad
 k(g)=(g^{-1})_{11}=(gv_1^*)(v_1).
\]
Under Peter--Weyl (Fact~\ref{fact:PW}), $f$ corresponds to
$v_1^*\otimes v_1\in V^\vee\otimes V$, and $k$ to
$v_1\otimes v_1^*\in V\otimes V^\vee$.
Calling $f$ a matrix coefficient of $V$ specifies the representation used
to evaluate $f$; the function $f$ belongs to the matrix-coefficient space
of $V$.

The product $fk$ is a matrix coefficient of the underlying tensor product
$W=V\otimes V^\vee$:
\[
 (fk)(g)=(v_1^*\otimes v_1)
             \bigl(g\cdot(v_1\otimes v_1^*)\bigr).
\]
Identify $W$ with $\End(V)$, with action $g\cdot T=gTg^{-1}$.
The input vector $v_1\otimes v_1^*$ becomes $E_{11}$ and the evaluating functional $v_1^*\otimes v_1$ becomes
$L(T)=T_{11}$. The decomposition
$\End(V)=\mathfrak{sl}_3\oplus\C I$ has scalar projection
$T\mapsto\frac{\tr T}{3}I$, so
\[
 \begin{gathered}
 E_{11} = (E_{11}- \frac{I}{3}) + \frac{I}{3} =: E_0 + \frac{I}{3},\quad \tr E_0=0,\\
 (fk)(g) = L(gE_0g^{-1}) + L(g\frac{I}{3}g^{-1}) = L(gE_0g^{-1}) + \frac{1}{3}.
 \end{gathered}
\]
Consequently $fk-\frac13$ is a matrix coefficient of the adjoint
representation $\mathfrak{sl}_3$, with input $E_0$ and functional
$L|_{\mathfrak{sl}_3}$; this representation is the Cartan component of $W$.
The constant $\frac13$ is a matrix coefficient of the trivial
representation $\C$: take the vector $I$ and the functional $zI\mapsto z/3$.

Let $F$ be the highest-weight filtration of $\C[H]$, assigning degree
$\ell(\lambda)$ to the entire
Peter--Weyl summand $V(\lambda)^\vee\otimes V(\lambda)$.
The representations $V$, $V^\vee$, $\mathfrak{sl}_3$ and $\C$ have
highest weights $\varpi_1$, $\varpi_2$, $\varpi_1+\varpi_2$ and $0$,
respectively. Since $2\rho^\vee=2\alpha_1^\vee+2\alpha_2^\vee$ in type
$A_2$, the definition $\ell(\lambda)=\langle\lambda,2\rho^\vee\rangle$ gives
\[
 \ell(\varpi_1)=\ell(\varpi_2)=2,\qquad
 \ell(\varpi_1+\varpi_2)=4,\qquad \ell(0)=0.
\]
Thus $f,k$ have degree $2$, $fk-\frac13$ has degree $4$, and the trivial
coefficient has degree $0$. Writing $[u]_d$ for the class of $u\in F_d$
in $F_d/F_{d-1}$, the associated-graded product is precisely
\[
 \inn_F(f)\inn_F(k)=[fk]_4=[fk-\tfrac13]_4.
\]
The constant contributes zero to this degree-$4$ class; constants remain
in degree $0$ of the associated graded ring.
\end{example}

To describe the associated graded ring, consider the \emph{three-flag
algebra}
\begin{equation}\label{eq:P3}
 R_H=\bigl(\C[H/U]^{\otimes3}\bigr)^H
 =\bigoplus_{\lambda_1,\lambda_2,\lambda_3\in X^+(T)}\bigl(V(\lambda_1)\otimes
 V(\lambda_2)\otimes V(\lambda_3)\bigr)^H ,
\end{equation}
graded by the \emph{boundary weights} $(\lambda_1,\lambda_2,\lambda_3)$; by
Fact~\ref{fact:basic-affine}(i), the multiplication of $R_H$ is the Cartan
product in each of the three slots, and $R_H$ is a domain. These labels
are the highest weights of the displayed $H$-modules; the corresponding
right-$T^3$ character is $(\lambda_1^*,\lambda_2^*,\lambda_3^*)$.
The two gradings agree on $(2\rho,2\rho,2\rho)$ because $(2\rho)^*=2\rho$.

\begin{proposition}[{cf.\ \cite[Props.~4.1, 4.11]{Manon2016}}]\label{prop:cartan-degeneration}
Filter $A_0$ by $\sum_e\ell(\lambda_e)$ in \eqref{eq:graph-tensors}. This filtration is a
positive filtration, and the associated graded ring is the subring of
$\bigotimes_{v\in V(\Gamma)}R_H$ spanned by the summands whose boundary
weights at the two ends of every edge are dual to each other (the three slots
at $v$ being the three half-edges at $v$).
\end{proposition}

\begin{proof}
For one edge, the discussion above shows that the associated graded ring of
$\C[H]$ for the filtration by $\ell(\lambda)$ is
$\bigoplus_\lambda V(\lambda)^\vee\otimes V(\lambda)$ with the Cartan product
in both factors. Taking the tensor product over all edges gives the associated
graded ring of $\C[H^{E(\Gamma)}]$. Since taking invariants under the
reductive group $H^{V(\Gamma)}$ is exact (Fact~\ref{fact:git}(i)), the
associated graded ring of the invariant ring is the invariant ring of the
associated graded ring, which is the ring described.
\end{proof}

\subsection{The local string valuation}\label{subsec:local-string}
The associated graded ring of Proposition~\ref{prop:cartan-degeneration} is
not yet a semigroup ring, because the spaces of invariants
$(V(\lambda_1)\otimes V(\lambda_2)\otimes V(\lambda_3))^H$ can have dimension
greater than one. A second filtration, by string parameters, distinguishes the
elements of these spaces.

The valuation has a concrete polynomial construction. Choose highest weight
vectors $v_\nu\in V(\nu)$, for dominant $\nu$, compatibly with Cartan
multiplication, and the lowest weight functionals $\ell_\nu\in V(\nu)^\vee$ normalized by
$\ell_\nu(v_\nu)=1$. An invariant of boundary weights
$(\mu,\lambda,\eta)$ determines a highest weight vector of weight $\mu^*$ in
$V(\lambda)\otimes V(\eta)$; contracting its first factor with $\ell_\lambda$
gives a vector $w\in V(\eta)$. Proposition~\ref{prop:poly-model} shows that
the assignment
\[
 T\longmapsto F_wL^\lambda M^\eta,\qquad
 F_w(n)=\ell_\eta(n^{-1}w)\quad(n\in U),
\]
defines an injective ring homomorphism
$R_H\hookrightarrow\C[U][L_1,\dots,L_r,M_1,\dots,M_r]$.
Here $L^\lambda=\prod_aL_a^{\langle\lambda,\alpha_a^\vee\rangle}$ and
$M^\eta=\prod_aM_a^{\langle\eta,\alpha_a^\vee\rangle}$ record the second and
third boundary weights. The conjugation weight of $F_w$ is
$\lambda+\eta-\mu^*$, so the first boundary weight is recoverable too.

Fix a reduced word $\mathbf i=(i_1,\dots,i_N)$ of $w_0$, choose raising
operators $e_j\in\Lie U$, and put $x_j(a)=\exp(ae_j)$. Pull back the
polynomial model along the birational map
\[
 \Phi_{\mathbf i}:\A^N\longrightarrow U,\qquad
 (a_1,\dots,a_N)\longmapsto
 x_{i_1}(a_1)\cdots x_{i_N}(a_N).
\]
Define $v_{\mathbf i}(T)$ to be the exponent of the lexicographically
largest monomial of the resulting polynomial, in the variable order
$(a_1,\dots,a_N,L_1,\dots,L_r,M_1,\dots,M_r)$. Leading monomials multiply,
so $v_{\mathbf i}(TT')=v_{\mathbf i}(T)+v_{\mathbf i}(T')$; addition gives
$v_{\mathbf i}(T+T')\le\max\{v_{\mathbf i}(T),v_{\mathbf i}(T')\}$ when
$T+T'\ne0$. This is the largest-term convention of
Section~\ref{subsec:bg-filtrations}. The coefficient of a leading monomial
also shows that the leaves are one-dimensional.

For $\SL_3$, take $e_1=E_{12}$, $e_2=E_{23}$ and write
$n=\left(\begin{smallmatrix}1&x&z\\0&1&y\\0&0&1\end{smallmatrix}\right)$.
Then $\Phi_{121}(a,b,c)$ gives $(x,y,z)=(a+c,b,ab)$, so $z$ and $xy-z$
pull back to $ab$ and $bc$, with distinct leading exponents $(1,1,0)$ and
$(0,1,1)$. Example~\ref{ex:sl3-good-basis} explains why these polynomials
belong to the compatible basis.

String parameters compute these leading exponents. Apply $e_{i_1}$ to $w$
as many times as possible without obtaining $0$, say $t_1$ times; then
apply $e_{i_2}$ as many times as possible, say $t_2$ times; and so on.
Lemma~\ref{lem:string-param} proves that the compatible basis used in
Section~\ref{app:string} has distinct resulting tuples
$(t_1,\dots,t_N)$, and that these are its leading exponents. Thus the value
of an element of boundary weights $(\mu,\lambda,\eta)$ has the form
$(\boldsymbol t;\lambda,\eta)$. For example, let $H=\SL_3$, $\mathbf i=(1,2,1)$, and let
$v_1,v_2,v_3$ be the standard basis of $\C^3$, with lowest weight vector $v_3$.
The raising operators are $e_1=E_{12}$ and $e_2=E_{23}$. We have
$e_1v_3=0$, then $e_2v_3=v_2$ and $e_2^2v_3=0$, then $e_1v_2=v_1$ and
$e_1^2v_2=0$; the string parameter of $v_3$ is $(0,1,1)$.

For a dominant weight $\lambda$ put
\[
 \beta_b^\vee=s_{i_1}\cdots s_{i_{b-1}}(\alpha_{i_b}^\vee),\qquad
 q_\lambda(\mathbf i)=\bigl(\langle\lambda,\beta_1^\vee\rangle,\dots,
 \langle\lambda,\beta_N^\vee\rangle\bigr)\in\N^N .
\]
By Lemma~\ref{lem:extremal}, $q_\lambda(\mathbf i)$ is the string parameter of a
lowest weight vector of $V(\lambda^*)$. For instance, for $\SL_3$ and
$\mathbf i=(1,2,1)$ we have $\beta^\vee_1=\alpha_1^\vee$,
$\beta^\vee_2=\alpha_1^\vee+\alpha_2^\vee$, $\beta^\vee_3=\alpha_2^\vee$, so
$q_{\varpi_1}=(1,1,0)$, $q_{\varpi_2}=(0,1,1)$ and $q_\rho=(1,2,1)$; the vector
$v_3$ above is a lowest weight vector of $\C^3=V(\varpi_2^*)$, and its string
parameter is $q_{\varpi_2}$.

The next theorem records the local semigroup, its boundary weights, and
the corner values. Their sum is an interior point, the key input to
Theorem~\ref{thm:corner}.

String cones originate in work of Littelmann
\cite{Littelmann1998}; Berenstein--Zelevinsky give polyhedral
tensor-product formulas \cite[Prop.~3.5, Thms.~2.3, 3.10]{BZ2001},
and Manon constructs the corresponding three-flag degeneration
\cite[Thm.~5.9]{Manon2016}. Section~\ref{app:string} proves the
formulation needed below from Lusztig's canonical-basis theory, including
the interior-corner calculation in part (iv).

\begin{theorem}\label{thm:local-cone}
There is a valuation $v_{\mathbf i}:R_H\setminus0\to\Z^{N+2r}$ with
one-dimensional leaves and the following properties.
\begin{enumerate}[label=\textup{(\roman*)}]
\item The value semigroup is $Q_{3,\mathbf i}=K_{3,\mathbf i}\cap\Z^{N+2r}$
for a full-dimensional pointed rational polyhedral cone $K_{3,\mathbf i}$.
Hence $\gr_{v_{\mathbf i}}R_H\simeq\C[Q_{3,\mathbf i}]$.
\item A homogeneous element of boundary weights $(\mu,\lambda,\eta)$ has
value $(\boldsymbol t;\lambda,\eta)$, where $\lambda$ and $\eta$ are written in
the basis of fundamental weights, $\boldsymbol t\in\N^N$, and
\begin{equation}\label{eq:omitted-weight}
 \mu^*=\lambda+\eta-\textstyle\sum_bt_b\alpha_{i_b}.
\end{equation}
In particular the boundary weights are given by a linear map
$\pi_\partial:Q_{3,\mathbf i}\to X^+(T)^3$, and
$\pi_\partial^{-1}(0)=\{0\}$.
\item For $1\le a\le r$, the spaces of invariants of boundary weights
$(0,\varpi_a,\varpi_a^*)$, $(\varpi_a^*,0,\varpi_a)$ and
$(\varpi_a,\varpi_a^*,0)$ are lines, spanned by the identity tensor of
$V(\varpi_a)$ placed in two of the three slots. We call these invariants the
\emph{corner pairings} $\kappa_{23,a}$, $\kappa_{31,a}$, $\kappa_{12,a}$. Their
values are
\[
 p_{23,a}=(q_{\varpi_a}(\mathbf i);\varpi_a,\varpi_a^*),\qquad
 p_{31,a}=(0;0,\varpi_a),\qquad
 p_{12,a}=(0;\varpi_a^*,0).
\]
\item The point $(q_\rho(\mathbf i);2\rho,2\rho)$ lies in the interior of
$K_{3,\mathbf i}$.
\end{enumerate}
\end{theorem}

The sum of all the corner-pairing values in (iii) is
$(q_\rho(\mathbf i);2\rho,2\rho)$. Part (iv), proved in
Proposition~\ref{prop:rho-interior}, therefore puts this sum in the
interior of the local cone. Theorem~\ref{thm:corner} identifies the
same sum with the canonical translator.

\begin{example}\label{ex:local-sl2-sl3}
(1) For $H=\SL_2$, Example~\ref{ex:basic-affine} gives
$\C[H/U]=\C[\C^2]$, so
$R_{\SL_2}=\C[(\C^2)^3]^{\SL_2}$. Let $v_1,v_2,v_3$ denote the three
vector coordinates, and put
$a=d_{12}$, $b=d_{23}$, $c=d_{31}$, where $d_{ij}=\det(v_i,v_j)$.
These determinants are invariant because $\det h=1$ for $h\in\SL_2$.
On $a\ne0$, let $g=(v_1,v_2/a)\in\SL_2$. Acting by $g^{-1}$ gives
the unique representative
\[
 v_1=(1,0),\qquad v_2=(0,a),\qquad v_3=(-b/a,-c).
\]
This change of basis and the displayed representative are regular on
$a\ne0$. Consequently the localized invariant ring is
$R_{\SL_2}[a^{-1}]=\C[a^{\pm1},b,c]$, and $a,b,c$ are algebraically
independent. The cyclically analogous chart $b\ne0$ shows that every
invariant belongs to
\[
 \C[a^{\pm1},b,c]\cap\C[a,b^{\pm1},c]=\C[a,b,c]
 \quad\text{inside }\C(a,b,c).
\]
Indeed, the first expression allows denominators only in $a$, and the
second only in the relatively prime variable $b$. Thus
$R_{\SL_2}=\C[d_{12},d_{23},d_{31}]$, with no relations.
The monomial $d_{12}^{n_3}d_{23}^{n_1}d_{31}^{n_2}$ has
boundary weights $(n_2+n_3,n_3+n_1,n_1+n_2)$. These weights determine
$n_1,n_2,n_3$ uniquely, so each fixed-boundary-weight space of invariants
has dimension at most one. Thus no string parameter is needed: $Q_3=\N^3$.
Theorem~\ref{thm:corner} is therefore immediate from Fact~\ref{fact:semigroup}:
\[
 d_{12}d_{23}d_{31}\,
 \frac{d(d_{12})}{d_{12}}\wedge
 \frac{d(d_{23})}{d_{23}}\wedge
 \frac{d(d_{31})}{d_{31}}
 =d(d_{12})\wedge d(d_{23})\wedge d(d_{31}).
\]
The form on the right is a regular canonical generator of the polynomial
ring. Thus the product of the three local corner pairings corresponds to
the interior exponent $(1,1,1)$, and the canonical translator is
$(1,1,1)$.

(2) For $H=\SL_3$, Example~\ref{ex:basic-affine} identifies $H/U$ with
nonzero incident pairs $(v,\phi)$ and shows that $\C[H/U]=\C[Z]$ for
$Z=\{\phi(v)=0\}$. Thus $R_{\SL_3}=\C[Z^3]^{\SL_3}$; zero vectors and
covectors may be included when computing regular functions.
Let $T$ have columns $v_1,v_2,v_3$, let $\Phi$ have rows
$\phi_1,\phi_2,\phi_3$, and put
\[
 Q=\Phi T=(q_{ij}),\qquad q_{ij}=\phi_i(v_j),\qquad
 E=\det T,\qquad D=\det\Phi .
\]
Under $h\in\SL_3$, the matrices transform as
$(T,\Phi)\mapsto(hT,\Phi h^{-1})$. Hence $q_{ij}$, $D$ and $E$ are
invariant, and the incidence conditions give $q_{ii}=0$.
Taking determinants of $Q=\Phi T$ gives
\[
 DE=P_++P_-,\qquad P_+=\phi_2(v_3)\phi_3(v_1)\phi_1(v_2),\quad
 P_-=\phi_3(v_2)\phi_1(v_3)\phi_2(v_1),
\]
because the only nonzero terms of the determinant are the two
three-cycles, both with positive sign.

To prove that these invariants give the full presentation, form
\[
 S=\C[q_{ij}\ (i\ne j),D,E]/(DE-P_+-P_-).
\]
Here $P_+=q_{23}q_{31}q_{12}$ and $P_-=q_{32}q_{13}q_{21}$.
On $E\ne0$, a change of basis in $\SL_3$ gives
$T=T_0=\operatorname{diag}(1,1,E)$ and $\Phi=QT_0^{-1}$.
This is a regular representative depending only on the six off-diagonal
entries of $Q$ and on $E^{\pm1}$, so the natural map $S\to R_{\SL_3}$
becomes an isomorphism
\[
 S[E^{-1}]\simeq R_{\SL_3}[E^{-1}]
       =\C[q_{ij}\ (i\ne j),E^{\pm1}].
\]
Using the covectors as a basis gives the analogous isomorphism after
inverting $D$. The ring $S$ is a domain: its defining polynomial is
primitive and linear in $E$. Moreover, $S$ is a hypersurface, smooth
where $D$ or $E$ is nonzero, and the locus $D=E=0$ has codimension $2$
(there it is the nonzero equation $P_++P_-=0$ in six variables).
Serre's criterion makes $S$ normal. Functions regular off that
codimension-two locus extend, so
$S=S[D^{-1}]\cap S[E^{-1}]$ in the fraction field
(Section~\ref{subsec:bg-algebra}). Every invariant lies in both
localized rings, hence $R_{\SL_3}=S$. This proves generation by the six
corner pairings and the two determinants, with exactly the displayed
relation.

The products
$P_+$ and $P_-$ have the same boundary weights
$(\varpi_1+\varpi_2,\varpi_1+\varpi_2,\varpi_1+\varpi_2)$, so the first
filtration does not separate $P_+$ from $P_-$. The string valuation for
$\mathbf i=(1,2,1)$ does: $v_{\mathbf i}(P_+)>v_{\mathbf i}(P_-)$, and
$\gr_{v_{\mathbf i}}R_{\SL_3}$ is the hypersurface ring $DE=P_+$ in eight
variables (Section~\ref{subsec:sl3-local}).
\end{example}

\subsection{The graph valuation}\label{subsec:graph-valuation}
At every vertex $v$ of $\Gamma$, label the three half-edges at $v$ in their
cyclic order, and choose a reduced word $\mathbf i(v)$. Let
\[
 Q_\Gamma\subset\textstyle\prod_{v\in V(\Gamma)}Q_{3,\mathbf i(v)}
\]
be the subsemigroup of tuples whose boundary weights at the two ends of every
edge are dual to each other. The semigroup $Q_\Gamma$ is the set of lattice
points of the cone obtained by intersecting $\prod_vK_{3,\mathbf i(v)}$ with
the linear subspace defined by these \emph{matching equations}; hence $Q_\Gamma$
is a finitely generated, saturated and pointed semigroup. Choose coordinates in
which $Q_\Gamma\subset\N^M$ and in which the first coordinate is the total
degree $\sum_e\ell(\lambda_e)$ (a linear function on $Q_\Gamma$).

\begin{theorem}[{cf.\ \cite[Thms.~1.1, 5.9, \S6]{Manon2016}}]\label{thm:graph-valuation}
There is a valuation $v_\Gamma:A_0\setminus0\to\N^M$ with one-dimensional
leaves and value semigroup $Q_\Gamma$, so that
$\gr_{v_\Gamma}A_0\simeq\C[Q_\Gamma]$. For $f\in A_0$, the value $v_\Gamma(f)$
is computed as follows: take the component $f_{\mathrm{top}}$ of $f$ of largest
total degree in \eqref{eq:graph-tensors}, regard $f_{\mathrm{top}}$ as an
element of $\bigotimes_vR_H$ by Proposition~\ref{prop:cartan-degeneration}, and
take the largest value of $f_{\mathrm{top}}$ under the product of the local
valuations $v_{\mathbf i(v)}$. Moreover
$\rank\gp(Q_\Gamma)=\dim A_0=(m-1)\dim H$.
\end{theorem}

\begin{proof}
By Theorem~\ref{thm:local-cone}, each $R_H$ has a basis of homogeneous
elements with distinct values (one element in each leaf), and the product of
two basis elements has, as leading term, a nonzero multiple of the basis
element whose value is the sum. Products over $v$ of these basis elements form
a basis of $\bigotimes_vR_H$ with the same property, and the products whose
boundary weights match form a basis of the ring of
Proposition~\ref{prop:cartan-degeneration}. The values of the matching basis
elements are exactly the points of $Q_\Gamma$. The rule in the statement
therefore defines a valuation (the leading term of a product is the product of
the leading terms, which is nonzero), with one-dimensional leaves (distinct
basis elements have distinct values) and with value semigroup $Q_\Gamma$.

Lemma~\ref{lem:scalarization} turns $v_\Gamma$ into a projective positive
filtration with associated graded ring $\C[Q_\Gamma]$, i.e.\ into a flat
degeneration of $A_0$ to $\C[Q_\Gamma]$. Hence
$\rank\gp(Q_\Gamma)=\dim\C[Q_\Gamma]=\dim A_0$.
For the final dimension formula, choose a regular semisimple $s\in T$ and
an element $h$ outside the finitely many proper parabolic subgroups containing
$T$. Every parabolic subgroup containing $s$ contains $Z_H(s)=T$, so
$(s,h)$ is $H$-irreducible. Since $m\ge2$, the irreducible locus in $H^m$
is therefore nonempty. The irreducible locus is open, because the reducible locus is a finite
union of closed incidence images over the projective varieties $H/Q$.
The orbits in the irreducible locus are closed with finite stabilizers by Lemma~\ref{lem:schur}. The generic quotient dimension is
$m\dim H-\dim H=(m-1)\dim H$.
\end{proof}

\begin{example}\label{ex:theta-sl2}
Let $H=\SL_2$ and let $\Gamma$ be the theta graph of Figure~\ref{fig:torus}.
For $\SL_2$ all representations are self-dual, and by
Example~\ref{ex:local-sl2-sl3}(1) a point of $Q_3$ is determined by its three
boundary weights. Hence $Q_\Gamma$ is the semigroup of triples
$(n_0,n_1,n_2)$ of edge weights that satisfy the triangle inequalities and
have even sum. This semigroup is free, generated by $(1,1,0)$, $(1,0,1)$ and
$(0,1,1)$, which are the edge weights of the trace functions $x=\tr A$,
$y=\tr B$ and $z=\tr AB$ (for instance, the loop representing $A$ traverses
$e_1$ and $e_0$ once each). Thus $v_\Gamma$ is the valuation given by the
exponents of $x,y,z$, compared first by total degree, and $Q_\Gamma\simeq\N^3$.
\end{example}

\subsection{Peripheral equations and their leading terms}\label{subsec:puncture-eq}
Fix semisimple elements $s_1,\dots,s_k\in H$, with no regularity
assumption, and write $a_i=\chi_H(s_i)$.  The equations below prescribe
these adjoint-quotient values; they allow every Jordan type in the fibers.
For a representation $W$ of $H$ and a loop $\gamma$, let
$\chi_W(\gamma)\in A_0$ be the function $\varphi\mapsto\tr_W\varphi(\gamma)$.
Let $\gamma_i$ be the loop around $p_i$, and put
\begin{equation}\label{eq:As}
 f_{i,a}=\chi_{\varpi_a}(\gamma_i)-\chi_{\varpi_a}(s_i),\qquad
 A_\bs=A_0/(f_{i,a}:1\le i\le k,\ 1\le a\le r).
\end{equation}

\begin{lemma}\label{lem:As}
$\Spec A_\bs$ is the scheme-theoretic fiber of the peripheral morphism
$\Spec A_0\to(H/\!/H)^k$ over $(\chi_H(s_i))_i$.
The reduction of this fiber is $\cM_B^{\mathrm{sp}}(H,\mathbf a)$; if all $s_i$ are
regular semisimple, then $\Spec A_\bs=\cM_B(H,\bcC)$ with
$\cC_i=H\cdot s_i$.
\end{lemma}

\begin{proof}
Using the surface relation to express $\gamma_k$ in terms of the free
generators, the variety $\Hom(F_m,H)=H^m$ parametrizes all solutions of
the surface relation without monodromy conditions.  Fact~\ref{fact:steinberg}(i)
identifies the prescribed invariant values with the equations
$f_{i,a}=0$.  They are invariant, so by Fact~\ref{fact:git}(i), taking
invariants commutes with imposing these equations.  Reduction commutes
with the quotient by the same exactness.  In the regular-semisimple case,
Fact~\ref{fact:steinberg}(ii) identifies each monodromy fiber with its
conjugacy class.  Proposition~\ref{prop:spectral-flatness} will show that
$A_\bs$ is already reduced for every parameter.
\end{proof}

To compute the leading terms of the functions $\chi_{\varpi_a}(\gamma_i)$, we
write them as tensor networks. In the gauge of Lemma~\ref{lem:gauge},
$\chi_W(\gamma_i)$ is obtained by placing the transport of $W$ on every edge
occurrence of the boundary walk of $p_i$, placing the identity tensor
$\id_W\in W^\vee\otimes W$ at every turn of the walk, and contracting all
indices. At a turn, the identity tensor occupies two of the three half-edges
at the vertex, and the trivial representation occupies the third half-edge:
the turn contributes a corner pairing.

Precisely, a corner of a triangle corresponds to a pair of half-edges at a
vertex $v$. For $W=V(\varpi_a)$, the incoming half-edge carries $W^\vee$
and the outgoing half-edge carries $W$, as in Figure~\ref{fig:vertex}.
Writing $\varpi_{a^*}=\varpi_a^*$, the boundary-oriented turns
$h_1\to h_2$, $h_2\to h_3$ and $h_3\to h_1$ therefore contribute
$\kappa_{12,a^*}$, $\kappa_{23,a^*}$ and $\kappa_{31,a^*}$, respectively,
in the notation of Theorem~\ref{thm:local-cone}(iii).
Let $p^v_{\iota,a}\in Q_{3,\mathbf i(v)}$ denote the value of the pairing
selected by the oriented occurrence $\iota$. Let $\mathcal I(i)$ be the set
of corner occurrences of the boundary walk of $p_i$, and put
\begin{equation}\label{eq:pia}
 p_{i,a}=\sum_{\iota\in\mathcal I(i)}p^{v(\iota)}_{\iota,a}
 \in\textstyle\prod_vQ_{3,\mathbf i(v)} .
\end{equation}
At the two ends of each traversal of an edge, the contributions to $p_{i,a}$
are $\varpi_a$ and $\varpi_a^*$, so $p_{i,a}$ satisfies the matching equations
and $p_{i,a}\in Q_\Gamma$.

\begin{proposition}\label{prop:peripheral-symbols}
For all $i,a$ there is $d_{i,a}\in\C^*$ such that
$\inn_{v_\Gamma}\bigl(\chi_{\varpi_a}(\gamma_i)\bigr)=d_{i,a}x^{p_{i,a}}$. In
particular $\inn_{v_\Gamma}(f_{i,a})=d_{i,a}x^{p_{i,a}}$.
\end{proposition}

\begin{proof}
If the boundary walk traverses an edge $n$ times, the tensor network above
carries $n$ copies of $V(\varpi_a)$ or its dual on that edge; by
Fact~\ref{fact:cartan}, the component of largest degree is obtained by
projecting the $n$-fold tensor product onto its Cartan component. Under the
identification of Proposition~\ref{prop:cartan-degeneration}, the image of
$\chi_{\varpi_a}(\gamma_i)$ in the associated graded ring is therefore the
tensor product over $v$ of the products, in $R_H$, of the corner pairings at
the occurrences at $v$. These products are nonzero because $R_H$ is a domain.
Since each $v_{\mathbf i(v)}$ is additive, the value of this element is
$p_{i,a}$, which proves the first assertion. The constant
$\chi_{\varpi_a}(s_i)$ has value $0<p_{i,a}$; the inequality holds because
$p_{i,a}$ has nonzero boundary weights.
\end{proof}

\begin{example}\label{ex:theta-sl2-symbol}
In Example~\ref{ex:theta-sl2}, the boundary walk traverses each edge twice, so
the leading part of $\tr[A,B]$ has edge weights $(2,2,2)$, which is the
exponent $(1,1,1)$ in the coordinates of $x,y,z$; indeed
$\inn(\tr[A,B])=-xyz$. The six corner occurrences of the walk contribute the
six corner pairings $d_{12},d_{23},d_{31}$ at $v$ and at $w$.
Each local product has boundary weights $(2,2,2)$, so the two products
match along all three edges. The glued element has edge weights $(2,2,2)$,
hence trace-coordinate exponent $(1,1,1)$: the global canonical monomial
is $xyz$.
\end{example}

\section{The local corner identity}\label{sec:corner}

We keep the notation of Section~\ref{sec:graph}. In
Example~\ref{ex:local-sl2-sl3}(1), the canonical translator of
$Q_3=\N^3$ is $(1,1,1)$, which is the sum of the values of the three corner
pairings $d_{12},d_{23},d_{31}$; equivalently, $d_{12}d_{23}d_{31}$ generates
the canonical module of $\C[Q_3]$. The main result of this section says that
the same holds for every simply connected semisimple group. Put
\[
 P_{\mathbf i}=\sum_{a=1}^r\bigl(p_{12,a}+p_{23,a}+p_{31,a}\bigr)
 =\bigl(q_\rho(\mathbf i);2\rho,2\rho\bigr),
\]
the sum of the values of all $3r$ corner pairings; the second equality follows
from Theorem~\ref{thm:local-cone}(iii), since $\sum_aq_{\varpi_a}=q_\rho$ and
$\sum_a(\varpi_a+\varpi_a^*)=2\rho$.

\begin{theorem}[Corner identity]\label{thm:corner}
For every reduced word $\mathbf i$ of $w_0$, the semigroup $Q_{3,\mathbf i}$
has a canonical translator $\vartheta_{\mathbf i}$, and
$\vartheta_{\mathbf i}=P_{\mathbf i}$. Equivalently, the product of the
leading forms of the $3r$ corner pairings generates the canonical module of
$\C[Q_{3,\mathbf i}]$.
\end{theorem}

\begin{idea}
The geometric volume form on the three-flag quotient has boundary weights
$\kappa=(2\rho,2\rho,2\rho)$
(Lemma~\ref{lem:basic-affine-form}, Proposition~\ref{prop:knop}).
Hilbert-series reciprocity transfers this degree to a canonical translator
$\vartheta_{\mathbf i}$ of $Q_{3,\mathbf i}$
(Lemma~\ref{lem:canonical-degree}). The point $P_{\mathbf i}$ is interior
and has the same boundary weights, so
$P_{\mathbf i}-\vartheta_{\mathbf i}\in Q_{3,\mathbf i}$ has zero boundary
weights and must vanish. The geometric degree computation is essential:
symmetry alone gives equal, self-dual boundary weights
\cite[Lem.~4.2]{LawtonManon2016}.
\end{idea}

\subsection{The canonical module of the three-flag algebra}\label{subsec:canonical-degree}
Let $S_H=\C[\mathcal Z_H^3]=\C[H/U]^{\otimes3}$, graded by $X^+(T)^3$ through
the right action of $T^3$, so that $R_H=S_H^H$ is a graded subring.

\begin{lemma}\label{lem:basic-affine-form}
The canonical sheaf of $\mathcal Z_H$ is generated by an $H$-invariant form
$\Omega$ on which right translation by $t\in T$ acts by $(2\rho)(t)$.
Consequently $\omega_{S_H}\simeq S_H(-\kappa)$, with $\kappa=(2\rho,2\rho,2\rho)$
and with an $H$-invariant generator.
\end{lemma}

\begin{proof}
At the base point of $H/U$, the tangent space is $\mathfrak h/\mathfrak u$, on
which $U$ acts with determinant $1$ because $U$ is unipotent. Hence a nonzero
top covector at the base point extends by left translation to a nowhere
vanishing $H$-invariant top form $\Omega$ on $H/U$. Right translation by $t$
acts on the tangent space at the base point by $\Ad(t^{-1})$ on
$\mathfrak h/\mathfrak u$, whose weights are $0$ and the negative roots; hence
$t$ multiplies $\Omega$ by $\prod_{\alpha>0}\alpha(t)=(2\rho)(t)$. By
Fact~\ref{fact:basic-affine}(ii), $\mathcal Z_H\setminus(H/U)$ has codimension
at least two, so $\Omega$ generates $\omega_{\mathcal Z_H}$. The exterior product
of three copies of $\Omega$ generates $\omega_{\mathcal Z_H^3}$.
\end{proof}

For instance, for $\SL_3$ the variety $\mathcal Z_H=\{(v,\phi):\phi(v)=0\}$ is
a hypersurface in $\C^6$; the coordinates of $v$ have weight $\varpi_1$, the
coordinates of $\phi$ have weight $\varpi_2$, and the equation has weight
$\varpi_1+\varpi_2$, so adjunction gives the weight
$3\varpi_1+3\varpi_2-(\varpi_1+\varpi_2)=2\rho$.

To obtain a canonical generator on the quotient, contraction along the
$H$-orbits must remain nonzero over every prime divisor. The next lemma
shows that stabilizers are finite away from codimension two and gives a
dense open product chart for the quotient.

\begin{lemma}\label{lem:principal}
For the diagonal action of $H$ on $\mathcal Z_H^3$:
\begin{enumerate}[label=\textup{(\roman*)}]
\item the set of points with positive-dimensional stabilizer has codimension
at least two;
\item there is a dense open $H$-stable subset of $\mathcal Z_H^3$, isomorphic to
$H\times T\times\mathcal Z_H$ with $H$ acting by left translation on the
first factor. In particular
$\C(\mathcal Z_H^3)^H=\Frac(R_H)$ and $\dim R_H=\dim\mathcal Z_H^3-\dim H$.
\end{enumerate}
\end{lemma}

\begin{proof}
(i) By Fact~\ref{fact:basic-affine}(ii) we may work on $(H/U)^3$. If the
stabilizer of $(g_1U,g_2U,g_3U)$ is positive-dimensional, then the
stabilizer contains a nontrivial unipotent one-parameter subgroup (the
stabilizer is contained in the unipotent group $g_1Ug_1^{-1}$), so
$\bigcap_j\Ad(g_j)\mathfrak u$ contains a nonzero nilpotent element $e$. We use
two standard facts (Fact~\ref{fact:nilpotent-cones}): there are finitely many
nilpotent $H$-orbits, and every nonzero nilpotent orbit $\mathcal O$ has even
dimension $d>0$. Scalar multiplication permutes the finitely many nilpotent
orbits, since it commutes with $H$. If $b$ is the order of the resulting
permutation group, every $z\in\C^*$ has a $b$-th root, so the permutation
induced by $z$ is a $b$-th power and hence the identity.
Thus each nonzero orbit is stable under scaling,
and $\PP\mathcal O$ has dimension $d-1$.
Let $\mathcal B$ be the flag variety and
$\mathcal B_e=\{\mathfrak b\in\mathcal B:e\in[\mathfrak b,\mathfrak b]\}$. The
variety of triples $(e,\mathfrak b',\mathfrak b'')$ with $e$ nilpotent and
$e\in[\mathfrak b',\mathfrak b']\cap[\mathfrak b'',\mathfrak b'']$ has
dimension $2N$: over the locus of pairs in relative position $w\in W$, which
has dimension $N+\ell(w)$, the fiber is the vector space
$[\mathfrak b',\mathfrak b']\cap[\mathfrak b'',\mathfrak b'']$ of dimension
$N-\ell(w)$. Restricting to $e\in\mathcal O$ gives $d+2\dim\mathcal B_e\le2N$,
so the set of $gU\in H/U$ with $e\in\Ad(g)\mathfrak u$ has dimension at most
$r+\dim\mathcal B_e\le r+N-d/2$. Consequently the set of pairs
($[e]\in\PP\mathcal O$, triple of points of $H/U$ compatible with $e$) has
dimension at most $(d-1)+3(r+N-d/2)=\dim(H/U)^3-1-d/2$, and the image in
$(H/U)^3$ has codimension at least $1+d/2\ge2$. Taking the union over the
finitely many nilpotent orbits proves (i).

(ii) Let $x_a\in\C[\mathcal Z_H^2]^H$ be the pairing of weights
$(\varpi_a,\varpi_a^*)$. A pair $(g_1U,g_2U)$ with all $x_a\ne0$ has
underlying flags in the open Bruhat orbit (opposite position). By the
Bruhat decomposition and the equality
$U\cap\dot w_0U\dot w_0^{-1}=\{1\}$, the map
$(h,t)\mapsto(hU,ht\dot w_0U)$ is an isomorphism from $H\times T$ onto the open
set $\{x_1\cdots x_r\ne0\}\subset\mathcal Z_H^2$. Adding a third point of
$\mathcal Z_H$ gives the open set $H\times T\times\mathcal Z_H\subset
\mathcal Z_H^3$, on which $H$ acts freely by left translation on the first
factor, with quotient $T\times\mathcal Z_H$. This open set is the
nonvanishing locus of the invariant function $x_1\cdots x_r$, which gives the
statements about rational invariants and dimension.
\end{proof}

The following contraction formula is closely related to Knop's
canonical-module theorem \cite{Knop1989}; we give a direct proof
under the hypotheses needed here. Contraction removes the orbit
directions from the volume form. The point to check is that the descended
form has no zero or pole along a divisor; factoriality ensures that a prime
divisor pulls back with multiplicity one.

\begin{proposition}\label{prop:knop}
Let a connected semisimple group $L$ act on a normal affine variety
$X=\Spec S$ such that $S$ is factorial with $S^*=\C^*$, and
$\omega_X=S\,\Omega$ for an $L$-invariant form $\Omega$. Put $R=S^L$, $Y=\Spec R$
and $\pi:X\to Y$, and assume that $\C(X)^L=\Frac R$, that
$\dim Y=\dim X-\dim L$, and that the points of $X$ with positive-dimensional
stabilizer form a set of codimension at least two. Let $\xi_1,\dots,\xi_n$ be
the vector fields generating the action of a basis of $\Lie L$. Then the
contraction $\iota_{\xi_n}\cdots\iota_{\xi_1}\Omega$ is the pullback
$\pi^*\beta$ of a rational top form $\beta$ on $Y$ which generates $\omega_Y$.
If a torus acts on $X$ commuting with $L$ and acts on $\Omega$ by a character
$\nu$, then the torus acts on $\beta$ by $\nu$.
\end{proposition}

\begin{proof}
\emph{The ring $R$ is factorial.} Let $f\in R$ and let $f=\prod_jf_j$ be the
prime factorization of $f$ in $S$. Since $L$ is connected and $S^*=\C^*$, each
$f_j$ spans an $L$-stable line, on which $L$ acts by a character; a semisimple
group has no nontrivial characters, so every $f_j$ is invariant. Hence $R$ is
factorial, and a prime element of $R$ remains prime in $S$.

\emph{Descent of the contraction.} At a general point of $X$ the vector fields
$\xi_j$ are linearly independent and span the kernel of $d\pi$ (by the
dimension assumption). The contraction
$\eta=\iota_{\xi_n}\cdots\iota_{\xi_1}\Omega$ is therefore a nonzero form which
vanishes on vertical vectors; $\eta$ is $L$-invariant, because the adjoint
representation of $L$ has determinant $1$. Choosing a transcendence basis
$y_1,\dots,y_d$ of $\C(Y)$, we can write
$\eta=b\cdot\pi^*(dy_1\wedge\dots\wedge dy_d)$ with $b\in\C(X)^L=\C(Y)$, so
$\eta=\pi^*\beta$ for a rational top form $\beta$ on $Y$.

\emph{The divisor of $\beta$ is zero.} Let $D\subset Y$ be a prime divisor and
$f\in R$ an equation of $D$. Then $E=V_X(f)$ is a prime divisor, $f$ vanishes
to order one along $E$, and $E$ dominates $D$ by Fact~\ref{fact:git}(i). At the
generic point of $E$, the variety $X$ is smooth, $\Omega$ is a nonvanishing
form, and the stabilizers are finite. Lift a transcendence basis
of $\C(D)$ to functions $z_2,\dots,z_d$ regular at the generic point of $D$;
then $df\wedge dz_2\wedge\dots\wedge dz_d$ generates $\omega_Y$ at $D$, and the
pullback of this form is nonzero at the generic point of $E$, because
$\pi^*f$ vanishes to order one along $E$ and $\C(E)/\C(D)$ is separable.
Hence $d\pi$ is surjective at the generic point of $E$, the vector fields
$\xi_j$ form a basis of $\ker d\pi$ there, and both $\eta$ and
$\pi^*(df\wedge dz_2\wedge\dots\wedge dz_d)$ are nonvanishing sections of the
line of top forms vanishing on $\ker d\pi$. So $\eta$ is a unit times the
pullback of a local generator of $\omega_Y$ at $D$. Thus
$\ord_D\beta=0$, and $\beta$ generates $\omega_Y$. A torus commuting with $L$
fixes the vector fields $\xi_j$, so the torus acts on $\eta$, hence on
$\beta$, by the character by which the torus acts on $\Omega$.
\end{proof}

\begin{lemma}\label{lem:canonical-degree}
The ring $R_H$ is Gorenstein with $\omega_{R_H}\simeq R_H(-\kappa)$,
$\kappa=(2\rho,2\rho,2\rho)$. For every reduced word $\mathbf i$, the semigroup
$Q_{3,\mathbf i}$ has a canonical translator $\vartheta_{\mathbf i}$, and
$\pi_\partial(\vartheta_{\mathbf i})=\kappa$.
\end{lemma}

\begin{proof}
(The summand $V(\lambda)$ of $\C[H/U]$ has right $T$-weight $\lambda^*$; since
$(2\rho)^*=2\rho$, this relabeling does not affect $\kappa$.) We compute the
canonical boundary degree geometrically, then transfer that degree to the
semigroup ring by Hilbert-series reciprocity. The positive Rees family
supplies the Cohen--Macaulay property needed for reciprocity. By
Fact~\ref{fact:basic-affine}(iii), Lemmas~\ref{lem:basic-affine-form}
and~\ref{lem:principal}, the action of $H$ on $\mathcal Z_H^3$ satisfies the
hypotheses of Proposition~\ref{prop:knop}; hence $\omega_{R_H}\simeq
R_H(-\kappa)$ as graded modules. Put $B=\C[Q_{3,\mathbf i}]$. By
Lemma~\ref{lem:scalarization}, there is a positive filtration of $R_H$ with
finitely generated Rees algebra $\cR$ and $\cR/t\cR=B$. The ring $B$ is
Cohen--Macaulay (Fact~\ref{fact:semigroup}), and $t$ is a nonzerodivisor, so
$\cR$ is Cohen--Macaulay at every prime containing $t$. The non-Cohen--Macaulay
locus of $\Spec\cR$ is closed and stable under the $\mathbb G_m$-action defined
by the grading; since the grading is positive, the closure of every orbit
contains the vertex, which is a Cohen--Macaulay point. Hence $\cR$ is
Cohen--Macaulay, and so is $R_H=\cR/(t-1)$; by the first part, $R_H$ is
Gorenstein.

The valuation $v_{\mathbf i}$ respects boundary weights, so $R_H$ and $B$ have
the same $X^+(T)^3$-graded Hilbert series. By
Fact~\ref{fact:canonical}(iii), applied to both rings,
\[
 H_{\omega_B}(\boldsymbol z)=(-1)^{\dim B}H_B(\boldsymbol z^{-1})
 =(-1)^{\dim R_H}H_{R_H}(\boldsymbol z^{-1})=H_{\omega_{R_H}}(\boldsymbol z)
 =\boldsymbol z^\kappa H_B(\boldsymbol z).
\]
By Theorem~\ref{thm:local-cone}(ii), $B$ has a one-dimensional component of
degree $0$, so $\omega_B$ has a one-dimensional component of degree $\kappa$,
spanned by an interior monomial $x^\vartheta$ (Fact~\ref{fact:semigroup}). The
submodule $x^\vartheta B\subset\omega_B$ has the same Hilbert series as
$\omega_B$, hence $\omega_B=x^\vartheta B$; so $\vartheta$ is the canonical
translator, and $\pi_\partial(\vartheta)=\kappa$.
\end{proof}

For $H=\SL_2$, Lemma~\ref{lem:canonical-degree} can be checked by hand: the
canonical module of $R_{\SL_2}=\C[d_{12},d_{23},d_{31}]$ is generated by
$d_{12}d_{23}d_{31}$ (with respect to the logarithmic form), whose boundary
weights are $(2,2,2)=(2\rho,2\rho,2\rho)$ (Example~\ref{ex:local-sl2-sl3}(1)).

\subsection{Proof of the corner identity}\label{subsec:corner-proof}
By Lemma~\ref{lem:canonical-degree}, $\vartheta_{\mathbf i}$ exists. By
Theorem~\ref{thm:local-cone}(iv), $P_{\mathbf i}=(q_\rho(\mathbf
i);2\rho,2\rho)$ lies in the interior of $K_{3,\mathbf i}$, so
$P_{\mathbf i}=\vartheta_{\mathbf i}+y$ for some $y\in Q_{3,\mathbf i}$. By
\eqref{eq:omitted-weight} and Lemma~\ref{lem:extremal}, the first boundary
weight $\mu$ of $P_{\mathbf i}$ satisfies
$\mu^*=4\rho-\sum_bq_{\rho,b}(\mathbf i)\alpha_{i_b}=4\rho-2\rho=2\rho$, so
$\pi_\partial(P_{\mathbf i})=\kappa$. By Lemma~\ref{lem:canonical-degree},
$\pi_\partial(y)=0$, and therefore $y=0$ by Theorem~\ref{thm:local-cone}(ii).
\qed

\begin{example}\label{ex:corner-sl3}
For $H=\SL_3$ and $\mathbf i=(1,2,1)$ (Example~\ref{ex:local-sl2-sl3}(2)), the
ring
\[
 \gr R_{\SL_3}=\C[A,B,C,\bar A,\bar B,\bar C,D,E]/(DE-ABC)
\]
is a hypersurface ring, where $A,\dots,\bar C$ are the six corner pairings. By adjunction
(Fact~\ref{fact:canonical}(i)), the canonical generator has degree
$v(A)+\dots+v(E)-v(DE)$, which is the sum of the values of the six corner
pairings, as Theorem~\ref{thm:corner} predicts. In coordinates,
$\vartheta=(1,2,1;2,2;2,2)$ and $q_\rho(121)=(1,2,1)$
(Section~\ref{subsec:sl3-local}).
\end{example}

\section{The special fiber}\label{sec:special}

Under the assumptions of Section~\ref{sec:graph}, the spectral ring
$A_\bs$ for arbitrary semisimple $s_i$ degenerates to a reduced
Gorenstein toric-face ring with canonical generator of degree zero
(Theorem~\ref{thm:special-fiber}). The model is $\C[x,y,z]/(xyz)$.

\subsection{The global corner identity}\label{subsec:global-corner}
By Theorem~\ref{thm:corner}, at every vertex $v$ the translator
$\vartheta_{\mathbf i(v)}=P_{\mathbf i(v)}$ has all three boundary weights
equal to $2\rho$, which is self-dual. Hence the tuple
\[
 \Omega_\Gamma=(P_{\mathbf i(v)})_{v\in V(\Gamma)}
\]
satisfies the matching equations, i.e.\ $\Omega_\Gamma\in Q_\Gamma$.

\begin{lemma}\label{lem:glue}
$\Omega_\Gamma$ is the canonical translator of $Q_\Gamma$.
\end{lemma}

\begin{proof}
Let $E$ be the subspace defined by the matching equations, so that the cone of
$Q_\Gamma$ is $E\cap\prod_vK_{3,\mathbf i(v)}$. The group
$\gp(Q_\Gamma)$ generated by $Q_\Gamma$ is the lattice of integral points of $E$ in the product of
the local value lattices: for any such $z$, the point $z+n\Omega_\Gamma$
belongs to $Q_\Gamma$ for all sufficiently large integers $n$, so
$z=(z+n\Omega_\Gamma)-n\Omega_\Gamma\in\gp(Q_\Gamma)$.
Since $\Omega_\Gamma\in E$ is an
interior point of the product cone, the interior of $E\cap\prod_vK_{3,\mathbf
i(v)}$ relative to $E$ is the intersection of $E$ with the interior of the
product cone. Hence a lattice point $q$ of $E$ lies in $\relint Q_\Gamma$ if
and only if $q-\Omega_\Gamma\in\prod_vQ_{3,\mathbf i(v)}$, which, as
$q-\Omega_\Gamma\in E$, means $q-\Omega_\Gamma\in Q_\Gamma$.
\end{proof}

\begin{proposition}\label{prop:global-corner}
The leading exponents $p_{i,a}$ of the monodromy equations
(Proposition~\ref{prop:peripheral-symbols}) satisfy
\begin{equation}\label{eq:global-corner}
 \sum_{i=1}^k\sum_{a=1}^rp_{i,a}=\Omega_\Gamma .
\end{equation}
\end{proposition}

\begin{proof}
Every corner of every triangle occurs exactly once among the boundary walks
(Section~\ref{subsec:ribbon}). The boundary-oriented corner pairings are
the three families in Theorem~\ref{thm:local-cone}(iii) with the index $a$
replaced by $a^*$ (Section~\ref{subsec:puncture-eq}). Since $a\mapsto a^*$
permutes the fundamental weights, summing over $a$ removes this difference.
Hence the left-hand side of \eqref{eq:global-corner} is the
tuple, indexed by the vertices $v$, of the sums of the values of all $3r$ corner
pairings at $v$, which is $(P_{\mathbf i(v)})_v$.
\end{proof}

In the running example, $\Omega_\Gamma=(1,1,1)$ in the coordinates of $x,y,z$,
and \eqref{eq:global-corner} says that the leading exponent of $\tr[A,B]$ is
$(1,1,1)$ (Example~\ref{ex:theta-sl2-symbol}).

\subsection{Canonical partitions}\label{subsec:canonical-partition}
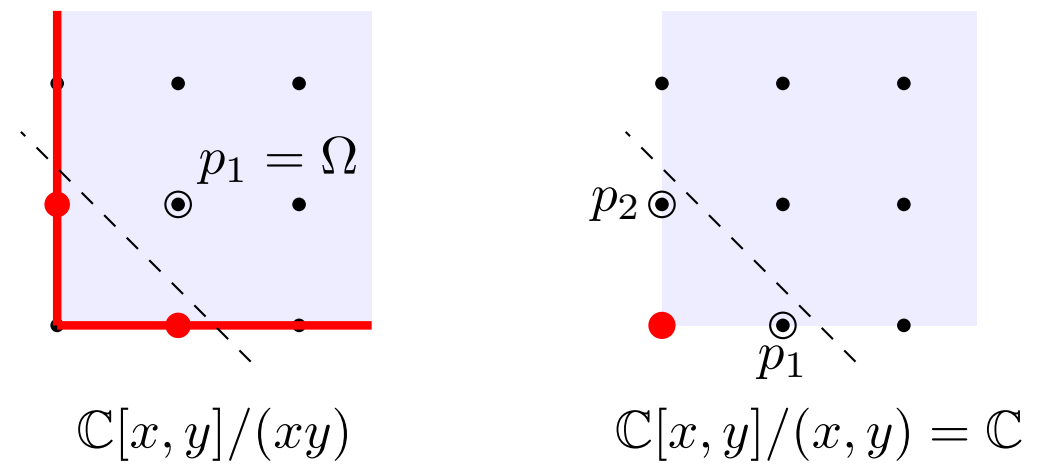
\begin{figure}[htbp]
\centering
\begin{tikzpicture}[scale=0.8,font=\small]
 \fill[blue!7] (0,0) rectangle (2.6,2.6);
 \foreach \x in {0,...,2}\foreach \y in {0,...,2}\fill (\x,\y) circle (1.6pt);
 \draw[red,ultra thick] (0,0)--(2.6,0);\draw[red,ultra thick] (0,0)--(0,2.6);
 \draw[dashed] (1.6,-0.3)--(-0.3,1.6);
 \fill[red] (1,0) circle (3pt);\fill[red] (0,1) circle (3pt);
 \node[above right] at (1,1) {$p_1=\Omega$};
 \draw (1,1) circle (3pt);
 \node at (1.3,-0.9) {$\C[x,y]/(xy)$};
 \begin{scope}[xshift=5cm]
 \fill[blue!7] (0,0) rectangle (2.6,2.6);
 \foreach \x in {0,...,2}\foreach \y in {0,...,2}\fill (\x,\y) circle (1.6pt);
 \fill[red] (0,0) circle (3.2pt);
 \draw[dashed] (1.6,-0.3)--(-0.3,1.6);
 \draw (1,0) circle (3pt) node[below] {$p_1$};
 \draw (0,1) circle (3pt) node[left] {$p_2$};
 \node at (1.3,-0.9) {$\C[x,y]/(x,y)=\C$};
 \end{scope}
\end{tikzpicture}
\caption{Canonical partitions of $\Omega=(1,1)$ in $\N^2$. The fan $\mathscr G$
of Proposition~\ref{prop:canonical-partition} is drawn in red; the dashed line
is the slice $x+y=1$, which meets $\mathscr G$ in a $0$-sphere on the left and
in the empty set on the right.}\label{fig:partition}
\end{figure}

The following proposition is the convex-geometric mechanism behind the
structure of the special fiber. To motivate the statement, consider
$Q=\N^2$, with translator $\Omega=(1,1)$ (Figure~\ref{fig:partition}). If
$\Omega=p_1$ is not split, then $x^{p_1}=xy$, the quotient $\C[x,y]/(xy)$ is the
toric-face ring of the fan consisting of the two coordinate rays, and the slice
$\{x+y=1\}$ of this fan is two points, a $0$-sphere. If
$\Omega=(1,0)+(0,1)$ is split into two parts, then $\C[x,y]/(x,y)=\C$ is the
toric-face ring of the fan $\{0\}$, whose slice is empty, the $(-1)$-sphere. In
both cases, every facet of the cone has height one over exactly one part $p_j$.
The same happens for $Q=\N^3$ and $\Omega=(1,1,1)=p_1$, which is the running
example: $\C[x,y,z]/(xyz)$ is the toric-face ring of the boundary of the octant,
whose slice is a triangle, a circle.

\begin{proposition}[Canonical partitions]\label{prop:canonical-partition}
Let $C\subset L_\R$ be a full-dimensional pointed rational cone of dimension
$d$, let $Q=C\cap L$, and suppose that the canonical translator $\Omega$ of
$Q$ is written as
\[
 \Omega=p_1+\dots+p_t,\qquad p_j\in Q\setminus\{0\}.
\]
For $1\le j\le t$, let $\mathcal E_j$ be the set of facets $F$ of $C$ with
$u_F(p_j)=1$.
\begin{enumerate}[label=\textup{(\roman*)}]
\item The sets $\mathcal E_1,\dots,\mathcal E_t$ are nonempty and form a
partition of the set of facets of $C$; for $F\in\mathcal E_j$ one has
$u_F(p_i)=\delta_{ij}$.
\item The monomials $x^{p_1},\dots,x^{p_t}$ form a regular sequence in
$\C[Q]$, and the ideal they generate is radical.
\item Let $\mathscr G$ be the fan consisting of the cones
$F_1\cap\dots\cap F_t$, with $F_j\in\mathcal E_j$, and all their faces. Then
$\C[Q]/(x^{p_1},\dots,x^{p_t})\simeq\C[\mathscr G]$ as $L$-graded rings.
\item Let $P$ be the span of $p_1,\dots,p_t$ and $\pi:L_\R\to V=L_\R/P$ the
projection. Then $\pi$ maps $|\mathscr G|$ homeomorphically onto $V$, maps each
cone of $\mathscr G$ linearly and injectively, and maps the fan $\mathscr G$ onto
a complete fan in $V$. In particular, for every linear form $h$ on $L_\R$
positive on $C\setminus0$, the slice $|\mathscr G|\cap h^{-1}(1)$ is homeomorphic
to a sphere of dimension $d-t-1$ (empty if $d=t$).
\end{enumerate}
\end{proposition}

\begin{idea}
By Lemma~\ref{lem:height-one}, $u_F(\Omega)=1$ for every facet $F$, so every
facet is ``charged'' by exactly one of the summands $p_j$. Divisibility by
$x^{p_j}$ is then tested only on the facets in $\mathcal E_j$, and the facet
sets for different $j$ do not interact; this independence gives (ii) and
(iii). For (iv), subtracting from a point $x$ the largest admissible multiple of
each $p_j$ moves $x$ into $|\mathscr G|$, and gives the inverse of $\pi$.
\end{idea}

\begin{proof}
(i) The integers $u_F(p_j)\ge0$ sum to $u_F(\Omega)=1$, so exactly one of them
is $1$. If $\mathcal E_j$ were empty, $p_j$ would lie on every facet, hence be
$0$, since $C$ is pointed. Note the consequence, for $q\in Q$:
\begin{equation}\label{eq:divisibility}
 q-p_j\in Q\iff u_F(q)\ge1\ \text{for all }F\in\mathcal E_j .
\end{equation}

(ii) All ideals below are monomial, since they are stable under the torus
$\Spec\C[L]$, whose weight spaces in $\C[Q]$ are lines. Suppose that
$x^{q+p_j}\in(x^{p_1},\dots,x^{p_{j-1}})$, i.e.\ $q+p_j-p_i\in Q$ for some
$i<j$. For $F\in\mathcal E_i$ we have $u_F(p_j)=0$ and $u_F(p_i)=1$, so
$u_F(q)=u_F(q+p_j-p_i)+1\ge1$; by \eqref{eq:divisibility}, $q-p_i\in Q$ and
$x^q\in(x^{p_1},\dots,x^{p_{j-1}})$. Thus $x^{p_j}$ is a nonzerodivisor modulo
$(x^{p_1},\dots,x^{p_{j-1}})$; the ideal is proper as all $p_j\ne0$. For
radicality, suppose $x^{nq}$ lies in the ideal, i.e.\ $nq-p_j\in Q$ for some
$j$. Then $nu_F(q)\ge1$, hence $u_F(q)\ge1$, for all $F\in\mathcal E_j$, and
$x^q$ lies in the ideal by \eqref{eq:divisibility}.

(iii) By \eqref{eq:divisibility}, the monomial $x^q$ survives in the quotient
if and only if, for each $j$, some $F_j\in\mathcal E_j$ satisfies
$u_{F_j}(q)=0$, i.e.\ if and only if $q\in F_1\cap\dots\cap F_t$ for some
choice of $F_j\in\mathcal E_j$. The product $x^qx^{q'}$ of surviving
monomials survives if and only if the same choice works for $q$ and $q'$,
because $u_F(q+q')=0$ if and only if $u_F(q)=u_F(q')=0$. This is the
multiplication of $\C[\mathscr G]$.

(iv) By (i), $u_F(p_i)=\delta_{ij}$ for $F\in\mathcal E_j$, so the $p_j$ are
linearly independent. For $x\in L_\R$ put $m_j(x)=\min_{F\in\mathcal
E_j}u_F(x)$ and $r(x)=x-\sum_jm_j(x)p_j$. For $F\in\mathcal E_j$ we get
$u_F(r(x))=u_F(x)-m_j(x)\ge0$, with equality for some $F\in\mathcal E_j$;
hence $r(x)\in|\mathscr G|$. Moreover $r(x+\sum_ja_jp_j)=r(x)$ for all real
$a_j$, and $r(q)=q$ for $q\in|\mathscr G|$ (all $m_j(q)$ vanish). Hence $r$
induces a continuous inverse $V\to|\mathscr G|$ of $\pi|_{|\mathscr G|}$. A
linear map which is injective on a cone is injective on the span of the cone, so
$\pi$ maps each cone of $\mathscr G$ isomorphically onto a cone in $V$, and
$\pi(\sigma)\cap\pi(\tau)=\pi(\sigma\cap\tau)$ is a common face; thus the images
form a complete fan in $V$. Finally, the cone over the slice
$|\mathscr G|\cap h^{-1}(1)$ is $|\mathscr G|\simeq V$, so radial projection
identifies the slice with the unit sphere of $V$, of dimension $d-t-1$.
\end{proof}

\subsection{The associated graded ring}\label{subsec:special-proof}

\begin{theorem}\label{thm:special-fiber}
There are a projective positive filtration $F$ of $A_\bs$, a rational fan
$\mathscr G_\Gamma$, and a linear form $h$ positive on
$|\mathscr G_\Gamma|\setminus0$, such that:
\begin{enumerate}[label=\textup{(\roman*)}]
\item $\gr_FA_\bs\simeq\C[\mathscr G_\Gamma]$, the grading of
$\C[\mathscr G_\Gamma]$ being given by $h$;
\item $\C[\mathscr G_\Gamma]$ is reduced and Gorenstein, and the canonical
module of $\C[\mathscr G_\Gamma]$ is generated in degree $0$;
\item the fan $\mathscr G_\Gamma$ maps, as in
Proposition~\ref{prop:canonical-partition}(iv), onto a complete fan in a vector
space of dimension $d_\bs=\dim A_\bs=(2g+k-2)\dim H-kr$.
\end{enumerate}
\end{theorem}

\begin{proof}
Apply Proposition~\ref{prop:canonical-partition} to $L=\gp(Q_\Gamma)$,
$C=\cone(Q_\Gamma)$ and the $t=kr$ points $p_{i,a}$; the hypotheses hold by
the saturation of $Q_\Gamma$, Lemma~\ref{lem:glue} and
Proposition~\ref{prop:global-corner}. Let $\mathscr G_\Gamma$ be the resulting
fan.

Next apply Lemma~\ref{lem:scalarization} to $A_0$, the valuation $v_\Gamma$, and
the elements $\chi_{\varpi_a}(\gamma_i)$; let $F$ be the resulting filtration
of $A_0$ and $h$ the resulting degree function. Since $h(p_{i,a})>0$ and
$\chi_{\varpi_a}(s_i)$ is a constant, the initial form of
$d_{i,a}^{-1}f_{i,a}$ is $x^{p_{i,a}}$. These initial forms form a regular
sequence by Proposition~\ref{prop:canonical-partition}(ii), so by
Lemma~\ref{lem:filtered-regular} the quotient filtration on $A_\bs$, again
denoted $F$, satisfies
\[
 \gr_FA_\bs\simeq\C[Q_\Gamma]/(x^{p_{i,a}})\simeq\C[\mathscr G_\Gamma],
\]
and $\cR_F(A_\bs)$ is a quotient of the finitely generated ring
$\cR_F(A_0)$; this proves (i).

Reducedness is Proposition~\ref{prop:canonical-partition}(ii). The ring
$\C[Q_\Gamma]$ is Gorenstein with canonical generator $x^{\Omega_\Gamma}$
(Lemma~\ref{lem:glue}), so by Fact~\ref{fact:canonical}(i) the quotient by
the regular sequence is Gorenstein with canonical generator of degree
$\Omega_\Gamma-\sum_{i,a}p_{i,a}=0$; this proves (ii). Finally
$\rank\gp(Q_\Gamma)=(m-1)\dim H=(2g+k-2)\dim H$ by
Theorem~\ref{thm:graph-valuation}, and $d_\bs=(2g+k-2)\dim H-kr$, so (iii)
is Proposition~\ref{prop:canonical-partition}(iv).
\end{proof}

\begin{example}[The cubic special fiber]\label{ex:cubic-special-fiber}
Put $P=A_0=\C[x,y,z]$, $\kappa=x^2+y^2+z^2-xyz-2$ and
$A_\tau=P/(\kappa-\tau)$. Let $F$ be the total-degree filtration on $P$
and $\bar F$ the quotient filtration on $A_\tau$.
Then $\inn_F(\kappa-\tau)=-xyz$ for every $\tau$.
Since $xyz$ is a nonzerodivisor
in $P$, Lemma~\ref{lem:filtered-regular} gives
$\gr_{\bar F}A_\tau=B=P/(xyz)$, as in Example~\ref{ex:rees-cubic}.
The ideal $(xyz)=(x)\cap(y)\cap(z)$ is radical, so $B$ is reduced.

The canonical module of $P$ is $P(-3)$, generated by
$dx\wedge dy\wedge dz$. On the dense coordinate torus
$\{xyz\ne0\}\simeq(\C^*)^3$ of $\Spec P$, put
$\eta_T=\frac{dx}{x}\wedge\frac{dy}{y}\wedge\frac{dz}{z}$.
Then $\omega_P=(xyz)P\eta_T$: multiplying $\eta_T$ by $xyz$ cancels the
simple poles along the coordinate planes and gives the regular volume
form $dx\wedge dy\wedge dz$. This realizes the interior-monomial ideal
$(xyz)$ as the canonical module (Fact~\ref{fact:semigroup}).
The coefficient $xyz$ of the reference torus form is, up to sign, the
initial puncture equation. This is the model for the global corner
identity \eqref{eq:global-corner}.
Graded adjunction
(Fact~\ref{fact:canonical}(i)) gives
\[
 \omega_B\simeq(\omega_P/(xyz)\omega_P)(3)\simeq B.
\]
By Lemma~\ref{lem:adjunction}, the degree-zero generator is the residue
$\eta=\operatorname{Res}_{\{xyz=0\}}\frac{dx\wedge dy\wedge dz}{xyz}$.
On the three coordinate planes $\{x=0\}$, $\{y=0\}$ and $\{z=0\}$,
the respective branch forms are
\[
 \frac{dy}{y}\wedge\frac{dz}{z},\qquad
 -\frac{dx}{x}\wedge\frac{dz}{z},\qquad
 \frac{dx}{x}\wedge\frac{dy}{y}.
\]
Their residues along each common line have opposite signs. The logarithmic
poles on the individual planes therefore describe a regular dualizing
section on the singular union. Each branch form is invariant under
simultaneous scaling, making the degree-zero claim explicit. Adjunction
divides the ambient canonical generator $xyz\eta_T$ by the initial
hypersurface equation before taking the residue, giving
$\operatorname{Res}\eta_T$ up to sign. The resulting degree-zero generator
is the input to the Rees calculation of Lemma~\ref{lem:rees-canonical},
which gives $K_{X_\tau}+D_\tau\sim0$ for normal fibers. This verifies
Theorem~\ref{paper:thm-degeneration}(i) in the running example.

Here $Q_\Gamma=\N^3$, $kr=1$ and $p_{1,1}=(1,1,1)$, so
$\mathscr G_\Gamma$ is the fan of proper faces of the octant.
The projection along $(1,1,1)$ identifies the support of this fan with
$\R^3/\R(1,1,1)\simeq\R^2$: every class of $q\in\R^3$ has the unique
representative
\[
 q-\min(q_1,q_2,q_3)(1,1,1)
\]
in the boundary of the nonnegative octant. The inverse is continuous and
linear on each cone, and the three coordinate faces give a complete fan
with three two-dimensional cones. The positive slice is
\[
 \{q\in\partial\R_{\ge0}^3:q_1+q_2+q_3=1\},
\]
the boundary of a triangle, hence $S^1$, as Theorem~\ref{paper:thm-degeneration}(ii) predicts for
$d=\dim A_\tau=2$. The projective boundary $\Proj B$ is a union of three
projective lines; the dual complex of $(X_\tau,D_\tau)$ is also a circle
by the incidence of those lines (\hyperref[subsec:intro-example]{the introductory model calculation} and
Example~\ref{ex:exceptional-cubics}).
\end{example}

\begin{remark}\label{rem:kim}
The canonical partition mechanism is not specific to our coordinates. In
Kim's skein-theoretic compactifications for $\SL_3$
\cite[Thm.~5.8, Prop.~4.6]{Kim2026}, the sum, over all punctures, of the
values of both orientations of the peripheral webs is again the interior
translator.
\end{remark}
\subsection{A uniform filtration and a flat peripheral map}
\label{subsec:spectral-filtration}
Prescribed values have degree zero and do not change the leading terms.
We use this uniformity to show that $A_0$ is free over the peripheral
character ring.

Let $H$ be simply connected and semisimple, of rank $r$, and assume
$2g+k\ge3$. In the coordinates of $H/\!\!/H$ given by the fundamental characters
(Fact~\ref{fact:steinberg}), write the prescribed values as
$\mathbf c=(c_{i,a})\in\A^{kr}$, and put
\[
 g_{i,a}=\chi_{\varpi_a}(\gamma_i),\qquad
 A_{\mathbf c}=A_0/(g_{i,a}-c_{i,a})_{i,a}.
\]
This is the ring $A_\bs$ of \eqref{eq:As} when $c_{i,a}=\chi_{\varpi_a}(s_i)$.
No regularity or genericity is assumed.

\begin{example}[Flatness of the trace family]\label{ex:cubic-flatness}
In the running example, $A_0=\C[x,y,z]$ and there is one peripheral character
$\kappa=\tr[A,B]=x^2+y^2+z^2-xyz-2$, with leading term $-xyz$. Every monomial
can be written uniquely as $(xyz)^n\cdot x^ay^bz^c$ with $\min(a,b,c)=0$, so the
monomials not divisible by $xyz$ form a basis of $\C[x,y,z]$ as a module over
$\C[xyz]$. Replacing $xyz$ by $-\kappa$ and inducting on the degree, the same
monomials span $\C[x,y,z]$ over $\C[\kappa]$. They are also linearly
independent: the leading monomials of
$\kappa^n x^ay^bz^c$, with $\min(a,b,c)=0$, have distinct exponents
$(n+a,n+b,n+c)$. Thus the map $\kappa:\A^3\to\A^1$ is faithfully flat,
and the same basis specializes to a vector-space basis in every fiber.

Every fiber is an integral hypersurface of dimension two. Indeed,
$\kappa-\tau$, viewed as a monic quadratic in $x$, has discriminant
\[
 \Delta=(z^2-4)y^2-4z^2+4(\tau+2).
\]
A polynomial square of degree two in $y$ would have a square leading
coefficient in $\C[z]$, whereas $z^2-4$ has two simple roots.
Unique factorization then shows that $\Delta$ is not a square in
$\C(y,z)$ either. The quadratic is therefore irreducible over
$\C(y,z)$, and Gauss's lemma gives irreducibility in $\C[x,y,z]$.
Consequently every fiber is reduced and connected; the hypersurface
structure gives the Gorenstein property. This verifies all assertions
of Theorem~\ref{paper:thm-degeneration}(iii) directly.
\end{example}

\begin{proposition}[Uniformity and freeness]\label{prop:spectral-flatness}
The filtration of $A_0$ in the proof of Theorem~\ref{thm:special-fiber} can be
chosen independently of $\mathbf c$ and stable under the group $Z(H)^m$ acting
by central translations of the free generators. For every $\mathbf c$ the
quotient filtration satisfies
\[
 \gr A_{\mathbf c}\simeq\C[Q_\Gamma]/(x^{p_{i,a}})\simeq\C[\mathscr G_\Gamma].
\]
Moreover, the peripheral characters $g_{i,a}$ are algebraically independent,
and $A_0$ is a free module over $\C[g_{i,a}]$. Consequently the peripheral
morphism
\[
 \Spec A_0\longrightarrow(H/\!\!/H)^k
\]
is faithfully flat, and every fiber $\Spec A_{\mathbf c}$ is nonempty, reduced,
connected and Gorenstein, of pure dimension $d=(2g+k-2)\dim H-kr$.
\end{proposition}

\begin{idea}
The canonical partition gives a monomial basis of $\C[Q_\Gamma]$ over
$\C[x^{p_{i,a}}]$. Lifting that basis and inducting on filtration degree
gives a free basis of $A_0$ over $\C[g_{i,a}]$.
\end{idea}

\begin{proof}
\emph{Uniformity.} In the tree gauge of Lemma~\ref{lem:gauge}, central
translations of the free generators translate the non-tree edge transports.
Each fixed-edge-weight summand of \eqref{eq:graph-tensors} is an eigenspace,
since the center acts on each $V(\lambda_e)$ by its central character.
Thus these translations preserve $v_\Gamma$ and act by characters on its
leaves. The transport along a boundary word is multiplied by a central
element, so each $g_{i,a}$ is an eigenvector as well. Apply the equivariant
part of Lemma~\ref{lem:scalarization} once for $Z(H)^m$, preserving these
initial forms. Constants have degree $0$, whereas every
$p_{i,a}$ has positive degree, so the initial form of $g_{i,a}-c_{i,a}$ is the
initial form of $g_{i,a}$, independently of $\mathbf c$.
Proposition~\ref{prop:canonical-partition} and
Lemma~\ref{lem:filtered-regular} give the displayed associated graded ring,
with the same grading and the same fan for every $\mathbf c$.

\emph{Freeness.} Enumerate the pairs $(i,a)$ as $j=1,\dots,t$, with $t=kr$, let
$\mathcal E_j$ be the facet blocks of Proposition~\ref{prop:canonical-partition},
and put $\mathscr B=|\mathscr G_\Gamma|\cap\gp(Q_\Gamma)$. Every $q\in Q_\Gamma$
has a unique expression
\begin{equation}\label{eq:spectral-free-decomposition}
 q=b+\sum_{j=1}^t n_jp_j,\qquad b\in\mathscr B,\quad n_j\in\N,
 \qquad n_j=\min_{F\in\mathcal E_j}u_F(q):
\end{equation}
indeed, the retraction $r$ in the proof of
Proposition~\ref{prop:canonical-partition}(iv) has exactly these integral
coefficients. Consequently
\[
 \C[Q_\Gamma]=\bigoplus_{b\in\mathscr B}\C[x^{p_1},\dots,x^{p_t}]\,x^b .
\]
Rescale the $g_j$ so that their initial forms are the monomials $x^{p_j}$, and
choose a lift $v_b\in A_0$ of each $x^b$, with $v_0=1$. Induction on the
filtration degree shows that the elements $v_b$ span $A_0$ over $\C[g_j]$. In a
nontrivial finite relation among the elements $g_1^{n_1}\cdots g_t^{n_t}v_b$,
the terms of largest degree have pairwise distinct initial monomials by
\eqref{eq:spectral-free-decomposition}, so they cannot cancel. This proves
freeness, including the algebraic independence of the $g_j$. A free module
containing $1$ in a basis is faithfully flat.

\emph{The fibers.} The fiber $A_{\mathbf c}$ is nonzero because its associated
graded ring is. If $A_{\mathbf c}$ had a nonzero nilpotent element, its initial
form would be a nonzero nilpotent of the reduced ring $\gr A_{\mathbf c}$; so
$A_{\mathbf c}$ is reduced. If an idempotent had positive degree, its initial
form would have square zero; so every idempotent lies in degree $0$, i.e.\ in
$\C$, and $\Spec A_{\mathbf c}$ is connected. For the Gorenstein property,
consider the Rees algebra of $A_0$: its special fiber $\C[Q_\Gamma]$ is
Gorenstein, so the Rees algebra is Gorenstein at the vertex; its
non-Gorenstein locus is closed and stable under the $\mathbb G_m$-action of the
grading, and the closure of every orbit contains the vertex, so this locus is
empty. Hence $A_0$ is Gorenstein. The elements $g_{i,a}-c_{i,a}$ form a regular
sequence by Lemma~\ref{lem:filtered-regular}; since $A_0$ is an equidimensional
affine domain, $A_{\mathbf c}$ is Gorenstein of pure dimension
$\dim A_0-kr=d$.
\end{proof}

\begin{remark}\label{rem:spectral-normality-scope}
Reduced and Gorenstein does not imply normal; indeed the toric-face ring
$\C[\mathscr G_\Gamma]$ is never normal when the fan has more than one maximal
cone. The normality needed for the compactification
(Theorem~\ref{thm:AB}) is provided by Proposition~\ref{prop:spectral-interior}
under intrinsic genericity, and by Fact~\ref{fact:gl-closure-normality} for all
parameters in type~$A$.
\end{remark}

\section{Log canonical compactifications from Rees algebras}\label{sec:rees}

Let $A$ be a normal affine domain with a projective positive filtration
$F$ as in Theorem~\ref{thm:special-fiber}. Rees adjunction gives
$K_X+D\sim0$ for $X=\Proj\cR_F(A)$ and $D=\Proj\gr_FA$
(Lemma~\ref{lem:rees-canonical}). To prove log canonicity, we use the
degeneration argument of Arg\"uz--Bousseau
\cite[Thm.~3.7, \S4]{ArguzBousseau2026}: degenerate the pair to a union
of projective toric varieties, compute its different on the normalized
components, and apply inversion of adjunction. A $\mathbb G_m$-action
then transports log canonicity to $(X,D)$
(Theorem~\ref{thm:AB}).

\begin{example}\label{ex:cubic-degeneration}
For the cubic surface of \hyperref[subsec:intro-example]{the introductory model calculation}, the degeneration
replaces the equation $t(x^2+y^2+z^2)-xyz-(\tau+2)t^3=0$ by $xyz=0$:
\[
 Y_0=\{xyz=0\}\subset\PP^3_{[x:y:z:t]},\qquad D_0=Y_0\cap\{t=0\}.
\]
The variety $Y_0$ is the union of three projective planes, glued pairwise along
lines through the point $[0:0:0:1]$; the divisor $D_0$ is a triangle of lines.
On each plane, the two gluing lines and the line $D_0\cap\PP^2$ form the toric
boundary (Figure~\ref{fig:ball}).
\end{example}

\begin{figure}[htbp]
\centering
\begin{tikzpicture}[scale=1.1,font=\small]
 \coordinate (O) at (0,0);
 \coordinate (A) at (90:1.8);\coordinate (B) at (210:1.8);\coordinate (C) at (330:1.8);
 \fill[blue!8] (A)--(B)--(C)--cycle;
 \draw[thick] (O)--(A);\draw[thick] (O)--(B);\draw[thick] (O)--(C);
 \draw[red,very thick] (A)--(B)--(C)--cycle;
 \fill (O) circle (2pt);
 \node at (150:0.75) {$\PP^2$};\node at (30:0.75) {$\PP^2$};\node at (270:0.62) {$\PP^2$};
 \node[right] at (2.3,0.1) {\parbox{6.3cm}{2-cells $\leftrightarrow$ components $\PP^2$ of $Y_0$\\
 interior edges $\leftrightarrow$ double lines of $Y_0$\\
 boundary edges (red) $\leftrightarrow$ lines of $D_0$\\
 central vertex $\leftrightarrow$ the point $[0:0:0:1]$}};
\end{tikzpicture}
\caption{The degenerate fiber $Y_0=\{xyz=0\}\subset\PP^3$ of
Example~\ref{ex:cubic-degeneration}, encoded by the polyhedral ball
$\{q\in|\mathscr G|:h(q)\le1\}$, where $\mathscr G$ is the boundary fan of the
octant and $h=x+y+z$. In general, the maximal cells of the ball correspond to
the toric components of $Y_0$, and the boundary sphere $h=1$ corresponds to
$D_0$.}\label{fig:ball}
\end{figure}

\subsection{The canonical module of the Rees algebra}

\begin{lemma}\label{lem:rees-canonical}
Let $A$ be a normal domain and $F$ a projective positive filtration of $A$ such
that $B=\gr_FA$ is reduced and Gorenstein with canonical module generated in
degree $0$. Put $R=\cR_F(A)$, $X=\Proj R$ and $D=\Proj B$. Then:
\begin{enumerate}[label=\textup{(\roman*)}]
\item $R$ is a normal Gorenstein domain, and $A$ is Gorenstein;
\item $\omega_R\simeq R(-1)$; more precisely, every element $\Theta\in\omega_R$
of degree one whose image under
$\omega_R\to\omega_R/t\omega_R\simeq\omega_B(-1)$ generates $\omega_B$ is a
generator of $\omega_R$;
\item $X$ is normal and Cohen--Macaulay, $D$ is a reduced Weil divisor, and
$\omega_X\simeq\mathcal O_X(-D)$ as rank-one reflexive sheaves. In particular
$K_X+D\sim0$, trivialized by the rational top form
\begin{equation}\label{eq:projective-form}
 \Omega_X=\iota_{\mathcal E}(\Theta/t),
\end{equation}
where $\mathcal E$ is the Euler vector field of the grading on $\Spec R$,
i.e.\ the vector field generating the $\mathbb G_m$-action defined by the
grading.
\end{enumerate}
\end{lemma}

\begin{proof}
(i) The ring $R[t^{-1}]=A[t^{\pm1}]$ is normal, $t$ is a nonzerodivisor, and
$R/tR=B$ is reduced and Cohen--Macaulay. Hence $R$ satisfies $S_2$: along
$V(t)$ because the depth of $R$ is one more than the depth of $R/tR$, and away
from $V(t)$ because $A[t^{\pm1}]$ is normal. If $\mathfrak p\ni t$ is a prime of
height one, then $(R/tR)_{\mathfrak p}$ is a reduced local ring of dimension
zero, i.e.\ a field, so $\mathfrak pR_{\mathfrak p}=tR_{\mathfrak p}$ and
$R_{\mathfrak p}$ is a discrete valuation ring. By Serre's criterion, $R$ is
normal. The ring $R$ is Gorenstein at every prime containing $t$, because
$R/tR$ is Gorenstein and $t$ is a nonzerodivisor. The non-Gorenstein locus of
$\Spec R$ is closed and stable under the $\mathbb G_m$-action defined by the
grading; since all degrees are positive, the closure of every orbit contains
the vertex, which is a Gorenstein point. Hence $R$ is Gorenstein, and so is
$A=R/(t-1)$.

(ii) By Fact~\ref{fact:canonical}(i), $\omega_R/t\omega_R\simeq\omega_B(-1)$.
Let $\Theta$ be as in the statement and $M$ the cokernel of
$R(-1)\to\omega_R$, $1\mapsto\Theta$. Then $M=tM$, and $M$ is graded and
bounded below, so $M=0$. The kernel is a torsion submodule of $R(-1)$, hence
zero. So $\omega_R=R\Theta$.

(iii) Normality and Cohen--Macaulayness of $X$ follow from (i), by passing
to a Veronese subring generated in degree one and taking its affine charts.
The quotient $R/tR=B$ is reduced, so its projective spectrum is reduced.
By
Fact~\ref{fact:veronese}, on a Veronese subring $R^{(e)}$ generated in degree
one, the canonical module is $\bigoplus_j\Theta R_{ej-1}$ and the ideal of $D$
is $\bigoplus_jtR_{ej-1}$. The isomorphism $ta\mapsto\Theta a$ of these modules
gives $\omega_X\simeq\mathcal O_X(-D)$ as rank-one reflexive sheaves.
Viewed as a rational top form on the cone, $\Theta/t$ is invariant under the
grading action; its contraction with $\mathcal E$ is horizontal and descends
to the form $\Omega_X$ in \eqref{eq:projective-form}. On $t\ne0$ this is
ordinary descent along the free grading action; the displayed module
isomorphism extends it across the boundary and gives
$\divi(\Omega_X)=-D$.
\end{proof}

\begin{remark}\label{rem:weighted-proj}
The reflexive formulation in (iii) is necessary: a Gorenstein graded ring
need not have Gorenstein projective spectrum. For example, filter
$A=\C[x,y,z]/(xyz-1)$ by weights $(1,1,2)$. Then
$\gr_FA=\C[x,y,z]/(xyz)$ has canonical generator of degree zero, and
$R=\C[t,X,Y,Z]/(XYZ-t^4)$, with weights $(1,1,1,2)$, satisfies
$\omega_R\simeq R(-1)$. The chart at $[0:0:0:1]$ on $\Proj R$ is the quotient of
$XY=t^4$ by $(t,X,Y)\mapsto(-t,-X,-Y)$. Writing
$(X,Y,t)=(a^4,b^4,ab)$ identifies it with the cyclic quotient
$\A^2/\mu_8$ of weights $(1,3)$. The toric character lattice is
$M=\{(i,j)\in\Z^2:i+3j\equiv0\pmod8\}$; since $(1,1)\notin M$ but
$(2,2)\in M$, the canonical index is two. Thus $D$ and $K_X$ need not be Cartier,
although $K_X+D\sim0$.
\end{remark}

For the cubic surface, $R=\C[t,X,Y,Z]/(\Phi)$ with $\Phi$ the cubic equation
(Example~\ref{ex:rees-cubic}), all variables of degree one, so
$\omega_R\simeq R(-4+3)=R(-1)$ by Fact~\ref{fact:canonical}(i); on the affine
part $t=1$, the form $\Omega_X$ is the residue of
$\frac{dx\wedge dy\wedge dz}{\Phi(1,x,y,z)}$, which has a simple pole along each
line of $D_\tau$.

\subsection{Projective cones over toric-face rings}
Let $\mathscr G$ be a rational fan in $L_\R$ and $h$ an integral linear form
positive on $|\mathscr G|\setminus0$. Suppose that the cones of $\mathscr G$
are mapped, by a homeomorphism $|\mathscr G|\to V$ which is linear on each
cone, onto the cones of a complete fan in a real vector space $V$ of dimension
$n$ (as in Proposition~\ref{prop:canonical-partition}(iv)). Two properties of
complete fans will be used: every maximal cone of $\mathscr G$ has dimension
$n$, and every cone of dimension $n-1$ is a face of exactly two maximal cones.
Put
\[
 B=\C[\mathscr G],\qquad Y_0=\Proj B[u],\qquad D_0=\{u=0\}=\Proj B,
\]
where $x^q$ has degree $h(q)$ and $u$ has degree $1$. For a maximal cone
$\sigma$, let $Y_\sigma=\Proj\C[\sigma\cap L][u]$; this is the normal
projective toric variety of the polytope $\delta_\sigma=\{q\in\sigma:h(q)\le1\}$
(Fact~\ref{fact:semigroup}). The polytopes $\delta_\sigma$ fit together into a
ball, as in Figure~\ref{fig:ball}.

\begin{lemma}\label{lem:stable-toric}
\begin{enumerate}[label=\textup{(\roman*)}]
\item The scheme $Y_0$ is reduced, its irreducible components are the
varieties $Y_\sigma$, and its normalization is $\coprod_\sigma Y_\sigma$.
\item Every toric divisor of a component $Y_\sigma$ is of one of two kinds:
either the divisor corresponds to a facet $\tau$ of $\sigma$, is contained
in exactly one other component $Y_{\sigma'}$ (with $\tau=\sigma\cap\sigma'$),
and $Y_0$ has a node $\{ab=0\}$ at the generic point of the divisor; or the
divisor is $\{u=0\}\subset Y_\sigma$, is a component of $D_0$, and $Y_0$ is
smooth at the generic point of the divisor.
\end{enumerate}
\end{lemma}

\begin{proof}
(i) For a monomial $x^qu^a$, let $J(q)$ be the set of maximal cones containing
$q$. Restriction of monomials embeds $B[u]$ into
$\prod_\sigma\C[\sigma\cap L][u]$: a monomial of $B[u]$ survives in the factor
indexed by any $\sigma\in J(q)\ne\varnothing$, and products are computed
factorwise. Hence $B[u]$ is reduced, with minimal primes the kernels of the
projections to the factors, which correspond to the components $Y_\sigma$. The
normalization of $Y_0$ is therefore the disjoint union of the normal varieties
$Y_\sigma$.

(ii) The toric divisors of $Y_\sigma$ correspond to the facets of the polytope
$\delta_\sigma$. The facets through $0$ are the truncations of the facets $\tau$
of $\sigma$; each such $\tau$ is a face of exactly one other maximal cone
$\sigma'$, and the corresponding divisor is the intersection
$Y_\sigma\cap Y_{\sigma'}=\Proj\C[\tau\cap L][u]$. Near the generic point of this
divisor, only the two components $Y_\sigma$ and $Y_{\sigma'}$ are present; both
are smooth there, and the local ring of $Y_0$ is the fiber product of their
local rings over the common residue field (by the description of $B[u]$ in
(i)), which after completion is $K[[a,b]]/(ab)$, a node. The remaining facet,
$\{h=1\}\cap\sigma$, corresponds to the divisor $\{u=0\}$ of $Y_\sigma$; near
its generic point only the component $Y_\sigma$ is present, and $Y_\sigma$ is
smooth there.
\end{proof}

In Example~\ref{ex:cubic-degeneration}, the facets of type $\tau$ give the three
double lines and the facets $\{h=1\}$ give the three lines of $D_0$.

\subsection{The criterion}

\begin{theorem}\label{thm:AB}
Let $A$ be a normal affine domain and $F$ a projective positive filtration of
$A$. Suppose that $\gr_FA\simeq\C[\mathscr G]$ as graded rings, where the fan
$\mathscr G$ and the grading $h$ are as in the previous subsection, and that
$\gr_FA$ is Gorenstein with canonical module generated in degree $0$. Then the
pair $(X,D)=(\Proj\cR_F(A),\Proj\gr_FA)$ of Lemma~\ref{lem:rees-canonical} is
log canonical.
\end{theorem}

\begin{proof}
Put $R=\cR_F(A)$, so that $R/tR=B=\C[\mathscr G]$. By
Lemma~\ref{lem:rees-canonical}, $R$ is a normal Gorenstein domain with
$\omega_R\simeq R(-1)$.

\emph{Step 1: the degeneration.} Let
\[
 \mathscr R=R[x,y]/(xy-t),
\]
graded by giving the elements of $R_d$ degree $d$, $x$ degree $0$ and $y$
degree $1$ (so that $xy-t$ is homogeneous), and let
$\rho:\mathcal Y=\Proj_{\C[x]}\mathscr R\to\A^1_x$ be the relative Proj over the
degree-zero part $\C[x]$; this is the projective family whose fiber over
$c\in\A^1$ is $\Proj\bigl(\mathscr R/(x-c)\bigr)$. Let $\mathcal D=\{y=0\}\subset\mathcal Y$. We record the
properties of this family.
\begin{itemize}
\item \emph{Fibers over $c\ne0$.} Setting $x=c$ gives $y=t/c$ and
$\mathscr R/(x-c)\simeq R$; so the fiber is $X$, and $\mathcal D$ restricts to
$\{t=0\}=D$. More precisely, $\mathscr R[x^{-1}]=R[x^{\pm1}]$, so
$\mathcal Y\setminus\rho^{-1}(0)\simeq X\times\mathbb G_m$ and
$\mathcal D\setminus\rho^{-1}(0)\simeq D\times\mathbb G_m$.
\item \emph{The fiber over $0$.} $\mathscr R/(x)=R[y]/(t)=B[y]$, so the fiber is
$Y_0=\Proj B[y]$, and $\mathcal D\cap Y_0=\{y=0\}=D_0$, with $u=y$ in the
notation of Lemma~\ref{lem:stable-toric}.
\item \emph{A $\mathbb G_m$-action.} For $\lambda\in\C^*$, the automorphism
$x\mapsto\lambda x$, $y\mapsto\lambda^{-1}y$ of $\mathscr R$ (fixing $R$)
preserves the grading, so defines an automorphism of $\mathcal Y$ which
preserves $\mathcal D$ and maps the fiber over $c$ to the fiber over $\lambda c$.
\item \emph{$\mathcal Y$ is normal and Cohen--Macaulay.} The ring
$\mathscr R$ is a hypersurface in the Gorenstein ring $R[x,y]$, hence
Gorenstein. The projective charts of $\mathcal Y$ are Cohen--Macaulay: for a homogeneous
element $a$ of positive degree $e$, the chart is the quotient of
$\Spec(\mathscr R/(a-1))$ by $\mu_e$; the hypersurface is Cohen--Macaulay,
and the finite quotient is Cohen--Macaulay by
Fact~\ref{fact:finite-quotient-cm}. Normality of $\mathscr R$, and hence of
$\mathcal Y$, follows from Serre's criterion exactly as in
Lemma~\ref{lem:rees-canonical}(i), using that $x$ is a nonzerodivisor,
$\mathscr R/x\mathscr R=B[y]$ is reduced and Cohen--Macaulay, and
$\mathscr R[x^{-1}]=R[x^{\pm1}]$ is normal.
The morphism $\rho$ is flat: $\mathscr R$ is torsion-free over $\C[x]$,
and the degree-zero parts of its homogeneous localizations are likewise
flat over $\C[x]$.
\item \emph{$K_{\mathcal Y}+\mathcal D+Y_0\sim0$.} Grade $R[x,y]$ by $\Z^2$,
giving $R_d$ degree $(d,d)$, $x$ degree $(1,0)$ and $y$ degree $(0,1)$; this grading
is positive. The canonical module of $R[x,y]$ is generated by
$\Theta\wedge dx\wedge dy$, of degree $(1,1)+(1,0)+(0,1)=(2,2)$, and $xy-t$ is a
nonzerodivisor of degree $(1,1)$; by Fact~\ref{fact:canonical}(i),
$\omega_{\mathscr R}$ is generated in degree $(1,1)$, i.e.\ in degree $1$ for the
grading defining $\mathcal Y$. Since $\mathscr R/y\mathscr R=B[x]$ is reduced,
the argument of Lemma~\ref{lem:rees-canonical}(iii) gives
$\omega_{\mathcal Y}\simeq\mathcal O_{\mathcal Y}(-\mathcal D)$ as
rank-one reflexive sheaves. The relative Veronese formula over $\C[x]$
follows from finite duality for the cyclic grading action: the invariant
part of the natural canonical module is the canonical module of the
invariant ring. Concretely, on a chart defined by a homogeneous element $a$
of degree $e$, use the cover $a=1$ and take $\mu_e$-invariants. The canonical
module and the $y$-boundary ideal have the same invariant pieces, as in
Lemma~\ref{lem:rees-canonical}(iii). Thus the boundary-ideal calculation
applies on the relative projective charts. Finally
$Y_0=\divi(x)$ is a principal Cartier divisor. Hence
$\mathcal O(K_{\mathcal Y}+\mathcal D+Y_0)\simeq\mathcal O_{\mathcal Y}$; let
$\beta$ be a generator, a rational top form with
$\divi(\beta)=-\mathcal D-Y_0$.
\end{itemize}

\emph{Step 2: the different on $Y_0$.} On
$\mathcal Y\setminus\mathcal D$, the preceding reflexive identity makes
$\omega_{\mathcal Y}$ invertible. Since $\mathcal Y$ is Cohen--Macaulay,
it is Gorenstein there, and Lemma~\ref{lem:adjunction} applies to the
Cartier divisor $Y_0$ on this open set. We claim that, on each component
$Y_\sigma$ of the normalization of $Y_0$, the residue $\operatorname{Res}_{Y_0}\beta$
has divisor exactly $-\partial Y_\sigma$, the full toric boundary. Since
$\operatorname{Res}_{Y_0}\beta$ is a rational form on $Y_\sigma$ whose only
possible zeros and poles lie on $\partial Y_\sigma$ (away from
$\mathcal D\cup\operatorname{Sing}Y_0$, the form $\beta$ is a generator of
$\omega_{\mathcal Y}(Y_0)$, whose residue generates $\omega_{Y_0}$ by
Lemma~\ref{lem:adjunction}), it suffices to compute the order along each toric
divisor, using Lemma~\ref{lem:stable-toric}(ii).
\begin{itemize}
\item At the generic point of a divisor of the first kind, $\mathcal D$ is
absent, $\beta$ generates $\omega_{\mathcal Y}(Y_0)$, and $Y_0$ is a node there.
By Lemma~\ref{lem:adjunction}, $\operatorname{Res}\beta$ generates
$\omega_{Y_0}$, and its restriction to the branch $Y_\sigma$ has a simple pole
along the divisor.
\item At the generic point of a divisor of the second kind, $Y_0$ is smooth,
hence so is $\mathcal Y$ (a regular local ring modulo a nonzerodivisor
forces the ambient local ring to be regular). At this codimension-two point
$\mathcal D$ is Cartier, and its equation modulo the equation of $Y_0$
is a parameter, since the scheme-theoretic intersection
$\mathcal D\cap Y_0=D_0$ is reduced. Thus $Y_0$ and $\mathcal D$ meet
transversally. In local coordinates with
$Y_0=\{x=0\}$ and $\mathcal D=\{v=0\}$ we have
$\beta=(\text{unit})\cdot\frac{dx}{x}\wedge\frac{dv}{v}\wedge\gamma$, whose residue
has a simple pole along $\{v=0\}$.
\end{itemize}
Hence the different of $(\mathcal Y,Y_0+\mathcal D)$ along $Y_0$ is
$\coprod_\sigma\partial Y_\sigma$. By Fact~\ref{fact:toric}, each toric pair
$(Y_\sigma,\partial Y_\sigma)$ is log canonical.

\emph{Step 3: conclusion.} By inversion of adjunction
(Fact~\ref{fact:inversion}), the pair $(\mathcal Y,Y_0+\mathcal D)$ is log
canonical on a neighborhood $\mathcal N$ of $Y_0$, and so is
$(\mathcal Y,\mathcal D)$, which has a smaller boundary. Since $\rho$ is proper,
$\rho(\mathcal Y\setminus\mathcal N)$ is closed and does not contain $0$.
The complement of $\rho(\mathcal Y\setminus\mathcal N)$ is an open
neighborhood $W$ of $0$ with
$\rho^{-1}(W)\subset\mathcal N$. Choose $c\in W\setminus\{0\}$.
Every point of $\mathcal Y\setminus Y_0=X\times\mathbb G_m$ is moved into
$\rho^{-1}(c)$ by the $\mathbb G_m$-action, which preserves
$\mathcal D$; hence $(X\times\mathbb G_m,D\times\mathbb G_m)$ is log canonical,
and therefore $(X,D)$ is log canonical (Section~\ref{subsec:bg-pairs}).
\end{proof}

\begin{remark}
Recall that a reduced ring is \emph{seminormal} if every element $f$ of its
total ring of fractions with $f^2$ and $f^3$ in the ring belongs to the ring.
For an equidimensional $S_2$ seminormal scheme $Y_0$ which is nodal in
codimension one, a pair $(Y_0,D_0)$ is \emph{semi-log canonical} (slc) if
$K_{Y_0}+D_0$ is $\Q$-Cartier, the boundary contains no component of the
double locus, and the normalization with the conductor divisor and
the pulled-back boundary is lc
\cite[Def.--Lem.~5.10]{KollarSingMMP2013}. Here $Y_0$ is Cohen--Macaulay, hence $S_2$.
Seminormality follows from the monomial description: $B[u]$ is the inverse
limit of its normal component rings over their common-face rings. An element
$f$ of the total quotient ring whose square and cube lie in $B[u]$ belongs to
each normal component ring, and the restrictions of $f$ agree on every reduced
overlap because their squares and cubes agree. Hence $f$ belongs to $B[u]$; this proves seminormality, which passes to the projective charts.
The boundary $D_0$ meets the double locus in smaller codimension.
Thus Step~2 says that the central fiber $(Y_0,D_0)$ is slc with
$K_{Y_0}+D_0\sim0$, and Theorem~\ref{thm:AB} is an instance of the principle
that slc log Calabi--Yau pairs deform to lc ones.
Arg\"uz--Bousseau \cite[Thm.~3.7]{ArguzBousseau2026} prove a more general
criterion using the Cohen--Macaulay property and trivial logarithmic
dualizing sheaf of the projective toric-face fiber, which they establish
from the topology of its cone complex. In our situation the graded
Gorenstein computation in Theorem~\ref{thm:special-fiber} supplies these
dualizing properties directly.
\end{remark}

In Example~\ref{ex:cubic-degeneration}, the family $\mathcal Y$ is the pencil
of cubic surfaces
\[
 \bigl\{\,s\,t(x^2+y^2+z^2)=xyz+(\tau+2)s^3t^3\,\bigr\}\subset\PP^3_{[x:y:z:t]},
 \qquad s\in\A^1,
\]
whose fiber over $s\ne0$ is isomorphic to $X_\tau$ by
$[x:y:z:t]\mapsto[x:y:z:st]$, and
whose fiber over $s=0$ is $\{xyz=0\}$. Step~2 is the observation that on each
plane of $\{xyz=0\}$, say $\{z=0\}$ with affine coordinates $x/t$ and $y/t$, the
residue of the logarithmic volume form is a multiple of
$\frac{d(x/t)}{x/t}\wedge\frac{d(y/t)}{y/t}$, which has simple poles along the
two double lines $\{x=0\}$, $\{y=0\}$ and along the line $\{t=0\}$ of $D_0$.

\section{Central covers and spectral compactifications}
\label{sec:proof}\label{sec:main-proofs}

We pass from a simply connected group to $G$ through the central cover
$Z(G)^\circ\times G_\der^{\mathrm{sc}}\to G$. Integral descent requires an
invariant logarithmic volume form. The corner identity supplies it;
Proposition~\ref{prop:spectral-normal-suffices} then needs only normality
of the lifted spectral rings, supplied by
Proposition~\ref{prop:spectral-interior} under genericity and by
Fact~\ref{fact:gl-closure-normality} in type~$A$.

\subsection{The simply connected construction}

\begin{proposition}\label{prop:sc-pair}
Let $H$ be simply connected semisimple, let $2g+k\ge3$, and suppose the
nonempty spectral fiber $\Spec A_\bs$ is normal. Let $F$ be the uniform
filtration of Proposition~\ref{prop:spectral-flatness}. Then $\Spec A_\bs$
is connected, and
\[
 X_H=\Proj\cR_F(A_\bs),\qquad D_H=\Proj\gr_FA_\bs
\]
is a normal Cohen--Macaulay projective lc pair with
$X_H\setminus D_H=\Spec A_\bs$, $D_H$ reduced, and
$\omega_{X_H}(D_H)\simeq\mathcal O_{X_H}$. Moreover
$\omega_{\cR_F(A_\bs)}\simeq\cR_F(A_\bs)(-1)$, and for every generator
$\Theta$ of degree one the form
$\Omega_H=\iota_{\mathcal E}(\Theta/t)$ of \eqref{eq:projective-form}
satisfies $\divi(\Omega_H)=-D_H$.
\end{proposition}

\begin{proof}
Proposition~\ref{prop:spectral-flatness} gives connectedness and
Gorensteinness. Since a connected normal affine variety is integral,
$A_\bs$ is a normal domain. The remaining assertions follow from
Theorems~\ref{thm:special-fiber} and~\ref{thm:AB} and
Lemma~\ref{lem:rees-canonical}.
\end{proof}

For the once-punctured torus and $H=\SL_2$, the pair is the cubic surface
with its triangle of lines (\hyperref[subsec:intro-example]{the introductory model calculation}).

\subsection{Lift sectors}\label{subsec:sectors}\label{subsec:spectral-descent}
Let $p:\widetilde G=T_0\times H\to G$ be the central isogeny of
Fact~\ref{fact:isogeny}, with $T_0=Z(G)^\circ$, $H=G_\der^{\mathrm{sc}}$
and finite central kernel $F_p$. Since $\pi_1(\Sigma^\circ)\simeq F_m$ is free,
every representation in $G$ lifts to $\widetilde G$, in $|F_p|^m$ ways. Each
prescribed point $a_i\in G/\!\!/G$ has finitely many preimages
$(t_i,b_i)\in T_0\times(H/\!\!/H)$; a tuple $\mathbf b=((t_i,b_i))_i$ of such
preimages is \emph{admissible} if $\prod_it_i=1$ and the lifted surface relation
has a solution. For such $\mathbf b$, let $A_{\mathbf b}$ be the spectral ring
of $H$ with adjoint data $(b_i)_i$ (Section~\ref{subsec:puncture-eq}). Fix spectral data $\mathbf a$
and write $R=(\cR_B^{\mathrm{sp}}(G,\mathbf a))_{\red}\subset G^m$.
The pullback $\widehat\cR=R\times_{G^m}\widetilde G^m$ is reduced because
$p$ is \'etale. The deck group $F_p^m$ of this covering multiplies free generators by
central elements, commutes with conjugation, and gives
\begin{equation}\label{eq:quotient-cover}
 \cM_B^{\mathrm{sp}}(G,\mathbf a)
 =\bigl(\widehat\cR/\!/\widetilde G\bigr)/F_p^m.
\end{equation}

\begin{lemma}[Lift sectors]\label{lem:sectors}\label{lem:spectral-sectors}
Assume $2g+k\ge3$. The quotient upstairs in \eqref{eq:quotient-cover}
is a finite disjoint union of connected spaces
\begin{equation}\label{eq:spectral-sectors}
 V_{\mathbf b}=T_0^{2g}\times\Spec A_{\mathbf b},
\end{equation}
indexed by the admissible tuples $\mathbf b$. These summands are called
\emph{lift sectors}. Every connected component downstairs is
the quotient of one sector by its stabilizer in $F_p^m$.
If $G_\der$ is simply connected, the spectral character variety is
connected whenever nonempty. For regular semisimple data, the lifted
puncture fibers are regular semisimple classes of $H$, so
$\Spec A_{\mathbf b}=\cM_B(H,\bcC')$ in each sector.
\end{lemma}

\begin{idea}
The surface relation splits in $T_0\times H$: puncture torus values
must multiply to one, while the $2g$ handle torus values remain free.
The choices of lifted puncture parameters distinguish the sectors.
\end{idea}

\begin{proof}
Conjugation invariants in the central isogeny give
\begin{equation}\label{eq:spectral-base-cover}
 G/\!/G\simeq\bigl(T_0\times(H/\!/H)\bigr)/F_p.
\end{equation}
If $(t,b)$ runs through the finitely many points over $a_i$, then
\[
 p^{-1}\bigl((\chi_G^{-1}(a_i))_{\red}\bigr)
 =\coprod_{(t,b)\mapsto a_i}\{t\}\times\chi_H^{-1}(b).
\]
This is an equality of reduced schemes: $p$ is \'etale, and the $H$-fibers
are reduced by Fact~\ref{fact:steinberg}(iii). In $T_0\times H$ the
surface relation splits. The puncture torus values must multiply to $1$;
the handle torus values remain free, and the $H$ relation prescribes the
adjoint data $(b_i)_i$. A tuple is admissible when this lifted relation has a solution.

Imposing the invariant equations and taking invariants commute by
Fact~\ref{fact:git}(i). Reduction does not change the quotient summand:
Proposition~\ref{prop:spectral-flatness} makes $A_{\mathbf b}$ reduced,
so the nilradical of the representation ring has no nonzero invariant.
Exactness of invariants gives \eqref{eq:spectral-sectors}, and the same
proposition gives connectedness. The deck group permutes the open and
closed sectors. The orbits of the deck group give connected components of the quotient;
each component is the quotient of one sector by its stabilizer.

For connectedness downstairs, put $K=F_p\cap(\{1\}\times H)$.
Multiplying the lifts of $\gamma_1,\dots,\gamma_{k-1}$ by
$f_1,\dots,f_{k-1}\in F_p$ multiplies the lift of $\gamma_k$ by
$(f_1\cdots f_{k-1})^{-1}$; handle deck transformations do not change
the puncture tuple. Two admissible tuples differ by $(f_i)\in F_p^k$,
and their torus relations force $\prod_i f_i\in K$.
If $G_\der$ is simply connected, then $K=1$. Thus $\prod_i f_i=1$,
and all admissible tuples lie in a single puncture-deck orbit.

Finally, for a regular semisimple class $\cC_i$, a lift
$\tilde s=(t,h)$ has the same centralizer Lie algebra as $p(\tilde s)$.
Hence $\dim Z_H(h)=\rank H$ and $h$ is regular semisimple.
The preimage $p^{-1}(\cC_i)$ is the union of the central translates of
$\{t\}\times H\cdot h$; the distinct such classes are closed and pairwise
disjoint. This proves the
regular-class specialization.
\end{proof}

\subsection{An invariant logarithmic volume form}\label{subsec:residue}
Let $V=T_0^{2g}\times V_H$ be a lift sector with stabilizer $\Gamma_V$, where
$V_H=\Spec A_\bs$, assumed normal. The group $\Gamma_V$ acts on $V_H$
through its image $\Delta$ in $Z(H)^m$, by multiplying the free generators by
central elements of $H$.

Since $K_{X_H}+D_H\sim0$, the group $\Delta$ acts on the line
$H^0(X_H,\omega_{X_H}(D_H))$ by a character. Integral descent requires
this character to be trivial. Even on $\PP^1$, the map $z\mapsto\zeta z$
fixes $dz/z$ but multiplies $dz$ by $\zeta$; logarithmic and regular
forms behave differently. The corner identity \eqref{eq:global-corner}
gives the required cancellation.

\begin{lemma}\label{lem:invariant-residue}
The uniform filtration $F$ of Proposition~\ref{prop:spectral-flatness}
is $\Delta$-stable, and then $\Delta$ fixes a generator $\Theta$ of degree one of
$\omega_{\cR_F(A_\bs)}$. Hence
$\Omega_H=\iota_{\mathcal E}(\Theta/t)$ is a $\Delta$-invariant rational top
form on $X_H$ with $\divi(\Omega_H)=-D_H$.
\end{lemma}

\begin{proof}
\emph{An invariant filtration.} In the gauge of Lemma~\ref{lem:gauge},
multiplying the free generators by central elements multiplies edge transports
by central elements, which act by scalars on each summand of
\eqref{eq:graph-tensors} and preserve the local bases of
Theorem~\ref{thm:local-cone}. So $\Delta$ preserves $v_\Gamma$ and acts on each
leaf through a character. The function $\chi_{\varpi_a}(\gamma_i)$ is an
eigenvector, say of character $\xi_{i,a}$; since $\Delta$ preserves the prescribed
values, $\xi_{i,a}(\delta)\chi_{\varpi_a}(s_i)=\chi_{\varpi_a}(s_i)$ for all
$\delta\in\Delta$, so $f_{i,a}$ is an eigenvector of character $\xi_{i,a}$ as
well, including when $\chi_{\varpi_a}(s_i)=0$. Proposition~\ref{prop:spectral-flatness} chooses the ambient filtration
$Z(H)^m$-stable, before any values are prescribed. The quotient filtration
on $A_\bs$ is therefore $\Delta$-stable and has the properties of
Theorem~\ref{thm:special-fiber}.

\emph{Cancellation of characters on the special fiber.} Let $S=\C[Q_\Gamma]$
and $I=(x^{p_{i,a}})$, so that $S/I=\gr_FA_\bs$. The group $\Delta$ acts on
$S$ diagonally on monomials, i.e.\ through the torus $\Spec\C[\gp Q_\Gamma]$,
which fixes the logarithmic volume form $\eta$ of that torus. Hence the
generator $x^{\Omega_\Gamma}\eta$ of $\omega_S$ (Fact~\ref{fact:semigroup},
Lemma~\ref{lem:glue}) transforms by the character of $x^{\Omega_\Gamma}$,
which is $\prod_{i,a}\xi_{i,a}$ by \eqref{eq:global-corner}. By adjunction,
$\omega_{S/I}\simeq\omega_S\otimes_S\det(I/I^2)^\vee$, equivariantly, and
$\det(I/I^2)$ is generated by the wedge product of the classes of the
$x^{p_{i,a}}$, which transforms by the same character
$\prod_{i,a}\xi_{i,a}$. Hence $\omega_{S/I}$ has a $\Delta$-invariant generator
of degree $0$.

\emph{Lifting.} The surjection $\omega_R\to\omega_{S/I}(-1)$ of
Lemma~\ref{lem:rees-canonical}, with $R=\cR_F(A_\bs)$, is
$\Delta$-equivariant, since $\Delta$ fixes $t$. Choose a homogeneous
degree-one lift of the invariant generator to $\omega_R$ and average over
$\Delta$; the result $\Theta$ is invariant and
still lifts the generator, hence generates $\omega_R$ by
Lemma~\ref{lem:rees-canonical}. The action commutes with the grading
and fixes $t$, so it also fixes the Euler contraction defining $\Omega_H$.
\end{proof}

\begin{lemma}\label{lem:descent}
Let a finite group $\Gamma$ act on a normal projective Cohen--Macaulay lc pair
$(Y,E)$ with $E$ reduced, fixing a rational top form $\Omega_Y$ with
$\divi(\Omega_Y)=-E$. Let $q:Y\to X=Y/\Gamma$ and $D=q(E)$. Then $X$ is normal,
projective and Cohen--Macaulay, $\Omega_Y=q^*\Omega_X$ for a rational top form
$\Omega_X$ with $\divi(\Omega_X)=-D$, and $(X,D)$ is lc. In particular
$K_X+D\sim0$.
\end{lemma}

\begin{proof}
The invariant form $\Omega_Y$ descends to a rational form $\Omega_X$ on $X$. Let
$F\subset Y$ be a prime divisor over a prime divisor $P\subset X$, with
ramification index $e$. In characteristic zero,
$\ord_F(q^*\Omega_X)=e\cdot\ord_P(\Omega_X)+e-1$. If $F\subset E$, the left side
is $-1$, which forces $\ord_P\Omega_X=-1$. If $F\not\subset E$, the left side is
$0$, which forces $e=1$ and $\ord_P\Omega_X=0$. Hence $\divi(\Omega_X)=-D$ and
$K_Y+E=q^*(K_X+D)$, and $(X,D)$ is lc by Fact~\ref{fact:pairs}(ii). The quotient
of a normal projective variety by a finite group is normal and projective. The
Cohen--Macaulay assertion follows from
Fact~\ref{fact:finite-quotient-cm}.
\end{proof}

\subsection{Compactification from normal lifted fibers}

\begin{proposition}[Normal spectral fibers suffice]
\label{prop:spectral-normal-suffices}
In the stable range $2g+k\ge3$, suppose all the nonempty lifted
$H$-character fibers $\Spec A_{\mathbf b}$ of Lemma~\ref{lem:spectral-sectors}
are normal. Then every connected component of the reduced spectral
$G$-character variety admits a normal Cohen--Macaulay projective lc
compactification $(X,D)$ with reduced boundary and $K_X+D\sim0$ integrally.
\end{proposition}

\begin{proof}
Write the component as $U=V/\Gamma_V$ with
$V=T_0^{2g}\times\Spec A_{\mathbf b}$. The $H$-part of $\Gamma_V$ acts
through a subgroup $\Delta\subset Z(H)^m$ preserving the prescribed values.
Proposition~\ref{prop:sc-pair} and Lemma~\ref{lem:invariant-residue}
give $(X_H,D_H,\Omega_H)$ with a $\Delta$-stable filtration and an invariant
form. Let $(Y_T,E_T)$ be $(\PP^1)^{2g\dim T_0}$ with its toric boundary
and logarithmic volume form $\Omega_T$, and put
\[
 (Y,E)=(Y_T\times X_H,\ E_T\times X_H+Y_T\times D_H),\qquad
 \Omega_Y=\Omega_T\wedge\Omega_H.
\]
Then $Y\setminus E=V$, the pair is normal, projective, Cohen--Macaulay
and lc by Fact~\ref{fact:pairs}(iii), and $\divi(\Omega_Y)=-E$.
The torus translations in $\Gamma_V$ extend to $Y_T$ and fix $\Omega_T$;
the action on $X_H$ fixes $\Omega_H$. Lemma~\ref{lem:descent} now gives
the required quotient pair, including integral linear equivalence.
\end{proof}

\begin{example}[$G=\GL_3$]\label{ex:gl3}
The central cover is $\C^*\times\SL_3\to\GL_3$, $(z,h)\mapsto zh$, with kernel
$\{(\zeta^{-1},\zeta I):\zeta^3=1\}$. Take $(g,k)=(1,1)$ and the puncture
eigenvalues $2,3,1/6$ of Section~\ref{subsec:sl3-torus}, and let $U_{\SL_3}$ be
the corresponding $\SL_3$ relative character variety. Write the $\GL_3$
holonomies as $z_1A$ and $z_2B$ with $A,B\in\SL_3$. The scalars cancel in the
commutator, so the lifted variety is $(\C^*)^2\times U_{\SL_3}$, and the $\GL_3$
character variety is the quotient by
$(\zeta,\eta)\cdot(z_1,z_2,A,B)=(\zeta^{-1}z_1,\eta^{-1}z_2,\zeta A,\eta B)$,
$\zeta^3=\eta^3=1$. The form $\frac{dz_1}{z_1}\wedge\frac{dz_2}{z_2}$ on
$(\PP^1)^2$ is invariant, Lemma~\ref{lem:invariant-residue} provides an
invariant logarithmic form on the compactification of $U_{\SL_3}$, and the
quotient therefore satisfies $K_X+D\sim0$.
\end{example}

\subsection{Proof of Theorem~\ref{thm:spectral-main}}\label{subsec:proofA}

The interior result supplies normality, $\Q$-factoriality and klt
singularities. The degeneration and descent supply a Cohen--Macaulay compactification
and integral $K_X+D\sim0$. The restriction of this relation to $U$ makes $U$ Gorenstein; the
integral discrepancies then upgrade klt to canonical, allowing the dlt
modification to preserve the interior.

\begin{proof}[Proof of Theorem~\ref{thm:spectral-main}]
\emph{The torus case.} If $G$ is a torus, compatibility identifies the spectral
character variety with $G^{2g}\simeq(\C^*)^{2g\dim G}$, and the scheme-theoretic
peripheral fiber is already reduced. Compactify by $(\PP^1)^{2g\dim G}$ with its
toric boundary; the product of the logarithmic forms $dz/z$ has divisor minus
the boundary, and this smooth dlt pair satisfies all the conclusions.

\emph{The stable range.} If $G$ is not a torus, generic nonempty data force
$2g+k\ge3$. Indeed, if $g=0$ and $k\le2$, a representation is generated by at
most one element, so the image lies in a Borel subgroup, and the surface
relation contradicts Lemma~\ref{lem:spectral-genericity}.

\emph{Lifting.} Let $p:T_0\times H\to G$ be the central isogeny of
Fact~\ref{fact:isogeny}. Every lifted representation is $H$-irreducible, since
a proper parabolic subgroup $Q_H\subset H$ gives a proper parabolic subgroup
$p(T_0\times Q_H)\subset G$. Moreover the $H$-data $(b_i)_i$ in every admissible
tuple $\mathbf b$ of lifted parameters (Lemma~\ref{lem:spectral-sectors}) are
intrinsically $H$-generic: if lifted elements $x_i\in(\chi_H^{-1}(b_i))_{\red}\cap Q_H$
satisfied $\prod_ix_i\in[Q_H,Q_H]$, then, since admissibility gives
$\prod_it_i=1$ for the torus components $t_i$, the product of the elements
$p(t_i,x_i)\in(\chi_G^{-1}(a_i))_{\red}$ would lie in the commutator
subgroup of $p(T_0\times Q_H)$, contrary to genericity of $\mathbf a$.

\emph{Compactification.} By Proposition~\ref{prop:spectral-interior}, every
lifted spectral character variety $\Spec A_{\mathbf b}$ is normal, so
Proposition~\ref{prop:spectral-normal-suffices} gives a normal Cohen--Macaulay
projective lc compactification $(X,D)$ of $U$ with reduced boundary and
$K_X+D\sim0$.

\emph{The interior.} By Proposition~\ref{prop:spectral-interior}, $U$ is
$\Q$-factorial and klt. Restricting the logarithmic volume form to $U$ gives
$K_U\sim0$, so the canonical sheaf of $U$ is invertible; since $U$ is open in
the Cohen--Macaulay variety $X$, the variety $U$ is Gorenstein. The
discrepancies of $U$ are therefore integers, and being greater than $-1$ they are
nonnegative: $U$ has canonical singularities. Corollary~\ref{cor:dlt-open} now
gives the projective $\Q$-factorial dlt modification which is an isomorphism over
$U$.

\emph{Simply connected derived group.} Connectedness is
Lemma~\ref{lem:spectral-sectors}. If $G_\der$ is simply connected, then
$F_p\cap(\{1\}\times H)=1$, so $F_p$ maps injectively to $T_0$ and acts freely
on $T_0\times(H/\!\!/H)$ by translations in the first factor. Hence the map
\eqref{eq:spectral-base-cover} of parameter spaces is finite \'etale, and the
scheme-theoretic peripheral conditions pull back to the conditions defining
disjoint unions of scheme-theoretic fibers of $\chi_H$, which are reduced by
Fact~\ref{fact:steinberg}(iii). By Lemma~\ref{lem:spectral-relation-local}
and Proposition~\ref{prop:spectral-interior}, the lifted representation schemes
are reduced; by finite \'etale descent and exactness of invariants, the
scheme-theoretic peripheral fiber of $G$ is reduced, i.e.\ equal to
$\cM_B^{\mathrm{sp}}(G,\mathbf a)$.
\end{proof}

\subsection{Regular semisimple classes and arbitrary spectra in type \texorpdfstring{$A$}{A}}
\label{subsec:proof-corollaries}

\begin{proof}[Proof of Corollary~\ref{cor:regular-semisimple}]
For regular semisimple $s_i$, the identity component of $Z_G(s_i)$ is a torus,
whose unipotent variety is a point. Hence
Lemma~\ref{lem:spectral-fiber-local} identifies each reduced adjoint fiber with
the class $\cC_i$, and $\cM_B(G,\bcC)=\cM_B^{\mathrm{sp}}(G,\mathbf a)$. All
conclusions of Theorem~\ref{thm:spectral-main} apply, and Corollary~\ref{cor:symplectic-coarse} gives
the symplectic local models.
\end{proof}

For groups of type $A$ we use a normality theorem which holds without any
genericity. A \emph{quiver} is a directed graph, and a \emph{multiplicative
quiver variety} is a quotient of collections of linear maps along the arrows
and reversed arrows of a quiver, satisfying prescribed matrix-product
relations, by changes of basis at the vertices \cite[\S1]{CrawleyBoeveyShaw2006} and \cite[Def.~3.1]{Yamakawa2008}; see also \cite[\S2]{SchedlerTirelli2022}.
Character varieties with prescribed class closures at the punctures are open
pieces of such varieties.

\begin{fact}[Normality for $\GL_n$ class closures]
\label{fact:gl-closure-normality}
For arbitrary conjugacy classes $C_1,\dots,C_k\subset\GL_n$, the reduced affine
variety
\[
 \Bigl\{(A_j,B_j,M_i)\in\GL_n^{2g}\times\textstyle\prod_i\overline{C_i}:
 \prod_j[A_j,B_j]\prod_iM_i=1\Bigr\}_{\red}/\!\!/\GL_n
\]
is normal.
\end{fact}

This fact combines the identification with an open subset of a multiplicative
quiver variety \cite[Thm.~3.6, Rem.~3.7]{SchedlerTirelli2022} and the normality
of all reduced multiplicative quiver varieties
\cite[Cor.~5.5]{KaplanSchedler2023}. The identification on coordinate rings uses
the quotient by the auxiliary groups along the legs of the quiver
\cite[Lem.~9.1]{CrawleyBoevey2003}, the change of coordinates in
\cite[Prop.~2]{CrawleyBoevey2013} and \cite[Lem.~3.5]{SchedlerTirelli2022}, and exactness of invariants, which allows
the surface relation to be imposed before the remaining $\GL_n$-quotient is
taken. The statement concerns reduced varieties; it says nothing about
non-reduced scheme-theoretic fibers.

\begin{corollary}[Arbitrary spectra in type $A$]
\label{cor:spectral-type-a}
Suppose that the simply connected cover of $G_\der$ is a product of special
linear groups. For every compatible $\mathbf a\in(G/\!\!/G)^k$, every connected
component $U$ of the nonempty spectral character variety admits a normal
Cohen--Macaulay projective lc compactification $(X,D)$ with reduced boundary and
$K_X+D\sim0$. If moreover $U$ is $\Q$-factorial and klt, then the dlt
conclusion of Theorem~\ref{thm:spectral-main} holds for $U$ as well.
\end{corollary}

\begin{idea}
Normality for $\GL_n$ follows from multiplicative quiver varieties.
A finite \'etale base change on handle determinants gives the
$\SL_n$ case, and Proposition~\ref{prop:spectral-normal-suffices}
then applies.
\end{idea}

\begin{proof}
\emph{$\GL_n$.} Lemma~\ref{lem:gl-spectral-closure} identifies each reduced
$\GL_n$ spectral fiber with the closure of a regular class, so
Fact~\ref{fact:gl-closure-normality} shows that the spectral $\GL_n$-character
variety is normal for all data.

\emph{$\SL_n$.} Consider the $\GL_n$ spectral character variety with the same
peripheral polynomials, which have constant term $(-1)^n$. The determinants of
the $2g$ handle matrices define a morphism to $(\mathbb G_m)^{2g}$, and the
pullback of this morphism along the coordinatewise $n$-th power map of
$(\mathbb G_m)^{2g}$ is
\begin{equation}\label{eq:spectral-sl-gl}
 (\mathbb G_m)^{2g}\times\cM_B^{\mathrm{sp}}(\SL_n,\mathbf a).
\end{equation}
Indeed, chosen $n$-th roots $\alpha_j,\beta_j$ of the handle determinants
replace $A_j,B_j$ by $\alpha_j^{-1}A_j,\beta_j^{-1}B_j$ and leave the puncture
matrices unchanged; the scalars cancel in each commutator. This identification
holds on invariant rings, since the base change is flat and $\GL_n$ and $\SL_n$
have the same conjugation invariants on these spaces. The $n$-th power map of
the torus is finite \'etale, so the pullback \eqref{eq:spectral-sl-gl} is
normal, and hence so is the factor $\cM_B^{\mathrm{sp}}(\SL_n,\mathbf a)$.
Products of special linear groups have products of these varieties as spectral
character varieties, which are again normal.

\emph{Conclusion.} In the stable range $2g+k\ge3$,
Proposition~\ref{prop:spectral-normal-suffices} applies to the central cover of
$G$. Outside the stable range, $\Sigma^\circ$ is a sphere with at most two
punctures, and the nonempty spectral character variety is a point, since the
closed orbits correspond to the unique semisimple class. The torus case is as in
Theorem~\ref{thm:spectral-main}. Finally, the constructed $U$ is Cohen--Macaulay with $K_U\sim0$; if
$U$ is klt, it is therefore canonical, and Corollary~\ref{cor:dlt-open}
applies.
\end{proof}

\begin{remark}\label{rem:spectral-final-scope}
For a general reductive group and nongeneric data, the argument needs exactly
one further input, isolated in Proposition~\ref{prop:spectral-normal-suffices}:
normality of the lifted spectral character varieties. The uniform degeneration
of Proposition~\ref{prop:spectral-flatness} does not supply normality by itself.
Without genericity we make no claim that the interior is $\Q$-factorial, nor
that a dlt modification can leave the whole interior unchanged.
\end{remark}

\part{The Poisson construction}\label{part:poisson}

\section{Poisson compactification and logarithmic nonvanishing}
\label{sec:poisson-bridge}
For a class closure smaller than a full adjoint fiber, restricting a spectral
compactification gives no control of $K_X+D$. This section builds a model
for arbitrary class closures, starting from the normalized interiors of
Section~\ref{new:sec-normalization}, and reduces log canonicity to a single
nonvanishing statement, in three steps. First, if $X$ carries a Poisson structure which is symplectic on $U$ and
preserves $D$, then one nonzero logarithmic pluricanonical form on a resolution
of $U$ forces $K_X+D\sim0$ and log canonicity
(Theorem~\ref{thm:poisson-nonvanishing}). Second, such Poisson compactifications
come from filtrations for which the bracket has degree zero, i.e.\
$\{F_a,F_b\}\subset F_{a+b}$. Third, matrix-coefficient filtrations of character
varieties have this property for every reductive group, and the conclusion
survives finite central covers. In the running example, Goldman's product formula
gives $\{x,y\}=c\,(2z-xy)$ for a nonzero constant $c$, an element of
$F_2=F_{1+1}$; so the total-degree filtration has degree zero.

\subsection{The logarithmic invariant used in the criterion}

Let $U$ be an irreducible complex variety.  Choose a projective resolution
$\rho:U'\to U$, and a smooth projective compactification $U'\subset Y$ with
simple-normal-crossings boundary $B=Y\setminus U'$.  Set
\[
 \overline P_m^{\mathrm{res}}(U)
   =h^0\bigl(Y,\mathcal O_Y(m(K_Y+B))\bigr),\qquad m\geq1.
\]
These numbers are independent of these choices: take a common resolution and
use invariance of logarithmic pluricanonical sections under proper birational
maps of smooth open varieties.  The corresponding Iitaka dimension is denoted
$\overline\kappa^{\mathrm{res}}(U)$.  For smooth $U$ this is its usual logarithmic
Kodaira dimension.  Divisors exceptional over the interior $U$ are \emph{not}
included in $B$.  Thus this convention does not artificially allow poles over
singular points of $U$.  If $U$ has canonical singularities and a trivialized
canonical divisor, its volume form has nonnegative order along those interior
exceptional divisors.

In the canonical, canonically trivial case one can equivalently use the usual
logarithmic plurigenus of $U_{\mathrm{reg}}$.  Indeed, a logarithmic pluriform
there is $f\Omega^{\otimes m}$.  At each prime divisor of $U$, regularity gives
$\operatorname{ord}(f)\geq0$, so normality extends $f$ to a regular function
on $U$.  Canonical singularities make the pullback of $\Omega$ regular on a
resolution.  Hence the pluriform has no poles along divisors exceptional over
$U$, even if those divisors were initially included in the boundary of a
compactification of $U_{\mathrm{reg}}$.  Comparison on a common resolution
proves equality of the two spaces of logarithmic pluriforms.

\begin{theorem}[Poisson compactification and nonvanishing]
\label{thm:poisson-nonvanishing}
Let $X$ be an irreducible normal projective complex variety of dimension
$d=2n>0$, let $U\subset X$ be a dense open subset, and suppose that
$D=(X\setminus U)_{\mathrm{red}}$ is a divisor.  Assume that:
\begin{enumerate}
\item $X$ carries a regular algebraic Poisson bracket;
\item the bracket is symplectic at every codimension-one point of $U$;
\item $D$ is a Poisson closed subset, meaning that its reduced ideal sheaf is
preserved by every Hamiltonian derivation.
\end{enumerate}
If $\overline P_m^{\mathrm{res}}(U)>0$ for one $m\geq1$, then
\[
 K_X+D\sim0
\]
integrally and $(X,D)$ is log canonical.  Moreover,
$\overline P_\ell^{\mathrm{res}}(U)=1$ for every $\ell\geq1$.
Without the nonvanishing hypothesis one still has
$\overline P_m^{\mathrm{res}}(U)\leq1$ for all $m$.  Consequently
\[
 \overline\kappa^{\mathrm{res}}(U)\in\{-\infty,0\}.
\]
\end{theorem}

\begin{idea}
Let $\pi$ be the Poisson bivector and $\Omega=(\pi^{\wedge n})^{-1}$ the dual
volume form. Since $D$ is Poisson, $\pi^{\wedge n}$ vanishes along every boundary
component, so $\Omega$ has a pole of order $e_i\ge1$ along each $D_i$. Every
logarithmic pluricanonical form is $f\,\Omega^{\otimes m}$ with $f$ regular in
codimension one, hence constant. A nonzero such form therefore forces $e_i=1$
for all $i$, i.e.\ $K_X+D\sim0$, and logarithmic regularity of
$\Omega^{\otimes m}$ on a log resolution is log canonicity. In the running
example, $\Omega$ restricts on $U_\tau$ to a constant multiple of
$dx\wedge dy/(2z-xy)$, which has a simple pole along each line of the triangle
$D_\tau$.
\end{idea}

\begin{proof}
Write $\pi$ for the Poisson bivector on the smooth locus of $X$, and let
$\Omega$ be its inverse top exterior power, normalized by
$\Omega(\pi^{\wedge n})=1$.  It is a rational top form on $X$.  Regularity of
$\pi$, together with its nondegeneracy in codimension one of $U$, gives
\[
 \operatorname{div}_X(\Omega)=-\sum_i e_iD_i,
 \qquad e_i\geq0,
\]
where $D_i$ are the prime components of $D$.  In fact $e_i\geq1$.  To see this,
work near the generic point of $D_i$, where $X$ and $D_i$ are smooth.  If $x=0$
defines $D_i$, the Poisson-ideal condition says that every coefficient of
$\pi$ involving the normal direction $\partial_x$ is divisible by $x$.
Every term in the top exterior power involves that normal direction, so
$\pi^{\wedge n}$ vanishes along $D_i$.

Choose a projective log resolution $p:Y\to X$ such that
$U'=p^{-1}(U)$ is a resolution of $U$ and
$B=(Y\setminus U')_{\mathrm{red}}$ is simple normal crossings.  The union of
$B$ with the exceptional locus can also be required to have simple normal
crossings.  Any nonzero logarithmic pluricanonical section can be written
uniquely as
\[
 \sigma=f\Omega^{\otimes m},\qquad f\in\C(X)^*.
\]
At the strict transform of $D_i$, logarithmic regularity implies
\[
 \operatorname{ord}_{D_i}(f)\geq m(e_i-1)\geq0.
\]
At the strict transform of every other prime divisor of $X$, it implies
$\operatorname{ord}(f)\geq0$.  Thus $f$ has no pole at any codimension-one
point of $X$.  Normality makes $f$ a global regular function on $X$, and
projectivity makes it constant.  This argument also proves that each space
of logarithmic pluricanonical sections has dimension at most one.

For the assumed nonzero section, the constant $f$ forces $e_i=1$ for all $i$.
Choosing the canonical divisor defined by $\Omega$, we therefore have the
equality of Weil divisors
\[
 K_X+D=\operatorname{div}_X(\Omega)+D=0.
\]
In particular this sum is Cartier; no prior $\Q$-Gorenstein hypothesis was
needed.  On $Y$, the fact that $\Omega^{\otimes m}$ is logarithmically regular
implies
\[
 \operatorname{ord}_E(\Omega)\geq-1
 \quad\text{if }p(E)\subset D,
 \qquad
 \operatorname{ord}_E(\Omega)\geq0
 \quad\text{if }p(E)\cap U\ne\varnothing.
\]
For an exceptional divisor $E$, these orders are precisely the discrepancies
of $(X,D)$, because $K_X+D=0$.  The log-resolution criterion proves log
canonicity.  The same inequalities show that $\Omega$ itself is a logarithmic
section.  The positive powers of $\Omega$ give the asserted one-dimensional spaces for
every $\ell\geq1$.
\end{proof}

\begin{proposition}[An ample boundary gives rational singularities]
\label{prop:poisson-fano-type}
In the nonvanishing case of Theorem~\ref{thm:poisson-nonvanishing}, suppose
additionally that $U$ is klt and that an effective ample Cartier divisor $Z$
has support $D$.  Then $X$ is of Fano type.  In particular, $X$ has rational
singularities and is Cohen--Macaulay.
\end{proposition}

\begin{proof}
Write $Z=\sum_i r_iD_i$ and choose a rational number
$0<\epsilon<\min_i r_i^{-1}$.  The divisor
$\Delta=D-\epsilon Z$ is effective, and
\[
 K_X+\Delta=-\epsilon Z
\]
is $\Q$-Cartier and anti-ample.  For every divisorial valuation whose center
lies in $D$, the discrepancy of $(X,\Delta)$ equals that of $(X,D)$ plus
$\epsilon\operatorname{ord}_E(Z)>0$.  It is therefore strictly greater than
$-1$.  Valuations with center outside $D$ satisfy the same strict inequality
because $U$ is klt.  Thus $(X,\Delta)$ is klt.  The rational-singularity
theorem for varieties underlying klt pairs gives the last assertion
\cite[Thm.~5.22]{KollarMori1998}.
\end{proof}

\begin{remark}
The theorem proves log canonicity and integral triviality of the logarithmic
canonical divisor.  It does not assert that $X$ is Cohen--Macaulay or that $D$
is Cartier.  Those stronger requirements, if imposed in the definition of a
preferred compactification, require additional arguments.  Proposition
\ref{prop:poisson-fano-type} supplies Cohen--Macaulayness for the Rees
compactifications below when the interior is canonical; Cartierness of the
reduced divisor $D$ is still not asserted.
\end{remark}

\subsection{Degree-zero Poisson Rees algebras}

\begin{proposition}[Poisson Rees compactification]
\label{prop:poisson-rees}
Let $A$ be a finitely generated integral Poisson $\C$-algebra.  Suppose that it
has an exhaustive increasing multiplicative filtration $F_\bullet A$, with
$F_0A=\C$ and finite-dimensional filtered pieces, such that
\[
 \{F_pA,F_qA\}\subset F_{p+q}A,
 \qquad
 R=\bigoplus_{m\geq0}F_mA\,t^m\quad\text{is finitely generated}.
\]
Let $\overline R$ be the normalization of $R$, and let
$\overline A$ be the normalization of $A$.  Then
\[
 X=\operatorname{Proj}\overline R
\]
is a normal projective Poisson compactification of
$U=\operatorname{Spec}\overline A$.  The reduced boundary of this compactification is a Poisson
divisor.  No reducedness, normality, or Gorenstein hypothesis on
$\operatorname{gr}_F A$ is needed for this assertion.
\end{proposition}

\begin{proof}
Give $R\subset A[t]$ the bracket with $t$ central.  The stated inequality
makes this a homogeneous Poisson bracket of degree zero.  Kaledin's
normalization theorem for Poisson algebras extends it uniquely to
$\overline R$; uniqueness also preserves its grading and the centrality of
$t$ (Lemma~\ref{paper:lem-poisson-normalization}).

The normalization is finite and positively graded, with degree-zero part
$\C$.  The bracket extends to every homogeneous localization and preserves
its degree-zero subring.  These subrings are the affine charts of
$\operatorname{Proj}\overline R$, so they define its regular Poisson bracket.
Normalization commutes with localization, whence
\[
 \overline R[t^{-1}]=\overline A[t,t^{-1}],
 \qquad D_+(t)=\operatorname{Spec}\overline A.
\]
The closed complement $V_+(t)$ is Poisson because $(t)$ is a homogeneous
Poisson ideal.  The reduction of $V_+(t)$ is Poisson as well: in characteristic zero
the radical of a Poisson ideal is a Poisson ideal
Lemma~\ref{paper:lem-poisson-normalization}.  Finally, after passing to a sufficiently
divisible Veronese grading, $t^N$ is a nonzero section of the ample line bundle
$\mathcal O_X(N)$.  The zero locus of $t^N$ has pure codimension one and has support
$V_+(t)$.  This proves the boundary assertion.
\end{proof}

\subsection{The Poisson compactification input for all reductive groups}

The filtration in the next proposition is an ordinary matrix-coefficient
filtration.  The proof of the proposition does not assert Poisson compatibility of the more
special string filtrations used for the spectral compactification theorem.

Poisson structures on projective character-variety compactifications have
precedents.  Biswas--Lawton--Ramras construct such structures, including
Poisson boundaries, on wonderful compactifications for semisimple adjoint
groups \cite[Thm.~6.5]{BiswasLawtonRamras2019}; their Question~6.7 asks how
those structures relate to Goldman's bracket.  The construction below
starts with the surface's Goldman bracket itself, works for every reductive
isogeny type, and is arranged to apply to its closed Poisson restrictions.

\begin{proposition}[A uniform matrix-coefficient construction]
\label{prop:matrix-poisson-compactification}
Let $G$ be a connected complex reductive group, and let $\Sigma$ be an
oriented surface of genus $g$ with $k\geq1$ punctures.  Every irreducible
component of every reduced closed Poisson subvariety of
\[
 \mathcal X_G(\Sigma)=G^{2g+k-1}\mathbin{/\!\!/}G
\]
has a normalization admitting the compactification of
Proposition~\ref{prop:poisson-rees}.
\end{proposition}

\begin{proof}
Put $r=2g+k-1$ and choose a ciliated ribbon-graph spine of $\Sigma$.  The
Fock--Rosly formula for graph connections expresses the bivector on $G^r$
as a finite sum, with constant coefficients, of wedges of left- and
right-invariant vector fields on its factors.  The invariant bracket of this bivector is the
Goldman bracket on $\C[G^r]^G$; see \cite[formula~(16), Proposition~3 and the
discussion following formula~(17), Prop.~5 and formula~(24)]{FockRosly}.  One may equivalently use the
canonical quasi-Poisson formula upstairs and its Poisson restriction to
invariants.  The degree estimate below uses only the displayed vector-field
description, not an equivariance assertion about the upstairs bivector.

Choose a finite-dimensional subspace $W\subset\C[G]$ containing $1$, generating
$\C[G]$ as an algebra, and stable under both translations.  Such a space
exists: take matrix coefficients of a faithful representation and its
contragredient, or take the finite-dimensional translation span of a finite
set of algebra generators.  On $\C[G^r]$, let $E_m$ be the span of products
of at most $m$ elements from the copies of $W$ in the individual factors.
Every left- or right-invariant vector field preserves $E_m$.  The
Fock--Rosly formula consequently gives
\[
 \{E_p,E_q\}\subset E_{p+q}.
\]
The pieces are finite dimensional and stable under simultaneous conjugation.
Their Rees algebra $S$ is generated by $t$ and the finitely many degree-one
matrix coefficients multiplied by $t$.  Reductivity gives finite generation
of
\[
 S^G=\bigoplus_{m\geq0}E_m^G t^m.
\]
The invariant pieces define the required degree-zero Poisson filtration on
$\C[G^r]^G$.

If $I$ is the radical ideal of the chosen closed Poisson subvariety, give
$A=\C[G^r]^G/I$ the quotient filtration.  The Rees algebra of $A$ is a quotient
of $S^G$, hence finitely generated, and its bracket has the same degree
estimate.  Minimal primes of a Poisson algebra in characteristic zero are
Poisson, so the same holds after selecting an irreducible component.
Proposition~\ref{prop:poisson-rees} now applies.  The zero-dimensional case
simply gives a point.
\end{proof}

\begin{lemma}[The relative locus is Poisson]
\label{lem:relative-poisson-ideal}
Prescribing arbitrary conjugacy-class closures at the punctures defines a
reduced closed Poisson subvariety of $\mathcal X_G(\Sigma)$.
\end{lemma}

\begin{proof}
Use the free presentation
$M=G^{2g}\times G^{k-1}$, with the last puncture monodromy equal to the inverse
of the total moment map $\mu$.  Equip it with its fused quasi-Poisson
bivector $P$.  On a free puncture factor the standard bivector $P_G$ is tangent
to every conjugacy class; the extra fusion terms involve infinitesimal
conjugation fields.  Thus invariant Hamiltonian derivations preserve every
prescribed class-closure ideal on these factors.  For the last factor the
moment-map identity is
\[
 P^\sharp d(\mu^*f)=(\mu^*\mathcal Df)_M,
\]
where the right side is a combination of infinitesimal $G$-action fields.
If $h$ is $G$-invariant, this gives $\{\mu^*f,h\}=0$.  The Hamiltonian
derivation therefore preserves the final class-closure condition as well.
These identities are algebraic; they extend from a compact real form to
its complexification, which is the given connected complex reductive group.
See
\cite[Definition~2.2, Proposition~3.4 and Proposition~5.1]{AKSMQuasiPoisson}.
Consequently the ideal of all the conditions is stable under invariant
Hamiltonian derivations.  The radical of this ideal remains stable in characteristic zero,
and contraction to $\C[M]^G$ is a Poisson ideal.  Exactness of reductive
invariants identifies its quotient with the coordinate algebra of the
reduced relative locus.
\end{proof}

In the generic situation, smooth boundary evaluation shows that every
component has the exact-class locus dense.  On the
normalization, the complement of that locus has codimension at least two;
the possible finite stabilizer quotients preserve its symplectic volume.
Thus the codimension-one nondegeneracy required by
Theorem~\ref{thm:poisson-nonvanishing} is supplied by
Theorem~\ref{new:thm-normalized-interior} in Section~\ref{new:sec-normalization}.

\begin{corollary}[A reduction of the normalized compactification problem]
\label{cor:normalized-nonvanishing-reduction}
Let $U$ be an irreducible component of the normalization of a generic
relative character variety for a connected complex reductive group and a
punctured surface.  Assume the interior conclusions established above
(Section~\ref{new:sec-normalization}),
in particular the codimension-one nondegeneracy of its Goldman Poisson
structure.  Then $U$ has a normal projective Poisson compactification with
Poisson boundary and
\[
 \overline\kappa^{\mathrm{res}}(U)\in\{-\infty,0\}.
\]
If $\overline\kappa^{\mathrm{res}}(U)\geq0$, that compactification is a log
canonical log Calabi--Yau pair with integral trivialization.  If $U$ is
canonical, the compactification is also of Fano type and Cohen--Macaulay.
An alternative sufficient condition in terms of intersection cohomology
is given in Section~\ref{paper:app-ih}.
\end{corollary}

\begin{corollary}[Equivalent nonvanishing criteria]
\label{cor:equivalent-nonvanishing}
In Corollary~\ref{cor:normalized-nonvanishing-reduction}, suppose that $U$ has
canonical singularities and $K_U\sim0$.  Let $U'\to U$ be a projective
resolution, and let $(X,D)$ be the Poisson Rees compactification constructed
above.  The following conditions are equivalent:
\begin{enumerate}
\item $U$ admits a normal projective log canonical compactification
$(X_0,D_0)$ with reduced complement and $K_{X_0}+D_0\sim_\Q0$;
\item $\overline\kappa^{\mathrm{res}}(U)\geq0$;
\item $\overline P_1^{\mathrm{res}}(U)>0$;
\item $F^{\dim U}H^{\dim U}(U',\C)\ne0$;
\item $(X,D)$ is log canonical and $K_X+D\sim0$ integrally.
\end{enumerate}
One may replace $\overline\kappa^{\mathrm{res}}(U)$ by
$\overline\kappa(U_{\mathrm{reg}})$ and $U'$ by $U_{\mathrm{reg}}$ in the
corresponding conditions.  The top-weight intersection-cohomology condition
of Section~\ref{paper:app-ih} is sufficient, but
is not claimed to be equivalent to these conditions.
\end{corollary}

\begin{proof}
Conditions (ii), (iii), and (v) are equivalent by
Theorem~\ref{thm:poisson-nonvanishing}; (iii) and (iv) are equivalent by
logarithmic Hodge--de Rham degeneration.  Clearly (v) implies (i).  For the
remaining implication take a log resolution of $(X_0,D_0)$ and denote by
$B$ only the reduced inverse image of $D_0$.  The divisor $K_Y+B$ is the
pullback of $K_{X_0}+D_0$ plus an effective exceptional divisor.  Over $D_0$
effectivity follows from log canonicity, since the coefficient is $a_E+1$;
over the interior it follows from canonicity of $U$, since the coefficient
is $a_E$.  A sufficiently divisible positive multiple therefore has a nonzero
section.  This proves (i)$\Rightarrow$(ii).  The assertions involving the
smooth locus follow from the comparison of pluriforms at the beginning of
this section and logarithmic Hodge--de Rham degeneration.
\end{proof}

\subsection{Finite central covers and volume forms}

\begin{lemma}[Valuative transfer of the volume form]
\label{new:gl-finite-volume}
Let $p:V_1\to V_2$ be a finite surjective morphism of irreducible normal
complex varieties.  Suppose both varieties have canonical singularities
and rational top forms $\Omega_i$ trivializing their canonical sheaves,
with $p^*\Omega_2=\Omega_1$.  Then $\Omega_1$ extends as a logarithmic
top form on a smooth compactification of a resolution of $V_1$ if and
only if $\Omega_2$ has the corresponding property for $V_2$.
\end{lemma}

\begin{proof}
For a divisorial valuation $v$ of $\C(V_2)$, choose an extension $w$ to
$\C(V_1)$, and let $e$ be its ramification index.  In characteristic zero,
the codimension-one ramification formula for rational top forms is
\[
 \operatorname{ord}_{w}(p^*\Omega_2)+1
   =e\bigl(\operatorname{ord}_{v}(\Omega_2)+1\bigr).
\]
The residue-field extension is separable, so the different contributes
$e-1$.  The formula can be checked on the normalization of a model on
which $v$ is a divisor.  Restrictions and extensions of divisorial
valuations under a finite function-field extension remain divisorial.

Moreover $w$ has a center on $V_1$ if and only if its restriction has a
center on $V_2$.  One direction follows by composition.  The other follows
from the valuative criterion for the finite, hence proper, morphism $p$.
Thus the displayed formula transfers the bound
$\operatorname{ord}(\Omega_i)\geq-1$ exactly between valuations at
infinity.  At valuations centered in the interiors, the stronger bound
$\operatorname{ord}(\Omega_i)\geq0$ follows from canonicity and the
chosen canonical trivialization.  These are precisely the bounds for
a logarithmic top form on a smooth compactification of a resolution.
\end{proof}

\begin{proposition}[Reduction to simply connected simple groups]
\label{new:gl-simple-reduction}
Fix $g\geq0$ and $k\geq1$.  Suppose that every nonempty irreducible
component $V$ of every normalized relative character variety of a
simply connected simple complex group, with these $g,k$ and compatible
intrinsically generic puncture conjugacy classes of arbitrary Jordan
type, satisfies
\[
 \overline\kappa^{\mathrm{res}}(V)\geq0.
\]
Then the same logarithmic nonvanishing holds for every connected
complex reductive group with such puncture data.  In particular,
Corollary~\ref{cor:normalized-nonvanishing-reduction} gives the
normalized log canonical log Calabi--Yau compactification, with
integral $K+D\sim0$ and its Fano-type consequences, for every such
component.  No stronger genericity condition on the given data is
required.
\end{proposition}

\begin{proof}
Let $G$ be connected reductive and choose its finite central covering
\[
 f:\widetilde G=T\times\prod_{a=1}^sG_a^{\mathrm{sc}}
          \longrightarrow G,\qquad K=\ker(f),
\]
where $T=Z(G)^\circ$ and the $G_a^{\mathrm{sc}}$ are simply connected
simple.  Let $Y_i=\overline C_i$, let $R$ be the reduced relation
variety for $G$, and put $U_0=R/\!/G$.  Use the \emph{full} preimages
$\widetilde Y_i=f^{-1}(Y_i)$ to define the reduced lifted relation
variety $\widetilde R$ with relation equal to $1$ in $\widetilde G$.

Eliminate the last puncture coordinate from the relation.  The map
\[
 \widetilde R\longrightarrow R
\]
is the base change of
$\widetilde G^{\,2g+k-1}\to G^{\,2g+k-1}$, hence is a finite
\'etale torsor for $\Gamma=K^{2g+k-1}$.
Indeed the eliminated last coordinate automatically lies in
$f^{-1}(Y_k)$ whenever its image lies in $Y_k$.
In terms of the coordinates, $\Gamma$ multiplies the handle and
first $k-1$ puncture lifts by arbitrary central elements, and
multiplies the last puncture lift by the inverse product of the
puncture multipliers.  This explains both surjectivity and the use
of $k\geq1$.  It also includes all possible lift sectors.

The preimage $f^{-1}(C_i)$ is a finite disjoint union of
$\widetilde G$-conjugacy classes.  To see finiteness, lift a
representative of $C_i$ and then multiply it by the finitely many
elements of $K$; lifting conjugating elements exhausts the preimage.
Since $f$ is open, $f^{-1}(C_i)$ is dense in $\widetilde Y_i$, so
$\widetilde Y_i$ is the union of the closures of these lifted
classes.  Distinct closures are allowed to meet along smaller
orbits.  Each lifted class has the form
\[
 \widetilde C_i=\{z_i\}\times\prod_a C_{i,a}.
\]

Every lifted choice is intrinsically generic.  In fact, if
semisimple elements in a proper parabolic $\widetilde P$ upstairs
had product in $[\widetilde P,\widetilde P]$, their images would
lie in the corresponding proper parabolic $P=f(\widetilde P)$
and have product in $[P,P]$, contradicting the original genericity.
This is only a downward implication for a hypothetical violation;
the converse genericity assertion is neither used nor asserted.
For a nonempty lifted relation sector, its central coordinates
satisfy $\prod_i z_i=1$.  Each factor tuple $(C_{i,a})_i$ is then
generic as well: a violation in a proper parabolic $P_a$ would give
a violation in
\[
 T\times P_a\times\prod_{b\ne a}G_b^{\mathrm{sc}},
\]
whose derived subgroup is
$\{1\}\times[P_a,P_a]\times\prod_{b\ne a}G_b^{\mathrm{sc}}$.
The equality of the central product is exactly what is needed here.

Apply Theorem~\ref{new:thm-normalization-commutes} to each lifted
choice.  It identifies its normalized relation variety with the
one obtained by normalizing the puncture closures before imposing
the relation, and its normalized quotient with the quotient of
that normal relation variety.  The exact-class locus is dense in
every component.  Since distinct lifted choices have disjoint
exact-class loci, normalization separates the choices, including
any intersections along smaller orbits.  It follows that
\[
 \widetilde U:=\widetilde R^{\nu}/\!/\widetilde G
   \simeq
 \coprod_{\text{nonempty lifted choices}}
       \left(T^{2g}\times\prod_a U_a^{\nu}\right),
\]
where $U_a$ is the relative character variety for the factor
$G_a^{\mathrm{sc}}$ with puncture closures $\overline C_{i,a}$.
This identity is understood componentwise when a factor is reducible.
The product formula follows by splitting the relation and the
conjugation action into their torus and simple factors.

The action of $\Gamma$ lifts to $\widetilde R^{\nu}$ and commutes
with conjugation; the conjugation action of $\widetilde G$ factors
through $G$.  Thus $\widetilde U/\Gamma$ is normal and finite over
$U_0$.  One can check finiteness at the coordinate rings: the
normalization is finite, and invariants of a finite equivariant
algebra over an affine reductive quotient are finite over the
invariant base.  Over the dense exact-class locus normalization is
an isomorphism, and commuting invariants give
\[
 \bigl(\C[\widetilde R]^ {\widetilde G}\bigr)^\Gamma
      =\bigl(\C[\widetilde R]^\Gamma\bigr)^G
      =\C[R]^G.
\]
Consequently the finite normal map is birational componentwise,
and
\[
 U_0^{\nu}\simeq\widetilde U/\Gamma.
\]
In particular every component $U$ of $U_0^\nu$ is the finite
surjective image of a component $\widetilde V$ of $\widetilde U$.
Only finiteness is claimed for this map of coarse spaces.

By the hypothesis and
Corollary~\ref{cor:normalized-nonvanishing-reduction}, each simple
factor component has a log canonical log Calabi--Yau Poisson
compactification.  The Poisson criterion identifies its rational
Goldman volume as a logarithmic top form there.  On $T^{2g}$ the
Goldman symplectic volume is a nonzero constant multiple of the
product of the coordinate $d\log$ forms, which is logarithmic on
$(\mathbb P^1)^{2g\dim T}$.  For $g=0$ this factor is a point.
Products of smooth compactifications of resolutions therefore make
the product volume logarithmic on $\widetilde V$.

Choose a nondegenerate invariant bilinear form on $\Lie G$ and use
its pullback upstairs.  Invariance makes its central and simple
summands orthogonal.  Hence the lifted Goldman form is the sum of
the factor forms, and its top exterior power is a nonzero constant
multiple of the product volume just considered.  On exact loci,
$df$ identifies the adjoint local systems and their parabolic
cup-product pairings.  Thus, for the finite map
$p:\widetilde V\to U$, the Goldman volumes obey
\[
 p^*\Omega_U=\Omega_{\widetilde V}
\]
as rational forms.  In particular, the central-twisting action
preserves the Goldman form and its volume, without a character
multiplier.

Both interiors are canonical, and these volumes trivialize their
canonical sheaves, by
Theorem~\ref{new:thm-normalized-interior}.
Lemma~\ref{new:gl-finite-volume} therefore descends the logarithmic
volume to $U$.  Explicitly, for any divisorial valuation $v$ of
$\C(U)$ and an extension $w$ of ramification index $e$,
\[
 \operatorname{ord}_{w}(p^*\Omega_U)+1
      =e\bigl(\operatorname{ord}_{v}(\Omega_U)+1\bigr).
\]
Finiteness preserves whether these valuations have a center in the
interior.  The bound of at most a simple pole at infinity thus
descends even if the coarse finite map ramifies.  It gives
$\overline P_1^{\mathrm{res}}(U)>0$, hence the required logarithmic
nonvanishing and, by
Corollary~\ref{cor:normalized-nonvanishing-reduction}, the stated
compactification.
\end{proof}

\section{The Goldman Hodge class and a concrete period test}
\label{sec:goldman-period-test}

There is a useful refinement of the nonvanishing problem in
Section~\ref{sec:poisson-bridge}.  The Goldman form belongs to the second
Deligne Hodge filtration even for nonsemisimple prescribed conjugacy classes.
This requires an argument: a regular algebraic differential form on an affine
variety need not represent a class in the corresponding Deligne Hodge
filtration.  For example, differentials of the second kind on affine curves
can represent classes outside $F^1$.

\emph{Idea.} The Goldman form is obtained from algebraic Maurer--Cartan forms by
an algebraic surface-transgression, and algebraic constructions respect the
Hodge filtration (Fact~\ref{paper:fact-hodge}). Following Mellit's argument for
$\GL_n$ \cite[\S4]{Mellit2025}, we make this transgression algebraic for every
reductive group. Then the top power of the Goldman form lies in
$F^{2n}H^{2n}$, which consists of logarithmic volume forms on a smooth
compactification; so a compact cycle with nonzero Goldman-volume period produces
the nonzero logarithmic form required by Theorem~\ref{thm:poisson-nonvanishing}.

\begin{proposition}[All-type adaptation of Mellit's Hodge-filtration argument]
\label{prop:goldman-f2-all-types}
Let $G$ be a connected complex reductive group, and prescribe arbitrary
conjugacy classes at the punctures of an oriented surface.  On the locus of
closed orbits in the relative character variety where the stabilizer is a finite extension
of $Z(G)$, the Goldman form associated with the trace pairing of a faithful
representation has cohomology class in $F^2H^2$ on its smooth locus.
\end{proposition}

\begin{proof}
We first make explicit the Borel-subgroup version of the argument of
\cite[Section~4, Propositions~4.4--4.5]{Mellit2025} for arbitrary elements
$C_1,\ldots,C_k\in\mathrm{GL}_N$.  Mellit's Section~4.4 already constructs
the parabolic character stack for arbitrary fixed matrices; the displayed
splitting-principle proof in Section~4.5 then assumes that they are diagonal.
The following argument supplies the replacement for that
diagonal assumption before passing to arbitrary reductive groups.
Put $H_i=Z_{\mathrm{GL}_N}(C_i)$.
Each $H_i$ is connected: it is the nonempty open set of invertible elements
in the vector space of matrices commuting with $C_i$.  The reductive quotient of $H_i$
is a product of general linear groups.  Choose a Borel subgroup $B_i$ of
$H_i$.  Since $C_i$ is central in $H_i$, it belongs to $B_i$.  After replacing
$C_i$ by a conjugate, arrange that $B_i$ is contained in the fixed
upper-triangular Borel $B$ of $\mathrm{GL}_N$.

For clarity, write $\mathcal B_H$ for the mapping-cone complex of varieties
used by Mellit:
\[
 \mathcal B_H=\operatorname{cone}\left(
   \bigoplus_i BH_i[1]\xrightarrow{\sum_i\sigma_{C_i}}
   B\mathrm{GL}_N\right).
\]
Here $\sigma_{C_i}$ is the bar-complex map which inserts the fixed commuting
element $C_i$.  The analogous complexes obtained from $B_i\subset B$ and
their torus quotients will be denoted $\mathcal B_B$ and $\mathcal B_T$.
The cohomology of $BH_i$ is that of its reductive quotient, so it is
concentrated in even degrees and is Tate.  In degree four there is a
commutative diagram of short exact sequences
\[
 0\longrightarrow\bigoplus_i H^2(BH_i)
 \longrightarrow H^4(\mathcal B_H)
 \longrightarrow H^4(B\mathrm{GL}_N)\longrightarrow0,
\]
with analogous rows for $\mathcal B_B$ and $\mathcal B_T$.  Restriction from
$H_i$ to $B_i$ is injective on cohomology by the usual splitting principle,
while passage from a connected solvable group to its torus quotient induces
an isomorphism.  Consequently
\[
 H^4(\mathcal B_H)\longrightarrow H^4(\mathcal B_B)
       \longleftarrow H^4(\mathcal B_T)
\]
has its first arrow injective and its second arrow an isomorphism, both in
mixed Hodge structures.

The universal trace cocycle is
\[
 \eta_2=\operatorname{Tr}(x^{-1}dx\wedge dy\,y^{-1}),\qquad
 \eta_3=\operatorname{Tr}((x^{-1}dx)^{\wedge3}),\qquad
 \eta=\eta_2+\tfrac13\eta_3.
\]
Mellit's cocycle identity applies unchanged: its extra centralizer boundary
term vanishes because $C_i$ is constant.  On upper-triangular matrices
$\eta_3=0$, and $\eta_2$ is pulled back from
$\sum_a d\log x_a\wedge d\log y_a$ on the diagonal torus.  It is therefore a
logarithmic two-form on the torus complex and represents an $F^2$ class.
Strictness for the injective morphism of mixed Hodge structures above implies
$[\eta]\in F^2H^4(\mathcal B_H)$.

The surface-transgression morphism of complexes is algebraic.  The pullback of the Goldman cocycle
is the usual triangle formula for the Goldman form on the relative character
stack, represented by a two-form on the representation-space atlas with no
higher simplicial component.  The identities $d\omega=0$ and
$\delta\omega=0$ make it a basic closed form, hence the pullback of the
Goldman form on the coarse smooth locus.  The class of the Goldman form therefore belongs to
$F^2$.  No semisimplicity of $C_i$ was used in
this extension.  The equality with the Goldman form is the cup-product
calculation in \cite[Corollary~4.7 and Theorem~4.13]{Mellit2025}, whose expression
uses only fixed monodromies and the invariant trace pairing.  Equivalently,
the triangle expression computes the pairing on parabolic first cohomology
constructed for arbitrary conjugacy classes by Guruprasad--Huebschmann--
Jeffrey--Weinstein \cite{GHJW1997}; it does not require diagonalization.

Now choose a faithful representation $G\hookrightarrow\mathrm{GL}_N$.
The resulting algebraic map of relative character stacks pulls the preceding
class back to the Goldman class for
$\langle a,b\rangle=\operatorname{Tr}(d\rho(a)d\rho(b))$.  This pairing is
nondegenerate: after making the representation unitary on a maximal compact
subgroup, it is negative definite on its real Lie algebra, and complexification
gives the assertion.  Finally the map from the character stack to its coarse
space induces an injection in degree-two rational cohomology on the stated
stable locus.  Indeed, factor this map as
$[R/G]\to[R/(G/Z(G))]\to R/\!/G$ on the locus of closed orbits.
The effective quotient stack has finite inertia, so its rational cohomology
agrees with that of its coarse space, as follows in local finite quotient
charts.  The first map is a $BZ(G)$-gerbe; since
$H^1(BZ(G),\mathbb Q)=0$, its Leray sequence gives an injection in degree two.
Both maps are algebraic and their pullbacks preserve mixed Hodge structures.
Strictness therefore descends the $F^2$ statement to the coarse smooth locus.
\end{proof}

\begin{lemma}[Regular extension preserves this $F^2$ class]
\label{lem:f2-resolution}
Let $U$ be normal, let $\rho:U'\to U$ be a projective resolution which is an
isomorphism over $U_{\mathrm{reg}}$, and let $\omega$ be a closed regular
two-form on $U_{\mathrm{reg}}$.  Suppose that it extends to a regular two-form
$\widetilde\omega$ on $U'$ and that $[\omega]\in F^2H^2(U_{\mathrm{reg}})$.
Then $[\widetilde\omega]\in F^2H^2(U')$.
\end{lemma}

\begin{proof}
The restriction map $r:H^2(U')\to H^2(U_{\mathrm{reg}})$ is a morphism of
mixed Hodge structures.  By strictness there is a class
$\alpha\in F^2H^2(U')$ with $r(\alpha)=[\omega]$.  The localization sequence
shows that the kernel of $r$ is spanned by the classes of exceptional prime
divisors.  Thus
\[
 [\widetilde\omega]-\alpha=\sum_j c_j[E_j].
\]
Both terms on the left pair to zero with every complete algebraic curve in
$U'$.  For $\widetilde\omega$ this follows because the pullback of a
holomorphic two-form to the normalization of a curve is zero; for $\alpha$
it follows from $F^2H^2(C)=0$ for a smooth proper curve.  The exceptional
complex divisor $\sum_j c_jE_j$ is therefore numerically trivial over $U$.
Apply the negativity lemma \cite[Lem.~3.39]{KollarMori1998} separately
to its real and imaginary parts and their negatives; the
intersection-theoretic proof of that lemma applies to real coefficients as well.
An exceptional real divisor which is both relatively nef
and relatively anti-nef is zero.  Hence all coefficients vanish, proving the
claim.
\end{proof}

For a normalized generic relative character variety $U$, let $U^\circ$
denote the inverse image of the locus where every puncture monodromy lies
in the prescribed conjugacy class itself, rather than a smaller class in
its closure.  The Goldman form is initially defined on the smooth part
of this exact-class locus.  The exact-class locus has complement of
codimension at least two in $U$, and the smooth part of $U^\circ$ has
complement of codimension at least two in $U_{\mathrm{reg}}$.
Removing such a subset from a smooth variety does not
change $H^2$, as a mixed Hodge structure.  Proposition
\ref{prop:goldman-f2-all-types} therefore gives the $F^2$ statement on all of
$U_{\mathrm{reg}}$.  Whenever the symplectic-singularity extension theorem
applies to $U$, Lemma~\ref{lem:f2-resolution} gives it on a projective
resolution as well.

\begin{corollary}[A Goldman-period criterion]
\label{cor:goldman-period}
In the setting of Corollary~\ref{cor:normalized-nonvanishing-reduction}, assume
that $U$ has canonical singularities and $K_U\sim0$.  Let
$\dim U=2n$, and choose its trace-pairing Goldman form as in
Proposition~\ref{prop:goldman-f2-all-types}.  The regular extension of the Goldman form
$\widetilde\omega$ to a projective resolution $U'$ exists by
Theorem~\ref{new:thm-normalized-interior}.  The following conditions are equivalent:
\begin{enumerate}
\item $U$ has a normal log canonical log Calabi--Yau compactification.
\item $[\widetilde\omega]^n\ne0$ in $H^{2n}(U',\mathbb C)$.
\item There is a finite rational singular $2n$-cycle $\Gamma$ in $U'$ such that
\[
 \int_\Gamma\widetilde\omega^{\wedge n}\ne0.
\]
\end{enumerate}
When these conditions hold, $U$ admits the specific compactification of
Section~\ref{sec:poisson-bridge}.  Under the canonical-interior hypotheses,
that compactification is of Fano type.
\end{corollary}

\begin{proof}
Conditions (2) and (3) are equivalent by the de Rham pairing and the universal
coefficient theorem; a finite singular cycle has compact support, and no
smooth representative is required.  By the preceding results and
multiplicativity of the Hodge filtration,
$[\widetilde\omega]^n\in F^{2n}H^{2n}(U')$.  Thus (2) implies (1) by
Corollary~\ref{cor:equivalent-nonvanishing}.

Conversely, if (1) holds, the Poisson criterion gives a log canonical
log Calabi--Yau pair $(X,D)$ for the specific compactification constructed
there.  On a common log resolution its symplectic volume
$\widetilde\omega^{\wedge n}$ is a nonzero logarithmic top form.  By the
birational invariance of the logarithmic top-form space proved in
Section~\ref{sec:poisson-bridge}, it remains logarithmic on a smooth projective
compactification of the fixed $U'$ with simple normal crossings boundary.
Deligne's $E_1$ degeneration identifies the space of these forms with
$F^{2n}H^{2n}(U')$.  In particular the cohomology class of this form is nonzero.
This proves (2).
\end{proof}

\begin{example}[The period test in the running example]
Let $G=\SL_2$, $(g,k)=(1,1)$ and $-2<\tau<2$. The characters of
$\mathrm{SU}(2)$-representations form a compact component of the real points of
$U_\tau=\{x^2+y^2+z^2-xyz=\tau+2\}$, a $2$-sphere $\Gamma$ contained in $[-2,2]^3$
\cite{Goldman2003}. The Goldman form is a constant multiple of
$\widetilde\omega=dx\wedge dy/(2z-xy)$; on $\Gamma$ it restricts to a real
nondegenerate two-form, so $\int_\Gamma\widetilde\omega\ne0$ and condition~(3) of
Corollary~\ref{cor:goldman-period} holds with $n=1$. For a general group and
nonsemisimple local monodromy no such unitary cycle is at hand; this is why
Section~\ref{r2:sec-dolbeault-closure} constructs the cycle on the Dolbeault
side.
\end{example}

\section{The Dolbeault closure and the compact Goldman period}
\label{r2:sec-dolbeault-closure}

The construction in this section uses the tame nonabelian Hodge
correspondence separately on each exact residue stratum.  No assertion
about its continuity across different nilpotent orbits is needed.
Throughout the proof $G$ is semisimple and simply connected; the central
isogeny and torus reductions are applied afterwards.  Fix an intrinsically
generic tuple $C_i=G\cdot(S_i u_i)$ and a nonempty component of its
exact-class character variety.

\emph{Plan.} We need a compact cycle with nonzero Goldman period. In the
running example the exact Dolbeault stratum is a moduli space of rank-two
parabolic Higgs bundles on an elliptic curve with one marked point; the Hitchin
map is the quadratic differential $\det\Phi$, the base is a line, and a general
fiber is an elliptic curve, a compact Lagrangian. In general, on the
Dolbeault side the natural candidates are fibers of the Hitchin map. The first
subsection matches the Betti and Dolbeault data at the punctures, including the
closure order of orbits. The second builds an algebraic closure of the exact
Dolbeault stratum, with projective Hitchin map, on which every object is stable.
The third shows that the symplectic class on the exact stratum is a complex
linear combination of Chern classes. The fourth extends this identity to a
resolution; Hodge theory then makes projective subvarieties isotropic, Poisson
commutativity of the Hitchin Hamiltonians gives the reverse bound, and a general
Hitchin fiber is a compact Lagrangian inside the exact stratum
(Theorem~\ref{r2:alg-hitchin-rank}). Hyperk\"ahler rotation turns this fiber
into a cycle with nonzero Goldman period.

\subsection{The local orbit dictionary preserves the closure order}

Choose a maximal compact subgroup and a common complex torus for the
commuting polar factors of $S_i$.  The local tame correspondence with
\emph{zero Betti weight} determines a real alcove weight $\alpha_i$ and
a semisimple Dolbeault residue $b_i$ with
\begin{equation}\label{r2:eq-dol-polar-data}
 b_i+\overline{b_i}=0,\qquad
 S_i=\exp\bigl(-2\pi\sqrt{-1}(\alpha_i+2b_i)\bigr).
\end{equation}
The bar in this expression is the conjugation on the complexified
cocharacter space used in the local correspondence. Equation
\eqref{r2:eq-dol-polar-data} is the zero-weight case of Simpson's table
\cite[table, p.~720]{Simpson1990}: writing $b_i=\sqrt{-1}\,t_i$ with $t_i$ real
in the compact form, $S_i$ is the commuting product of the unitary factor
$\exp(-2\pi\sqrt{-1}\alpha_i)$, recorded by the parahoric weight, and the
positive factor $\exp(4\pi t_i)$, recorded by the residue.  Let $L_i$ be the
reductive quotient of the parahoric group at $\alpha_i$.  The usual
root-space identification identifies $L_i$ with
$Z_G(\exp(-2\pi\sqrt{-1}\alpha_i))$.  The latter is connected because
$G$ is simply connected.  Uniqueness of the commuting compact and positive
polar factors gives
\begin{equation}\label{r2:eq-dol-residue-centralizer}
 H_i:=Z_{L_i}(b_i)=Z_G(S_i).
\end{equation}
Let $\mathcal N_i$ be the $H_i$-orbit of $\log u_i$, and set
\[
 \mathcal O_i=L_i\cdot(b_i+\mathcal N_i)
       \subset\mathfrak l_i,\qquad
 Y_i=\overline{\mathcal O_i}.
\]
Here and below an invariant nondegenerate bilinear form identifies a
reductive Lie algebra with its dual.

\begin{lemma}[Exact local types and their boundary types]
\label{r2:lem-local-dol-dictionary}
With these fixed weights and semisimple residues, the exact graded-residue
condition $\mathcal O_i$ corresponds to the ordinary Betti condition
$C_i$.  More generally, the $L_i$-orbits in $Y_i$ correspond bijectively
to the $G$-orbits in $\overline C_i$, with the same nilpotent-orbit
closure order in $H_i$.  This assertion concerns exact orbit labels and
does not assert a homeomorphism of either pair of orbit closures.
\end{lemma}

\begin{proof}
The adapted model and local table of
\cite[\S5.1, (5.4)--(5.8), Prop.~6.8 and Table~1]{BiquardGarciaPradaMundet2020},
with the parahoric conventions of \cite[\S4.5]{HKSZ2025}, give, for a nilpotent
graded residue $f$ completed to an $\mathfrak{sl}_2$-triple $(e,h,f)$,
the unipotent monodromy factor
\[
 \exp\bigl(2\pi\sqrt{-1}(h+e-f)\bigr).
\]
The triple lies in $\mathfrak h_i$.  In $\mathfrak{sl}_2$ the element
$h+e-f$ is a nonzero nilpotent whenever $f\ne0$, and is conjugate to
$f$; every nonzero scalar multiple of either is conjugate to it as
well.  Applying the associated $\mathrm{SL}_2$-homomorphism proves
that the nilpotent monodromy logarithm belongs to precisely the same
complex $H_i$-orbit as $f$.  The assertion for $f=0$ is immediate.
The Betti weight is $-(b_i+\overline{b_i})=0$ throughout, so there is
no unspecified unipotent radical of a nontrivial Betti filtration.

The semisimple-part morphism identifies
\[
 Y_i\simeq L_i\mathbin{\times^{H_i}}
                (b_i+\overline{\mathcal N_i}),\qquad
 \overline C_i\simeq G\mathbin{\times^{H_i}}
                (S_i\exp\overline{\mathcal N_i}).
\]
For the Lie-algebra formula, regularity of the semisimple-part map
follows by polynomial interpolation on the fixed finite set of
eigenvalues; the group formula is proved in the same manner.
Exponential and logarithm are algebraic inverse isomorphisms on the
nilpotent and unipotent varieties.  Both orbit decompositions therefore
have the asserted common nilpotent-orbit poset.
\end{proof}

\begin{lemma}[The smooth exact locus with finite stabilizers]
\label{r2:lem-finite-stabilizer-hk}
Let $G$ be complex semisimple.  Fix degree-zero tame data
$(\alpha,\mathbf L)$ and $(0,\mathbf C)$ related by the local
correspondence of \cite[Proposition~6.8, Table~1 and Theorem~7.10]{BiquardGarciaPradaMundet2020}.
Here the graded residues lie in the \emph{exact} orbits $\mathbf L$
and the monodromies in the \emph{exact} classes $\mathbf C$.
Restrict to the stable locus with finite automorphism groups, and let
$O_I$ and $O_B$ be the smooth loci of its Dolbeault and Betti coarse
spaces.  The exact tame correspondence restricts to a diffeomorphism
\[
 \Phi:O_I\longrightarrow O_B.
\]
Choose the invariant bilinear form on $\Lie G$ to be negative definite
on a compact real form, and use it in both symplectic constructions.
Then $O_I$ has the hyperk\"ahler structures $(I,J,K,g)$ coming from
the tame equations, with
\[
 \Omega_{\mathrm{Dol}}=\omega_J+\sqrt{-1}\omega_K,
 \qquad
 \Phi^*\omega_{\mathrm{Goldman}}
      =c(\omega_K+\sqrt{-1}\omega_I),
 \qquad c\in\C^*.
\]
The first equality uses the algebraic symplectic form on the exact
Dolbeault stratum.  The constant $c$ only records the fixed conventions
for the flat connection and the Goldman pairing.
Consequently, if $\dim_{\C}O_I=2n$ and $F\subset O_I$ is a nonempty
smooth projective $I$-complex Lagrangian submanifold of dimension $n$,
then
\[
 \int_{\Phi(F)}\omega_{\mathrm{Goldman}}^{\wedge n}
   =c^n\sqrt{-1}^{\,n}\int_F\omega_I^{\wedge n}\ne0,
\]
where $F$ has its complex orientation and $\Phi(F)$ its transported
orientation.
\end{lemma}

\begin{proof}
We retain the finite stabilizer in the local construction, rather than
identifying stability with simplicity.  Use the adapted tame model and
weighted spaces of \cite[Section~5.1, equations~(5.4)--(5.8), and
Section~5.4]{BiquardGarciaPradaMundet2020}.
For fixed exact local type, the deformation complex is the parabolic
deformation complex with variations tangent to the prescribed residue
orbits.  The invariant pairing and Serre duality identify its obstruction
space with the dual of its infinitesimal automorphism space.  Equivalently,
one may include the residue coadjoint orbits in the Hamiltonian
description; their Kirillov--Kostant forms supply exactly the boundary
pairing.  Finiteness of the automorphism group therefore makes both
spaces zero.  This is the complex-group vanishing used immediately
before \cite[Proposition~7.1]{BiquardGarciaPradaMundet2020}.

The elliptic slice and implicit-function construction of
\cite[Section~7.1]{BiquardGarciaPradaMundet2020} consequently gives a
smooth hyperk\"ahler uniformizing chart $V$ for the space of tame
solutions, with the remaining finite automorphism group $H$ acting
by hyperk\"ahler isometries.  No step of this construction requires
$H=Z(G)$: that stronger condition only removes effective finite
isotropy.  This is the stable nonsimple orbifold phenomenon described
in \cite[Remark~7.6]{BiquardGarciaPradaMundet2020}.  The Dolbeault and
flat descriptions are the $I$- and $J$-complex structures on this same
solution chart, with the same $H$-action.  Existence and uniqueness of
the adapted harmonic reduction, followed by the Riemann--Hilbert
correspondence, give the exact correspondence of
\cite[Proposition~7.9 and Theorem~7.10]{BiquardGarciaPradaMundet2020} on these quotients.
Thus this argument also covers finite stabilizers not contained in
the center; it does not use the generic-weight diffeomorphism statement
of Theorem~7.11 there.

Here is why the coarse smooth loci correspond.  At a point of $V$,
divide its stabilizer by the kernel of its action on the germ.
Finite holomorphic linearization identifies the resulting quotient
germ with the quotient of its complex tangent space.  In either
complex structure $I$ or $J$ the effective stabilizer preserves a
nondegenerate complex symplectic form.  A nonidentity symplectic
linear transformation cannot fix a hyperplane, since its fixed
subspace has even codimension.  The Chevalley--Shephard--Todd theorem
\cite[\S5]{ShephardTodd1954} therefore says that this quotient is smooth precisely when the
effective stabilizer is trivial.  Triviality on the germ is the same
condition for $I$ and $J$, since the action on the underlying real
chart is the same.  The hyperk\"ahler tensors descend on that open
locus, and the common charts give the asserted diffeomorphism.

Finally the tensors are the required tensors, not just abstract
symplectic forms.  On the $I$-chart the holomorphic symplectic pairing
is the Serre-duality pairing of the fixed-exact-type deformation
complex.  This is also the pairing obtained by algebraic cotangent
reduction with the residue coadjoint orbits; hence it equals
$\Omega_{\mathrm{Dol}}$.  On the flat chart the same invariant
bilinear form gives the de Rham pairing on parabolic cohomology.
The identification of this pairing with the group-cohomology Goldman pairing is
\cite[Theorems~8.3 and~10.5]{GHJW1997}; those pairing identities apply
to the complex adjoint local system as well.  The hyperk\"ahler
identities now give the two displayed formulas, with one fixed
normalization constant.  On an $I$-complex Lagrangian $F$, the real
and imaginary parts of $\Omega_{\mathrm{Dol}}|_F=0$ give
$\omega_J|_F=\omega_K|_F=0$.  The remaining form $\omega_I|_F$
is K\"ahler, so its top integral is positive.  This proves the last
assertion.
\end{proof}

\subsection{Closed algebraic moduli and stability on every boundary stratum}

Let $M_0$ be the reduced closed subspace of the moduli space of
semistable weak parahoric Higgs torsors of weight $\alpha=(\alpha_i)$
on which the graded residues lie in $Y_i$.  Use the degree-zero
topological condition of tame nonabelian Hodge theory.  Let
$M=M_0^\nu$, and let $M^\circ\subset M$ denote the inverse image of
the exact graded-residue locus.  When discussing one component, take
the component whose exact locus corresponds to the chosen Betti
component.

The projectivity input uses the same parabolic-degree (Ramanathan)
stability as the tame correspondence.  The hypotheses of this input are a smooth
projective curve with finitely many marked points, a semisimple complex
group, fixed real parahoric weights and topological type, and closed
conjugacy-invariant conditions on the graded residues.  These are
exactly the conditions here: the degree-zero topological type is fixed,
and each $Y_i$ is a graded-residue orbit closure.  Under these
hypotheses Lemma~\ref{r2:lem-proper-parahoric} gives a projective
Hitchin morphism on $M_0$ before genericity is used to prove stability
of every object.  Finite normalization then gives projectivity on $M$.

The proof in Section~\ref{paper:app-properness} first replaces the
real weights, for the algebraic moduli construction only, by rational
weights in the same parahoric facet and stability cell, preserving
all stability equalities.  A finite Galois cover then identifies
parahoric Higgs torsors with equivariant Higgs bundles, with twist
$\pi^*K_X(D)$; trivial inertia representations at auxiliary branch
points prevent additional Higgs poles.  Uniqueness of the
Harder--Narasimhan reduction makes it equivariant, so the resulting
$\Theta$-stratification has precisely the Ramanathan semistable open.
The valuative existence statement for the full weak parahoric Higgs
stack and semistable reduction give properness of the Hitchin
morphism on the separated quasiprojective moduli space.  Since its
base is affine, this morphism is projective.  Closed residue
conditions preserve that conclusion.  Thus the projectivity used
below applies to the actual stability condition and to the entire
algebraic closure, including its boundary strata.

\begin{proposition}[Algebraic closure and its strata]
\label{r2:prop-dol-closure}
Every point of $M_0$ is stable and has finite automorphism group.  The
space $M$ has symplectic singularities and carries a projective Hitchin
morphism to an affine variety.  Each component has a dense exact locus,
and the complement of that locus has codimension at least two.  There
is a finite refinement of the exact graded-residue decomposition into
smooth symplectic strata.  The Hitchin functions Poisson-commute on
each stratum with respect to its own symplectic form.

The smooth exact locus has the usual hyperk\"ahler identification with
the corresponding smooth ordinary Betti exact locus.  No identification
of the whole spaces $M$ and $U$ is used.
\end{proposition}

\begin{proof}
The graded-residue condition is algebraic and closed.  Lemma
\ref{r2:lem-proper-parahoric}, whose proof applies to any closed
conjugacy-invariant residue condition, supplies a projective Hitchin
morphism on $M_0$.  Normalization is finite and hence preserves
projectivity of this morphism.

Apply the exact-data correspondence of
\cite[Theorem~7.10]{BiquardGarciaPradaMundet2020} to each possible tuple
of boundary residue orbits.  Lemma
\ref{r2:lem-local-dol-dictionary} identifies its Betti local conditions
with a tuple of classes contained in $\overline C_i$.  Intrinsic
genericity says that all these representations are irreducible.  Their
effective centralizers are finite.  Since $G$ is semisimple, their full
centralizers are finite as well.  The exact correspondence preserves
automorphisms and polystability; thus every polystable object of $M_0$
is stable with finite automorphisms.  This proves the assertion about
all points, including every boundary stratum, using only exact-data
correspondences.  It also excludes strictly semistable objects: the
polystable graded object of a strictly semistable object has a proper
Levi reduction and cannot be stable.  Thus the entire semistable
stack in question is stable, not just its chosen closed-orbit
representatives.

We give the local algebraic and differential-form details.  Let
$\mathfrak B^{\mathrm{fr}}$ be the stack of parahoric bundles with
framings of their reductive quotient fibers at the marked points, and
let $L=\prod_iL_i$.  The cotangent stack of
$\mathfrak B^{\mathrm{fr}}$ parametrizes weak Higgs fields with these
framings; change of framing has moment map equal to the graded residue
in $\bigoplus_i\mathfrak l_i^*$.  This also follows directly from
local duality.  If $\mathfrak p_\alpha$ is the parahoric lattice and
$\mathfrak p_\alpha^+$ its pro-unipotent radical, then
\[
 (\mathfrak p_\alpha^+)^\perp
      =\mathfrak p_\alpha
\]
for the residue pairing $\operatorname{Res}B(a,b)\,dz/z$.
Indeed, for a root $r$, the two relevant exponents are
$\lceil-r(\alpha)\rceil$ and $\lfloor r(\alpha)\rfloor+1$;
their sum is one.  The toral exponents are zero and one.
Consequently the cotangent vectors to bundles with reductive framings
are exactly sections of the weak parahoric Higgs lattice, and its
quotient by the pro-unipotent lattice is the stated moment map.

This stack description can be made finite dimensional near any
point.  Choose a bounded smooth quotient presentation
$[Q/A]$ of the framed bundle stack, using a frame of sufficiently
many global sections of a faithful associated vector bundle.  It may
be chosen $L$-equivariant by retaining the independent puncture
framings.  Here $Q$ is smooth and $A$ is a general linear group.
The classical cotangent presentation is the zero-level reduction of
$T^*Q$ by $A$.  Thus the relative Higgs stack near the point is the
reduction of a smooth symplectic variety $P=T^*Q$ by $K=A\times L$,
with moment constraint
\begin{equation}\label{r2:eq-dol-moment-constraint}
 \mu:P\longrightarrow\mathfrak a^*\oplus
                \bigoplus_i\mathfrak l_i^*,\qquad
 \mu\in\{0\}\times\prod_iY_i.
\end{equation}
The stabilizer of this point is precisely the automorphism group of the Higgs
object.  It is finite at every point under consideration.  The moment
identity identifies the transpose of $d\mu$ with the infinitesimal
$K$-action, followed by the symplectic isomorphism.  Hence $d\mu$ is
surjective.  Shrinking about the point, $\mu$ is a smooth morphism.

 The normalizations $Y_i^\nu$ are symplectic varieties.  Here is an
 explicit form-extension argument, in addition to the associated-bundle
 description.  Take an $H_i$-equivariant resolution
 $Z\to\overline{\mathcal N_i}^{\nu}$, write $\nu_Z$ for its map to
 $\mathfrak h_i^*$, and let $\omega_Z$ be the regular extension of the
 nilpotent orbit form.  On $L_i\times Z$ the coupling form
 \[
 \omega_Z+d\langle b_i+\nu_Z,g^{-1}dg\rangle
 \]
 is basic for $(g,z)\cdot h=(gh,h^{-1}z)$.  Indeed the contractions
 of its two terms with a vertical vector are respectively
 $d\langle\nu_Z,\xi\rangle$ and
 $-d\langle\nu_Z,\xi\rangle$.  It descends to
 $L_i\times^{H_i}Z$ and restricts to the coadjoint orbit form, with
 the convention $\iota_{\xi^\#}\omega=-d\langle\mu,\xi\rangle$,
 consistent with the cotangent convention $\omega=d\lambda$.
 The nilpotent-orbit normalization theorem used for the Betti interior
 therefore supplies the required regular extension also for $Y_i^\nu$.
 Smooth base change gives the normalization
of the constrained space in \eqref{r2:eq-dol-moment-constraint} by replacing
$Y_i$ with $Y_i^\nu$.  Taking its finite-stabilizer quotient gives
$M$ locally.  In particular each component has dense exact locus.
The codimension of a boundary exact stratum is
\begin{equation}\label{r2:eq-dol-boundary-codimension}
 \sum_i\bigl(\dim\mathcal O_i-\dim\mathcal O'_i\bigr),
\end{equation}
so every proper boundary stratum has codimension at least two.

For clarity, extension of the symplectic form does not follow merely
from this statement about singularity types.  Choose $L_i$-equivariant
projective resolutions $Z_i\to Y_i^\nu$, with the regular extensions
$\omega_i$ of the orbit forms.  Such equivariant resolutions exist
by functorial resolution; regularity of these extensions is the
symplectic-singularity property of $Y_i^\nu$.  Smooth base change by
$\mu$ gives a smooth variety
\[
 \widetilde P=P\mathbin{\times}_{\mathfrak a^*\oplus
                    \bigoplus_i\mathfrak l_i^*}
                \left(\{0\}\times\prod_iZ_i\right).
\]
On it the regular closed form
\begin{equation}\label{r2:eq-dol-resolution-form}
 \omega_P-\sum_i\omega_i
\end{equation}
is $K$-invariant and horizontal.  Horizontality follows from the
moment identities on the dense orbit locus and then on all of $Z_i$
by regularity.  The stabilizers on $\widetilde P$ are subgroups of
the already finite stabilizers on $P$.  The quotient is therefore a
smooth Deligne--Mumford stack, proper and birational over the local
normalized moduli space; its coarse space has only finite quotient
singularities.  The displayed form descends as a regular orbifold
form.  Lemma~\ref{paper:lem-finite-quotient-forms} extends the regular
orbifold form over a resolution of this coarse finite quotient.
A common resolution consequently extends the exact-locus
symplectic form over a resolution of $M$.

On the smooth locus of $M$ the form first extends across the
codimension-two complement of the exact locus.  It stays nondegenerate:
the zero locus of its top exterior power, if nonempty, would be a
divisor, whereas it is contained in a subset of codimension at least
two.  Thus $M$ has symplectic singularities in the usual sense.

There are finitely many nilpotent orbits in each $H_i$.  The inverse
image in the normalization of an exact residue stratum is a finite
cover of that stratum; on local orbit charts it is a finite equivariant
cover of a homogeneous orbit, hence is smooth before the finite
stabilizer quotient.  Each such stratum carries its own symplectic
form by ordinary Hamiltonian reduction of the corresponding product
of coadjoint orbits.  Refining by the finitely many local stabilizer
types gives smooth symplectic strata.  These are the strata relevant
here even when normalization makes two ambient symplectic leaves join.

The Hitchin functions commute on each exact stratum.  In the cotangent
description they are invariant polynomials in the Higgs field.  If $p$
is invariant, the invariant pairing identifies $dp(\phi)$ with a section
centralizing $\phi$, since
$B(dp(x),[\xi,x])=0$.  Pairing $p(\phi)$ with a Cech--Serre-dual
coefficient gives a Hamiltonian variation of the form
$(\delta E=c\,dp(\phi),\delta\phi=0)$.  The cocycle condition for this variation is
$[\phi,dp(\phi)]=0$; two such variations pair to zero in the cotangent
form and have zero graded-residue variation.  Thus their pairwise
brackets vanish.  Hamiltonian
reduction preserves this identity, as do the finite covers and the
fixed-stabilizer strata.  Equivalently this is the usual local
commutativity calculation for the principal Hitchin system, which
does not depend on the choice of the coadjoint residue orbit.
Finally, Lemma~\ref{r2:lem-finite-stabilizer-hk} identifies the smooth
coarse exact loci and their symplectic forms, including points with
finite noncentral automorphisms. This proves the asserted
hyperk\"ahler identification.
This completes the proof.
\end{proof}

\begin{remark}[What has and has not been glued]
\label{r2:rem-no-closure-nah}
The closure in this section is taken entirely in the algebraic
Dolbeault moduli space.  The only use of a boundary Betti moduli space
was to check, stratum by stratum, that its points have finite
automorphisms.  In particular none of the conclusions presupposes a
continuous extension of the exact nonabelian Hodge homeomorphism to
the whole orbit closure.  A compact cycle lying in the smooth exact
Dolbeault locus may be transported by the exact correspondence alone.
\end{remark}

\subsection{The exact Dolbeault symplectic class is algebraic in degree two}
\label{r2:sec-dolbeault-class}

The next calculation concerns the complex de Rham class of the
Dolbeault symplectic form.  The coefficients of this class need not be rational or
integral.  The assertion is that the class belongs to the complex span
of first Chern classes of algebraic line bundles.  This distinction is
essential when the semisimple residues have arbitrary complex
eigenvalues.  We prove the assertion on the entire stable exact-residue
stack, including points whose underlying bundle is unstable.

The characteristic-class description has classical antecedents in
cotangent reduction.  In the logahoric setting it also appears
explicitly in \cite[Theorem~6.9 and Section~6.4]{KydonakisZhao2026}
over the regularly stable underlying-torsor locus.  The argument below
supplies the scope needed here on the entire Higgs-stable stack,
including its finite stabilizers; the basic class phenomenon is a
prior input rather than a claimed new discovery.

This scope does not require a nonempty regularly stable underlying-bundle
locus.  Together with extension over the algebraic residue closure, it
makes the formula applicable uniformly in the genus and residue orbit.

\begin{lemma}[Degree two of an algebraic classifying space]
\label{r2:lem-classifying-character-span}
Let $H$ be a complex affine algebraic group, not necessarily connected.
Then
\[
 H^2(BH,\mathbb C)
   =\operatorname{span}_{\mathbb C}
       \{c_1(\mathbb C_\chi):\chi\in X^*(H)\}.
\]
Consequently, if $L$ is a complex reductive group and
$\mathcal O=L/H$ is a coadjoint orbit, then every class in
$H^2_L(\mathcal O,\mathbb C)$ is a complex linear combination of
equivariant first Chern classes of algebraic $L$-linearized line
bundles on $\mathcal O$.
\end{lemma}

\begin{proof}
Suppose first that $H$ is connected.  The unipotent radical of $H$ is
contractible as a complex Lie group, so passage to the reductive
quotient does not change its classifying-space cohomology.  A connected
reductive group is, up to a finite central isogeny, the product of its
connected central torus and a semisimple group.  Finite kernels have
no positive-degree cohomology with complex coefficients.  The
semisimple factor has zero degree-two classifying-space cohomology,
whereas for a torus $(\mathbb C^*)^a$ the first Chern classes of its
characters form a basis of $H^2$.  Character groups also become
isomorphic after tensoring with $\mathbb C$ under these finite
isogenies.  This proves the connected case.

For general $H$, put $\Pi=H/H^0$.  The classifying-space spectral
sequence and the vanishing of positive-degree cohomology of the
finite group $\Pi$ give
\[
 H^2(BH,\mathbb C)=H^2(BH^0,\mathbb C)^\Pi.
\]
The invariant subspace of $X^*(H^0)\otimes\mathbb C$ is spanned by
invariant integral characters: average an integral basis over the
finite group.  If $\chi$ is an invariant character of $H^0$, the
algebraic representation $\operatorname{Ind}_{H^0}^{H}\chi$ has
dimension $|\Pi|$, and the restriction of its determinant to $H^0$
is $\chi^{|\Pi|}$.  Thus a positive multiple of every such invariant
character extends to $H$.  This proves the stated spanning assertion
also in the disconnected case.

Finally, the Borel construction of a homogeneous space gives
$H^*_L(L/H,\mathbb C)=H^*(BH,\mathbb C)$.  The character $\chi$
corresponds under this identification to the algebraic homogeneous
line bundle $L\times^H\mathbb C_\chi$, with its natural
$L$-linearization.
\end{proof}

\begin{proposition}[Chern-class formula for an exact Dolbeault leaf]
\label{r2:prop-dolbeault-algebraic-class}
Let $G$ be semisimple, and consider a finite-type moduli space of
stable weak parahoric $G$-Higgs torsors with fixed exact graded-residue
coadjoint orbits $\mathcal O_i\subset\mathfrak l_i^*$.
Assume the Higgs automorphism groups are finite.  Let $O$ be the
smooth locus of this exact coarse moduli space, with its Dolbeault
symplectic form $\Omega$.  There are algebraic line bundles
$L_1,\ldots,L_N$ on $O$ and coefficients $a_j\in\mathbb C$ such that
\begin{equation}\label{r2:eq-dol-chern-class}
 [\Omega]=\sum_{j=1}^N a_j c_1(L_j)
                 \quad\text{in }H^2(O,\mathbb C).
\end{equation}
The statement applies to arbitrary nilpotent parts of the residues
and to arbitrary real parahoric weights whose algebraic moduli
problem is represented by the fixed facet and stability condition.
\end{proposition}

\begin{proof}
We first make the cotangent description sufficiently global for a
cohomology calculation.  Let $\mathfrak B^{\mathrm{fr}}$ be the
stack of parahoric bundles with a frame of the reductive quotient
fiber at every puncture, and put $L=\prod_iL_i$.  The group $L$
changes these frames.  The canonical cotangent one-form is intrinsic
to the cotangent stack of $\mathfrak B^{\mathrm{fr}}$; its contraction
with an infinitesimal change of frame is the corresponding graded
residue.  The local residue-pairing calculation in
Proposition~\ref{r2:prop-dol-closure} identifies this cotangent stack
with framed weak Higgs fields.  Thus its moment map for $L$ is the
graded-residue map.  This is the stack version of the framed
cotangent and Atiyah-algebroid constructions
\cite[Sections~6--8]{Markman1994}, \cite[Sections~3.1--3.2]{LogaresMartens2010}.

No infinite-dimensional quotient is needed.  The bundles underlying
the finite-type exact Higgs moduli form a bounded family.  A
sufficiently positive twist of a faithful associated vector bundle
is generated by its global sections, with vanishing first
cohomology, throughout this family.  Choosing a basis of those
sections gives a finite-dimensional quotient presentation $[Q/A]$
of a bounded open of $\mathfrak B^{\mathrm{fr}}$, with $Q$ smooth
and $A$ a general linear group, or a product of such groups in an
equivariant-cover presentation.  Smoothness follows from
$H^2$-vanishing for bundle deformations on a curve.  The puncture
frames are independent of the section bases, so the $L$-action on
$Q$ commutes with $A$.  For parahoric torsors one may construct the
same presentation on a rational-weight equivariant cover and then
descend; the fixed-facet replacement of real weights does not change
this algebraic presentation.  The construction is also compatible
with restriction to smaller bounded opens.  Equivalently, its
form identities may be read on the finite-dimensional simplicial
nerve of these quotient charts.

Put $K=A\times L$ and $P=T^*Q$.  Let $\lambda$ be the canonical
Liouville one-form on $P$ and use the convention
$\omega_P=d\lambda$, so that
$\iota_{\xi_P}\omega_P=-d\langle\mu_P,\xi\rangle$.  The cotangent action preserves $\lambda$,
and its moment map satisfies
\begin{equation}\label{r2:eq-dol-cotangent-class}
 \langle\mu_P,\xi\rangle=\iota_{\xi_P}\lambda,
 \qquad
 d_K\lambda=\omega_P-\mu_P.
\end{equation}
Here $d_K=d-\iota$ denotes the Cartan differential, with the
linear moment-map term in equivariant degree two.  One may restrict
to a maximal compact subgroup of $K$ when computing equivariant
cohomology; all differential forms and coefficients are complex.
The displayed equation is a global identity of equivariant forms,
not just an assertion that ordinary local primitives exist.

Give each $\mathcal O_i$ its Kirillov--Kostant--Souriau form
$\omega_i$, with moment map the inclusion
$\iota_i:\mathcal O_i\hookrightarrow\mathfrak l_i^*$.
The equivariantly closed extension of $\omega_i$ is $\omega_i-\iota_i$ in the
same convention.  If $H_i$ is the stabilizer of a chosen point of
$\mathcal O_i$, Lemma~\ref{r2:lem-classifying-character-span}
expresses its class as
\begin{equation}\label{r2:eq-dol-orbit-class}
 [\omega_i-\iota_i]
   =\sum_j b_{ij}\,
       c_1^{L_i}(L_i\times^{H_i}\mathbb C_{\chi_{ij}})
       \quad\text{in }H^2_{L_i}(\mathcal O_i,\mathbb C).
\end{equation}
In particular no integrality assumption on the residue is being
made.  The stabilizer $H_i$ may be disconnected or have a nontrivial
unipotent radical.

Use the shifting description of reduction at coadjoint orbits.
On $P\times\prod_i\mathcal O_i$ take the form
$\omega_P-\sum_i\omega_i$ and the moment map
$\mu_P-\sum_i\iota_i$, where the orbit terms have zero
$\mathfrak a^*$-component.  The zero locus $Z$ of this moment map is precisely
\[
 \mu_A=0,\qquad \mu_{L_i}\in\mathcal O_i,
\]
with the chosen residues recorded.  Restrict to the stable open in
$Z$.  The quotient stack $\mathcal O=[Z/K]$ is the exact stable
Higgs stack, and its symplectic form is the reduction of the
displayed form.  Finiteness of the stabilizers makes the moment
constraint regular and $\mathcal O$ a smooth Deligne--Mumford
stack, as in Proposition~\ref{r2:prop-dol-closure}.

The restriction map in equivariant cohomology sends the equivariant
symplectic class of the product to the reduced symplectic class:
on $Z$ the total moment map is zero, and the remaining two-form is
invariant and horizontal.  By
\eqref{r2:eq-dol-cotangent-class}, the cotangent
contribution $\omega_P-\mu_P$ is equivariantly exact.  Hence
\begin{equation}\label{r2:eq-dol-reduced-class}
 [\Omega_{\mathcal O}]
       =-\sum_{i,j}b_{ij}\,c_1(\mathcal L_{ij}),
\end{equation}
where $\mathcal L_{ij}$ is obtained by restricting the indicated
$K$-linearized homogeneous line bundle to $Z$ and descending to
$[Z/K]$.  These are actual algebraic line bundles on the stack.
They can also be defined intrinsically: fixing the exact residue
reduces its reductive frame torsor to the stabilizer $H_i$, and
$\chi_{ij}$ defines the associated line bundle.  This intrinsic
description shows that the calculation and its line bundles glue
under changes of the bounded quotient presentation.

Lemma~\ref{r2:alg-finite-inertia} descends the Chern-class expression
\eqref{r2:eq-dol-reduced-class} to the smooth coarse locus,
proving \eqref{r2:eq-dol-chern-class}.
\end{proof}

The proposition supplies the algebraic-class hypothesis in
Theorem~\ref{r2:alg-hitchin-rank}.  The proof of the proposition does not assert that the
Dolbeault form is exact, and it does not infer a pure Hodge type on
the nonproper moduli space.  The pure Hodge-type argument is applied
only after restriction to a smooth projective variety in
Lemma~\ref{r2:alg-projective-isotropy}.

\subsection{Algebraic symplectic classes and projective isotropy}
\label{r2:alg-isotropy-section}

For a smooth complex algebraic variety $X$, write
\[
 \operatorname{Alg}^2(X)
   :=\operatorname{span}_{\C}
       \{c_1(L):L\in\operatorname{Pic}(X)\}
       \subset H^2(X,\C).
\]
The coefficients in this span may be arbitrary complex numbers.  The
following argument converts an algebraic description of the symplectic
class on the exact open stratum into a constraint on every projective
subvariety of a resolution.  It uses neither an identification of the
closures under nonabelian Hodge theory nor a vanishing theorem for their
ordinary or intersection cohomology.

\begin{lemma}[Extension of an algebraic degree-two class]
\label{r2:alg-class-extension}
Let $X$ be a smooth irreducible complex quasi-projective variety,
$O\subset X$ a dense Zariski-open subset, and $\xi\in H^2(X,\C)$.
If $\xi|_O\in\operatorname{Alg}^2(O)$, then
$\xi\in\operatorname{Alg}^2(X)$.
In particular this applies to the de Rham class of any closed regular
two-form on $X$.
\end{lemma}

\begin{proof}
Every line bundle on $O$ extends to $X$.  Indeed, choose a nonzero
rational section of the line bundle, take the closure in $X$ of its
Cartier divisor on $O$, and use the fact that every Weil divisor on the
smooth variety $X$ is Cartier.  Thus
$\operatorname{Pic}(X)\to\operatorname{Pic}(O)$ is surjective.

Write $Z=X\setminus O$, put $d=\dim X$, and let
$D_1,\ldots,D_s$ be its irreducible components of codimension one.
The localization sequence gives
\[
 H^2_Z(X,\C)\longrightarrow H^2(X,\C)
                  \longrightarrow H^2(O,\C).
\]
Purity in degree two, equivalently Poincar\'e duality with supports
followed by the description of top Borel--Moore homology, identifies
$H^2_Z(X,\C)$ with the span of the fundamental classes of the $D_j$.
More explicitly, ignoring the Tate twist which is immaterial for the
underlying complex vector space,
\[
 H^2_Z(X,\C)\simeq H^{\mathrm{BM}}_{2d-2}(Z,\C)
       =\bigoplus_{j=1}^s\C[D_j].
\]
Components of $Z$ of higher codimension contribute nothing in this
degree.  The image of $[D_j]$ in $H^2(X,\C)$ is
$c_1(\mathcal O_X(D_j))$.

Express $\xi|_O$ as a finite linear combination of Chern classes and
extend the corresponding line bundles to $X$.  Subtract that linear
combination from $\xi$.  The difference restricts to zero on $O$, so
the localization sequence expresses it as a linear combination of the
boundary divisor classes.  This proves the assertion.
\end{proof}

\begin{lemma}[Projective isotropy]
\label{r2:alg-projective-isotropy}
Let $X$ be smooth and quasi-projective, and let $\widetilde\Omega$ be a
closed regular two-form with
$[\widetilde\Omega]\in\operatorname{Alg}^2(X)$.
For every morphism $f:P\to X$ from a smooth projective variety,
\[
 f^*\widetilde\Omega=0.
\]
Consequently the pullback of $\widetilde\Omega$ to a projective
resolution of any projective subvariety of $X$ is zero.
\end{lemma}

\begin{proof}
The class $[f^*\widetilde\Omega]$ is a complex linear combination of
first Chern classes of algebraic line bundles on $P$, so it belongs to
$H^{1,1}(P)$.  Since $f^*\widetilde\Omega$ is a holomorphic two-form,
the same class belongs to $H^{2,0}(P)$.  The pure Hodge decomposition
of the smooth projective variety $P$ gives
$H^{1,1}(P)\cap H^{2,0}(P)=0$.  Finally the natural map from
holomorphic two-forms to de Rham cohomology on a smooth projective
variety is injective, so the form itself vanishes.  For a singular
projective subvariety, apply this argument after a projective
resolution.  These are the usual consequences of projective Hodge
theory; see also \cite{Deligne1971}.
\end{proof}

\begin{remark}[Why the projective restriction matters]
\label{r2:alg-open-hodge-caution}
The conclusion does not assert that $\widetilde\Omega$ is exact on
$X$.  On a nonproper smooth variety the de Rham class of a closed
regular two-form need not have pure Hodge type $(2,0)$.  The two
incompatible pure Hodge types used in the proof occur after pullback
to the smooth projective variety $P$.
\end{remark}

\subsubsection*{Passing from the stable stack to its coarse space}

\begin{lemma}[Finite-inertia descent of Chern-class spans]
\label{r2:alg-finite-inertia}
Let $q:\mathcal O\to O$ be the coarse, hence good, moduli morphism of
a finite-type complex Deligne--Mumford stack with finite inertia.
Let $O^{\circ}\subset O$ be a smooth open subset and
$\mathcal O^{\circ}=q^{-1}(O^{\circ})$.
Suppose a closed regular two-form $\Omega$ on $O^{\circ}$ satisfies
\[
 q^*[\Omega]
     \in\operatorname{span}_{\C}
          \{c_1(\mathcal L):\mathcal L\in
                        \operatorname{Pic}(\mathcal O^{\circ})\}.
\]
Then $[\Omega]\in\operatorname{Alg}^2(O^{\circ})$.
\end{lemma}

\begin{proof}
The orders of geometric stabilizers are bounded on the finite-type
stack with finite inertia.  Choose a positive integer $N$ divisible
by all those orders.  For every line bundle $\mathcal L$, stabilizers
act trivially on the fibers of $\mathcal L^{\otimes N}$.
The descent criterion for good moduli spaces
\cite[Theorem~10.3]{AlperGoodModuli2013} therefore gives a line bundle
$L$ on $O^{\circ}$ with
\[
 \mathcal L^{\otimes N}\simeq q^*L.
\]
Since finite groups have no positive-degree cohomology over $\C$,
the coarse morphism induces an isomorphism
\[
 q^*:H^2(O^{\circ},\C)\xrightarrow{\ \sim\ }
                       H^2(\mathcal O^{\circ},\C).
\]
This can be seen on finite-quotient charts and then from the Leray
spectral sequence.  Replacing each stack line bundle in the assumed
finite expression by its $N$th power and dividing its coefficient by
$N$ proves the claim.
\end{proof}

For generic stable objects with semisimple structure group, the center
is finite and the full stabilizers are finite.  Thus the preceding
descent applies on the finite-type stable moduli stack.  For a
reductive group with positive-dimensional center, finite
\emph{effective} stabilizers alone do not suffice to descend arbitrary
line bundles: a nontrivial character of the central torus is not
annihilated by taking a tensor power.  One must rigidify and verify the
central weights, or apply the established reduction to semisimple
simply connected factors and torus factors.  The latter reduction is
sufficient for the compactification application here.
In particular, Proposition~\ref{r2:prop-dolbeault-algebraic-class}
supplies the required algebraic class on the exact Dolbeault stratum;
Lemma~\ref{r2:alg-finite-inertia} transfers that statement from the
stable stack to its smooth coarse open subset.

\subsubsection*{Forcing the Hitchin rank and a compact smooth fiber}

\begin{theorem}[Algebraic-class criterion for a Lagrangian fibration]
\label{r2:alg-hitchin-rank}
Let $M$ be an irreducible normal complex quasi-projective variety of
dimension $2n$.  Suppose:
\begin{enumerate}
\item the smooth locus carries a closed symplectic form $\Omega$,
and on a projective resolution $\pi:X\to M$, chosen to be an
isomorphism over $M_{\mathrm{reg}}$, its pullback extends to a regular
two-form $\widetilde\Omega$;
\item there is a dense smooth open subset $O\subset M_{\mathrm{reg}}$
on which $[\Omega|_O]\in\operatorname{Alg}^2(O)$;
\item $h:M\to A$ is projective, $A$ is affine, and the pullbacks of
regular functions on $A$ Poisson-commute on $M_{\mathrm{reg}}$.
\end{enumerate}
Then the reduced image $B=h(M)_{\mathrm{red}}$ has dimension $n$, and
$h$ has generic rank $n$.  A general fiber is generically Lagrangian
where it meets $M_{\mathrm{reg}}$.

Assume further that $M$ has a finite stratification by smooth locally
closed symplectic varieties
\[
 M=O\sqcup\coprod_{\alpha\ne0}S_\alpha,
 \qquad \dim S_\alpha<2n,
\]
such that the restrictions of the pulled-back functions from $A$
Poisson-commute with respect to the symplectic form of \emph{each}
stratum.  Then a dense open subset $B^{\circ}\subset B$ has the
property that every fiber over $B^{\circ}$ is a nonempty smooth
projective Lagrangian variety contained entirely in $O$.
\end{theorem}

\begin{proof}
The form $\widetilde\Omega$ is closed: its exterior derivative
vanishes on a dense open subset and is a regular form.  View $O$ as
an open subset of $X$.  Lemma~\ref{r2:alg-class-extension} gives
\[
 [\widetilde\Omega]\in\operatorname{Alg}^2(X).
\]
Let $r=\dim B$, which equals the generic differential rank in
characteristic zero.  The Poisson-commuting hypothesis implies
$r\leq n$: the Hamiltonian vectors of any $r$ generically independent
pullbacks span an isotropic subspace in a symplectic vector space of
dimension $2n$.

Choose a general point of $M_{\mathrm{reg}}$ and an irreducible
component $F$ of its fiber.  This component has dimension $2n-r$ and
is projective because $h$ is projective.  The strict transform of $F$ in $X$
is projective over $F$, hence projective over $\C$.  Resolve that
strict transform projectively.  Lemma~\ref{r2:alg-projective-isotropy}
makes the pullback of $\widetilde\Omega$ zero.  Since $\pi$ is an
isomorphism over $M_{\mathrm{reg}}$, the smooth part of $F$ there is
generically isotropic.  Therefore
\[
 2n-r=\dim F\leq n.
\]
Combining the two inequalities gives $r=n$.  The same argument for
the components of a general fiber proves the first assertion about
Lagrangian fibers.

For the final assertion, the commuting restrictions on $S_\alpha$
give
\[
 \dim\overline{h(S_\alpha)}
       \leq\tfrac12\dim S_\alpha<n.
\]
Here the dimension of the image is its generic differential rank on
each irreducible component of the stratum.  The forms on these
strata need not be obtained by restricting the ambient form across
singularities; their own commuting-function hypothesis is what is
used.  Remove the finitely many proper closed subsets
$\overline{h(S_\alpha)}$ from $B$, and then shrink to a smooth dense
open over which $h|_O$ is smooth.  This is possible by generic
smoothness in characteristic zero.  Every fiber over the resulting
open $B^{\circ}$ is contained in $O$, is smooth of dimension $n$, and
is projective.  The isotropy already proved makes every connected
component Lagrangian.
\end{proof}

\begin{corollary}[The compact fiber supplies a Goldman period]
\label{r2:alg-goldman-cycle}
Under the final hypotheses of
Theorem~\ref{r2:alg-hitchin-rank}, suppose the exact open stratum $O$
has a hyperk\"ahler identification with a smooth open subset of the
normalized Betti moduli space.  Let $I$ be the Dolbeault complex
structure, and choose the hyperk\"ahler conventions
\[
 \Omega_I=\omega_J+i\omega_K,\qquad
 \Omega_J=\omega_K+i\omega_I,
\]
where $\Omega_J$ corresponds to the Goldman form, up to a fixed
nonzero scalar.  Then the top power of the Goldman form has a
nonzero period on a compact real $2n$-cycle in the exact Betti open
stratum.
\end{corollary}

\begin{proof}
Take a connected component $F$ of a fiber over $B^{\circ}$.  It is a
compact complex Lagrangian submanifold for $\Omega_I$.
Consequently $\omega_J|_F=\omega_K|_F=0$, and
\[
 \int_F\Omega_J^n=i^n\int_F\omega_I^n\ne0,
\]
because $\omega_I|_F$ is a K\"ahler form.  A fixed nonzero scalar in
the Goldman identification multiplies this period by its $n$th
power and cannot make it zero.  Transporting the fundamental cycle
of $F$ under the exact-stratum identification proves the assertion.
For $n=0$, the same statement is the nonzero degree of a point.
\end{proof}

Only the exact open stratum needs the hyperk\"ahler identification in
this argument.  The algebraic closure of this stratum provides the proper Hitchin
map and the extension of the symplectic form.  Thus the remaining
geometric inputs for applying the Goldman-period compactification
criterion are explicit: the exact-stratum Chern-class formula,
symplectic extension on the normalized Dolbeault closure, projective
Hitchin properness, and the finite symplectic stratification with
commuting restrictions.

\section{Compactifications for conjugacy-class closures}
\label{r2:sec-full-goal}

We now combine the Poisson criterion with the compact cycle constructed
inside the exact Dolbeault stratum. The algebraic closure supplies a
projective Hitchin morphism and an extension of the two-form; the
exact-stratum correspondence transports the compact cycle to the
Betti side.

\begin{theorem}[All generic normalized relative character varieties]
\label{r2:thm-full-normalized-logcy}
Let $G$ be any connected complex reductive algebraic group, and let
$\Sigma$ be a compact Riemann surface of any genus with a nonempty
finite set of punctures.  Let $\mathbf C=(C_1,\ldots,C_k)$ be a
compatible intrinsically generic tuple of conjugacy classes, and put
\[
 U_0=
 \left\{(A_j,B_j,M_i)\in G^{2g}\times\prod_i\overline C_i:
       \prod_{j=1}^g[A_j,B_j]\prod_{i=1}^kM_i=1\right\}_{\red}/\!/G.
\]
For every nonempty irreducible component $U$ of $U_0^\nu$, there is a
normal projective compactification $U\subset X$ whose reduced boundary
$D=(X\setminus U)_{\red}$ is a divisor and such that
\[
 (X,D)\ \text{is log canonical},\qquad
 K_X+D\sim0
 \quad\text{integrally}.
\]
One may use the Poisson compactification constructed in
Section~\ref{sec:poisson-bridge}.  The boundary of this compactification is a Poisson closed
subset and supports an effective ample Cartier divisor.  The variety
$X$ is of Fano type, hence has rational singularities and is
Cohen--Macaulay.  If $U_0$ is normal, the conclusion applies to its
components without changing the interior.
\end{theorem}

\begin{idea}
By Theorem~\ref{thm:poisson-nonvanishing}, it suffices to construct a Poisson
compactification of $U$ (Section~\ref{sec:poisson-bridge}) and one nonzero
logarithmic pluricanonical form on a resolution of $U$. The form comes from the
Goldman class (Proposition~\ref{prop:goldman-f2-all-types}) once a compact cycle
with nonzero Goldman-volume period is available, and
Section~\ref{r2:sec-dolbeault-closure} supplies such a cycle: a compact
Lagrangian Hitchin fiber in the exact Dolbeault stratum, transported to the
Betti side by hyperk\"ahler rotation. Central isogenies and the torus factor are
treated by Proposition~\ref{new:gl-simple-reduction}.
\end{idea}

\begin{proof}
First suppose that $G$ is semisimple and simply connected, and write
$\dim U=2n$.  If $n=0$, the component is a point and the assertion is
immediate, with empty boundary.  Suppose $n>0$.

The inverse image of the prescribed exact-class locus is dense in $U$.
The intersection of this inverse image with $U_{\mathrm{reg}}$ is a nonempty
irreducible smooth open subset. The exact-data tame correspondence identifies
this intersection with a smooth exact Dolbeault locus $O$.  Take the component of the algebraic normalized
Dolbeault closure containing $O$, as in
Proposition~\ref{r2:prop-dol-closure}, and denote it by $M$.
This proposition establishes all of the following algebraically:
$M$ is quasiprojective with symplectic singularities; its Hitchin
morphism is projective over an affine base; its complement of $O$
has a finite refinement by lower-dimensional symplectic strata; and
the Hitchin functions commute on every stratum.  The dimension of
$M$ is $2n$, because its exact open subset is identified with the
Betti exact open subset.

Proposition~\ref{r2:prop-dolbeault-algebraic-class}, including its
finite-inertia descent, gives
\[
 [\Omega_I|_O]\in\operatorname{Alg}^2(O).
\]
Theorem~\ref{r2:alg-hitchin-rank} therefore produces a smooth
projective Lagrangian fiber, entirely contained in $O$.  Let $F$ be
one of its connected components.  By
Corollary~\ref{r2:alg-goldman-cycle}, its image under the
\emph{exact-stratum} hyperk\"ahler correspondence is a compact real
$2n$-cycle in the smooth exact Betti locus with
\[
 \int_F\omega_{\mathrm{Goldman}}^{\wedge n}
      =c^n\sqrt{-1}^{\,n}\int_F\omega_I^{\wedge n}\ne0
\]
for a fixed nonzero normalization constant $c$.

Choose a projective resolution $\widetilde U\to U$ that is an
isomorphism over $U_{\mathrm{reg}}$.  The compact cycle lies in that
open subset and hence lifts unchanged to $\widetilde U$.  Thus
\[
 [\widetilde\omega_{\mathrm{Goldman}}]^n
          \ne0\quad\text{in }H^{2n}(\widetilde U,\C).
\]
The normalized Betti interior has the canonical symplectic
singularities and integral trivialization of its canonical divisor
already established in Theorem~\ref{new:thm-normalized-interior}.
The all-type $F^2$ statement of
Proposition~\ref{prop:goldman-f2-all-types} and its resolution extension
therefore apply.  Corollary~\ref{cor:goldman-period} now gives the
claimed compactification.  The effective ample boundary support and
Proposition~\ref{prop:poisson-fano-type} give the assertion about Fano
type and its consequences.

For general connected reductive $G$, take the finite central covering
\[
 T\times\prod_a G_a^{\mathrm{sc}}\longrightarrow G
\]
by a torus and simply connected simple groups.  A punctured surface
has free fundamental group, so its representations lift after choosing
lifts of free generators.  Retain all the lifted puncture classes and
normalize componentwise.  Intrinsic genericity passes to the lifted
data, and the simply connected result applies to each semisimple
factor. The torus factor has its invariant logarithmic volume form.
The central-twisting action preserves the product Goldman volume.
Proposition~\ref{new:gl-simple-reduction} and its finite-map valuation
argument show that this invariant volume descends with at most
logarithmic poles. Thus logarithmic nonvanishing descends to $U$;
applying the Poisson criterion there gives the compactification.
This proves the statement
for every connected reductive group.
\end{proof}

\begin{remark}[The two cohomological inputs]
The Dolbeault form has a class generated by algebraic line bundles;
this is used only after restriction to a smooth projective variety,
where Hodge decomposition forces isotropy.  The Betti Goldman form
has its independent Deligne $F^2$ property; together with the nonzero
period this provides logarithmic nonvanishing at infinity.  The
hyperk\"ahler correspondence transports the compact cycle and the
actual differential forms on the exact open strata.  It is not
asserted to preserve the Betti mixed Hodge structure.

The stronger demand for a $\Q$-factorial dlt compactification keeping
the whole original interior remains impossible in general, for the
normal $\SL_3$ example constructed in
Theorem~\ref{new:sl3-main}.  The theorem above concerns the normal
projective log canonical compactification with integral $K_X+D\sim0$.
\end{remark}

\part{Normality and examples}\label{part:normality}

\section{An exact normality criterion and the positive-genus classification}
\label{r2:norm-section}

Theorem~\ref{paper:thm-normalized} applies to normalized interiors.
The classification below is independent of both compactification
constructions; its application is given in
Section~\ref{subsec:normality-compactification}.
Finite stabilizer invariants can remove a normalization defect. We
determine exactly when this occurs and prove that it cannot occur for
intrinsically generic data in positive genus.

Keep the compatible intrinsically generic data of
Section~\ref{new:sec-base-change}, put $G_{\ad}=G/Z(G)$, and write
$q:R_{\mathbf Y}\to U_{\mathbf Y}$ for the geometric quotient.
At a tuple
\[
 r=(A_1,B_1,\ldots,A_g,B_g,M_1,\ldots,M_k)
\]
put $p=q(r)$ and $F=(G_{\ad})_r$, a finite group.  The
\emph{realized exact classes} are $\mathcal O_i=G\cdot M_i$; they may
be smaller than the originally prescribed open classes $\cC_i$.
Their relative representation stratum is smooth.  The quotient of this stratum has an
\'etale-local presentation as the quotient of a smooth germ by $F$.
Let $W_r$ denote the tangent representation of this smooth germ, and put
\[
 K_r=\ker\bigl(F\longrightarrow\operatorname{GL}(W_r)\bigr).
\]
Normality always refers to the given reduced variety.

\emph{Strategy.} Normality is a local property, so we work in the completed
local ring at $p$. A computation in the group directions gives
$\widehat{\mathcal O}_{U_{\mathbf Y},p}\simeq(A[[W_r^*]])^F$, where $A$ is the
completed product of the puncture-transversal rings
(Proposition~\ref{r2:norm-equivariant-model}). If $F$ acts faithfully on $W_r$,
polynomial invariants of $F$ detect every normalization defect of $A$, so
normality at $p$ is equivalent to normality of every $Y_i$ at $M_i$; in general
only the kernel $K_r$ matters (Theorem~\ref{r2:norm-kernel-criterion}). In
genus at least two the action on $W_r$ is faithful. In genus one we instead
show that the kernel cannot remove a puncture-normalization defect
(Lemma~\ref{r2:norm-genus-one-kernel}). On the sphere it can do so;
Section~\ref{r2:sec-sphere-masking} gives an explicit generic example.

\subsection{An equivariant formal local model}

\begin{proposition}[Puncture transversals and the deformation representation]
\label{r2:norm-equivariant-model}
For each $i$, there is an $F$-invariant formal transversal $T_i$ to
$\mathcal O_i$ in $Y_i$ at $M_i$.  If
\[
 A=\widehat\bigotimes_{i=1}^k\widehat{\mathcal O}_{T_i,M_i},
\]
then a formal Luna slice at $r$ has coordinate ring
\[
 A[[W_r^*]]
\]
with the diagonal $F$-action: $F$ acts on $A$ through the puncture
transversals and linearly on $W_r$.  Consequently
\begin{equation}\label{r2:norm-local-invariants}
 \widehat{\mathcal O}_{U_{\mathbf Y},p}
       \simeq\bigl(A[[W_r^*]]\bigr)^F.
\end{equation}
The ring $A$ is reduced and excellent, and it is normal if and only if
every $Y_i$ is normal at $M_i$.
\end{proposition}

\begin{proof}
Write $\mathfrak g=\Lie G$, $\mathfrak g_{\ad}=\Lie G_{\ad}$ and
$E_i=T_{M_i}\mathcal O_i\subset T_{M_i}G$. All coordinate changes below are
first chosen analytically and then completed, so normality of analytic germs and
of their completions may be compared freely.

\emph{Step 1: puncture transversals.} Choose $F$-stable complements
$N_i'\oplus E_i=T_{M_i}G$ and $\mathfrak q_i\oplus\Lie\bigl((G_{\ad})_{M_i}\bigr)
=\mathfrak g_{\ad}$. Averaged coordinates turn $N_i'$ into an $F$-stable smooth
formal germ $N_i\subset\widehat G_{M_i}$ through $M_i$, and the map
\[
 \widehat{\mathfrak q_i}\times N_i\longrightarrow\widehat G_{M_i},\qquad
 (u,x)\longmapsto\exp(u)\,x\,\exp(-u),
\]
is $F$-equivariant with invertible differential, hence an isomorphism of formal
germs. Since $Y_i$ is conjugation invariant, the inverse image of
$\widehat{(Y_i)}_{M_i}$ is $\widehat{\mathfrak q_i}\times T_i$ with
$T_i=Y_i\cap N_i$; this holds scheme-theoretically, because conjugation and its
inverse both preserve the ideal of $Y_i$.

\emph{Step 2: the relation.} Exponential coordinates on the handles give the
ambient germ $\operatorname{Spf}A[[L^*]]$ with
$L=\mathfrak g^{2g}\oplus\bigoplus_i\mathfrak q_i$, on which $F$ acts linearly
on $L$. The surface word takes values in $G_\der$, and by
Lemma~\ref{lem:differential} and irreducibility its differential in the
$L$-directions is an $F$-equivariant surjection $D:L\to\mathfrak g_\der$. With an
$F$-stable right inverse of $D$, the formal implicit function theorem over the
complete base $A$ gives an $F$-equivariant isomorphism
\[
 \widehat{\mathcal O}_{R_{\mathbf Y},r}\simeq A[[(\ker D)^*]];
\]
equivariance follows from uniqueness of the implicit solution.

\emph{Step 3: the slice.} Since the effective stabilizer is finite, the tangent
space of the conjugation orbit is a copy $J\simeq\mathfrak g_{\ad}$ inside
$\ker D$. Choose an $F$-stable complement $\ker D=J\oplus W_r$ and set the $J$
coordinates to zero. For a formal point $x$ there is a unique formal group element
$u$ such that the $J$ coordinates of $u^{-1}x$ vanish, again by the implicit
function theorem, because the derivative in $u$ is invertible. Uniqueness gives
both the product decomposition along the orbit and $F$-equivariance, so
$A[[W_r^*]]$ is the coordinate ring of a formal Luna slice \cite{Luna1973}. Taking
$F$-invariants gives \eqref{r2:norm-local-invariants}.

\emph{Step 4: properties of $A$.} Setting the puncture transversals to their base
points gives the realized exact-class stratum, whose tangent representation is
therefore $W_r$. Each $\widehat{(Y_i)}_{M_i}$ is a smooth germ times $T_i$, so
reducedness, excellence and normality of $A$ are equivalent to the corresponding
properties of the $Y_i$ at the $M_i$.
\end{proof}

\subsection{Which normalization defects survive finite invariants}

The following example shows both phenomena behind
Theorem~\ref{r2:norm-kernel-criterion}: a finite quotient can hide a
normalization defect, and a faithful smooth factor makes the defect visible.

\begin{example}
Let $A=\C[[t^2,t^3]]\subset\C[[t]]$ be the completed cusp, and let $\mu_2$ act by
$t\mapsto-t$. Then $A^{\mu_2}=\C[[t^2]]$ is normal although $A$ is not. Now add
a smooth direction $w$ on which $\mu_2$ acts faithfully, $w\mapsto-w$. The ring
$(A[[w]])^{\mu_2}$ contains $t^2$ and $t^3w$, hence its fraction field contains
$tw=t^3w/t^2$, which is integral because $(tw)^2=t^2w^2$. But $tw\notin
(A[[w]])^{\mu_2}$, since $A$ has no term $t$; so $(A[[w]])^{\mu_2}$ is not
normal. In the notation above, $K_r=\mu_2$ in the first case and $K_r=1$ in the
second.
\end{example}

\begin{lemma}[Detection after a faithful smooth factor]
\label{r2:norm-faithful-factor}
Let a finite group $F$ act on a reduced excellent complete local
$\C$-algebra $A$ and on a finite-dimensional vector space $W$.
Put $K=\ker(F\to\operatorname{GL}(W))$, let $\overline A$ be the
normalization of $A$, and set $M=\overline A/A$.
Then
\begin{equation}\label{r2:norm-defect-invariants}
 \overline{\bigl(A[[W^*]]\bigr)^F}\,
       \big/\bigl(A[[W^*]]\bigr)^F
       \simeq\bigl(M[[W^*]]\bigr)^F.
\end{equation}
In particular,
\begin{equation}\label{r2:norm-kernel-test}
 \bigl(A[[W^*]]\bigr)^F\text{ is normal}
       \quad\Longleftrightarrow\quad A^K\text{ is normal}.
\end{equation}
If $Q=F/K$, the right side is equivalently the vanishing of the finite
collection of covariant spaces
\begin{equation}\label{r2:norm-finite-degree-test}
 \bigoplus_{j=0}^{|Q|-1}
       \bigl(M^K\otimes_\C\operatorname{Sym}^jW^*\bigr)^Q.
\end{equation}
Thus a normalization defect which survives $K$ is already detected in
polynomial degree at most $|Q|-1$ in the smooth deformation variables.
\end{lemma}

\begin{proof}
\emph{The defect module.} Since $A$ is excellent, $\overline A$ is finite over
$A$, and $\overline A[[W^*]]$ is the normalization of $A[[W^*]]$. Normalization
commutes with taking invariants of a finite group in characteristic zero, so the
normalization of $(A[[W^*]])^F$ is $(\overline A[[W^*]])^F$, and exactness of
$F$-invariants gives \eqref{r2:norm-defect-invariants}. In the same way
$\overline A^K$ is the normalization of $A^K$, with defect $M^K$.

\emph{Reduction to $Q$.} Since $K$ acts trivially on $W$, the defect equals
$(M^K[[W^*]])^Q$ with $Q=F/K$ acting faithfully on $W$. If $M^K=0$ this
vanishes. Conversely, let $M^K\ne0$. The faithful action of $Q$ has a free orbit
$Qw\subset W$, because the fixed spaces of the finitely many nonidentity elements
are proper. Products of at most $|Q|-1$ affine linear functions interpolate the
delta functions on $Qw$, so restriction to the orbit gives a $Q$-equivariant
surjection
\[
 \textstyle\bigoplus_{j=0}^{|Q|-1}\operatorname{Sym}^jW^*\twoheadrightarrow\C[Qw]\simeq\C[Q].
\]
Every irreducible $Q$-representation occurs in the regular representation
$\C[Q]$, hence, by semisimplicity, in the left side. A nonzero element of $M^K$
spans a finite-dimensional $Q$-stable subspace, which contains an irreducible
$V$; pairing $V$ with a copy of $V^*$ in polynomials of degree at most
$|Q|-1$ gives a nonzero element of \eqref{r2:norm-finite-degree-test}, and
polynomial tensors embed in $M^K[[W^*]]$. This proves \eqref{r2:norm-kernel-test}
and the finite-degree test.
\end{proof}

\begin{theorem}[Exact local normality test in every genus]
\label{r2:norm-kernel-criterion}
For arbitrary $g\geq0$, the generic relative character variety
$U_{\mathbf Y}$ is normal at $p=q(r)$ if and only if
\[
 A^{K_r}\text{ is normal},
\]
where $A$ is the completed product of the puncture transversals in
Proposition~\ref{r2:norm-equivariant-model} and $K_r$ is computed from
the realized exact-class deformation representation $W_r$.
Equivalently, the finite-degree test
\eqref{r2:norm-finite-degree-test} vanishes for this puncture ring.
If the action of $F$ on $W_r$ is faithful, normality at $p$ is equivalent
to normality of every $Y_i$ at $M_i$.
\end{theorem}

\begin{proof}
Apply Lemma~\ref{r2:norm-faithful-factor} to
\eqref{r2:norm-local-invariants}.  Excellence identifies normality of
the local algebraic ring with normality of its completion.  For a
faithful action $K_r=1$, and the last assertion is the puncture
normality equivalence in Proposition~\ref{r2:norm-equivariant-model}.
\end{proof}

\subsection{The kernel in positive genus}

\begin{lemma}[The moving deformation space]
\label{r2:norm-moving-space}
For $f\in F$, put
\[
 V_f=\operatorname{im}(\operatorname{Ad}(f)-1)\subset\mathfrak g.
\]
This space is preserved by the tuple.  The sum of the nonidentity
eigenspaces of $f$ on $W_r$ has dimension
\begin{equation}\label{r2:norm-moving-dimension}
 (2g-2)\dim V_f+
       \sum_{i=1}^k\operatorname{rank}
          \bigl((\operatorname{Ad}(M_i)-1)|_{V_f}\bigr).
\end{equation}
Consequently the $F$-action on $W_r$ is faithful when $g\geq2$.
When $g=1$ and $f\in K_r$, every $M_i$ acts trivially on each simple
ideal of $\mathfrak g_{\der}$ on which $f$ acts nontrivially.
\end{lemma}

\begin{proof}
The deformation complex for the realized exact classes is
\[
 \mathfrak g\xrightarrow{\delta}
 \mathfrak g^{2g}\oplus\bigoplus_i
          \operatorname{im}(\operatorname{Ad}(M_i)-1)
 \xrightarrow{D}\mathfrak g_{\der}.
\]
The middle cohomology of this complex is $W_r$.  The kernel of $\delta$ is
$\mathfrak z(\mathfrak g)$, and $D$ is surjective by irreducibility.
All maps commute with $f$.  On the nonidentity eigenspaces, the end
spaces both identify with $V_f$, the first map is injective, and the
last is surjective.  Taking dimensions proves
\eqref{r2:norm-moving-dimension}.  Every nonidentity $f\in G_{\ad}$ acts
nontrivially on $\mathfrak g$, so the first term is positive when
$g\geq2$.

For $g=1$ and $f\in K_r$, all the nonnegative ranks in
\eqref{r2:norm-moving-dimension} vanish.  Hence every $M_i$ centralizes
$V_f$ pointwise.  The Lie subalgebra generated by $V_f$ is an ideal in
$\mathfrak g$: it is stable under $\mathfrak g^f$, since
$[\mathfrak g^f,V_f]\subset V_f$, and under $V_f$ itself.
The Lie subalgebra generated by $V_f$ is therefore the sum of precisely
the simple ideals on which $f$ acts nontrivially.  Every $M_i$ acts trivially on that ideal, proving
the final assertion.
\end{proof}

\begin{lemma}[Generic central data exclude non-A factors]
\label{r2:norm-central-data}
Suppose all the puncture elements have trivial image in a fixed simple
factor of $G/Z(G)$.  Intrinsic genericity implies that this factor has
type A.
\end{lemma}

\begin{proof}
For every irreducible root system other than type A there is a
fundamental weight $\varpi$ in the root lattice.  One can take
$\varpi_1$ in type $B_\ell$, $\varpi_2$ in types $C_\ell$ and
$D_\ell$, the adjoint fundamental weight in types $E_6$ and $E_7$,
and any fundamental weight in the center-free types $E_8,F_4,G_2$.
The usual low-rank coincidences are understood as type A where
appropriate.  These are the standard root- and weight-lattice
calculations \cite{HumphreysLie1972}.

Take the maximal parabolic corresponding to $\varpi$ in the selected
factor and its inverse image $P$ in $G$.  To see exactly which
characters occur, use a finite central isogeny
\[
 \widetilde G=Z(G)^\circ\times\prod_a G_a^{\mathrm{sc}}
       \longrightarrow G.
\]
The characters of the inverse image of $P$ are a character of the
central torus together with an integral multiple of $\varpi$.
Since $\varpi$ is in the root lattice, it is trivial on the center of
the selected simply connected factor.  It therefore descends through
every central kernel, independently of the central-torus character.
It follows that
\[
 X^*(P)=X^*(G)\oplus\Z\varpi.
\]
Choose representatives of the semisimple puncture classes in a common
maximal torus of $P$.  Compatibility makes the product of their values
under every character from $X^*(G)$ equal to one.  Their images in the
selected adjoint factor are trivial, so their values under $\varpi$
are individually one.  Their product consequently has trivial image
in $P/[P,P]$, contradicting intrinsic genericity.  Thus the selected
simple factor must be of type A.
\end{proof}

\begin{lemma}[The genus-one kernel cannot remove a puncture defect]
\label{r2:norm-genus-one-kernel}
In genus one, with $A$ and $K_r$ as above,
\[
 A^{K_r}\text{ is normal}\quad\Longleftrightarrow\quad A\text{ is normal}.
\]
\end{lemma}

\begin{proof}
Write $K=K_r$, and let $I$ be the set of simple factors of $G/Z(G)$ on
which some element of $K$ acts nontrivially.  By
Lemma~\ref{r2:norm-moving-space}, all punctures have trivial image in
each of these adjoint factors.  Lemma~\ref{r2:norm-central-data} makes
every factor in $I$ type A.  The action of $K$ on all other simple
factors and on the connected center is trivial.

We explain the local factorization through a central isogeny, including
the possible intersection of different lifted conjugacy closures.
Choose lifts $\widetilde M_i$ under
$\widetilde G\to G$.  The inverse image of one conjugacy class is the
union of the central translates of one $\widetilde G$-conjugacy class.
The reduced closure of this preimage is the union of the corresponding translated
closures.  In a factor indexed by $I$, the chosen point
$\widetilde M_i$ is central.  Every translated closure passing through
that point must have the same central semisimple value in this factor.
The central multiplier of the branch at that point is therefore fixed.  All such branches
have the identical unipotent orbit closure in that type-A factor;
central multiplication does not change its unipotent Jordan orbit.
The reduced union thus factors locally as
\[
  \left(\text{a product of type-A unipotent orbit closures}\right)
       \times\left(\text{a reduced union in the other factors}\right).
\]
The first factor is normal by Kraft--Procesi
\cite{KraftProcesi1979}.  Finite central isogenies are \'etale over
$\C$, so the same statement identifies the completed puncture germs.

There is a small equivariance point in using these lifts.  An arbitrary
$f\in F$ can act on a chosen lifted germ by conjugation followed by a
central translation restoring the chosen point.  For $k\in K$, this
translation is the identity: on the factors in $I$ the point is central,
and on the other factors the lift of $k$ is central.  Hence $k$ commutes
with every chosen $\widetilde M_i$.  It acts by conjugation on the
displayed type-A factor and identically on the other factor.

Choose the formal orbit transversals factor by factor.  There are no
orbit directions at the central points in the factors indexed by $I$.
Inner automorphisms do not permute the simple factors, so the
transversals can also be chosen $F$-invariant.  Taking their completed
product yields
\[
 A\simeq A_0\widehat\otimes_\C A_1,
 \qquad
 A^K\simeq A_0^K\widehat\otimes_\C A_1,
\]
where $A_0$ is normal and $K$ acts trivially on $A_1$.
The second equality follows from the Reynolds projector, which commutes
with these tensor products and their adic inverse limits.  The ring
$A_0^K$ is normal.  A product of complex algebraic germs is normal
exactly when both factors are normal; the same statement holds for
their excellent completions.  Consequently both $A$ and $A^K$ are
normal exactly when $A_1$ is normal.  This proves the lemma.
\end{proof}

\begin{theorem}[Normality in every positive genus]
\label{r2:norm-positive-genus}
Let $G$ be any connected complex reductive group, and let
$\mathbf Y$ be compatible intrinsically generic conjugacy-class
closures.  If $g\geq1$, then for every tuple $r$ and $p=q(r)$,
\begin{equation}\label{r2:norm-pointwise-equivalence}
 \begin{split}
 U_{\mathbf Y}\text{ is normal at }p
 &\quad\Longleftrightarrow\quad R_{\mathbf Y}\text{ is normal at }r\\
 &\quad\Longleftrightarrow\quad
       Y_i\text{ is normal at }M_i\text{ for every }i.
 \end{split}
\end{equation}
Moreover $R_{\mathbf Y}\ne\varnothing$, and
\begin{equation}\label{r2:norm-global-equivalence}
 U_{\mathbf Y}\text{ is normal}
       \quad\Longleftrightarrow\quad
       \text{every }Y_i\text{ is normal}.
\end{equation}
\end{theorem}

\begin{proof}
For $g\geq2$, Lemma~\ref{r2:norm-moving-space} makes the deformation
representation faithful, so Theorem~\ref{r2:norm-kernel-criterion}
gives the equivalence with puncture normality.
For $g=1$, combine that theorem with
Lemma~\ref{r2:norm-genus-one-kernel}.  The representation-space
equivalence is Proposition~\ref{new:prop-local-transfer}.

For the global assertion, choose arbitrary $M_i\in Y_i$.
Compatibility implies that $(\prod_i M_i)^{-1}\in G_{\der}$.
By the commutator theorem of Pasiencier--Wang and Ree, every element of a
connected complex semisimple algebraic group is a commutator
\cite{ReeCommutators1964}; see \cite[Thm.~1.1]{MalleBourbaki1069}.
Choose $A_1,B_1\in G_{\der}$ with
$[A_1,B_1]=(\prod_iM_i)^{-1}$, and set all remaining handles equal to
the identity.  This realizes the chosen puncture tuple.  Intrinsic
genericity makes the resulting representation irreducible.
Every puncture point is therefore realized, and applying
\eqref{r2:norm-pointwise-equivalence} proves
\eqref{r2:norm-global-equivalence} and nonemptiness.
\end{proof}

\subsection{The conductor locus, its codimension, and Cohen--Macaulayness}

\begin{corollary}[The nonnormal locus is determined by the punctures]
\label{r2:norm-conductor-locus}
Assume $g\geq1$.  Let $N_i\subset Y_i$ be the reduced nonnormal
locus, and set
\[
 E=\left(\bigcup_i\operatorname{ev}_i^{-1}(N_i)\right)_{\red}
       \subset R_{\mathbf Y}.
\]
Then
\[
 \operatorname{Nonnormal}(U_{\mathbf Y})=q(E),\qquad
 q^{-1}\bigl(U_{\mathbf Y}^{\mathrm{normal}}\bigr)
       =R_{\mathbf Y}\setminus E.
\]
If $\mathfrak c_U$ is the normalization conductor of
$\C[U_{\mathbf Y}]\subset\C[U_{\mathbf Y}^{\nu}]$, then
\begin{equation}\label{r2:norm-radical-conductor}
 \sqrt{\mathfrak c_U}=I(E)^G_{\ad}
       \subset\C[R_{\mathbf Y}]^G_{\ad}=\C[U_{\mathbf Y}].
\end{equation}
Writing the codimension of an empty locus as $+\infty$, one has
\begin{equation}\label{r2:norm-codimension}
 \operatorname{codim}_{U_{\mathbf Y}}
       \operatorname{Nonnormal}(U_{\mathbf Y})
   =\min_i\operatorname{codim}_{Y_i}N_i.
\end{equation}
For a chosen connected component, the corresponding normality test uses
exactly those puncture germs realized by that component.
\end{corollary}

\begin{proof}
The set equality and the inverse-image equality are
Theorem~\ref{r2:norm-positive-genus}.  The set $E$ is invariant and
closed.  A good quotient sends it onto the closed subset with ideal
$I(E)^G_{\ad}$, and this ideal is radical.  The support of the finite module
$\C[U_{\mathbf Y}^{\nu}]/\C[U_{\mathbf Y}]$ is precisely the nonnormal
locus.  The annihilator of this quotient is the conductor, giving
\eqref{r2:norm-radical-conductor}.

Each $N_i$ is a union of finitely many conjugacy orbits.  For a realized
tuple of orbits $\mathcal O_i\subset Y_i$, the relative quotient
stratum has dimension
\[
 2g\dim G+\sum_i\dim\mathcal O_i-2\dim G_{\ad}.
\]
The codimension of this stratum in $U_{\mathbf Y}$ is therefore
$\sum_i(\dim Y_i-\dim\mathcal O_i)$.
All choices of these orbits are realized in positive genus by the
commutator argument.  A nonnormal stratum has at least one
$\mathcal O_i\subset N_i$, and its codimension is consequently at
least the right side of \eqref{r2:norm-codimension}.
Equality is attained by taking a largest-dimensional orbit in an $N_i$
of minimal codimension and the prescribed open classes at all other
punctures.  This proves the formula.
\end{proof}

\begin{corollary}[Serre's condition and the Cohen--Macaulay locus]
\label{r2:norm-cm}
Every generic relative character variety considered here satisfies
Serre's condition $R_1$.  The normal locus of such a variety equals its Cohen--Macaulay
locus.  In positive genus the following global conditions are
equivalent:
\begin{enumerate}
\item every prescribed puncture closure is normal;
\item $U_{\mathbf Y}$ is normal;
\item $U_{\mathbf Y}$ satisfies Serre's condition $S_2$;
\item $U_{\mathbf Y}$ is Cohen--Macaulay.
\end{enumerate}
When they hold, the original interior has the rational Gorenstein
symplectic singularities of its normalization.
\end{corollary}

\begin{proof}
The complement of the originally prescribed exact-class locus has
codimension at least two.  On that exact-class locus the quotient is
\'etale-locally a finite quotient of a smooth variety and is normal.
Every codimension-one point is therefore regular, proving $R_1$.
At a normal point, the finite normalization is an isomorphism on a
neighborhood; Theorem~\ref{new:thm-normalized-interior} makes that
neighborhood Cohen--Macaulay.  Conversely, if the local ring at a point
is Cohen--Macaulay, that local ring and its localizations satisfy $S_2$.
Together with $R_1$, Serre's criterion makes the local ring normal.
This proves the equality of loci.  Serre's criterion, followed by
Theorem~\ref{r2:norm-positive-genus}, gives the global equivalences.
\end{proof}

\subsection{Normality on the sphere}

For $g=0$, the kernel $K_r$ in
Theorem~\ref{r2:norm-kernel-criterion} can be nontrivial.  The following
condition makes the action faithful.

\begin{corollary}[Two strongly regular punctures on the sphere]
\label{r2:norm-sphere-two-tori}
Suppose $g=0$ and two of the realized punctures have centralizers in
$G$ which are connected maximal tori.  Then the $F$-action on $W_r$
is faithful.  In particular,
\[
 U_{\mathbf Y}\text{ is normal at }q(r)
 \quad\Longleftrightarrow\quad
 Y_i\text{ is normal at }M_i\text{ for every }i.
\]
If two prescribed puncture classes have connected maximal-torus
centralizers, this equivalence holds at every point of the spherical
character variety.  The nonnormal locus of the spherical character variety is then $q(E)$, with $E$ as
in Corollary~\ref{r2:norm-conductor-locus}; no assertion that every
puncture tuple is realized on the sphere is needed.
\end{corollary}

\begin{proof}
Renumber the punctures so that $Z_G(M_i)=T_i$ is a connected maximal
torus for $i=1,2$, and suppose $1\ne f\in K_r$.
Any lift of $f$ to $G$ centralizes every tuple element and therefore
belongs to both $T_i$.  Consequently
\[
 \ker(\operatorname{Ad}(M_i)-1)=\operatorname{Lie}T_i
       \subset\mathfrak g^f,
 \qquad i=1,2.
\]
The restriction of $\operatorname{Ad}(M_i)-1$ to $V_f$ is thus
injective for both punctures.  Substituting their two ranks
$\dim V_f$ in \eqref{r2:norm-moving-dimension} with $g=0$ shows that
all remaining ranks vanish.  Each remaining $M_j$ therefore acts
trivially on every simple ideal moved by $f$, by the Lie-ideal argument
in Lemma~\ref{r2:norm-moving-space}.

Choose one such simple adjoint factor.  In that factor the surface
relation reduces to $M_1M_2=1$, and all other punctures are the
identity.  The projected tuple is generated by one semisimple element
and hence lies in a maximal torus.  The original tuple consequently
lies in the inverse image of a proper parabolic of this simple factor,
contrary to irreducibility.  Thus $K_r=1$.
The normality assertion follows from
Theorem~\ref{r2:norm-kernel-criterion}.  For the global locus statement,
semisimple conjugacy classes are closed, so the two prescribed
punctures retain their torus centralizers at every point.  Apply the
pointwise assertion and take the invariant closed image, as in
Corollary~\ref{r2:norm-conductor-locus}.
\end{proof}

The connectedness condition concerns the full centralizer, not only
its identity component.  It cannot be omitted from this argument:
in $\operatorname{PGL}_2$ a regular semisimple involution has a
disconnected centralizer, and a spherical triple of the three
nonidentity elements of a Klein four subgroup has a nontrivial
effective stabilizer acting trivially on its zero-dimensional
deformation space.

The algebraic obstruction to dropping it is already visible in
\[
 A=\C\oplus(x,y)^2\subset B=\C[x,y],
 \qquad F=\{\pm1\},\qquad (x,y)\longmapsto(-x,-y).
\]
The ring $A$ is nonnormal only at the origin, while
\[
 A^F=B^F=\C[x^2,xy,y^2]
\]
is normal. The spherical example in Section~\ref{r2:sec-sphere-masking}
realizes the same loss of a normalization defect. We first construct
the puncture closure used in that example.

\subsection{A nonnormal conjugacy-class closure in \texorpdfstring{$\PGL_4$}{PGL4}}\label{paper:app-pgl4-raw}

The following explicit construction shows that the normalization in Theorem~\ref{paper:thm-normalized}
is needed already in type $A$, when the centralizer of a semisimple element is
disconnected.

\begin{proposition}[A raw type-A closure can already be nonnormal]
\label{new:gl-pgl4-nonnormal}
There is a conjugacy-class closure in $\PGL_4$ which is nonnormal and
fails $S_2$.  Consequently there are intrinsically generic relative
$\PGL_4$ character varieties, already for a twice-punctured torus,
that are nonnormal and fail $S_2$.  The normalization qualification in
Theorem~\ref{r2:thm-full-normalized-logcy} is essential even within type A.
\end{proposition}

\begin{proof}
Let $p:\SL_4\to\PGL_4$ be the central isogeny, and set
\[
 J=\begin{pmatrix}1&1\\0&1\end{pmatrix},\qquad
 c=\operatorname{diag}(J,-I_2),\qquad
 s=\operatorname{diag}(I_2,-I_2).
\]
Let $C$ be the $\PGL_4$ conjugacy class of $p(c)$ and put $Y=\overline C$.
The inverse image $p^{-1}(C)$ is the union of the $\SL_4$ classes of
$\zeta c$, $\zeta\in\mu_4$.  Because $p$ is finite etale, and hence
open, taking the inverse image of $Y$ gives the union of their closures.

The classes of $c$ and $-c$ in $\SL_4$ are distinct: at eigenvalue $1$
the former has a size-two Jordan block, while the latter has two
size-one blocks.  Their closures are distinct irreducible components
of $p^{-1}(Y)$ and have the same dimension.  Both contain $s$.
For $c$, conjugation by $\operatorname{diag}(t,t^{-1},1,1)$ and the
limit $t\to0$ give $s$.  For $-c$, first interchange the two
two-dimensional blocks; this permutation matrix has determinant one.
The same degeneration in the second block then gives $s$.
The other two components have eigenvalues $\{i,-i\}$ and do not meet
$s$.  Thus $p^{-1}(Y)$ has two distinct local irreducible components
through $s$, and is not normal there.  If $Y$ were normal at $p(s)$,
its etale inverse image would be normal at $s$, a contradiction.

There are exactly two conjugacy orbits in $Y$: $C$ of dimension $10$,
and the semisimple class of $p(s)$ of dimension $8$.  Indeed the sole
proper degeneration of the regular unipotent class in a $2$-by-$2$
block is the identity; the two choices of active block upstairs are
identified after projection.  The dimensions follow respectively
from centralizer dimensions $5$ and $7$ in the $15$-dimensional group
$\SL_4$.  It follows that $Y$ is smooth in codimension one and satisfies
$R_1$.  Since it is nonnormal at $p(s)$, Serre's criterion shows that
its local ring there fails $S_2$.

Theorem~\ref{new:thm-universal-realization}, applied to this conjugacy
closure and the point $p(s)$, realizes its germ up to smooth factors
at a point with trivial effective stabilizer in a generic relative
$\PGL_4$ character variety of a twice-punctured torus.  Normality and
$S_2$ ascend and descend under these smooth factors.  The resulting
relative germ therefore has the same two failures.
\end{proof}

\subsection{A generic sphere quotient that removes a puncture defect}
\label{r2:sec-sphere-masking}

The following example shows that the positive-genus hypothesis cannot
be removed from Theorem~\ref{r2:norm-positive-genus}.

\begin{proposition}[Normality can be restored on a three-punctured sphere]
\label{r2:prop-sphere-masking}
There is a compatible intrinsically generic tuple of three conjugacy
classes in $\PGL_4$ and a point of its relative character variety
whose first puncture lies in the nonnormal locus of its prescribed
class closure, although the character variety is normal at that
point.  More precisely, the completed local ring of the character
variety at the point is
\[
 \C[[u,v,w]]/(uv-w^2).
\]
The effective stabilizer has order two, its exact-stratum deformation
representation is zero, and it interchanges the two branches of the
puncture transversal.
\end{proposition}

\begin{proof}
All displayed matrices below represent projective classes.  Set
\[
 S=\begin{pmatrix}0&I_2\\ I_2&0\end{pmatrix},\qquad
 B=\begin{pmatrix}1&0\\0&2\end{pmatrix},\qquad
 C=\begin{pmatrix}-1&2\\-4&5\end{pmatrix},\qquad
 T=\operatorname{diag}(B,C),\qquad R=(ST)^{-1}.
\]
The eigenvalues of $S$ are $1,1,-1,-1$; those of $T$ are
$1,1,2,3$.  The matrix $C$ is conjugate to
$\operatorname{diag}(1,3)$.  Write
\[
 \lambda_\pm=\frac{9\pm\sqrt{57}}2.
\]
These are the distinct eigenvalues of $BC$, whose trace is $9$ and
determinant is $6$.  Since
$(ST)^2=\operatorname{diag}(CB,BC)$, the eigenvalues of $R$ are
\[
 \lambda_+^{-1/2},\ -\lambda_+^{-1/2},\
 \lambda_-^{-1/2},\ -\lambda_-^{-1/2}.
\]
In particular $R$ is regular semisimple.

Let $C_1$ be the class with semisimple part $S$, regular unipotent
part on the $+1$ eigenspace of $S$, and trivial unipotent part on
its $-1$ eigenspace.  This is the nonnormal $\PGL_4$ class closure
of Proposition~\ref{new:gl-pgl4-nonnormal}, after conjugating its
displayed semisimple representative.  Let $C_2$ and $C_3$ be the
closed semisimple classes of $T$ and $R$.  The tuple
\[
 r=(S,T,R),\qquad STR=1,
\]
belongs to the prescribed closure representation space.

We verify intrinsic genericity, including the character lattice of
$\PGL_4$.  Choose determinant-one lifts
\[
 S,\qquad T'=6^{-1/4}T,\qquad R'=6^{1/4}R
\]
to $\SL_4$, with the same choice of fourth root.  They still have
product one.  For the maximal parabolic corresponding to an
$a$-dimensional subspace, the character group of the $\PGL_4$
parabolic pulls back to the subgroup generated by
\[
 (\det_a)^{m_a},\qquad m_a=\frac4{\gcd(4,a)},
 \qquad a=1,2,3.
\]
Indeed the central scalar $\zeta I_4$ acts on $\det_a$ by
$\zeta^a$.  Thus it suffices to check that for every choice of
$a$ eigenvalues at each puncture, the product of their determinants
has nontrivial $m_a$th power.  Checking all such choices is stronger
than checking those compatible with the specified unipotent part.

For $a=1$ the scalar factors in $T'$ and $R'$ cancel.  The possible
products are $\pm e/\sqrt{\lambda_\pm}$, with
$e\in\{1,2,3\}$; their fourth powers are never one, since
$\lambda_\pm\ne\pm e^2$.  For $a=3$ the determinant-one condition
reduces the calculation to the complementary eigenvalue and hence
to $a=1$.

For $a=2$, the products of two eigenvalues of $T$, with multiplicities,
are
\[
 p\in\{1,2,3,6\}.
\]
The corresponding products for $R$ are
$-1/\lambda_+$, $-1/\lambda_-$ and $\pm1/\sqrt6$;
the products for $S$ are $\pm1$.  Their total squares therefore
cannot be one: the first two alternatives would require
$\lambda_\pm=\pm p$, whereas the last would require $p^2=6$.
All the excluded equalities are false.  Every proper parabolic is
contained in a maximal parabolic, so these checks establish
intrinsic genericity.  In particular the tuple $r$ is irreducible.

The effective stabilizer of the tuple $r$ is precisely the group of order two
generated by
\[
 f=\operatorname{diag}(I_2,-I_2).
\]
Indeed $fTf^{-1}=T$ and $fSf^{-1}=-S$, so $f$ centralizes the
projective tuple.  Conversely, if a matrix $P$ centralizes the
projective class of $T$, then $PTP^{-1}=cT$.  The unique eigenvalue
of $T$ of multiplicity two is $1$, forcing $c=1$.  If $P$ also
centralizes the projective class of $S$, then
$PSP^{-1}=\varepsilon S$ with $\varepsilon\in\{1,-1\}$.
When $\varepsilon=1$, irreducibility and Schur's lemma imply that
$P$ is scalar.  When $\varepsilon=-1$, the same argument applies
to $f^{-1}P$.  This proves the stabilizer assertion.

The realized exact conjugacy classes at $r$ have dimensions
$8$, $10$ and $12$: their centralizer dimensions in $\PGL_4$
are respectively $7$, $5$ and $3$.  Since $\dim\PGL_4=15$,
the exact-stratum slice has dimension
\[
 \dim W=-2\cdot15+8+10+12=0.
\]
Consequently the entire effective stabilizer is the kernel
$K_r=F_r\simeq\boldsymbol\mu_2$ of its action on $W$.

We describe the puncture transversal explicitly.  In an eigenbasis
of $S$, its two eigenspaces have dimension two.  Pulling the
closure $\overline C_1$ back along the finite etale isogeny
$\SL_4\to\PGL_4$ gives, near the chosen lift of $S$, two
branches.  The first allows a nilpotent element in the first
$\mathfrak{sl}_2$ block and sets the second block to zero; the
second reverses these roles.  The fixed semisimple eigenvalues
exclude central directions.  Their transverse reduced union is
therefore
\[
 (\mathcal N_{\mathfrak{sl}_2}\times\{0\})
       \ \cup\
 (\{0\}\times\mathcal N_{\mathfrak{sl}_2}),
\]
whose two branches meet scheme-theoretically at the reduced origin.
Writing
\[
 B_0=\C[[u,v,w]]/(uv-w^2),
\]
its completed local ring is
\[
 A=B_0\mathbin{\times_{\C}}B_0.
\]
In particular this puncture germ is nonnormal.

The element $f$ interchanges the $+1$ and $-1$ eigenspaces of
$S$; in the local $\SL_4$ chart its action is conjugation followed
by multiplication by $-1$ to restore the chosen lift of $S$.
This element consequently interchanges the two branches.  Identifying the
two nilpotent cones using this map gives
$f(a,b)=(b,a)$.  The other two prescribed closures are smooth, so
their transverse rings are just $\C$.  The equivariant local
model of Proposition~\ref{r2:norm-equivariant-model} now gives,
for the quotient point $p=q(r)$,
\[
 \widehat{\mathcal O}_{U,p}
     \simeq A^{\boldsymbol\mu_2}
     \simeq B_0,
 \qquad a\longmapsto(a,a).
\]
The $A_1$ surface ring $B_0$ is normal: it is a hypersurface with
its only singularity at the codimension-two origin, hence is
$S_2$ and $R_1$.  Excellence then implies that $U$ is normal at
$p$.  This proves every assertion.
\end{proof}

\subsection{Consequences for compactification}\label{subsec:normality-compactification}
Combining the positive compactification theorem with the normality criterion
now determines when the given reduced interior can be retained.

\begin{corollary}[The original interior and its normality]
\label{r2:cor-original-interior}
For a nonempty relative variety $U_0$ with compatible intrinsically generic
conjugacy-class data as in Theorem~\ref{r2:thm-full-normalized-logcy}, a normal
projective log canonical log Calabi--Yau compactification preserving
the given reduced interior exists if and only if $U_0$ is normal.
In positive genus this is equivalent to normality of every prescribed
conjugacy-class closure.  It is also equivalent to Serre's condition
$S_2$, and to Cohen--Macaulayness of $U_0$.

In genus zero, normality is determined at every point by
Theorem~\ref{r2:norm-kernel-criterion}: the completed puncture-transversal
ring $A$ must have normal invariant ring $A^{K_r}$.  With two prescribed
punctures having connected maximal-torus centralizers, this simplifies
pointwise to normality of all the realized puncture germs.
\end{corollary}

\begin{proof}
An open subscheme of a normal variety is normal.  The converse is
Theorem~\ref{r2:thm-full-normalized-logcy}, since normalization then
does not change the interior; take the disjoint union of the
component compactifications when there is more than one component.
The remaining assertions are
Theorem~\ref{r2:norm-positive-genus} and
Corollaries~\ref{r2:norm-cm} and~\ref{r2:norm-sphere-two-tori}.
Normal components alone do not make an intersecting reduced union
normal, so this conclusion concerns the given reduced interior as
a scheme.
\end{proof}

\section{Examples}\label{sec:examples}

The examples exhibit singular interiors, the need for reduced spectral
fibers, explicit local and global $\SL_3$ constructions, and obstructions
to preserving the original interior.

\subsection{The exceptional \texorpdfstring{$\SL_2$}{SL2} cubic surfaces}
\begin{example}[Direct verification for $\tau=\pm2$]
\label{ex:exceptional-cubics}
Recall the compactifications from \hyperref[subsec:intro-example]{the introductory model calculation}:
\[
 X_\tau=\{t(x^2+y^2+z^2)-xyz=(\tau+2)t^3\}\subset\PP^3_{[x:y:z:t]},\qquad
 D_\tau=X_\tau\cap\{t=0\},\qquad U_\tau=X_\tau\setminus D_\tau.
\]
We check the conclusions directly at $\tau=\pm2$. The affine singularity
equations are $2x=yz$, $2y=xz$, $2z=xy$. For $\tau=-2$ the only singular
point is the origin; for $\tau=2$ there are four:
\[
 (x,y,z)=(2\varepsilon_1,2\varepsilon_2,2\varepsilon_3),\qquad
 \varepsilon_i\in\{1,-1\},\quad \varepsilon_1\varepsilon_2\varepsilon_3=1.
\]
The Hessian determinants are respectively $8$ and $-32$, so all these
singularities are ordinary double points, analytically the quadric cone
$\{uw=v^2\}$. This cone is the quotient of $\C^2$ by
$(s,r)\mapsto(-s,-r)$, and is canonical and $\Q$-factorial
(Fact~\ref{fact:symplectic-quotient}).

The projective surfaces are smooth along $D_\tau$: at $t=0$, vanishing of
all four derivatives would require
$xy=xz=yz=x^2+y^2+z^2=0$, which is impossible in projective space.
Near the vertex $[1:0:0:0]$, the surface equation expresses $t$ as $yz$
times a unit, so $D_\tau$ has local equation $yz=0$; the other vertices are
the same. Thus the boundary remains a simple normal crossings triangle.
The isolated singularities also imply that the cubics are irreducible
(distinct projective components would meet along a singular curve) and
normal. As hypersurfaces, $X_\tau$ are Cohen--Macaulay and Gorenstein.
Adjunction still gives $K_{X_\tau}+D_\tau\sim0$. Since the canonical nodes
lie away from the boundary, both pairs are log canonical, indeed
$\Q$-factorial dlt. The identity morphism therefore supplies the required
dlt modification without changing the interior. The dual complex of
$(X_\tau,D_\tau)$ is still a circle.

Finally, the scheme-theoretic peripheral fiber has coordinate ring
$\C[x,y,z]/(x^2+y^2+z^2-xyz-\tau-2)$, an integral domain, so both fibers
are connected and reduced. These computations verify the conclusions
independently of the general theorems: the intrinsically generic case $U_{-2}$
illustrates Theorem~\ref{thm:spectral-main}, and the nongeneric case $U_2$
illustrates Corollary~\ref{cor:spectral-type-a} (see Example~\ref{ex:generic-sl2} for genericity).
\end{example}

\subsection{A nilpotent-cone singularity}
\begin{example}[Repeated eigenvalues and a nilpotent-cone germ]
\label{ex:clock-shift-spectral}
Let $G=\SL_n$ with $n\ge2$, let $(g,k)=(1,1)$, let $\zeta$ be a primitive $n$-th
root of unity, and prescribe the characteristic polynomial $(t-\zeta^{-1})^n$
at the puncture. The reduced fiber is $S=\zeta^{-1}\mathcal U(\SL_n)$, the
scalar $\zeta^{-1}$ times the unipotent variety, which is singular at
$\zeta^{-1}I$. The data are intrinsically generic: on an invariant subspace of
dimension $d$ with $0<d<n$, the determinant of the puncture monodromy is
$\zeta^{-d}\ne1$.

On a basis $e_0,\dots,e_{n-1}$ of $\C^n$, let $A$ be the clock matrix
$Ae_j=\zeta^je_j$ and $B$ the cyclic shift $Be_j=e_{j+1}$ (indices modulo $n$),
both rescaled to have determinant one. Then $AB=\zeta BA$, so $[A,B]=\zeta I$ and
$M=[A,B]^{-1}=\zeta^{-1}I\in S$. The eigenlines of $A$ are distinct and $B$
permutes them transitively, so the pair $(A,B)$ is irreducible and its
centralizer is the center of $\SL_n$.

We claim that the germ of the spectral character variety at the class $u$ of
$(A,B)$ is the germ of $\mathcal N(\mathfrak{sl}_n)$ at $0$. The spectral
representation variety is $\{(A',B'):[A',B']\in\zeta\,\mathcal U(\SL_n)\}$. By
Lemma~\ref{lem:differential}, the commutator map $c:\SL_n^2\to\SL_n$ has
surjective differential at $(A,B)$, so the kernel of the differential has
dimension $\dim\SL_n$. Since $c$ is equivariant for conjugation and
$c(A,B)=\zeta I$ is central, this kernel contains the tangent space to the
conjugation orbit, which also has dimension $\dim\SL_n$ because the stabilizer
in $\PGL_n$ is trivial. Take a Luna slice in the smooth ambient variety
$\SL_n^2$ through $(A,B)$ (Fact~\ref{fact:git}(iii)). The slice is transverse
to the orbit, so the restriction of $c$ to the slice is \'etale at $(A,B)$.
The inverse image of $\zeta\,\mathcal U(\SL_n)$ in this slice is a local
model for the spectral character variety at $u$. The germ of the spectral
character variety at $u$ is therefore the
germ of $\zeta\,\mathcal U(\SL_n)$ at $\zeta I$, which the logarithm identifies with
$(\mathcal N(\mathfrak{sl}_n),0)$. For $n=2$ we recover the node of the cubic
surface $U_{-2}=\{x^2+y^2+z^2-xyz=0\}$ at the origin.
\end{example}

\subsection{A nonreduced peripheral fiber}
\begin{example}[A nonreduced spectral scheme for $\PGL_2$]
\label{ex:pgl2-spectral-thickening}
Let $G=\PGL_2$, $(g,k)=(1,1)$, and choose the coordinate $q$ on
$G/\!/G=\A^1$ whose pullback to $\SL_2$ is $(\tr M)^2$.
The reduced fiber $\{q=0\}$ is the regular semisimple class of the
nontrivial involutions, the image of $\diag(i,-i)\in\SL_2$. The datum is
intrinsically generic: an involution in a Borel subgroup has nontrivial image
in the torus quotient of that Borel subgroup, whereas the commutator subgroup
has trivial image.

Nevertheless the scheme-theoretic peripheral fiber is not reduced.  On the free $\SL_2$ character variety write
$x=\tr A$, $y=\tr B$, $z=\tr AB$, so that
\[
 \C[\SL_2^2]^{\SL_2}=\C[x,y,z],\qquad
 \kappa=\tr[A,B]=x^2+y^2+z^2-xyz-2
\]
\cite[\S2.1, (2.1.4)]{Goldman2003}.  Changing the signs of the two lifts
acts by $(x,y,z)\mapsto(\epsilon x,\eta y,\epsilon\eta z)$.
Consequently the free $\PGL_2$ character ring is
\[
 B=\C[X,Y,Z,W]/(W^2-XYZ),
 \qquad (X,Y,Z,W)=(x^2,y^2,z^2,xyz),
\]
and $\kappa=X+Y+Z-W-2$ belongs to $B$.  The puncture inverse has
the same trace, so its prescribed invariant is $q=\kappa^2$.
The scheme fiber at $q=0$ is therefore
\begin{equation}\label{eq:pgl2-double-fiber}
 \Spec\frac{\C[X,Y,Z,W]}
 {(W^2-XYZ,\ (X+Y+Z-W-2)^2)}.
\end{equation}
The class of $\kappa$ in this ring is a nonzero nilpotent element. Indeed $B$ is
a domain, and $\kappa$ is neither zero nor a unit, since $\kappa$ vanishes at the
point $(X,Y,Z,W)=(0,0,2,0)$ of $\Spec B$; hence $\kappa\notin(\kappa^2)$. The same
point shows that the fiber is nonempty. A nonreduced scheme is not an open
subscheme of a normal variety, so the conclusions of Theorem~\ref{thm:spectral-main} fail for this
scheme-theoretic fiber, although the datum is generic and regular semisimple.
This is consistent with the last assertion of Theorem~\ref{thm:spectral-main}, since the derived
group $\PGL_2$ is not simply connected. The reduced fiber $\Spec B/(\kappa)$ is
the relative character variety, to which Corollary~\ref{cor:regular-semisimple}
applies.
\end{example}

\subsection{The local \texorpdfstring{$\SL_3$}{SL3} algebra}\label{subsec:sl3-local}
We make Example~\ref{ex:local-sl2-sl3}(2) explicit. On the open set where the
four corner pairings $L_1=\bar C=\phi_2(v_1)$, $L_2=C=\phi_1(v_2)$,
$M_1=B=\phi_3(v_1)$ and $M_2=\bar B=\phi_1(v_3)$ are invertible, the ring
$R_{\SL_3}$ embeds into $\C[x,y,z,L_1,L_2,M_1,M_2]$ via
\[
 A=\phi_2(v_3)=zL_1M_2,\quad\bar A=\phi_3(v_2)=(xy-z)L_2M_1,\quad
 D=xL_1M_1,\quad E=yL_2M_2 ,
\]
where $(x,y,z)$ are the coordinates of the upper unitriangular matrix
$\left(\begin{smallmatrix}1&x&z\\0&1&y\\0&0&1\end{smallmatrix}\right)$
(Proposition~\ref{prop:poly-model}). The substitution $(x,y,z)=(a+c,b,ab)$
corresponds to the reduced word $\mathbf i=(1,2,1)$, and taking the
lexicographically largest exponent in $(a,b,c;L_1,L_2;M_1,M_2)$ gives the
valuation $v_{\mathbf i}$. The values of the generators are:
\[
\begin{array}{@{}c|cccccccc@{}}
 &A&B&C&\bar A&\bar B&\bar C&D&E\\\hline
 \boldsymbol t&(1,1,0)&0&0&(0,1,1)&0&0&(1,0,0)&(0,1,0)\\
 \lambda&(1,0)&0&(0,1)&(0,1)&0&(1,0)&(1,0)&(0,1)\\
 \eta&(0,1)&(1,0)&0&(1,0)&(0,1)&0&(1,0)&(0,1)
\end{array}
\]
Since $DE=(ab+bc)L_1L_2M_1M_2$, $ABC=ab\,L_1L_2M_1M_2$ and
$\bar A\bar B\bar C=bc\,L_1L_2M_1M_2$, we get
$v(DE)=v(ABC)>v(\bar A\bar B\bar C)$, and
\[
 \gr R_{\SL_3}=\C[A,\dots,E]/(DE-ABC).
\]
Writing a point of $\Z^7$ as $(t_1,t_2,t_3;l_1,l_2;m_1,m_2)$ and eliminating the
coefficients of the eight generators, one finds that $Q_{3,\mathbf i}$ is cut
out by nine inequalities:
\[
\begin{array}{lll}
 t_1\ge0,&t_3\ge0,&t_2-t_3\ge0,\\
 l_1-t_1\ge0,&l_2-t_3\ge0,&l_2+t_1-t_2\ge0,\\
 m_1-t_3\ge0,&m_2-t_2+t_3\ge0,&m_1-t_1+t_2-2t_3\ge0 .
\end{array}
\]
Each of the nine linear forms takes the value $1$ at the sum of the values of the
six corner pairings, $(1,2,1;2,2;2,2)$. By Lemma~\ref{lem:height-one}, this sum
is the canonical translator, as asserted by Theorem~\ref{thm:corner}.

\subsection{The once-punctured torus for \texorpdfstring{$\SL_3$}{SL3}}\label{subsec:sl3-torus}
Let $(g,k)=(1,1)$ and $M=[A,B]^{-1}$. For $a^3=3/2$, the matrices
\[
 A=\diag(a,a/3,2a),\qquad B=\begin{pmatrix}0&0&1\\1&0&0\\0&1&0\end{pmatrix}
\]
lie in $\SL_3$ and satisfy $M=\diag(2,3,1/6)$. No nonempty proper subset of
$\{2,3,1/6\}$ has product $1$, so the class of $M$ is regular semisimple and
generic, and the relative character variety $U$ is nonempty, connected and of
dimension $6$. The two monodromy equations are $\tr M=31/6$ and
$\tr M^{-1}=41/6$. On the theta graph, $Q_\Gamma$ is cut out from two copies of
$Q_{3,121}$ by six matching equations, so $\rank\gp(Q_\Gamma)=2\cdot7-6=8$.
The boundary walk has six corner occurrences, and the leading exponents $p_1,p_2$
of the two equations satisfy
$p_1+p_2=\Omega_\Gamma=\bigl((1,2,1;2,2;2,2),(1,2,1;2,2;2,2)\bigr)$. Hence
$\gr\C[U]=\C[Q_\Gamma]/(x^{p_1},x^{p_2})$ is a reduced Gorenstein toric-face
ring, its cone complex maps onto a complete fan in $\R^6$, and the Rees
construction gives a log canonical compactification $(X,D)$ of $U$ with
$K_X+D\sim0$.

\subsection{Two obstructions for prescribed conjugacy closures}
\label{subsec:closure-obstructions}

The following examples separate the need to normalize an interior from
the stronger obstruction to preserving it in a $\Q$-factorial model.
Both use intrinsically generic data on a genus-one surface with two
punctures.  The complete constructions and proofs are in
Appendix~\ref{new:sec-counterexample}.

For $G$ of type $G_2$, prescribe the unipotent closure
$\exp(\overline{\mathcal O_8})$, where $\mathcal O_8$ is the
eight-dimensional nilpotent orbit of type $\widetilde A_1$, and the
regular semisimple class of an element $t$ with
$\alpha(t)=1/12$ and $\beta(t)=24$, for $\alpha$ short and $\beta$
long.  Theorem~\ref{new:thm-g2-counterexample} gives a point with trivial
effective stabilizer whose analytic germ is
\[
 (\overline{\mathcal O_8},0)\times(\C^{12},0).
\]
This twenty-dimensional relative variety is nonnormal and fails
$S_2$ at that point.  The component through the point therefore
cannot remain an open subvariety of either a normal log-canonical
or a semi-log-canonical compactification: normality and $S_2$,
respectively, are necessary already on the interior.

For $G=\SL_3$, instead prescribe the unipotent closure of Jordan
type $(2,1)$ and the regular semisimple class of
$\diag(2,3,1/6)$.  Theorem~\ref{new:sl3-main} gives a normal
ten-dimensional relative variety which is not $\Q$-factorial.
The algebraic obstruction is a projective birational morphism from a
smooth variety with no exceptional divisors and with a fiber $\PP^2$:
$\Q$-factoriality of the
target would force a positive multiple of a relatively ample divisor
to descend, contradicting its positive degree on a line in that
fiber.  Hence no $\Q$-factorial compactification, and in
particular no $\Q$-factorial dlt model, can preserve the whole
interior.  Normalization changes nothing in this example; nevertheless
Theorem~\ref{r2:thm-full-normalized-logcy} supplies a normal
log-canonical log Calabi--Yau compactification.

\appendix
\addtocontents{toc}{\protect\setcounter{tocdepth}{2}}
\makeatletter
\def\toclevel@section{0}
\def\toclevel@subsection{1}
\def\toclevel@subsubsection{2}
\makeatother
\addtocontents{toc}{\protect\addvspace{6pt}}

\section{Inputs to the spectral construction}\label{app:spectral-inputs}
The arguments in this appendix support the spectral construction in
Part~\ref{part:spectral} and the spectral interior analysis in
Section~\ref{sec:interior}.
Section~\ref{app:string} proves the local string-valuation theorem;
Section~\ref{app:scalar} converts valuations into positive integer
filtrations; Section~\ref{app:geometric-facts} supplies the additional
geometric facts used in the canonical-module and singularity arguments.

\subsection{The local string valuation}\label{app:string}

In this subsection we prove Theorem~\ref{thm:local-cone}. The proof has three
parts. First, a tensor invariant is determined by a single polynomial function
on $U$ (Proposition~\ref{prop:poly-model}); this gives a concrete valuation, by
leading exponents in coordinates attached to the reduced word $\mathbf i$.
Second, Lusztig's canonical basis provides a basis of these polynomial functions
with distinct leading exponents, the string parameters, which are computed by a
piecewise-linear procedure; this identifies the value semigroup with the
lattice points of a cone (Theorem~\ref{thm:app-local-cone}). Third, the point
$q_\rho(\mathbf i)$ is the string parameter of the basis element with all
Lusztig coordinates equal to $1$, which is interior for a simple topological
reason (Proposition~\ref{prop:rho-interior}).

Let $H$ be simply connected semisimple, with $T\subset B=TU$ as before.
Index the simple roots by $I=\{1,\dots,r\}$, choose root vectors
$e_i\in\Lie U$, and put $x_i(s)=\exp(se_i)$. For $\lambda\in X^+(T)$
write $\lambda_i=\langle\lambda,\alpha_i^\vee\rangle$, let $v_\lambda\in
V(\lambda)$ be a highest weight vector, and let $\ell_\lambda\in V(\lambda)^\vee$
be the lowest weight vector with $\ell_\lambda(v_\lambda)=1$. Let
$\mathcal U=U(\Lie U)$ be the enveloping algebra, graded by
$\deg e_i=\alpha_i$. We use the following standard presentation
\cite[\S21.4]{HumphreysLie1972}:
\begin{equation}\label{eq:presentation}
 V(\lambda)^\vee=\mathcal U\cdot\ell_\lambda,\qquad
 \operatorname{Ann}_{\mathcal U}(\ell_\lambda)=\textstyle\sum_i\mathcal U\,e_i^{\lambda_i+1}.
\end{equation}

\subsubsection{What is needed, and why a special basis}\label{subsec:why-basis}
We shall see (Proposition~\ref{prop:poly-model}) that each graded piece of the
three-flag algebra $R_H$ is identified with a space of polynomial functions on
$U$ of the form
\[
 S_{\lambda,\eta}=\{F\in\C[U]:\mathcal L_j^{\lambda_j+1}F=0\text{ and }
 \mathcal R_j^{\eta_j+1}F=0\text{ for all }j\},
\]
where $\mathcal L_j$ and $\mathcal R_j$ are the derivatives of $F$ along left and
right multiplication by the root subgroup $x_j(s)$. The valuation $v_{\mathbf i}$
takes the lexicographically largest exponent of $F$ in coordinates
$(a_1,\dots,a_N)$ on $U$ attached to $\mathbf i$. To describe the set of values
of $v_{\mathbf i}$ on $S_{\lambda,\eta}$, the natural strategy is to find a basis
of $S_{\lambda,\eta}$ whose elements have pairwise distinct leading exponents:
then the leading exponent of any nonzero linear combination is the largest
leading exponent of a basis element occurring in it, so the value set is exactly
the set of leading exponents of the basis.

We need such bases for all $(\lambda,\eta)$ at once, compatibly. The cleanest way
is a single basis $\{F_b\}$ of $\C[U]$ such that every kernel
$\ker\mathcal L_j^{n}$ and every kernel $\ker\mathcal R_j^n$ is spanned by a
subset of the basis (then every intersection $S_{\lambda,\eta}$ is spanned by a
subset too), and whose leading exponents are distinct and computable. A basis
compatible with one of these filtrations is easy to find; a basis compatible
with all of them simultaneously is not, as the following example shows.

\begin{example}\label{ex:sl3-good-basis}
Let $H=\SL_3$ and write $n=\left(\begin{smallmatrix}1&x&z\\0&1&y\\0&0&1
\end{smallmatrix}\right)\in U$. Left and right multiplication by
$x_1(s)=I+sE_{12}$ and $x_2(s)=I+sE_{23}$ give
\[
 \mathcal L_1=\partial_x+y\,\partial_z,\quad \mathcal L_2=\partial_y,\qquad
 \mathcal R_1=\partial_x,\quad \mathcal R_2=\partial_y+x\,\partial_z .
\]
In the degree $\alpha_1+\alpha_2$ (for the conjugation action of $T$), the
polynomial functions are spanned by $z$ and $xy$, and
\[
 \ker\mathcal L_1=\C(xy-z),\quad\ker\mathcal L_2=\C z,\quad
 \ker\mathcal R_1=\C z,\quad\ker\mathcal R_2=\C(xy-z).
\]
So the monomial basis $\{z,xy\}$ fails to be compatible with
$\ker\mathcal L_1$ and $\ker\mathcal R_2$, whereas $\{z,xy-z\}$ is compatible
with all four kernels. For the reduced word
$\mathbf i=(1,2,1)$, the coordinates on $U$ come from
$x_1(a)x_2(b)x_1(c)=\left(\begin{smallmatrix}1&a+c&ab\\0&1&b\\0&0&1
\end{smallmatrix}\right)$, i.e.\ $(x,y,z)=(a+c,b,ab)$. Then $z=ab$ and
$xy-z=bc$ have the distinct leading exponents $(1,1,0)$ and $(0,1,1)$, whereas
$xy=ab+bc$ has the same leading exponent as $z$. The basis $\{z,xy-z\}$ is the
part of degree $\alpha_1+\alpha_2$ of the dual canonical basis described below.
Similarly, the representation $\C^3=V(\varpi_1)$ corresponds to
$\{F:\mathcal R_1^2F=\mathcal R_2F=0\}=\Span\{1,\,x,\,xy-z\}$, spanned by three
elements of the same basis, with leading exponents $0$, $(1,0,0)$ and $(0,1,1)$.
\end{example}

For a general group, a basis with these simultaneous compatibility
properties is supplied by the canonical-basis theory of Lusztig
\cite[Thms.~14.3.2, 14.4.3]{Lusztig1993} and, independently,
Kashiwara \cite[Thms.~6--7]{Kashiwara1991}. We only use the
combinatorial consequences listed in Section~\ref{subsec:lusztig}.

\subsubsection{Input from the theory of canonical bases}\label{subsec:lusztig}
The pairing $\langle D,F\rangle=(\mathcal R_DF)(1)$, where $\mathcal R_D$ is the
differential operator on $U$ obtained by differentiating right multiplication
and defined by $D\in\mathcal U$,
identifies $\C[U]$ with the graded dual of $\mathcal U$. Let $\mathcal L_i$ and
$\mathcal R_i$ be the derivations of $\C[U]$ given by
$\mathcal L_iF(n)=\frac{d}{ds}F(x_i(s)n)|_{s=0}$ and
$\mathcal R_iF(n)=\frac{d}{ds}F(nx_i(s))|_{s=0}$. Then
$\langle D,\mathcal L_i^nF\rangle=\langle e_i^nD,F\rangle$ and
$\langle D,\mathcal R_i^nF\rangle=\langle De_i^n,F\rangle$.

Let $\mathcal B$ be Lusztig's canonical basis of $\mathcal U$ (specialized at
$q=1$), and let $\{F_b\}_{b\in\mathcal B}\subset\C[U]$ be the dual basis. For
$b\in\mathcal B$ and $j\in I$, put
$\epsilon_j(b)=\max\{n:b\in e_j^n\mathcal U\}$ and
$\epsilon^R_j(b)=\max\{n:b\in\mathcal Ue_j^n\}$. We use the following
properties \cite[Thms.~14.3.2, 14.4.3, 14.4.9, 14.4.11; \S17.3; Thm.~18.3.8;
\S\S20.1, 42.1]{Lusztig1993}.

\begin{enumerate}[label=(B\arabic*)]
\item The ideals $e_j^n\mathcal U$ and $\mathcal Ue_j^n$ are spanned by subsets
of $\mathcal B$. For $b\in\mathcal B$ with $n=\epsilon_j(b)$, there is a unique
$\tilde e_j^{\max}b\in\mathcal B$ with $\epsilon_j(\tilde e_j^{\max}b)=0$ and
$\frac{e_j^n}{n!}\,\tilde e_j^{\max}b\equiv b\pmod{e_j^{n+1}\mathcal U}$.
Conversely, for $b_0$ with $\epsilon_j(b_0)=0$ and $n\ge0$ there is a unique
$\tilde f_j^nb_0$ with $\epsilon_j=n$ and $\tilde e_j^{\max}\tilde f_j^nb_0=b_0$.
The element $1$ is the only $b$ with $\epsilon_j(b)=0$ for all $j$.
\item For $\eta\in X^+(T)$, the elements $b$ with $\epsilon^R_j(b)\le\eta_j$
for all $j$ map to a basis of
\[
 \mathcal U/\sum_j\mathcal Ue_j^{\eta_j+1}\simeq V(\eta)^\vee,
\]
and this set of elements is stable under $\tilde e_j^{\max}$. If $b_0$ is in this set,
$\epsilon_j(b_0)=0$, and $b_0$ has degree $\gamma$, then $\tilde f_j^nb_0$ is in
this set if and only if $n\le\langle\eta-\gamma,\alpha_j^\vee\rangle$.
\item (Lusztig's parametrization.) For each reduced word $\mathbf h$ of $w_0$
there is a bijection $b_{\mathbf h}:\N^N\to\mathcal B$ such that, if $h_1=j$,
then $\epsilon_j(b_{\mathbf h}(c))=c_1$ and
$\tilde e_j^{\max}b_{\mathbf h}(c)=b_{\mathbf h}(0,c_2,\dots,c_N)$. For two
reduced words, the transition map $b_{\mathbf h'}^{-1}\circ b_{\mathbf h}$ is
the restriction to $\N^N$ of a piecewise-linear bijection of
$\R^N_{\ge0}$, given by finitely many integral linear formulas on the cones of
a rational polyhedral subdivision. In simply laced type, if $\mathbf h$ and
$\mathbf h'$ differ by a braid move $iji\to jij$ in positions $k,k+1,k+2$, the
transition changes only these three coordinates, by
$(a,b,c)\mapsto(b+c-\min(a,c),\,\min(a,c),\,a+b-\min(a,c))$; for a move
$ij\to ji$ the transition swaps two coordinates. In non-simply-laced type we use the admissible orbit Cartan datum described
below. The basis labels of this datum identify with the fixed labels in simply laced type,
and the maximal label operator of an orbit is the product of the commuting maximal label
operators in the corresponding orbit.
\item The element $b_{\mathbf 1}=b_{\mathbf h}(1,\dots,1)$ does not depend on
$\mathbf h$, has degree $2\rho$, and satisfies $\epsilon_j(b_{\mathbf 1})=1$
for every $j$.
\end{enumerate}

Let us explain these statements. (B1) says that the canonical basis is
compatible with left divisibility by powers of $e_j$, and that removing the
largest possible power of $e_j$ from a basis element gives, up to lower order
terms, another basis element; the map $b\mapsto\tilde e_j^{\max}b$ on labels is
the combinatorial shadow of ``dividing by $e_j$ as often as possible''. In
Example~\ref{ex:sl3-good-basis}, the canonical basis in degree
$\alpha_1+\alpha_2$ is $\{e_1e_2,e_2e_1\}$, whose dual basis under the pairing
$\langle D,F\rangle=(\mathcal R_DF)(1)$ is $\{z,xy-z\}$. One has
$e_1e_2\in e_1\mathcal U$ with $\tilde e_1^{\max}(e_1e_2)=e_2$, whereas
$e_2e_1\notin e_1\mathcal U$ (in degree $\alpha_1+\alpha_2$ the space
$e_1\mathcal U$ is the line $\C e_1e_2$), so $\epsilon_1(e_2e_1)=0$. (B2) says
that the canonical basis is also compatible with the quotient maps onto the
representations $V(\eta)^\vee$, and describes how long the strings of the
operators $\tilde f_j$ inside $V(\eta)^\vee$ are, exactly as for $\SL_2$-strings
in a representation. (B3) provides coordinates in which all these operations are
piecewise linear: Lusztig labels the canonical basis by the exponents of
Poincar\'e--Birkhoff--Witt monomials in root vectors (one coordinate for each
positive root, ordered by the reduced word), and the transition between the
labelings for two reduced words is given by formulas involving minima.

For $H=\SL_2$ all of this is elementary: $\mathcal U=\C[e]$, the canonical basis
is $\{e^n/n!\}$, $\epsilon(e^n/n!)=n$, the dual basis of $\C[U]=\C[x]$ is
$\{x^n\}$, Lusztig's parametrization is $n\mapsto e^n/n!$, and
$b_{\mathbf 1}=e$ has degree $\alpha=2\rho$. For $\SL_3$, the canonical basis of
$\mathcal U$ is not a basis of monomials in $e_1,e_2$, and the coordinate changes
of (B3) are genuinely piecewise linear (Example~\ref{ex:sl3-string}).

In simply laced type, assertion (B4) follows from (B3): every elementary
braid formula fixes the all-ones tuple, and a word beginning with $j$ gives
$\epsilon_j(b_{\mathbf1})=1$.  The degree of $b_{\mathbf1}$ is the sum of the positive roots,
namely $2\rho$.

Here are the folding details needed for (B3) and (B4) in the remaining
types.  Choose a simply laced diagram with an automorphism $\sigma$ such
that no two vertices in one orbit are adjacent.  The admissible orbit
Cartan data from $A_{2n-1}$, $D_{n+1}$, $E_6$ and $D_4$ give respectively
$B_n$, $C_n$, $F_4$ and $G_2$ (with the orbit Cartan convention).
If $J$ is an orbit, its simple root maps to
$\sum_{j\in J}\widehat\alpha_j$, and its reflection maps to
$\prod_{j\in J}\widehat s_j$; these reflections commute.  This is the
orbit-root normalization of \cite[\S\S14.1.1, 14.1.6, Thm.~14.4.9]{Lusztig1993},
not an identification of the enveloping algebras.  The compatibility of
this folding with piecewise-linear parameters is made explicit in
\cite[\S\S1.5, 1.8, 1.12]{Lusztig2011}.

Expand a reduced word downstairs by replacing each orbit letter with its
commuting letters upstairs. For example, for $A_3$ with $\sigma$ exchanging the
vertices $1$ and $3$, the orbits are $\{1,3\}$ and $\{2\}$, the downstairs
diagram has rank two with four positive roots, and the downstairs reduced word
$(\{1,3\},\{2\},\{1,3\},\{2\})$ expands to $(1,3,2,1,3,2)$, a reduced word of
length six, the number of positive roots of $A_3$.  The result is a reduced word of the upstairs
longest element: the inversion roots in each block are permuted by
$\sigma$, have the same sign, and sum to the positive downstairs inversion
root, so the expansion is reduced. The product of the expanded word sends every simple-root
orbit sum, hence every simple root, to a negative root, so it is the longest
element. By the cited folding theorem, a fixed basis label has
PBW coordinates constant on each such block, and the common block entries
are its downstairs coordinates.  Thus a change of downstairs coordinates
is obtained by duplicating block entries, applying the simply laced
piecewise-linear changes, and retaining one entry from each output block.
The same procedure computes the maximal label operators: in a word
beginning with the orbit, set every entry of the initial block to zero.
These maps and their inverses have finitely many integral linear pieces;
their identities extend from integral points to real points by homogeneity
and continuity.

The downstairs all-ones tuple expands to the upstairs all-ones tuple.
Consequently the label $b_{\mathbf1}$ is independent of the word and has
$\epsilon_j(b_{\mathbf1})=1$ for every downstairs $j$.  The degree of $b_{\mathbf1}$ is
$\sum_{\beta>0}\beta=2\rho$ downstairs: under the root map above,
the roots in each expanded inversion block sum to the corresponding
downstairs inversion root.  Products of irreducible root systems are
handled factor by factor.  This also justifies the folding step in the
termination proof below: maximal raising along an expanded longest word
ends at the label $1$, so maximal raising along the downstairs longest word
does the same. The corresponding right-label statements follow by applying
the canonical-basis-preserving antiautomorphism which exchanges left and
right multiplication.

By (B1) and the pairing,
\begin{equation}\label{eq:kernels}
 \ker\mathcal L_j^{n+1}=\Span\{F_b:\epsilon_j(b)\le n\},\qquad
 \ker\mathcal R_j^{n+1}=\Span\{F_b:\epsilon^R_j(b)\le n\},
\end{equation}
and $\mathcal L_j^nF_b=n!\,F_{\tilde e_j^{\max}b}$ if $n=\epsilon_j(b)$.

\subsubsection{A polynomial model of the three-flag algebra}
For $w\in V(\eta)$ put $F_w(n)=\ell_\eta(n^{-1}w)$. By \eqref{eq:presentation},
the map $w\mapsto F_w$ is an isomorphism of $V(\eta)$ onto
$\mathcal V_\eta=\{F:\mathcal R_j^{\eta_j+1}F=0\ \forall j\}$, and
$\mathcal L_jF_w=-F_{e_jw}$.

\begin{proposition}\label{prop:poly-model}
There is an injective ring homomorphism
\[
 R_H\hookrightarrow\C[U][L_1,\dots,L_r,M_1,\dots,M_r]
\]
whose image is
$\bigoplus_{\lambda,\eta}\{F\in\C[U]:\mathcal L_j^{\lambda_j+1}F=
\mathcal R_j^{\eta_j+1}F=0\ \forall j\}\,L^\lambda M^\eta$, where
$L^\lambda=\prod_aL_a^{\lambda_a}$ and $M^\eta=\prod_aM_a^{\eta_a}$. The summand
$(V(\mu)\otimes V(\lambda)\otimes V(\eta))^H$ maps into the summand indexed by
$(\lambda,\eta)$, and $\mu$ is recovered as follows: if the coefficient $F$ has
weight $\gamma$ for the conjugation action of $T$ on $U$, then
$\mu^*=\lambda+\eta-\gamma$. The corner pairings of types
$(\varpi_a^*,0,\varpi_a)$ and $(\varpi_a,\varpi_a^*,0)$ map to $M_a$ and
$L_{a^*}$ respectively, where $\varpi_{a^*}=\varpi_a^*$.
\end{proposition}

\begin{proof}
An invariant $T\in V(\mu)\otimes V(\lambda)\otimes V(\eta)$ is the same as a
highest weight vector of weight $\mu^*$ in $V(\lambda)\otimes V(\eta)$. Viewing
such a vector as a linear map $V(\lambda)^\vee\to V(\eta)$, the condition
$e_jT=0$ says that the map is $\mathcal U$-linear. By \eqref{eq:presentation},
a $\mathcal U$-linear map is determined by the image $w$ of $\ell_\lambda$, and
$w$ can be any vector with $e_j^{\lambda_j+1}w=0$ for all $j$; the weight of
$w$ is $\mu^*-\lambda$. Send the invariant to $F_wL^\lambda M^\eta$. Since
$\mathcal L_jF_w=-F_{e_jw}$, the conditions on $w$ become
$\mathcal L_j^{\lambda_j+1}F_w=0$, which gives the image.

For multiplicativity, recall that the product in $R_H$ is the Cartan product
in each slot. In the second slot, the Cartan product of the two invariants has
leading component $v_{\lambda+\lambda'}\otimes(w\cdot w')$, where $w\cdot w'$ is
the Cartan product of $w$ and $w'$ in the third slot; and the functional
$\ell_{\eta+\eta'}$ restricts to $\ell_\eta\otimes\ell_{\eta'}$ on the Cartan
component, so $F_{w\cdot w'}=F_wF_{w'}$.

The weight formula follows from
$F_w(tnt^{-1})=\eta(t)\gamma_w(t)^{-1}F_w(n)$, where $\gamma_w=\mu^*-\lambda$ is
the weight of $w$, so that $F_w$ has conjugation weight $\eta-\mu^*+\lambda$.
Finally, for the type $(\varpi_a^*,0,\varpi_a)$, the vector $w$ is a highest
weight vector of $V(\varpi_a)$, so $F_w=1$ and the image is $M_a$; the type
$(\varpi_a,\varpi_a^*,0)$ is similar.
\end{proof}

\subsubsection{The valuation and string parameters}
Fix a reduced word $\mathbf i=(i_1,\dots,i_N)$ of $w_0$. The map
\[
 \Phi_{\mathbf i}:\A^N\to U,\qquad
 (a_1,\dots,a_N)\mapsto x_{i_1}(a_1)\cdots x_{i_N}(a_N),
\]
is birational \cite[Thm.~1.2]{FominZelevinsky1999}; see also \cite[Prop.~4.5, with $u=e$, $v=w_0$]{BZ2001}. Composing
Proposition~\ref{prop:poly-model} with $\Phi_{\mathbf i}^*$ and taking the
lexicographically largest monomial in $\bigl(a_1,\dots,a_N,\allowbreak L_1,\dots,L_r,\allowbreak M_1,\dots,M_r\bigr)$
defines a valuation $v_{\mathbf i}$ on $R_H$. Since conjugation by $t\in T$
multiplies $a_b$ by $\alpha_{i_b}(t)$, a monomial $a^{\boldsymbol t}$ has
conjugation weight $\sum_bt_b\alpha_{i_b}$, and Proposition~\ref{prop:poly-model}
gives \eqref{eq:omitted-weight}.

\begin{lemma}\label{lem:string-param}
For $b\in\mathcal B$, define $b^{(0)}=b$ and, for $1\le k\le N$,
$t_k=\epsilon_{i_k}(b^{(k-1)})$ and $b^{(k)}=\tilde e^{\max}_{i_k}b^{(k-1)}$.
Then $b^{(N)}=1$, and $\boldsymbol t(b)=(t_1,\dots,t_N)$, the \emph{string
parameter} of $b$, is the lexicographically largest exponent of
$\Phi_{\mathbf i}^*F_b$, with coefficient $1$. In Lusztig's
coordinates, the map $S:c\mapsto\boldsymbol t(b_{\mathbf h}(c))$ extends to an
injective continuous piecewise-linear map $S:\R^N_{\ge0}\to\R^N_{\ge0}$.
\end{lemma}

\begin{proof}
\emph{Termination.} Write $E_j=\tilde e^{\max}_j$, viewed as a map on labels.
Each $E_j$ is idempotent, and $E_jE_k=E_kE_j$ if $\langle\alpha_j,\alpha_k^\vee
\rangle=0$. In simply laced type, if $\langle\alpha_j,\alpha_k^\vee\rangle=-1$,
compute in a reduced word beginning with $jkj$ using (B3):
\[
 (a,b,c)_{jkj}\xrightarrow{E_j}(0,b,c)_{jkj}=(b+c,0,b)_{kjk}
 \xrightarrow{E_k}(0,0,b)_{kjk}=(b,0,0)_{jkj}\xrightarrow{E_j}(0,0,0),
\]
with the remaining coordinates unchanged; by symmetry $E_kE_jE_k$ has the same
effect. So the $E_j$ satisfy the braid relations, and by Matsumoto's theorem
(Section~\ref{subsec:bg-groups}) the composite $E_{i_N}\cdots E_{i_1}$ does not depend on the reduced word $\mathbf i$
of $w_0$. Choosing a reduced word ending with $j$ shows that
$E_{i_N}\cdots E_{i_1}(b)$ has $\epsilon_j=0$, for every $j$; by (B1) this label
equals $1$. In non-simply-laced type, the same holds by folding (B3).

\emph{Leading exponent.} For $1\le k\le N+1$ put
$G_k(a_k,\dots,a_N)=F_{b^{(k-1)}}(x_{i_k}(a_k)\cdots x_{i_N}(a_N))$, so that
$G_1=\Phi^*_{\mathbf i}F_b$ and $G_{N+1}=F_{b^{(N)}}=F_1=1$. Taylor expansion in
$a_k$ gives
\[
 G_k=\sum_{n\le t_k}\frac{a_k^n}{n!}\,(\mathcal L_{i_k}^nF_{b^{(k-1)}})
 (x_{i_{k+1}}(a_{k+1})\cdots),
\]
whose top term is $a_k^{t_k}G_{k+1}$ by \eqref{eq:kernels}: the factor
$t_k!$ from the maximal derivative cancels the Taylor denominator.
By descending induction on $k$, the lexicographically largest exponent of
$G_k$ is $(t_k,\dots,t_N)$, with coefficient $1$.

\emph{Piecewise linearity.} The procedure defining $\boldsymbol t$ consists of
changes of reduced word, readings of a first coordinate and settings of a first
coordinate to $0$; by (B3), all these are piecewise-linear and continuous on
$\R^N_{\ge0}$. The reverse procedure (starting from $0$ and applying
$\tilde f^{t_N}_{i_N},\dots,\tilde f^{t_1}_{i_1}$, i.e.\ adding $t_k$ to the
first coordinate in a word beginning with $i_k$) is also piecewise linear and
inverts $S$ on integral, hence by homogeneity on rational, and by continuity on
all points of $\R^N_{\ge0}$. Hence $S$ is injective.
\end{proof}

\begin{example}\label{ex:sl3-string}
For $\SL_3$ and $\mathbf h=\mathbf i=(1,2,1)$, the procedure gives
$S(a,b,c)=(a,b+c,b)$; the image is the \emph{string cone}
$\{t_1\ge0,\ t_2\ge t_3\ge0\}$, and $S(1,1,1)=(1,2,1)=q_\rho(121)$. For the
coordinate $z$ of $U$ (the matrix entry in the upper right corner),
$\Phi^*_{\mathbf i}z=ab$, with exponent $(1,1,0)$; for $xy-z$,
$\Phi^*_{\mathbf i}(xy-z)=bc$, with exponent $(0,1,1)$.
\end{example}

By (B2), $F_b\in\mathcal V_\eta$ if and only if the string parameter satisfies
\begin{equation}\label{eq:string-bounds}
 \beta_k(\boldsymbol t):=t_k+\textstyle\sum_{l>k}t_l\langle\alpha_{i_l},\alpha_{i_k}^\vee\rangle
 \le\eta_{i_k}\qquad(1\le k\le N).
\end{equation}
Indeed, $b=\tilde f^{t_1}_{i_1}\cdots\tilde f^{t_N}_{i_N}1$, the label
$b^{(k)}$ has degree $\gamma_k=\sum_{l>k}t_l\alpha_{i_l}$ and satisfies
$\epsilon_{i_k}(b^{(k)})=0$, and by (B2) the step from $b^{(k)}$ to $b^{(k-1)}$
stays inside the set of (B2) if and only if
$t_k\le\langle\eta-\gamma_k,\alpha_{i_k}^\vee\rangle$, which is
\eqref{eq:string-bounds}.

\subsubsection{The local cone}
Fix the Lusztig coordinate word $\mathbf h$ used in $S$ and write
$\epsilon_j(c)=\epsilon_j(b_{\mathbf h}(c))$ on integral points,
extended piecewise linearly to $\R^N_{\ge0}$ by (B3).
The coordinates of $\lambda$ and $\eta$ in the inequalities below are
their fundamental-weight coordinates.

\begin{theorem}\label{thm:app-local-cone}
The elements $F_bL^\lambda M^\eta$ with $\epsilon_j(b)\le\lambda_j$ and
$F_b\in\mathcal V_\eta$ form a basis of $R_H$, with pairwise distinct values
$(\boldsymbol t(b);\lambda,\eta)$. The value semigroup is
$Q_{3,\mathbf i}=K_{3,\mathbf i}\cap\Z^{N+2r}$, where
\[
 \begin{aligned}
 K_{3,\mathbf i}&=\bigl\{(S(c);\lambda,\eta):(c,\lambda,\eta)\in
 \mathcal P\bigr\},\\
 \mathcal P&=\bigl\{c\ge0,\ \lambda\ge0,\ \eta\ge0,\
 \epsilon_j(c)\le\lambda_j\ \forall j,\
 \beta_k(S(c))\le\eta_{i_k}\ \forall k\bigr\}.
 \end{aligned}
\]
$K_{3,\mathbf i}$ is a full-dimensional pointed rational polyhedral cone. The product of two
basis elements has as leading term a nonzero multiple of the basis element whose
value is the sum of their values.
\end{theorem}

\begin{proof}
By Proposition~\ref{prop:poly-model}, \eqref{eq:kernels} and
\eqref{eq:string-bounds}, the listed elements form a basis, indexed by the
integral points of $\mathcal P$. Their values are distinct by
Lemma~\ref{lem:string-param}, so the value of a linear combination is the
largest value of a basis element occurring in it, and the value semigroup is the
set of integral points of $K_{3,\mathbf i}$ coming from integral points of
$\mathcal P$. Every integral point of $K_{3,\mathbf i}$ arises in this way,
because the inverse of $S$ maps integral points to integral points.

The functions $\epsilon_j$ and $\beta_k\circ S$ are piecewise linear, so
$\mathcal P$, and hence its image $K_{3,\mathbf i}$ under the piecewise-linear
map $(S,\id,\id)$, is a finite union of rational polyhedral cones. The value
semigroup is closed under addition, because $R_H$ is a domain. This forces
$K_{3,\mathbf i}$ to be convex: for rational points $x,y\in K_{3,\mathbf i}$
there is $m>0$ with $mx,my$ integral, so $m(x+y)\in K_{3,\mathbf i}$ and
$x+y\in K_{3,\mathbf i}$; since rational points are dense in a finite union of
rational polyhedral cones and $K_{3,\mathbf i}$ is closed, $K_{3,\mathbf i}$ is
convex, hence a rational polyhedral cone. The cone $K_{3,\mathbf i}$ is pointed
because $K_{3,\mathbf i}$ lies in the positive orthant, and full-dimensional because the image of $R_H$ generates the
fraction field $\C(a_1,\dots,a_N,L,M)$, as every $F\in\C[U]$ lies in some
kernel space. The last assertion follows from the additivity of leading
exponents.
\end{proof}

This proves Theorem~\ref{thm:local-cone}(i) and (ii); the assertion
$\pi_\partial^{-1}(0)=\{0\}$ holds because the invariants of the trivial
representation are the constants.

\subsubsection{Corner values}

\begin{lemma}\label{lem:extremal}
Let $w$ span the lowest weight line of $V(\lambda^*)$. Then the string parameter
of $F_w$ is $q_\lambda(\mathbf i)$, and
$\sum_bq_{\lambda,b}(\mathbf i)\alpha_{i_b}=\lambda+\lambda^*$.
\end{lemma}

\begin{proof}
An extremal vector spans a weight line whose weight belongs to the Weyl
orbit of the highest weight. By $\mathcal L_jF_w=-F_{e_jw}$, the string parameter of $F_w$ records the maximal
powers of $e_{i_1},e_{i_2},\dots$ that successively apply to $w$. The lowest
weight of $V(\lambda^*)$ is $w_0\lambda^*=-\lambda$. Suppose that after $b-1$
steps the vector is an extremal vector of weight
$-s_{i_{b-1}}\cdots s_{i_1}\lambda$. Since the word is reduced, the coroot
$\beta_b^\vee=s_{i_1}\cdots s_{i_{b-1}}\alpha_{i_b}^\vee$ is positive, so this
extremal vector lies at the bottom of its $\alpha_{i_b}$-string, of length
$\langle s_{i_{b-1}}\cdots s_{i_1}\lambda,\alpha_{i_b}^\vee\rangle
=\langle\lambda,\beta_b^\vee\rangle$; maximal raising leads to the extremal
vector of weight $-s_{i_b}\cdots s_{i_1}\lambda$. After $N$ steps the weight is
$-w_0\lambda=\lambda^*$, which gives the last formula.
\end{proof}

By Proposition~\ref{prop:poly-model} and Lemma~\ref{lem:extremal}, the corner
pairing of type $(0,\varpi_a,\varpi_a^*)$ is $F_wL_aM_{a^*}$ with $w$ in the
lowest weight line of $V(\varpi_a^*)$, of value $(q_{\varpi_a};\varpi_a,
\varpi_a^*)$; the other two types have values $(0;0,\varpi_a)$ and
$(0;\varpi_a^*,0)$. This proves Theorem~\ref{thm:local-cone}(iii).

\subsubsection{Interiority}

\begin{proposition}\label{prop:rho-interior}
For every reduced word $\mathbf i$, the point $(q_\rho(\mathbf i);2\rho,2\rho)$
lies in the interior of $K_{3,\mathbf i}$.
\end{proposition}

\begin{proof}
\emph{The string parameter of $b_{\mathbf 1}$ is $q_\rho(\mathbf i)$.} By (B4)
and \eqref{eq:kernels}, $F_{b_{\mathbf 1}}$ lies in $\bigcap_j\ker\mathcal L_j^2$
and has degree $2\rho$. Let $w$ span the lowest weight line of $V(\rho)$. Then
$F_w$ also has degree $2\rho$ (its conjugation weight is $\rho-(-\rho)$), and
lies in $\bigcap_j\ker\mathcal L_j^2$, because the $\alpha_j$-string through $w$
has length $\langle\rho,\alpha_j^\vee\rangle=1$. The inversion $n\mapsto n^{-1}$
of $U$ exchanges the operators $\mathcal L_j$ and $-\mathcal R_j$ and preserves
degrees, so the degree-$2\rho$ part of $\bigcap_j\ker\mathcal L_j^2$ has the
dimension of the degree-$2\rho$ part of $\mathcal V_\rho\simeq V(\rho)$, which is
the lowest weight line. Hence $F_{b_{\mathbf 1}}$ is a multiple of $F_w$, and by
Lemma~\ref{lem:extremal} (with $\rho^*=\rho$), $S(\mathbf 1)=q_\rho(\mathbf i)$.
In particular $F_{b_{\mathbf 1}}\in\mathcal V_\rho$.

\emph{The point is interior.} The point $(\mathbf 1;2\rho,2\rho)$ lies in the
interior of $\mathcal P$: its coordinates are positive, and the inequalities
defining $\mathcal P$ are strict there, since $\epsilon_j(b_{\mathbf 1})=1<2$
and $\beta_k(S(\mathbf1))\le\rho_{i_k}=1<2$ (because $F_{b_{\mathbf
1}}\in\mathcal V_\rho$). The map $(S,\id,\id)$ is continuous and injective on
the interior of $\mathcal P$, an open subset of $\R^{N+2r}$. By Brouwer's
invariance of domain (an injective continuous map from an open subset of
$\R^n$ to $\R^n$ is open), this map is open, so the image of a neighborhood of
$(\mathbf 1;2\rho,2\rho)$ is a neighborhood of
$(q_\rho(\mathbf i);2\rho,2\rho)$ contained in $K_{3,\mathbf i}$.
\end{proof}

This proves Theorem~\ref{thm:local-cone}(iv). The argument is uniform in the
root type and the reduced word, and needs no description of the facets of the
string cone.

\subsection{Scalarization of valuations}\label{app:scalar}

We prove Lemma~\ref{lem:scalarization}, following \cite[Thm.~2.17, Lem.~8.9, Prop.~8.10]{KavehManon2019}
but keeping track of the coefficients of the prescribed initial forms and of a
finite group action. The idea is that a lexicographic order on $\N^M$ agrees
with the linear form $\ell_C(a)=\sum_{k=1}^MC^{M-k}a_k$ on any finite set of
comparisons, if $C$ is large; and, by a Gr\"obner basis argument, finitely many
comparisons control the whole associated graded ring.

\begin{proof}[Proof of Lemma~\ref{lem:scalarization}]
\emph{Step 1: a presentation adapted to $v$.} Choose elements of $A$ whose
initial forms are the monomials $x^q$ for $q$ in a finite generating set of $Q$.
For each $i$ choose a monomial $m_i$ in these elements with
$\inn_v(m_i)=x^{q_i}$, and put $h_i=g_i-c_im_i$; discard $h_i$ if it is zero. Otherwise $h_i$ is
constant or $v(h_i)<q_i$. Let $b_1,\dots,b_s$ be the nonconstant elements among the chosen
elements, a set of algebra generators of $A$, the $g_i$ and the $h_i$. Let
$\pi:P=\C[y_1,\dots,y_s]\to A$ send $y_j$ to $b_j$, let $I=\ker\pi$, and let
$I_v$ be the kernel of $P\to\gr_vA$, $y_j\mapsto\inn_v(b_j)$; this map is
surjective, so $P/I_v\simeq\gr_vA\simeq\C[Q]$. Give $y_j$ the weight
$v(b_j)\in\N^M$. Then the initial ideal of $I$ for this weight (with respect to
the lexicographic order) equals $I_v$: an element of $I$ has an initial form
lying in $I_v$, because the leading terms must cancel in $\gr_vA$; conversely,
a weight-homogeneous element $G$ of $I_v$ has image of value smaller than its
weight, and subtracting from $G$ an expression of $\pi(G)$ as a polynomial in
the $b_j$ of smaller weights (which exists by repeatedly subtracting leading
terms, a process that terminates because $\N^M$ is well ordered) gives an
element of $I$ with initial form $G$. Note that $I$ contains the relations
$r_i=y_{\nu(i)}-c_iy^{a_i}-\tilde h_i$, where $b_{\nu(i)}=g_i$, the monomial
$y^{a_i}$ represents $m_i$, and $\tilde h_i$ represents $h_i$; the initial form
of $r_i$ is $y_{\nu(i)}-c_iy^{a_i}$.

\emph{Step 2: scalarization.} Recall that, for a weight on the monomials of a polynomial ring, the
\emph{initial ideal} of an ideal is spanned by the initial forms (the sums of
the terms of largest weight) of its elements, and a \emph{Gr\"obner basis} is a
finite subset of the ideal whose initial forms generate the initial ideal
\cite[\S15.2, pp.~327--329]{Eisenbud1995}. We homogenize so that every standard-degree
piece is finite dimensional; inclusion of initial ideals with the same
Hilbert function will then force equality. Let $I^h\subset P[z]$ be the homogenization of
$I$ with respect to the standard grading, with $z$ of weight $0$. Choose a
finite set $\mathcal G\subset I^h$ of homogeneous elements whose initial forms
generate the initial ideal of $I^h$ (a Gr\"obner basis), and add the
homogenized relations $r_i^h$. Only finitely many lexicographic comparisons
between weights of terms occur in the elements of $\mathcal G$; for $C\gg0$,
the linear form $\ell_C$ preserves all of them, so every $G\in\mathcal G$ has
the same initial form for the weight $\ell_C\circ v$ as for the lexicographic
weight. Hence the lexicographic initial ideal of $I^h$ is contained in the
$\ell_C$-initial ideal; both are homogeneous with the Hilbert function of
$P[z]/I^h$, so they are equal. Setting $z=1$ shows that the $\ell_C$-initial
ideal of $I$ is $I_v$.

\emph{Step 3: the filtration.} Give $y_j$ the degree $\ell_C(v(b_j))$, which
is a positive integer, and let $F$ be the image in $A$ of the corresponding
degree filtration of $P$. By Step~2, $\gr_FA\simeq P/I_v\simeq\C[Q]$, with $x^q$
of degree $h(q)=\ell_C(q)$, and $\inn_F(g_i)=c_ix^{q_i}$, because the
$\ell_C$-initial form of $r_i$ is $y_{\nu(i)}-c_iy^{a_i}$. The filtration is
positive because $\C[Q]_0=\C$, and projective because $\cR_F(A)$ is generated
by $t$ and the elements $b_jt^{\deg y_j}$.

\emph{Equivariance.} Suppose that $\Delta$ acts as in the statement. Replace
each chosen element by its projection to the eigenspace of the character by
which $\Delta$ acts on its leaf; the initial form is unchanged. Add the
eigencomponents of a set of algebra generators. Then $m_i$ and $g_i$ have the
same character, namely the character of the leaf of $q_i$, so $h_i$ is an
eigenvector, and all $b_j$ are eigenvectors. Let $\Delta$ act diagonally on
$P$; then $I$ is $\Delta$-stable, we may enlarge $\mathcal G$ by its
$\Delta$-translates, and the filtration of Step~3 is $\Delta$-stable.
\end{proof}

\subsection{Auxiliary geometric facts}\label{app:geometric-facts}
\label{app:background-proofs}
We justify the refinements concerning the basic affine space, scheme
fibers of the adjoint quotient and products of nilpotent cones recorded
in Section~\ref{sec:background}.

\paragraph{The basic affine space (Fact~\ref{fact:basic-affine})}\label{app:basic-affine-proof}
For parts (ii) and (iii), the geometric argument is as follows. The variety $\mathcal Z_H$ is the
image of the total space of the vector bundle
$\bigoplus_a\mathcal O(-\varpi_a)$ on $H/B$ under the proper map to
$\bigoplus_aV(\varpi_a)$ given by the sections; the locus where exactly the
coordinates indexed by a nonempty $J\subset\{1,\dots,r\}$ vanish has
codimension $|J|+|\Phi_J^+|\ge2$. For the class-group assertion, write
$\mathcal E=\bigoplus_a\mathcal O(-\varpi_a)$ and let
$\mathcal E^\circ$ be the complement of its coordinate zero divisors.
The torus bundle $\mathcal E^\circ\to H/B$ is $H/U\to H/B$.
Pullback identifies $\Pic(\mathcal E)$ with $\Pic(H/B)$, and the
localization sequence identifies $\Pic(\mathcal E^\circ)$ with the quotient
by the classes of those divisors. Their classes correspond to
$\mathcal O(-\varpi_a)$ and generate $\Pic(H/B)$, so $\Pic(H/U)=0$.
The codimension-two complement then gives
$\Cl(\mathcal Z_H)=\Pic(H/U)=0$. Units are constants because
$\C[H/U]$ is graded by $X^+(T)$ with degree-zero part $\C$.

\paragraph{Scheme fibers of the adjoint quotient (Fact~\ref{fact:steinberg})}\label{app:adjoint-fiber-proof}
For (iii), the dimension and orbit assertions are
\cite[Thm.~6.11 and Rem.~6.15(c)]{Steinberg1965}, and the differentials of
the fundamental characters are independent at every regular element by
\cite[Thm.~8.1]{Steinberg1965}.  The $r$ equations in the smooth group
$H$ have codimension $r$, so they define a Cohen--Macaulay complete
intersection.  The dense regular orbit in the fiber is smooth, hence the
fiber is reduced.  The nonregular complement has codimension at least two, so
Serre's criterion proves normality.

\paragraph{Products of nilpotent cones (Fact~\ref{fact:nilpotent-cones}(iv))}\label{app:nilpotent-products}
The $\Q$-factoriality in part~(iv) uses the proof, not just the statement,
of \cite[Prop.~4.1]{LehnNamikawaSorger2012}. Use the product orbit
stratification and the smooth factor. Every positive stratum codimension is
even, and a codimension-two transverse singularity is a surface quotient
$\C^2/F$ with $F\subset\SL_2$ finite. Choose a common multiple $m$ of the
finitely many orders of their local divisor class groups. We prove that
$mD$ is Cartier for every local analytic Weil divisor $D$, by induction on
the stratum codimension. Smooth and codimension-two points give the initial
cases. At a deeper stratum, take a sufficiently small contractible Stein
neighborhood $U$ and let $B$ be its intersection with that stratum; exclude
strata not incident to the chosen point. Then $\operatorname{codim}_U B\ge4$,
and $mD$ is Cartier on $U\setminus B$ by induction.

For a local complete intersection, the homological-depth theorem used in
the cited proof gives
\[
 H_2(U,U\setminus B;\Z)=H_3(U,U\setminus B;\Z)=0.
\]
The universal coefficient theorem and the exact sequence of the pair give
$H^2(U\setminus B,\Z)=0$. Cohen--Macaulayness and the codimension bound give
$H^1(U\setminus B,\mathcal O)=H^2(U\setminus B,\mathcal O)=0$.
The exponential sequence yields $\Pic_{\mathrm{an}}(U\setminus B)=0$.
The rank-one reflexive sheaf associated with $mD$ therefore extends as
$\mathcal O_U$ by normality, proving the induction. Canonicity follows from
the product of the Springer resolutions in (ii), under which canonical
divisors pull back.

\section{Inputs to the Poisson construction}\label{app:poisson-inputs}
This appendix collects supporting arguments for Part~\ref{part:poisson}.
Section~\ref{app:poisson-algebra} supplies the algebraic extension results
used for Poisson compactifications and the Dolbeault form.
Section~\ref{paper:app-properness} proves the Hitchin-projectivity input
used to construct the compact Goldman period.
Section~\ref{paper:app-ih} gives an alternative sufficient nonvanishing
criterion; the main period argument does not use this criterion.

\subsection{Differential forms and Poisson normalization}\label{app:poisson-algebra}

\subsubsection{Forms on finite quotients}\label{app:finite-quotient-forms}
\begin{lemma}[Extension of finite-quotient forms]
\label{paper:lem-finite-quotient-forms}
Let $V$ be a smooth affine complex variety and $F$ a finite group acting on $V$.
Every $F$-invariant regular $p$-form, $p\geq1$, descends rationally to $Z=V/F$,
and its pullback to every resolution of $Z$ is regular.
\end{lemma}
\begin{proof}
Work componentwise, replacing $F$ by the stabilizer of a component of $V$.
Rational descent is finite Galois descent. For a resolution $Y\to Z$,
normalize a dominant component of $V\times_ZY$, obtaining a finite surjection $W\to Y$ and a
morphism $W\to V$. The descended form $\tau$ therefore pulls back
regularly at every codimension-one point of $W$.

Fix prime divisors $E\subset Y$ and $E'\subset W$ above it. Put $d=\dim Y$.
Choose a uniformizer $t$ at $E$ and lifts $u_1,\ldots,u_{d-1}$ of a separating
transcendence basis of $\mathbb C(E)/\mathbb C$. For a uniformizer $s$
at $E'$, write $t=s^ev$, with $v$ a unit. Since $W$ is smooth near the
generic point of $E'$, the differentials $(dt,du_i)$ and $(ds,du_i)$
are bases over the respective discrete valuation rings, and
\[
 dt=s^{e-1}a\,ds+s^e\sum_i b_i\,du_i,
 \qquad a\in\mathcal O_{W,E'}^\times,\quad b_i\in\mathcal O_{W,E'}.
\]
Indeed $a\equiv ev\pmod{s}$, which is invertible in characteristic zero.
Write $\tau=dt\wedge\alpha+\beta$, using only the $du_i$ in $\alpha$
and $\beta$. A coefficient of order $m<0$ in $\alpha$ pulls back to
order $em+e-1<0$ in the $ds$ component, impossibly. Thus $\alpha$ is
regular. Subtract the regular pullback of $dt\wedge\alpha$. A coefficient of order
$m<0$ in $\beta$ then pulls back to order $em<0$, again impossibly.
Thus $\tau$ has no pole along $E$. Smoothness of $Y$ now gives
regularity everywhere.
\end{proof}

\subsubsection{Poisson normalization (Fact~\ref{paper:fact-poisson-forms}(i))}\label{app:poisson-normalization}
\begin{lemma}[Poisson normalization and reduction]
\label{paper:lem-poisson-normalization}
Let $A$ be a finitely generated complex Poisson algebra.
The radical of a Poisson ideal is Poisson, and the minimal primes of
$A$ are Poisson. If $A$ is a domain, its Poisson bracket extends
uniquely to its normalization $B$.
\end{lemma}

\begin{proof}
For a derivation $\delta$ of a complex algebra, the formal exponential
$\exp(t\delta)$ is an automorphism of $A[[t]]$. This automorphism sends nilpotent
elements to nilpotent elements. Since $A_{\mathrm{red}}[[t]]$ is reduced,
every coefficient of $\exp(t\delta)(a)$ is nilpotent when $a$ is.
Thus $\delta$ preserves the nilradical. Apply this to the Hamiltonian
derivations on a quotient by a Poisson ideal to obtain the radical
assertion. On a reduced ring, localization at a minimal prime is a
field and its kernel is that prime. Derivations extend to localizations,
so they preserve this kernel. Applying this after reduction proves
the minimal-prime assertion.

Now assume that $A$ is a domain, and put $K=\operatorname{Frac}(A)$.
The bracket extends uniquely to $K$ by the quotient rule. The ring $B$
is finite over $A$. Seidenberg's derivation theorem
\cite[\S3, pp.~168--169]{Seidenberg1966} implies that
$\{a,B\}\subset B$ for every $a\in A$. Consequently the bracket modulo
$B$ induces an alternating $B$-bilinear map
\[
 \bigwedge\nolimits_B^2\Omega^1_{B/A}\longrightarrow K/B,
 \qquad db\wedge dc\longmapsto\{b,c\}\bmod B.
\]
Let $\mathfrak p$ be a height-one prime of $B$, and put
$\mathfrak q=\mathfrak p\cap A$. The ring $B_{\mathfrak p}$ is a
discrete valuation ring, and $\kappa(\mathfrak p)/\kappa(\mathfrak q)$
is finite separable. The cotangent exact sequence therefore shows that
\[
 \Omega^1_{B/A}\otimes_B\kappa(\mathfrak p)
\]
is a quotient of the one-dimensional cotangent space of
$B_{\mathfrak p}$. Since $\Omega^1_{B/A}$ is finite, Nakayama's lemma
makes its localization at $\mathfrak p$ cyclic. The second exterior
power vanishes. The displayed map, localized at $\mathfrak p$, thus
gives $\{B,B\}\subset B_{\mathfrak p}$. Intersecting these discrete
valuation rings yields $\{B,B\}\subset B$. If $\dim B=0$, then
$A=B=\C$ and the assertion is immediate. The Poisson identities and
uniqueness are inherited from $K$.
\end{proof}

\subsection{Projectivity of the parahoric Hitchin morphism}\label{paper:app-properness}

This subsection proves that the Hitchin morphism on the algebraic Dolbeault
closure of Section~\ref{r2:sec-dolbeault-closure} is projective for the
parahoric stability condition used there. The method is to pass to a finite
cover on which parahoric Higgs torsors become equivariant Higgs bundles
(Lemma~\ref{r2:lem-auxiliary-parahoric-cover}), and then to apply the
semistable reduction theorem for $\Theta$-stratified stacks
\cite[Cor.~6.12]{AHLH2023}.

\begin{lemma}[An auxiliary cover in every genus]
\label{r2:lem-auxiliary-parahoric-cover}
Let $X$ be a smooth connected projective complex curve,
$D=\sum_{i=1}^k p_i$ with distinct $p_i$, and let
$\alpha_i\in X_*(T)\otimes_{\Z}\Q$ be rational parahoric weights.
Write $d_i$ for the least positive integer such that
$d_i\alpha_i\in X_*(T)$.  There is a connected finite Galois cover
$\pi:Y\to X$ whose inertia at every point above $p_i$ has order
$d_i$.  Additional branch points can be given trivial inertia
representation in $G$.  With the prescribed inertia representation
above $p_i$, equivariant principal $G$-bundles on $Y$ with Higgs
field twisted by
\[
 L_Y=\pi^*K_X(D)
\]
are equivalent to the original weak logarithmic parahoric Higgs
torsors on $X$.  The equivalence preserves the parabolic-degree
(Ramanathan) stability condition.  No restriction on the genus of
$X$ is needed.
The cover may also be chosen with $g(Y)\ge2$.
\end{lemma}

\begin{proof}
Choose $q\in X\setminus D$.  For each $i$ with $d_i>1$, Riemann--Roch,
applied to $\mathcal O_X(Nq)$ for sufficiently large $N$, supplies a
meromorphic function $f_i$ with a simple zero at $p_i$ and value $1$
at every $p_j$ with $j\ne i$.  The only poles of $f_i$ lie at $q$; its other
zeros are allowed to be auxiliary branch points.  Let $Y$ be the
normalization of $X$ in the compositum, taken in a common algebraic closure,
\[
 \C(X)\bigl(f_i^{1/d_i}:d_i>1\bigr).
\]
This is a field and a finite Galois extension because $\C$ contains
all roots of unity.  The normalization of the curve in this field is a smooth connected
projective curve.  At $p_i$, every $f_j$ with $j\ne i$ is a unit
and has a $d_j$-th root in $\C[[z]]$; the extension obtained from
$f_i$ has ramification index exactly $d_i$, since its valuation is
one.  Thus the total inertia has order exactly $d_i$.  If all
$d_i=1$, take the identity cover.

To arrange $g(Y)\ge2$, take a further double cover of $X$ defined by
$w^2=f_0$, where $f_0(p_i)=1$ and $f_0$ has a pole of even order
$2N$ at $q$ and $2N$ distinct simple zeros away from $D$.
For sufficiently large $N$, a generic function with these properties
exists by Riemann--Roch. Riemann--Hurwitz gives genus $2g(X)-1+N\ge2$.
Replace $Y$ by the normalization of the compositum with this cover.
The puncture inertia is unchanged, and the new curve dominates a curve
of genus at least two.

Locally write $z=w^{d_i}$ and let $\zeta\in\mu_{d_i}$ act by
$w\mapsto\zeta w$.  Give the fiber the action
\[
 \rho_i(\zeta)=(d_i\alpha_i)(\zeta)\in T.
\]
The gauge $\Delta_i(w)=w^{d_i\alpha_i}$ obeys
$\Delta_i(\zeta w)=\rho_i(\zeta)\Delta_i(w)$; conjugation by
$\Delta_i$ is exactly the local equivariant--parahoric construction
of \cite[arXiv v3, Section~2.3]{KSZ2024}.  At every auxiliary branch point use
the trivial representation, whose descended group scheme is the
ordinary constant group scheme.  These local equivalences glue
over the unramified complement, as in
\cite[arXiv v3, Theorem~2.10]{KSZ2024}.

The twist must keep track of the auxiliary ramification.  If $R_D$
is the reduced inverse image of $D$ and
$R_{\mathrm{aux}}^{\mathrm{ram}}=\sum_{y\notin\pi^{-1}(D)}(e_y-1)y$,
the differential identifies the naturally linearized line bundles
\[
 \pi^*K_X(D)
   \simeq K_Y\bigl(R_D-R_{\mathrm{aux}}^{\mathrm{ram}}\bigr).
\]
Above $p_i$, this sends $dz/z$ to $d_i\,dw/w$ and gives the
logarithmic Higgs correspondence of \cite[arXiv v3, Theorem~3.7]{KSZ2024}.
At an auxiliary point, $z=w^e$ and $\pi^*dz=e w^{e-1}dw$.
Since the fiber inertia there is trivial, an invariant section
$a(w)\pi^*dz$ has $a(w)\in\mathfrak g[[w^e]]$ and descends to a
regular Higgs field.  The use of $L_Y$, rather than a logarithmic
twist at every branch point, is therefore essential.

Equivariant parabolic reductions descend by the same construction.
For every character used in the stability test, the degree of the
associated line bundle upstairs is $\deg(\pi)$ times its parahoric
degree downstairs.  The local correction at $p_i$ is the pairing with $\alpha_i$
after transporting the character by the relative position of the reduction; at auxiliary points it is zero because the inertia
representation is trivial.  This is the degree calculation in
\cite[arXiv v3, Lemma~4.11 and Theorem~4.13]{KSZ2024}, unchanged by the
auxiliary points.  Higgs-invariance of a reduction is preserved by
the local identifications just proved.  The stability inequalities
are consequently equivalent.
\end{proof}

\subsubsection{Properness for the actual parahoric stability condition}

The properness input must use the same stability condition as tame
nonabelian Hodge theory.  It is insufficient simply to identify two
different uses of the word ``semistable.''  The following argument avoids
a comparison between arbitrary parabolic weights and auxiliary
determinant-character numerical invariants.

\begin{lemma}[Projective weak parahoric Hitchin morphism]
\label{r2:lem-proper-parahoric}
Fix real parahoric weights on a smooth projective curve with finitely
many marked points, a semisimple complex structure group, and a
topological type.  The moduli space of semistable logarithmic parahoric
Higgs torsors for the parabolic-degree (Ramanathan) stability condition
has a projective Hitchin morphism.  Prescribing closed conjugacy-invariant conditions, including arbitrary
orbit closures, for the graded residues gives a closed subspace with a projective
restricted Hitchin morphism.
\end{lemma}

\begin{proof}
Rational weights first suffice. Lemma~\ref{r2:lem-auxiliary-parahoric-cover}
identifies the required logarithmic parahoric Higgs torsors with
finite-deck-equivariant principal Higgs bundles on a smooth cover
$\pi:Y\to X$, twisted by $L_Y=\pi^*K_X(D)$. The equivalence preserves
parabolic-degree stability and holds in every genus. Auxiliary branch
points have trivial inertia representation and introduce no Higgs poles.
Choose the cover with $g(Y)\ge2$ when using the ordinary GIT
construction upstairs.

Let $\mathfrak H$ be the stack of principal $G$-Higgs bundles on $Y$
twisted by this line bundle.  The Harder--Narasimhan stratification is a
well-ordered $\Theta$-stratification in characteristic zero
\cite[Lemma~6.19]{AHLH2023}.  It is preserved by the finite deck group
$\Gamma$.  An unstable $\Gamma$-equivariant Higgs bundle has a unique
ordinary Harder--Narasimhan reduction.  Uniqueness makes the reduction
$\Gamma$-invariant, including the cocycle condition for its
linearization.  Hence equivariant semistability is equivalent to
ordinary semistability.

The assertion about families is equally precise.  Write
$\mathfrak H^{h\Gamma}$ for the stack of equivariant objects.  There is
a natural equivalence
\[
 \operatorname{Map}(\Theta,\mathfrak H^{h\Gamma})
    \simeq
 \operatorname{Map}(\Theta,\mathfrak H)^{h\Gamma}.
\]
For an invariant Harder--Narasimhan stratum
$\mathfrak S\subset\operatorname{Map}(\Theta,\mathfrak H)$, its
evaluation morphism is a locally closed monomorphism, and uniqueness of
the canonical filtration gives
\[
 \mathfrak S^{h\Gamma}
       \simeq\mathfrak S\times_{\mathfrak H}\mathfrak H^{h\Gamma}.
\]
Thus its evaluation is the base change of the stratum embedding.
The open substacks and the well-ordered index set defining the
$\Theta$-stratification are inherited.  Prescribing the finite inertia
representation types is an open and closed condition and is constant
along $\Theta$; restriction to those types preserves the
stratification.  Under descent this is the required well-ordered
$\Theta$-stratification of the parahoric Higgs stack, with exactly the
Ramanathan semistable open.

Let $\mathfrak P\to\mathcal A$ denote this weak parahoric Higgs stack
of fixed topological and inertia type over its affine Hitchin base.
It is algebraic and locally of finite type over the noetherian scheme
$\mathcal A$.  Indeed, spaces of Higgs fields over bundle charts are cut out by
finite systems of linear equations; the finite collection of deck
linearizations is described by isomorphisms satisfying closed cocycle
equations.  The prescribed inertia types are open and closed.
The diagonal is affine: bundle isomorphisms on the proper curve form
an affine scheme, and compatibility with the Higgs fields and deck
linearizations imposes closed equations.  Thus $\mathfrak P$ is also
quasi-separated with affine stabilizers.

The bounded GIT construction of \cite[arXiv v3, Thm.~6.1]{KSZ2024}
gives a finite-type atlas for $\mathfrak P^{\mathrm{ss}}$, so this
semistable open is quasi-compact.  To apply that construction on the
auxiliary cover, regard
$L_Y=K_Y(R_D-R_{\mathrm{aux}}^{\mathrm{ram}})$ as a subsheaf of
$K_Y(R_D+R_{\mathrm{aux}})$, where $R_{\mathrm{aux}}$ is the reduced
auxiliary ramification divisor.  Factoring the Higgs field through
$L_Y$ is a closed, deck-invariant linear condition on the logarithmic
parameter space.  It leaves Higgs-compatible reductions, hence
semistability, unchanged.  The fixed inertia conditions and the
choice $g(Y)\ge2$ therefore give exactly the bounded presentation
needed here.

The entire stack $\mathfrak P$, before imposing stability, satisfies
the existence part of the valuative criterion over $\mathcal A$
\cite[Lemma~5.25]{Rega2025}.  This assertion is independent of the
numerical stability condition.  Together with the well-ordered
$\Theta$-stratification constructed above and the preceding stack
hypotheses, \cite[Cor.~6.12]{AHLH2023} implies that
$\mathfrak P^{\mathrm{ss}}\to\mathcal A$ is universally closed.
The good-moduli morphism
$\mathfrak P^{\mathrm{ss}}\to M$ is universally surjective, so
$M\to\mathcal A$ is universally closed as well.  Since $M$ is a
separated quasiprojective variety of finite type, this morphism is
proper, hence projective over the affine base.  The argument does
not identify Ramanathan stability with the determinant-character
stability subject to the admissibility inequalities in
\cite[Section~3.2.5]{Rega2025}.

For real weights choose rational weights in the same parahoric facet
and the same stability cell, retaining every wall through the given
weights. The relevant walls form a finite rational affine arrangement:
in a bounded alcove there are finitely many fundamental-character and
relative-position types, and a bounded linear weight contribution can
equal only finitely many integral degrees. Each cell has a rational point
in its relative interior. Thus the replacement preserves the parahoric
group scheme and all signs and equalities in the stability tests; compare
\cite[arXiv v3, Rem.~6.11]{KSZ2024}. Only the algebraic moduli construction uses
this replacement. The actual real weights remain in the analytic
nonabelian Hodge correspondence.

Finally, the graded residue is an algebraic morphism to a quotient
stack for the reductive quotient of the parahoric group at a puncture.
Any closed conjugacy-invariant residue condition is closed there;
this includes arbitrary coadjoint-orbit closures.  The intersection of such a condition with the semistable locus has closed image
in the good moduli space, and the restricted Hitchin map remains
projective.
\end{proof}

\begin{remark}[Why this is a properness proof for the required stability]
The proof uses the existence lemma for the full stack from
\cite{Rega2025}, not an identification of the semistable locus in that
thesis with an unspecified parabolic-degree condition.  The latter
identification would require an additional argument.  The equivariant
Harder--Narasimhan construction supplies that condition directly and
then invokes the stability-independent semistable-reduction theorem.
\end{remark}

\subsection{An optional intersection-cohomology criterion}\label{paper:app-ih}
The Poisson proof constructs a Goldman period directly. The following optional
criterion records another sufficient input to the Poisson construction.

We use the convention
$IC_U=j_{!*}\Q_{U_{\mathrm{reg}}}[d]$ and
$IH^r(U)=\mathbb H^{r-d}(U,IC_U)$ for a $d$-dimensional irreducible variety.

\begin{proposition}[Top weight supplies a logarithmic volume form]
\label{prop:top-weight-log-form}
Let $U$ be an irreducible $d$-dimensional complex variety.  If
\[
 \operatorname{Gr}_{2d}^{W}IH^d(U,\Q)\ne0,
\]
then $\overline P_1^{\mathrm{res}}(U)>0$.
\end{proposition}

\begin{proof}
Let $\rho:U'\to U$ be a projective resolution.  The decomposition theorem in
the category of polarizable Hodge modules \cite{Saito1990} gives $IC_U$ as a direct summand of
$R\rho_*\Q_{U'}[d]$.  Hence $IH^d(U)$ is a direct summand of $H^d(U')$ as a
mixed Hodge structure.  This is the Hodge-module form of the decomposition
theorem; a splitting of rational constructible complexes alone would not
justify the assertion about weights.

It follows that $\operatorname{Gr}_{2d}^WH^d(U')\ne0$.  For a smooth
$d$-dimensional variety the latter graded piece is of Hodge type $(d,d)$.
Strictness of the Hodge filtration therefore implies $F^dH^d(U',\C)\ne0$.
For a smooth projective simple-normal-crossings compactification $(Y,B)$ of
$U'$, the logarithmic Hodge--de Rham spectral sequence degenerates at $E_1$
\cite[\S3.2]{Deligne1971}.
Since its complex has no forms of degree greater than $d$, this gives
\[
 F^dH^d(U',\C)
   \cong H^0\bigl(Y,\Omega_Y^d(\log B)\bigr)
   =H^0\bigl(Y,\mathcal O_Y(K_Y+B)\bigr).
\]
The last space is nonzero, as required.
\end{proof}

\begin{remark}[The amount of curious Lefschetz that is needed]
For $d=2n$, one does not need the entire curious-hard-Lefschetz package to
apply Proposition~\ref{prop:top-weight-log-form}.  By intersection Poincare
duality, it suffices to show
$\operatorname{Gr}_0^WIH_c^d(U)\ne0$.  In particular, a surjection
\[
 L^n:\operatorname{Gr}_0^WIH_c^d(U)
       \longrightarrow\operatorname{Gr}_{2d}^WIH_c^{2d}(U)
\]
is sufficient: the target is the nonzero fundamental-class space
$\Q(-d)$ when $U$ is irreducible.  Here $L$ has cohomological degree two and
weight four, as in curious hard Lefschetz.  The nonvanishing condition is
weaker than full intersection-cohomological curious hard Lefschetz, and the
existence of a logarithmic pluriform is weaker still.  None of these
nonvanishing assertions is inferred merely from the existence of a Springer
summand.
\end{remark}

\section{Realization and obstructions for conjugacy-class closures}
\label{new:sec-counterexample}

The extension from full adjoint fibers to arbitrary conjugacy-class closures
has an obstruction in the interior.  Genericity eliminates representations
contained in proper parabolic subgroups, but it does not eliminate the
singularities of the prescribed orbit closures.  This appendix proves a
realization theorem making that distinction precise and gives an explicit
generic relative character variety which is not normal.  The argument is
independent of the spectral compactification theorem. The final subsection
gives a different obstruction: a normal $\SL_3$ interior that is not
$\Q$-factorial, so that a $\Q$-factorial model cannot preserve that interior.

For a tuple $\mathbf X=(X_1,\ldots,X_k)$ of reduced conjugacy-class closures
in a connected reductive group $G$, write
\[
 R(\mathbf X)=
 \left\{(A_j,B_j,M_i)\in G^{2g}\times\prod_iX_i:
       \prod_j[A_j,B_j]\prod_iM_i=1\right\},
 \qquad
 M(\mathbf X)=R(\mathbf X)_{\red}/\!\!/G.
\]
The genericity condition is the intrinsic condition
\eqref{eq:intro-generic} on the semisimple parts.  We call a representation
\emph{principal} when its effective stabilizer in $G_{\ad}$ is trivial.
Equivalently, its centralizer in $G$ is exactly $Z(G)$.

\subsection{Transversality without normality}

\begin{lemma}[The relation with arbitrary invariant puncture conditions]
\label{new:lem-closure-transversality}
Let $X_i$ be reduced $G$-stable closed subvarieties of reduced adjoint
fibers, with compatible intrinsically generic semisimple data.  Put
$P=G^{2g}\times\prod_iX_i$.  At every point $p\in R(\mathbf X)$ there is
an \'etale morphism
\[
 Z\times R(\mathbf X)\longrightarrow P
\]
at a point $(e,p)$ mapping to $p$, where $Z$ is smooth of dimension
$\dim G_{\der}$.  In particular $R(\mathbf X)$ is reduced.  The same
construction applies after replacing any $X_i$ by an invariant locally
closed subvariety containing the corresponding entry of $p$.
\end{lemma}

\begin{proof}
By Lemma~\ref{lem:spectral-genericity}, whose proof uses only the semisimple
parts of the puncture entries, every point of $R(\mathbf X)$ is $G$-irreducible.
The proof of Lemma~\ref{lem:spectral-relation-local} uses only the irreducibility
of $p$ and the conjugation invariance of the puncture conditions, so it applies
verbatim to $P$ and gives the \'etale morphism. Since $P$ is reduced, so is the
germ of $Z\times R(\mathbf X)$ at $(e,p)$, and hence $R(\mathbf X)$ is reduced at
$p$. The construction preserves every invariant puncture condition, which proves
the last assertion.
\end{proof}

Thus the relation removes smooth directions; it does not repair a
nonnormal puncture germ.  At a principal point, the quotient also removes
smooth directions, by the Luna slice theorem.  The next result shows that
principal points realizing a prescribed puncture germ can always be
constructed after introducing one additional puncture.

\subsection{A realization theorem on a twice-punctured torus}

\begin{lemma}[A dominant commutator map on the principal locus]
\label{new:lem-principal-commutator}
Let $G$ be a connected reductive group.  There is a nonempty open subset
$V\subset G^2$ consisting of pairs with closed orbit and centralizer
$Z(G)$.  The commutator morphism
\[
 c:V\longrightarrow G_{\der},\qquad (A,B)\longmapsto[A,B]
\]
is smooth and dominant.
\end{lemma}

\begin{proof}
Fix a maximal torus $T$.  Choose $A\in T$ with full centralizer
$Z_G(A)=T$.  Such choices form a nonempty open subset of $T$: one may
exclude the root kernels and the fixed loci of the nonidentity Weyl-group
elements.  This use of the \emph{full} centralizer avoids a possible
disconnected centralizer in a non-simply-connected group.

There are only finitely many parabolic subgroups containing $T$.  Choose
$B$ outside all the proper ones.  We can simultaneously require $B$ to
belong to the big cell $U^-TU$ and to have nonzero simple-root
coordinates in its $U$-factor.  These are nonempty open conditions.
If a parabolic contains $A$, a maximal torus of that parabolic containing
$A$ must equal $T$; hence the pair $(A,B)$ is irreducible.  The orbit of $(A,B)$ is
closed by Lemma~\ref{lem:schur}.  Any element centralizing $A$ belongs
to $T$.  If it also centralizes $B$, uniqueness of the big-cell
factorization and the nonzero simple-root coordinates imply that every
simple root takes value $1$ on it.  It therefore belongs to $Z(G)$.
The Luna slice theorem at this pair, applied to the effective action of
$G_{\ad}$, gives a nonempty open principal stable locus $V$.

The differential of $c$ is surjective at every pair in $V$, even when
$[A,B]\ne1$: by the handle part of the computation in the proof of
Lemma~\ref{lem:differential}, the annihilator of its image is
$\Lie Z_G(A,B)=\mathfrak z(\mathfrak g)$, so the image is the orthogonal
complement $\mathfrak g_\der$ of the center. Since source and target are smooth,
$c|_V$ is smooth, and its image is a nonempty open, hence dense, subset of the
connected group $G_\der$.
\end{proof}

\begin{lemma}[An additional generic regular semisimple puncture]
\label{new:lem-add-generic-puncture}
Let $G$ be connected reductive, let $q:G\to G_{\ab}$ be the
abelianization, and let $s\in G$ be semisimple.  There is a dense open
subset
\[
 W_s\subset q^{-1}(q(s)^{-1})
\]
consisting of regular semisimple elements $t$ such that the pair of
semisimple classes $(G\cdot s,G\cdot t)$ is compatible and
intrinsically $G$-generic.
\end{lemma}

\begin{proof}
Conjugate $s$ into a maximal torus $T$ and put
$T_{\der}=T\cap G_{\der}$, a maximal torus of the connected group
$G_{\der}$.  The torus fiber with the required abelianization is the
irreducible coset $T_s=s^{-1}T_{\der}$.
For each conjugacy class of maximal parabolics choose a representative
$P_j$ containing $T$ and an algebraic character $\chi_j$ of $P_j$
which is trivial on $Z(G)^\circ$ and has nonzero restriction to
$T_{\der}$.  Such a character exists: on $P_j\cap G_{\der}$ take a
sufficiently divisible positive multiple of the corresponding
fundamental weight, trivialize it on $Z(G)^\circ$, and descend through
the finite central isogeny
$Z(G)^\circ\times G_{\der}\to G$.  Using an actual character, and
allowing a divisible multiple, is necessary for a general isogeny type.
In $T_s$, exclude the root kernels and the finitely many sets
\begin{equation}\label{new:eq-generic-character-avoidance}
 (w\chi_j)(t)=(v\chi_j)(s)^{-1},
 \qquad v,w\in W(G,T).
\end{equation}
Every excluded set is proper, since each root and each $w\chi_j$ has
nonzero restriction to $T_{\der}$.  Their complement is a nonempty
Weyl-invariant open subset of the regular semisimple locus of $T_s$.
The conjugation saturation of this set is a dense open subset of
$q^{-1}(q(s)^{-1})$: this fiber is an irreducible translate of
$G_{\der}$, and the regular semisimple conjugacy quotient identifies
the relevant open set with the Weyl quotient of $T_s$.
If $G$ is a torus, there are no root or parabolic exclusions and
$T_s=\{s^{-1}\}$, so the same statement holds.

To check genericity, suppose that a proper parabolic $Q$ contains
conjugates $s'$ of $s$ and $t'$ of $t$ with $s't'\in[Q,Q]$.
Contain $Q$ in a maximal parabolic, and conjugate the latter to some
$P_j$.  Since $[Q,Q]\subset[P_j,P_j]$, we have
$\chi_j(s')\chi_j(t')=1$.  Each semisimple element of $P_j$ is
$P_j$-conjugate into $T$.  The image of such an element under $\chi_j$ is therefore the
value on $s$ or $t$ of a Weyl translate of $\chi_j$.  This gives one of
the excluded equations \eqref{new:eq-generic-character-avoidance}, a
contradiction.
\end{proof}

\begin{theorem}[Realization of every conjugacy-closure germ]
\label{new:thm-universal-realization}
Let $G$ be connected reductive, put $d=\dim G$, $z=\dim Z(G)$,
$h=d-z$, and $r=\rank G$, and
let $X=\overline{\cC}$ be a reduced conjugacy-class closure in $G$.
For every point $x\in X$, there exist a regular semisimple conjugacy
class $\cD$ and a principal point
\[
 p=(A,B,x,y)\in R(X,\cD),\qquad y\in\cD,
\]
for a genus-one surface with two punctures, such that the semisimple data
are compatible and intrinsically $G$-generic.  If $u$ is the image of $p$ in
$M(X,\cD)$, there is an isomorphism of analytic germs
\begin{equation}\label{new:eq-universal-stable-germ}
 (M(X,\cD),u)\times(\C^h,0)
 \simeq (X,x)\times(\C^{2d+z-r},0).
\end{equation}
In particular, a nonnormal conjugacy-closure germ is realized in a
generic relative character variety at a principal point.
\end{theorem}

\begin{proof}
The semisimple parts of all points of $X$ belong to one semisimple class;
choose a representative $s$.  Let $V$ be the principal open subset of
Lemma~\ref{new:lem-principal-commutator}.  The morphism
\[
 f_x:V\longrightarrow q^{-1}(q(s)^{-1}),
 \qquad (A,B)\longmapsto([A,B]x)^{-1}
\]
is dominant, where $q:G\to G_{\ab}$ is the abelianization and
$q(x)=q(s)$ because the unipotent part maps trivially to the torus
$G_{\ab}$.  Therefore its image meets the open subset $W_s$ of
Lemma~\ref{new:lem-add-generic-puncture}.  Choose $(A,B)$ with
$y=f_x(A,B)\in W_s$, and put $\cD=G\cdot y$.  Then
$[A,B]xy=1$.  The data are compatible and generic, and the centralizer of the full
tuple is $Z(G)$, since the handle pair already has that centralizer.

The commutator map is a submersion at $(A,B)$.  The analytic implicit
function theorem, with the possibly singular germ $(X,x)$ as a
parameter, gives
\[
 (R(X,\cD),p)
 \simeq (X,x)\times(\cD,y)\times(\C^{2d-h},0)
 \simeq (X,x)\times(\C^{2d+z-r},0).
\]
For the first isomorphism, use local coordinates identifying the
commutator map with projection onto the smooth $h$-dimensional group
$G_{\der}$, then substitute $[A,B]=(xy)^{-1}$.  This substitution is
valid for all nearby parameters in $X\times\cD$, since their
abelianization product is constantly $1$.  This also shows directly that the
relation scheme is reduced near $p$.

The orbit of $p$ is closed by genericity, and its effective stabilizer
is trivial.  The Luna slice theorem gives an analytic local product
\[
 (R(X,\cD),p)
 \simeq(M(X,\cD),u)\times(\C^h,0).
\]
Combining the two products proves
\eqref{new:eq-universal-stable-germ}.  Normality is invariant under
\'etale localization and adding smooth factors, so nonnormality
passes from $(X,x)$ to $(M(X,\cD),u)$.  No cancellation of smooth
factors is used.
\end{proof}

\begin{corollary}[A normality criterion for the universal assertion]
\label{new:cor-universal-normality}
For a fixed connected reductive group $G$, the following are equivalent:
\begin{enumerate}[label=\textup{(\roman*)}]
\item every conjugacy-class closure in $G$ is normal;
\item every intrinsically generic relative character variety with
      prescribed conjugacy-class closures is normal;
\item every such variety for a genus-one surface with two punctures,
      with a regular semisimple class at the second puncture, is normal.
\end{enumerate}
\end{corollary}

\begin{proof}
If the puncture closures are normal, their product with the handle
factors is normal.  Lemma~\ref{new:lem-closure-transversality} implies
normality of the representation variety, and reductive invariants of a
normal ring are normal.  This proves (i)$\Rightarrow$(ii).
The implication (ii)$\Rightarrow$(iii) is immediate, while
Theorem~\ref{new:thm-universal-realization} proves
(iii)$\Rightarrow$(i).
\end{proof}

The criterion is a deduction from local transversality and realization;
it does not assert that normality alone gives a log Calabi--Yau
compactification.  It identifies a necessary interior condition for a
universal compactification statement.

\subsection{An explicit generic obstruction in type \texorpdfstring{$G_2$}{G2}}

\begin{fact}[The eight-dimensional orbit in $\mathfrak g_2$]
\label{new:fact-g2-orbit}
Let $\mathcal O_8$ be the eight-dimensional nilpotent orbit of type
$\widetilde A_1$ in $\mathfrak g_2$.  The closure of $\mathcal O_8$ is nonnormal, and
\[
 \overline{\mathcal O_8}\setminus\mathcal O_8
   =\mathcal O_6\cup\{0\},
\]
where $\mathcal O_6$ is the six-dimensional minimal nonzero orbit.
In particular $\overline{\mathcal O_8}$ is regular in codimension one.
These facts are classical; see Kraft~\cite[Thm.~1.2(b), \S2.1]{Kraft1989},
Levasseur--Smith~\cite[Thm.~3.13]{LevasseurSmith1988},
Wong~\cite[Introduction]{Wong2016}, and
Fu--Juteau--Levy--Sommers~\cite[\S1.8.4]{FuJuteauLevySommers2017}.
\end{fact}

\begin{lemma}\label{new:lem-g2-origin}
The germ $(\overline{\mathcal O_8},0)$ is nonnormal, is not
Cohen--Macaulay, and its local ring does not satisfy Serre's condition
$S_2$.
\end{lemma}

\begin{proof}
The nonnormal locus of an affine variety over $\C$ is closed.  Here
it is nonempty by Fact~\ref{new:fact-g2-orbit} and is stable under
scalar multiplication.  Indeed scalar multiplication preserves each
nilpotent orbit, for example by the cocharacter of an
$\mathfrak{sl}_2$-triple.  The closure of every scalar orbit contains
$0$, so the nonnormal locus contains $0$.

The variety is regular in codimension one by the orbit dimensions in
Fact~\ref{new:fact-g2-orbit}.  If its local ring at $0$ satisfied $S_2$,
Serre's criterion would make that local ring normal.  This contradicts
the first assertion.  A Cohen--Macaulay local ring satisfies $S_2$,
so the local ring is not Cohen--Macaulay either.
\end{proof}

Let now $G$ be the connected complex group of type $G_2$.  The simply
connected and adjoint forms of $G_2$ coincide, and the character lattice of a
maximal torus is the root lattice.  Fix simple roots $\alpha$ short and
$\beta$ long.  The exponential and logarithm polynomials on nilpotent
and unipotent elements give mutually inverse $G$-equivariant
isomorphisms.  Thus
\[
 X=\exp(\overline{\mathcal O_8})
\]
is a reduced unipotent conjugacy-class closure, nonnormal at $1$.

\begin{theorem}[A generic $G_2$ relative variety without a normal lc compactification]
\label{new:thm-g2-counterexample}
There are intrinsically generic conjugacy-class closures $(X,\cD)$ for
a genus-one surface with two punctures such that $M(X,\cD)$ is
twenty-dimensional and is nonnormal at a principal point $u$.
More precisely, $X=\exp(\overline{\mathcal O_8})$ and one can choose
\[
 \cD=G\cdot t,\qquad
 \alpha(t)=\frac1{12},\qquad \beta(t)=24.
\]
At the indicated point there is an analytic isomorphism
\begin{equation}\label{new:eq-g2-exact-germ}
 (M(X,\cD),u)
 \simeq (\overline{\mathcal O_8},0)\times(\C^{12},0).
\end{equation}
The relative variety is not Cohen--Macaulay and does not satisfy $S_2$
there.  Consequently its connected component containing $u$ admits
neither a normal log-canonical compactification nor a semi-log-canonical
compactification with that component as an unchanged open subvariety.
\end{theorem}

\begin{proof}
Let $w=s_\alpha s_\beta$.  In the basis $(\alpha,\beta)$ of the
character lattice, with columns giving the images of the basis vectors,
\[
 s_\alpha=\begin{pmatrix}-1&3\\0&1\end{pmatrix},\qquad
 s_\beta=\begin{pmatrix}1&0\\1&-1\end{pmatrix},\qquad
 w=\begin{pmatrix}2&-3\\1&-1\end{pmatrix}.
\]
In particular $w^2-w+I=0$ and $\det(I-w)=1$.
Choose $a\in T$ by
$\alpha(a)=2$ and $\beta(a)=3$, let $n\in N_G(T)$ represent $w$,
and put
\[
 \delta(a)=a\,w(a)^{-1},\qquad t=\delta(a)^{-1}.
\]
We use the character-action convention
$(w\lambda)(b)=\lambda(n^{-1}bn)$.  Therefore
$\delta^*=I-w^{-1}=w$, and direct evaluation gives
\[
 \alpha(\delta(a))=2^2\cdot3=12,
 \qquad
 \beta(\delta(a))=2^{-3}\cdot3^{-1}=\frac1{24}.
\]
These are the asserted coordinates of $t$.

Every nonzero character of $T$ is nontrivial on $a$, by unique prime
factorization of $2^m3^n$.  The same is true of $t$: the character
$m\alpha+n\beta$ takes value
\[
 (1/12)^m24^n=2^{-2m+3n}3^{-m+n},
\]
and the integral exponent transformation has determinant $1$.
Both $a$ and $t$ are therefore regular semisimple and have full
centralizer $T$.  Moreover all the parabolic-character tests in
Lemma~\ref{new:lem-add-generic-puncture}, with $s=1$, are avoided.
Thus $(1,t)$ is intrinsically generic.

The tuple
\[
 p=(a,n,1,t)
\]
satisfies the relation because $[a,n]=\delta(a)=t^{-1}$.
An element centralizing $a$ and $n$ lies in $T^w$.
The map $b\mapsto b\,w(b)^{-1}$ is a torus automorphism, since its
character map $I-w^{-1}$ is unimodular.  Hence $T^w=\{1\}$, and $p$
is principal.  Genericity implies that its orbit is closed.  Since
$t=[a,n]^{-1}$, the handle pair itself is irreducible and has closed
orbit as well.

We prove the exact germ statement without canceling smooth factors.
Take an analytic Luna slice $S$ for the handle pair in the smooth
variety $G^2$.  The effective stabilizer is trivial, so $S$ is smooth
of dimension $14$, and locally the quotient of
$G^2\times X\times\cD$ by simultaneous conjugation is
$S\times X\times\cD$.  The relative variety is locally cut out there
by the surface relation.  At $p$, the differential of that relation
in the $G^2$ handle directions is surjective by
Lemma~\ref{new:lem-principal-commutator}.  Passing from handle
directions to the slice discards the directions of simultaneous
conjugation of the handle pair.  Each discarded contribution can be
replaced by the opposite contribution from conjugating $t$ in $\cD$:
simultaneous conjugation of the full tuple has zero derivative on the
relation, and conjugating the entry $1\in X$ contributes zero.
Consequently the differential in the smooth $S\times\cD$ directions
is still surjective.  The implicit function theorem with $X$ as a
parameter gives
\[
 (M(X,\cD),u)\simeq (X,1)\times
 (\C^{\dim S+\dim\cD-\dim G},0)
 =(X,1)\times(\C^{12},0).
\]
Polynomial logarithm identifies $(X,1)$ with
$(\overline{\mathcal O_8},0)$, proving
\eqref{new:eq-g2-exact-germ}.  The dimension is $8+12=20$.

Lemma~\ref{new:lem-g2-origin} and invariance under smooth extension
give nonnormality and failure of Cohen--Macaulay\-ness and $S_2$.
An open subvariety of a normal variety is normal, so no normal
log-canonical compactification can preserve this interior.
A semi-log-canonical variety is $S_2$, and $S_2$ is preserved by
passing to open subvarieties, so that weaker category does not remove
the obstruction either.
\end{proof}

\begin{remark}[Why normalization is necessary]
\label{new:rem-counterexample-scope}
The full reduced adjoint fiber of $1\in G$ is the unipotent variety,
of dimension $12$.  Replacing $X$ by this fiber gives the generic
spectral character variety of dimension $24$, to which
Theorem~\ref{thm:spectral-main} applies.  The closed relative
subvariety in Theorem~\ref{new:thm-g2-counterexample} has a defective
interior even though the containing spectral variety has the asserted
log-canonical compactification.  Thus taking a closed puncture
restriction cannot automatically preserve the compactification theorem.

Theorem~\ref{r2:thm-full-normalized-logcy} compactifies the normalization
of this relative variety. The example shows why that normalization
is necessary: the original reduced interior cannot be an open
subvariety of a normal compactification.
\end{remark}

\subsection{A normal generic interior without a boundary-preserving
\texorpdfstring{$\Q$}{Q}-factorial model}
\label{new:sl3-section}

Normalization removes the obstruction exhibited by the nonnormal
$G_2$ example, but it does not recover every interior assertion of the
spectral theorem.  In particular, the following example separates the
existence of a normal log canonical compactification from the stronger
requirement that a $\Q$-factorial dlt model leave the entire interior
unchanged.  The argument proves algebraic local non-$\Q$-factoriality;
an analytic class-group calculation alone would not suffice.

\begin{theorem}\label{new:sl3-main}
Let $G=\SL_3$, let the surface have genus one and two punctures, and
prescribe the closures of the following conjugacy classes:
\[
 \cC_1=\{\text{unipotent elements of Jordan type }(2,1)\},
 \qquad
 \cC_2=G\cdot t,\qquad t=\diag(2,3,1/6).
\]
These data are compatible and intrinsically generic.  The relative
character variety
\[
 X=\Bigl\{(A,B,M_1,M_2)\in
 G^2\times\overline{\cC_1}\times\cC_2:
 [A,B]M_1M_2=1\Bigr\}/\!\!/G
\]
is reduced, normal and pure of dimension ten, but is not
$\Q$-factorial at a point $p$ described below.  There is a smooth
variety $\widetilde X$ and a projective morphism
$\pi:\widetilde X\to X$ which is birational on each component,
is an isomorphism when $M_1\ne I$, and has a fiber $\PP^2$ over $p$.
The exceptional locus of this morphism has codimension two in $\widetilde X$.

Consequently no $\Q$-factorial compactification can contain the whole
of this interior as an open subvariety.  In particular, no
$\Q$-factorial dlt model of a compactification can be isomorphic over
the whole interior.
\end{theorem}

\begin{proof}
Put $V=\C^3$ and write
\[
 Y=\{N\in\End(V):\rank N\le1,\ \tr N=0\}.
\]
For a rank-one endomorphism $N=u\otimes v$ one has
$N^2=v(u)N=(\tr N)N$.  Thus every element of $Y$ satisfies $N^2=0$,
and $N\mapsto I+N$ identifies $Y$ with
$\overline{\cC_1}$.  The complement of the origin in $Y$ is the single
nilpotent orbit of type $(2,1)$ and has dimension four.

\emph{Normality of the puncture closure.}
Consider the affine quadric
\[
 Q=\{(u,v)\in V\times V^*:v(u)=0\}.
\]
This variety is an irreducible hypersurface of dimension five, smooth away from
the origin, and hence normal by Serre's criterion.  The action
$z\cdot(u,v)=(zu,z^{-1}v)$ of $\C^*$ has invariant functions
generated by the $u_i v_j$.  Their relations are the rank-one
relations together with $\sum_i u_i v_i=0$.  Indeed, before imposing
the last equation the invariant ring is the coordinate ring of the
rank-at-most-one matrix variety; exactness of $\C^*$-invariants
commutes with quotienting by the invariant equation $v(u)$.  Therefore
$Y\simeq Q/\!\!/\C^*$ is normal.

\emph{Genericity and a specified point.}
The semisimple puncture classes are represented by $I$ and $t$.
For $\SL_3$, the proper parabolic test is the product of the
eigenvalues on a proper invariant subspace.  No nonempty proper
subset of $\{2,3,1/6\}$ has product one: the one-element products
are not one, and the two-element products are their complementary
inverses.  Hence the data are intrinsically generic; compatibility
follows from $\det t=1$.

Choose $c\in\C$ with $c^3=3/2$, and set
\[
 a=\diag(c,c/3,2c),\qquad
 n=\begin{pmatrix}0&0&1\\1&0&0\\0&1&0\end{pmatrix}.
\]
Both matrices lie in $\SL_3$, and direct multiplication gives
\[
 [a,n]=\diag(1/2,1/3,6)=t^{-1}.
\]
Thus $r=(a,n,I,t)$ is a point of the representation variety $R$,
and we set $p=[r]\in X$.  The eigenvalues of $a$ are distinct, so
its centralizer is the diagonal torus; a diagonal matrix commuting
with the three-cycle $n$ is scalar.  Hence
$Z_G(a,n)=Z(G)=\mu_3$.

In fact, intrinsic genericity makes every tuple of $R$ irreducible,
and Schur's lemma gives stabilizer exactly $Z(G)$ for every tuple.
The group-direction argument of
Lemma~\ref{lem:spectral-relation-local} applies verbatim to
$G^2\times Y\times\cC_2$: it uses conjugation invariance of the
puncture spaces, rather than their being full adjoint fibers.
The local product lemma identifies germs of that product with germs of $R$ times a smooth
factor of dimension eight.  Thus $R$ is reduced and normal.
All orbits are closed, and the effective group $H=\PGL_3$ acts
freely.  Luna's slice theorem makes $R\to X$ an $H$-torsor in the
\'etale topology.  In particular $X$ is normal.  The dimensions are
\[
 \dim R=2\cdot8+4+6-8=18,
 \qquad \dim X=18-8=10.
\]

\emph{The flag morphism.}
Define
\[
 \widetilde Y=
 \{(N,\ell)\in Y\times\PP(V):
       N(V)\subseteq\ell\subseteq\ker N\}.
\]
Over a line $\ell$, its fiber is
$\Hom(V/\ell,\ell)$, so $\widetilde Y\simeq T^*\PP^2$.
The projection $\widetilde Y\to Y$ is projective, since the
displayed incidence conditions are closed.  It is an isomorphism
over $N\ne0$, where $\ell=\operatorname{im}N$ is uniquely
determined, and its fiber over zero is the zero section $\PP^2$.
This exceptional locus has codimension two in the smooth
four-dimensional $\widetilde Y$.

Let $\widetilde R$ be the pullback of this flag morphism through
$R\to Y$, using $N=M_1-I$.  Equivalently, impose the surface
relation on the smooth variety
$G^2\times\widetilde Y\times\cC_2$.  The group directions include
simultaneous conjugation of $N$ and its line $\ell$.  Their word
differential is onto $\mathfrak{sl}_3$ at every point, since the
underlying tuple is irreducible.  Hence $\widetilde R$ is smooth of
dimension eighteen.  The effective $H$-action remains free.

The $H$-torsor $R\to X$ gives the associated
Brauer--Severi bundle
\[
 \mathcal P=R\times^H\PP^2\longrightarrow X.
\]
This variety is projective: $\mathcal O_{\PP^2}(3)$ has its natural
$\PGL_3$-linearization and descends to a relatively ample line
bundle.  The third power is relevant here;
$\mathcal O_{\PP^2}(1)$ itself does not have a
$\PGL_3$-linearization.  The invariant incidence conditions define
a closed subscheme $\widetilde X\subset\mathcal P$, whose pullback
to $R$ is $\widetilde R$.  This constructs the projective morphism
$\pi:\widetilde X\to X$.  Smoothness descends from
$\widetilde R$, so $\widetilde X$ is smooth of dimension ten.

The morphism is an isomorphism over $M_1\ne I$.  This open set is
dense in every component by the same group-direction local product,
which preserves puncture orbit strata.  Over a tuple with $M_1=I$
the fiber is $\PP^2$, because the effective stabilizer is trivial.
In particular $\pi^{-1}(p)=\PP^2$.

The exceptional locus has codimension two.  Indeed, before imposing
the surface relation it is the zero section in $\widetilde Y$,
which has codimension two and is preserved by the group directions.
The word differential restricted to this locus is still onto:
the same underlying tuples are irreducible.  The relation fiber of this locus
has dimension
\[
 2\cdot8+2+6-8=16,
\]
and quotienting by $H$ gives an eight-dimensional locus in the
ten-dimensional $\widetilde X$.  These are all exceptional points.
Thus $\pi$ has no exceptional divisors and is a nontrivial small
projective birational morphism on the component through $p$.

\emph{Algebraic non-$\Q$-factoriality.}
Suppose that $X$ were $\Q$-factorial on a neighborhood $W$ of $p$.
Choose a relatively ample Cartier divisor $L$ on
$\pi^{-1}(W)$ and let $D=\pi_*L$ be its Weil pushforward.
Some positive multiple $mD$ is Cartier after shrinking $W$.
The divisors $mL$ and $\pi^*(mD)$ agree outside the exceptional
locus.  Their difference is supported on exceptional divisors,
of which there are none, so
\[
 mL=\pi^*(mD).
\]
On a line in the fiber $\pi^{-1}(p)=\PP^2$, the left side has
positive degree and the right side has degree zero, a contradiction.
Thus $X$ is not $\Q$-factorial at $p$.  Since $\Q$-factoriality is
preserved on open subvarieties, the assertions about compactifications
and interior-preserving $\Q$-factorial dlt models follow.
\end{proof}

\begin{remark}\label{new:sl3-scope}
This obstruction does not disprove the existence of a normal log
canonical log Calabi--Yau compactification of this $X$.
This example shows that the additional $\Q$-factorial and
interior-preserving dlt conclusions of the spectral theorem require
revision for smaller conjugacy-class closures, even when those
closures are already normal.  Passing to normalization changes
nothing in the present example.
\end{remark}

\subsection*{Acknowledgments \& AI disclosure}
The author is grateful to Davesh Maulik for an encouraging conversation
around the time of ICBS 2025, when work on this project began in earnest,
and to Zhe Sun for suggesting a generalization of the results of
``On compactifications of the $\mathrm{SL}(2,\mathbb{C})$ character
varieties of punctured surfaces'' by Mohammad Farajzadeh-Tehrani
and Charles Frohman. Sun's suggestion motivated the use of ideal
triangulations to compactify character varieties for arbitrary connected
reductive groups $G$. The author also thanks Farajzadeh-Tehrani for
helpful correspondence concerning the log Calabi--Yau claim in an
earlier version of that paper and for sharing the accepted version.

The project has been developed with substantial AI assistance since 11 August 2026. 
The main results were obtained after weeks of discussion with AI, in which the author supplied the main ideas and directed the work
while the AI served as a tool to check feasibility and help develop technical arguments. 
The source draft records the author's independent checking and substantial revision of the arguments and their references.

\addtocontents{toc}{\protect\setcounter{tocdepth}{1}}

\end{document}